\pdfoutput=1
\documentclass[a4paper]{amsart}

\usepackage[T1]{fontenc}
\usepackage{amssymb}
\usepackage{mathtools}
\usepackage{a4wide}
\usepackage{enumitem}
\usepackage[dvipsnames]{xcolor}
\usepackage[all,cmtip]{xy}
\usepackage{amsaddr}
\usepackage[labelformat=simple]{subcaption}
\usepackage{hyperref}
\usepackage[abbrev, msc-links, nobysame]{amsrefs}
\usepackage{aliascnt}
\usepackage{cleveref}
\usepackage[overload]{textcase}

\setlist[itemize,1]{left = 0pt}
\setlist[enumerate,1]{left = 0pt}

\numberwithin{equation}{section}
\counterwithout{figure}{section} %

\RenewDocumentCommand{\paragraph}{m}{\par\noindent\textit{#1.}\ }

\hypersetup{
	colorlinks=true,
	linkcolor=black,
	citecolor=black
}

\BibSpec{arxiv}{%
  +{}{\PrintAuthors}{author}
  +{,}{ \textit}{title}
  +{}{ \parenthesize}{date}
  +{,}{ arXiv:}{eprint}
}

\theoremstyle{plain}
\newtheorem{theorem}{Theorem}
\newtheorem{prop}{Proposition}[section]
\newaliascnt{lemma}{prop}\newtheorem{lemma}[lemma]{Lemma}\aliascntresetthe{lemma}
\newaliascnt{cor}{prop}\newtheorem{cor}[cor]{Corollary}\aliascntresetthe{cor}
\newtheorem*{claim}{Claim}

\theoremstyle{definition}
\newaliascnt{defn}{prop}\newtheorem{defn}[defn]{Definition}\aliascntresetthe{defn}
\newaliascnt{remark}{prop}\newtheorem{remark}[remark]{Remark}\aliascntresetthe{remark}
\newaliascnt{example}{prop}\newtheorem{example}[example]{Example}\aliascntresetthe{example}

\let\Sec\S %

\crefformat{section}{\Sec\!\:#2#1#3}
\crefrangeformat{section}{\Sec\Sec\!\:#3#1#4~to~#5#2#6}
\crefmultiformat{section}{\Sec\Sec\!\:#2#1#3}{ and~#2#1#3}{, #2#1#3}{ and~#2#1#3}
\crefformat{enumi}{(#2#1#3)}
\crefrangeformat{enumi}{(#3#1#4)--(#5#2#6)}
\crefmultiformat{enumi}{(#2#1#3)}{ and~(#2#1#3)}{, (#2#1#3)}{ and~(#2#1#3)}
\crefrangemultiformat{enumi}{(#3#1#4)--(#5#2#6)}{ and~(#3#1#4)--(#5#2#6)}{, (#3#1#4)--(#5#2#6)}{ and~(#3#1#4)--(#5#2#6)}
\crefname{theorem}{Theorem}{Theorems}
\crefname{prop}{Proposition}{Propositions}
\crefname{lemma}{Lemma}{Lemmas}
\crefname{cor}{Corollary}{Corollaries}
\crefname{defn}{Definition}{Definitions}
\crefname{remark}{Remark}{Remarks}
\crefname{example}{Example}{Examples}
\crefname{figure}{Fig.}{Figs.}
\crefname{equation}{}{}

\DeclareFontFamily{U}{matha}{\hyphenchar\font45}
\DeclareFontShape{U}{matha}{m}{n}{<-6> matha5 <6-7> matha6 <7-8> matha7
  <8-9> matha8 <9-10> matha9 <10-12> matha10 <12-> matha12}{}
\DeclareSymbolFont{matha}{U}{matha}{m}{n}
\DeclareMathSymbol{\orightsymbol}{\mathbin}{matha}{"69}
\newcommand{\oright}[1]{\mathbin{\underset{#1}{\orightsymbol}}}

\DeclareFontFamily{U}{skakfig}{}
\DeclareFontShape{U}{skakfig}{m}{n}{<-> SkakNew-Figurine}{}
\DeclareSymbolFont{skakfig}{U}{skakfig}{m}{n}
\DeclareMathSymbol{\withsymbol}{\mathbin}{skakfig}{"76}
\newcommand{\with}[1]{\mathbin{\underset{#1}{\withsymbol}}}

\let\Re\relax
\let\Im\relax
\DeclareMathOperator{\Area}{Area}
\DeclareMathOperator{\colim}{colim}
\DeclareMathOperator{\Cone}{Cone}
\DeclareMathOperator{\diag}{diag}
\DeclareMathOperator{\Diff}{Diff}
\DeclareMathOperator{\ex}{ex}
\DeclareMathOperator{\Ext}{Ext}
\DeclareMathOperator{\Fun}{Fun}
\DeclareMathOperator{\GL}{GL}
\DeclareMathOperator{\gr}{gr}
\DeclareMathOperator{\Ham}{Ham}
\DeclareMathOperator{\Hom}{Hom}
\DeclareMathOperator{\id}{id}
\DeclareMathOperator{\Im}{Im}
\DeclareMathOperator{\Imm}{Imm}
\DeclareMathOperator{\Int}{Int}
\DeclareMathOperator{\Lag}{Lag}
\DeclareMathOperator{\MCG}{MCG}
\DeclareMathOperator{\Ob}{Ob}
\DeclareMathOperator{\Pin}{Pin}
\DeclareMathOperator{\pr}{pr}
\DeclareMathOperator{\Re}{Re}
\DeclareMathOperator{\rk}{rk}
\DeclareMathOperator{\SO}{SO}
\DeclareMathOperator{\so}{so}
\DeclareMathOperator{\Spec}{Spec}
\DeclareMathOperator{\Symp}{Symp}
\DeclareMathOperator{\Tw}{Tw}
\DeclareMathOperator{\vdim}{vdim}

\DeclareMathOperator{\dgcat}{dgcat}
\newcommand{\Mod}{\mathrm{Mod}}
\DeclareMathOperator{\Dp}{D_{perf}}
\DeclareMathOperator{\Dpdg}{D_{perf}^{dg}}
\DeclareMathOperator{\Shp}{Sh_{perf}}
\DeclareMathOperator{\Twp}{Tw_{perf}}
\DeclareMathOperator{\Uproj}{\U\text{-}proj}
\DeclareMathOperator{\Wgp}{\mathcal W^{\mathrm{geom}}_{\mathrm{pre}}}
\DeclareMathOperator{\Wgsu}{\mathcal W^{\mathrm{geom}}_{\mathrm{su}}}
\DeclareMathOperator{\Wt}{\mathcal W^{\mathrm{top}}}

\let\AA\relax
\newcommand{\AA}{\mathbb{A}}
\newcommand{\CC}{\mathbb{C}}
\newcommand{\FF}{\mathbb{F}}
\newcommand{\LL}{\mathbb{L}}
\newcommand{\PP}{\mathbb{P}}
\newcommand{\RR}{\mathbb{R}}
\renewcommand{\SS}{\mathbb{S}}
\newcommand{\TT}{\mathbb{T}}
\newcommand{\ZZ}{\mathbb{Z}}

\let\S\relax
\let\H\relax
\newcommand{\A}{\mathcal A}
\newcommand{\B}{\mathcal B}
\newcommand{\C}{\mathcal C}
\newcommand{\D}{\mathcal D}
\newcommand{\E}{\mathcal E}
\newcommand{\F}{\mathcal F}
\newcommand{\G}{\mathcal G}
\newcommand{\H}{\mathcal H}
\newcommand{\I}{\mathcal I}
\newcommand{\J}{\mathcal J}
\renewcommand{\L}{\mathcal L}
\newcommand{\M}{\mathcal M}
\renewcommand{\O}{\mathcal O}
\renewcommand{\P}{\mathcal P}
\newcommand{\R}{\mathcal R}
\newcommand{\S}{\mathcal S}
\newcommand{\T}{\mathcal T}
\newcommand{\U}{\mathcal U}
\newcommand{\V}{\mathcal V}
\newcommand{\W}{\mathcal W}
\newcommand{\LS}{\mathcal{LS}}
\newcommand{\SU}{\mathcal{SU}}

\renewcommand{\k}{\mathrm{k}}
\newcommand{\w}{\mathrm{w}}

\newcommand{\be}{\mathbf{e}}
\newcommand{\bm}{\mathbf{m}}
\newcommand{\bt}{\mathbf{t}}
\newcommand{\bu}{\mathbf{u}}
\newcommand{\bC}{\mathbf{C}}
\newcommand{\bL}{\mathbf{L}}
\newcommand{\bS}{\mathbf{S}}
\newcommand{\bT}{\mathbf{T}}
\newcommand{\bX}{\mathbf{X}}
\newcommand{\bY}{\mathbf{Y}}
\newcommand{\bmu}{\boldsymbol{\mu}}
\newcommand{\bph}{\boldsymbol{\varphi}}
\newcommand{\bPhi}{\boldsymbol{\Phi}}
\newcommand{\bGamma}{{\boldsymbol{\Gamma}}}

\newcommand{\eps}{\varepsilon}
\newcommand{\ph}{\varphi}
\newcommand{\pd}{\partial}
\newcommand{\res}[2]{\left.{{#1}}\right|_{{#2}}}

\title{Fukaya categories of surfaces with DG coefficients}
\author{Daniil Mamaev}

\address{Centre for Geometry and Physics, Uppsala University, 75106 Uppsala, Sweden}
\email{daniil.mamaev@math.uu.se}

\begin{document}

\begin{abstract}
Given a compact surface $S$ admitting a Liouville form, a finite set $Q \subset S$, a line field on $S$, and a local system of dg-categories on $S \setminus Q$ together with a section of the monodromy functor around each interior $q \in Q$, we construct a Fukaya category.
We give examples where such Fukaya categories are quasi-equivalent to derived categories of projective completions of line bundles and of iterated blow-ups of the affine plane; and to wrapped Fukaya categories of cotangent bundles of spheres.

We define two versions of the Fukaya category, a topological, which is constructed combinatorially from a full arc system, and a geometric, which includes all Lagrangian branes supported on immersed curves.
We prove that the topological category embeds into the geometric category in an essentially unique way.
To define Lagrangian branes we develop a generalisation of the Chekanov--Eliashberg algebra incorporating dg-coefficients and wrapping along the boundary.
\end{abstract}

\maketitle

\vspace*{-1\baselineskip}
{\small\tableofcontents}
\removelastskip\newpage

\addtocontents{toc}{\protect\setcounter{tocdepth}{1}}
\section*{Introduction}
\label{s: Introduction}

Let $S$ be a compact oriented surface admitting a Liouville form (i.e.\ each component of $S$ has non-empty boundary) and $Q \subset S$ a finite set, called the \emph{marking}.
The \emph{directed relative path 2-category} $\P^+_2(S, Q)$ is the (strict) 2-category whose objects are $S \setminus Q$, 1-morphisms are lazy paths in $S \setminus Q$ up to thin homotopy, and 2-morphisms are isotopy classes of homotopies of lazy paths in $S$ that intersect $Q$ positively, see \cref{d: directed relative path 2-category} for details.

Let $\k$ be a commutative ring and $\dgcat_\k$ the 2-category whose objects are small dg-categories over $\k$, 1-morphisms are dg-functors, and 2-morphisms are closed degree zero morphisms of dg-functors.
A \emph{weighted system of dg-categories} on $(S, Q)$ is a strict 2-functor $\C\colon \P^+_2(S, Q) \to \dgcat_\k$.
Informally one can think of $\C$ as a local system of dg-categories on $S \setminus Q$ (it is flattened by the coherent isomorphisms of dg-functors given by homotopies in $S \setminus Q$) together with a morphism of dg-functors from $\id$ to the monodromy around $q$ for each interior $q \in Q$.
The points of $Q \cap \pd S$ are stops, i.e.\ boundary paths are not allowed to pass through them.

Given a line field $\eta$ on $S$ and a weighted system of dg-categories $\C\colon \P^+_2(S, Q) \to \dgcat_\k$ we call $(\S = (S, Q, \eta); \C)$ a \emph{weighted surface}.

\medskip
\begin{center}\emph{
	In this paper we define and study Fukaya categories of weighted surfaces.
}\end{center}
\medskip

There is a related approach that associates the $\infty$-category of global sections to a perverse schober on a ribbon graph, see Kapranov--Schechtman~\cite{KS18} and Christ--Haiden--Qiu~\cite{CHQ26}.
We expect that \cite{CHQ26}*{Proposition 7.7} can be extended to associate a perverse schober $\G$ on the skeleton of any weighted surface $(\S; \C)$ in such a way that the global sections of $\G$ are equivalent to the $\infty$-categorical Ind-completion of the Fukaya category of $(\S; \C)$.
We will address this elsewhere.

\subsection*{Topological Fukaya category}
Let $(\S = (S, Q, \eta); \C)$ be a weighted surface. An \emph{arc brane} in $(\S; \C)$ is a properly embedded segment $\gamma\colon [0, 1] \to S \setminus Q$ (an \emph{arc}) together with a grading (i.e.\ a homotopy class of paths from $\gamma^*\eta$ to the tangent line $\dot \gamma$) and an object $M$ in $\C_{\gamma(0)}$.
A set of arc branes $\bGamma$ is called an \emph{arc brane system} if the underlying curves form a no more than countable set $\Gamma$ and arcs in $\Gamma$ are pairwise disjoint.
An arc brane system $\bGamma$ is called \emph{full} if $\Gamma$ cuts $S$ into topological discs such that each disc has at most one point from $Q$ and for each arc $\gamma \in \Gamma$ the objects in $\C_{\gamma(0)}$ obtained from arc branes supported on $\gamma$ generate $\C_{\gamma(0)}$ (i.e.\ the smallest full subcategory of $\Tw \C_{\gamma(0)}$ containing these objects and closed under taking shifts, cones, and homotopy equivalent objects contains $\C_{\gamma(0)}$).

Generalising the construction of Haiden--Katzarkov--Kontsevich~\cite{HKK17}*{\Sec 3} to our setting we associate an $A_\infty$-category $\V_\bGamma$ to any arc brane system $\bGamma$.
More concretely, $\Ob \V_\bGamma = \bGamma$, morphism spaces come from boundary paths in $\pd S \setminus Q$ between the ends of arcs in $\bGamma$, and $A_\infty$-operations come from the differential and composition in fibres of $\C$, concatenation of boundary paths, and $\C$ on 2-morphisms obtained from immersed polygons in $S$ bounded by arcs from $\bGamma$, see \cref{subsubs: A_infty-category associated to an arc system} for details.
Our first result is the following.

\begin{theorem} \label{t: A_infty-categories of arc systems}
Let $(\S = (S, Q, \eta); \C)$ be a weighted surface.
Then
\begin{enumerate}
	\item $\V_\bGamma$ is a strictly unital $\ZZ$-graded $\k$-linear $A_\infty$-category for any arc brane system $\bGamma$ in $(\S; \C)$.
	\item For a full arc brane system $\bGamma$ the triangulated envelope $\Tw \V_\bGamma$ does not depend on the choice of $\bGamma$ up to homotopy equivalence of strictly unital $A_\infty$-categories.
\end{enumerate}
\end{theorem}

By \emph{triangulated} in \cref{t: A_infty-categories of arc systems} and everywhere in the paper we mean closed under shifts and cones, not necessarily idempotent complete.
We define the \emph{topological Fukaya category} $\Wt(\S; \C)$ of $\S$ with coefficients in $\C$ as (the homotopy equivalence class of) $\Tw \V_\bGamma$ for a full arc brane system $\bGamma$ in $(\S; \C)$.
This definition generalises topological Fukaya categories by Haiden--Katzarkov--Kontsevich~\cite{HKK17} and relative wrapped Fukaya categories of surfaces by Evans--Lekili~\cite{EL26}.
\cref{t: A_infty-categories of arc systems} follows from \cref{t: topological Fukaya category inside geometric} below.

\subsection*{Homological mirror symmetry}
We call any geometrically motivated quasi-equivalence between a Fukaya category (A-side) and a derived category of algebro-geometric origin (B-side) or representation-theoretic origin (D-side) homological mirror symmetry (HMS).
The subject starts with Kontsevich's ICM address~\cite{K95}, by now it is a vast area of research and we will not attempt to give an overview.

In the topological Fukaya category $\Wt(\S; \C)$ we have a geometric input given by $\S$ and an algebraic input given by fibres of $\C$ (with the rest of the structure of $\C$ thought of as relating the two), and so we can consider it on the A- or D-side.
The theorem below is intended to illustrate the kind of HMS statements involving $\Wt(\S; \C)$ one can obtain and is not at all exhaustive.
$\Wt(\S; \C)$ is on the A-side in (1)\&(2) and on the D-side in (3).

\begin{theorem} \label{t: HMS}
	Each of the following is quasi-equivalent to $\Wt(\S; \C)$ for the indicated $\S$ and a suitable $\C$.
	\begin{enumerate}
		\item $\Dpdg \PP(\O_X \oplus \L)$ where $X$ is a separated Noetherian scheme and $\L$ a line bundle on $X$.
		The graded marked surface in this case is
		${\S = (S = \{z \in \CC \colon 1 \le |z| \le 2\}, Q = \{1, -2\}, \eta \text{ radial})}$.
		The fibres of $\C$ are $\Dpdg X$.
		\item $\Dpdg X$ where $X$ is a bamboo forest blow-up of the affine plane $\AA^2_\k$ over an infinite field $\k$.
		The graded marked surface in this case is
		$\S = (S = \{z \in \CC \colon 1 \le |z| \le 2\}, Q, \eta \text{ radial})$
		with $Q \cap \pd S = \{1\}$.
		The fibres of $\C$ are $\k[t]$ viewed as a dg-category with one object.
		\item The wrapped Fukaya category $\W(T^*S^n; \k)$ of $T^*S^n$ with coefficients in a commutative ring $\k$ for any $n \ge 2$.
		The graded marked surface in this case is
		${\S = (S = \{z \in \CC \colon |z| \le 2\}, Q = \{1, -1\}, \eta \text{ vertical})}$.
		For $n \ge 3$ the fibres of $\C$ are quasi-equivalent to $\W(T^*S^{n - 1}; \k)$ and for $n = 2$ to $\Tw \k[z]$.
	\end{enumerate}
\end{theorem}

$\Dpdg$ above is a dg-enhancement of the derived category of perfect complexes, it is known to be unique up to quasi-equivalence in the setting of \cref{t: HMS}\,(1)\&(2) by Canonaco--Neeman--Stellari~\cite{CNS22}.

Taking $X = \Spec \k$ in \cref{t: HMS}\,(1) we recover the classical homological mirror symmetry for $\PP^1$ going back to Kontsevich~\cite{K09}*{Pictures}.
Our result is therefore a relative version of this over $X$.
It can probably be extended to any rank two vector bundle $\V$ on $X$ in place of $\O_X \oplus \L$ by representing $\V$ as a twisted complex of line bundles.
We prove \cref{t: HMS}\,(1) in \cref{subs: Total spaces of projective completions of line bundles}.

The \emph{bamboo forest blow-up} in \cref{t: HMS}\,(2) is a successive blow-up of $\AA^2_{\k}$ at $\k$-points such that the exceptional locus is a disjoint union of chains of smooth rational curves with self-intersections $-2, \ldots, -2, -1$.
In the proof of \cref{t: HMS}\,(2) we use a projection $X \to \AA^1$ with regular fibre $\AA^1$ and singular fibres chains $\AA^1-\PP^1-\ldots-\PP^1$.
Our result is therefore a relative version of the mirror symmetry, proven by Lekili--Polishchuk~\cite{LP18}, between $\S$ as in \cref{t: HMS}\,(2) and the chain $\AA^1-\PP^1-\ldots-\PP^1$ with $|Q \cap \Int S|$ projective lines.
A similar relative statement for versal deformations of the chain $\AA^1-\PP^1-\ldots-\PP^1-\AA^1$ is proven by Evans--Lekili~\cite{EL26}.
We prove \cref{t: HMS}\,(2) in \cref{subs: Bamboo forest blow-ups of affine plane}.

\cref{t: HMS}\,(3) for $n = 1$ holds trivially, but with $S$ an annulus.
For $n = 2$ it is proven by Lekili--Segal~\cite{LS26} as the starting point for their conjectures on algebraic torus fibrations.
For $n \ge 3$ our construction is different from their framework as the regular fibre of the fibration we use is $T^*S^{n - 1}$, not a torus $(\CC^\times)^{n - 1}$.
More precisely, the standard Lefschetz fibration $T^*S^n \to \CC$ is singular over $Q = \{1, -1\}$ and has $T^*S^{n - 1}$ as regular fibre.
Moreover, $(\S; \C)$ can be viewed as a singular symplectic reduction of $(T^*S^n, \k)$
with respect to a Hamiltonian action of $\SO(n)$.
We prove \cref{t: HMS}\,(3) in \cref{subs: Cotangent bundles of spheres}.

\subsection*{Geometric Fukaya \texorpdfstring{$A_\infty$}{A∞}-pre-category}
\label{subs: Geometric Fukaya A_infty-pre-category}
We give a second, geometric construction of a Fukaya category associated to a graded marked surface $(S, Q, \eta)$ with coefficients in $\C\colon \P^+_2(S, Q) \to \dgcat_\k$.
To this end we choose a Liouville form $\lambda$ on $S$ \emph{compatible} with $Q$ (see \cref{d: Liouville structure on marked surface}).
We call $(\S = (S, Q, \eta, \lambda); \C)$ a \emph{weighted Liouville surface} and define the \emph{geometric Fukaya $A_\infty$-pre-category} $\Wgp = \Wgp(\S; \C)$ of $\S$ with coefficients in $\C$.
The reason we introduce $\Wgp$ is to have arbitrary Lagrangian branes supported on immersed curves in $S \setminus Q$ as objects.
On a more technical level, we use the geometric category to prove \cref{t: A_infty-categories of arc systems} about the topological category.

$A_\infty$-pre-categories were introduced by Kontsevich--Soibelman~\cite{KS01}.
An $A_\infty$-pre-category $\W$ is a set of objects $\Ob \W$, a collection of transverse sequences in $(\Ob \W)^n$ for each $n \ge 1$ such that all sequences of length one are transverse and subsequences of transverse sequences are transverse, a graded module of morphisms $\W(X_0, X_1)$ for each transverse sequence $(X_0, X_1)$, and an $A_\infty$-operation
$$
	\bmu^n\colon \W(X_{n - 1}, X_n) \otimes \ldots \otimes \W(X_0, X_1) \to \W(X_0, X_n)
$$
of degree $2 - n$ for each $n \ge 1$ and transverse sequence $(X_0, \ldots, X_n)$ such that the $A_\infty$-relations are satisfied, see \cref{subsubs: A_infty-pre-categories} for details.

The objects of $\Wgp$ are Lagrangian branes, which are certain exact immersed curves $\gamma$ in $S \setminus Q$ with additional brane data on them.
It takes some work to define what a Lagrangian brane is.
Namely, a \emph{Lagrangian pre-brane} $\L$ in $(\S; \C)$ is a compact proper generic immersion
${\gamma \colon (I, \pd I) \to (S \setminus Q, \pd S \setminus Q)}$
of a possibly disconnected $1$-dimensional manifold $I$
together with a generic \emph{potential} $P\colon I \to \RR$ (i.e.\ $\gamma^*\lambda + dP = 0$), a grading, and a local system of objects over $\gamma^*\C$, see \cref{d: Lagrangian pre-brane}.
We define a curved $A_\infty$-algebra $\A_\L$, the \emph{obstruction algebra} of the Lagrangian pre-brane $\L$.
Note that the pair $(\gamma, P)$ defines an embedded Legendrian (the \emph{Legendrian lift}) in the contactisation $((S \setminus Q) \times \RR_z, \lambda + dz)$.
The bar construction on our $\A_\L$ is a generalisation of the Chekanov--Eliashberg algebra of a Legendrian link introduced by Chekanov~\cite{C02} and Eliashberg--Givental--Hofer~\cite{EGH00}.
The generalisation goes in two directions, one is that $(\gamma, P)$ can go to the boundary, cf.\ Sivek~\cite{S11}, and another is that the coefficients are in $\C$ instead of a commutative ring.

The elements of $\A_\L$ are given by double points of $\gamma$ and Reeb chords of positive action in $\pd S \setminus Q$ between the ends of $\gamma$ (each such point or chord gives a graded $\k$-module of maps between respective fibres of the local systems).
The $A_\infty$-operations $\mu^n, \, n \ge 0$ in $\A_\L$ come from the differential and composition in fibres of $\C$, concatenation of Reeb chords, and $\C$ on 2-morphisms obtained from immersed polygons in $S$ bounded by $\gamma$, see \cref{subsubs: Obstruction algebra} for details.

A \emph{Lagrangian brane} is a pair $\LL = (\L, \delta)$ where $\L$ is a Lagrangian pre-brane and $\delta \in \A_{\L}$ a bounding cochain (i.e.\ a solution to the Maurer--Cartan equation
${\mu^0(1_\k) + \mu^1(\delta) + \mu^2(\delta, \delta) + \ldots = 0}$, the sum is finite since $\A_\L$ is nilpotent).
The objects of $\Wgp$ are the Lagrangian branes.
A sequence of Lagrangian branes is transverse if the underlying curves are in general position, see \cref{subsubs: Transverse Lagrangians} for details.
The morphism space $\Wgp(\LL_0, \LL_1)$ is given by the intersection points between $\gamma_0$ and $\gamma_1$ and the Reeb chords in $\pd S \setminus Q$ from an end of $\gamma_0$ to an end of $\gamma_1$.
The $A_\infty$-structure $\bmu$ is defined similarly to that in the obstruction algebra, with bounding cochains deforming $\mu$ so that $\bmu^0 = 0$, see \cref{subsubs: Fukaya category} for details.

The geometric Fukaya $A_\infty$-pre-category defined above enjoys the following properties.
\begin{theorem} \label{t: geometric Fukaya A_infty-pre-category is good}
Let $(\S = (S, Q, \eta, \lambda); \C)$ be a weighted Liouville surface.
Then
\begin{enumerate}
	\item The geometric Fukaya $A_\infty$-pre-category $\Wgp = \Wgp(\S; \C)$ is a small strongly triangulated $\ZZ$-graded $\k$-linear $A_\infty$-pre-category with enough homotopy equivalences.
	\item If $\LL_0 = (L_0 = (I, \gamma_0, P_0), \gr_0, \M_0, \delta_0)$ is a Lagrangian brane in $(\S; \C)$ and $L_1 = (I, \gamma_1, P_1)$ is isotopic to $L_0$, then there is a Lagrangian brane $\LL_1$ supported on $L_1$ such that $\LL_0$ and $\LL_1$ are homotopy equivalent in $\Wgp$.
	\item The strict equivalence class of $\Wgp$ is invariant under isotopy of $\eta$ or $\C$ and does not depend on the choice of a compatible Liouville form $\lambda$.
\end{enumerate}
\end{theorem}

A homotopy equivalence in an $A_\infty$-pre-category is a closed morphism of degree zero invertible up to homotopy, see \cref{d: homotopy inverse}, it is stronger than being a quasi-isomorphism in the sense of \cite{KS01}*{Definition 9}.
Having enough homotopy equivalences is an appropriate notion of unitality for $A_\infty$-pre-categories in our context.
Another option would be to use the homotopy unit framework of Fukaya--Oh--Ohta--Ono~\cite{FOOO}*{\Sec 3.3}, but it is significantly more involved.
By strongly triangulated we mean that there is a strict totalisation functor from the triangulated envelope $\Tw \Wgp$ to $\Wgp$.

In \cref{t: geometric Fukaya A_infty-pre-category is good}\,(2) two objects in an $A_\infty$-pre-category are homotopy equivalent if they can be connected by a chain of homotopy equivalences.  By `isotopic' we mean connected by a path in the space of generic $(I, \gamma, P)$ satisfying $\gamma^*\lambda + dP = 0$ with the $C^\infty$-topology.
This notion is weaker than Legendrian isotopy of the lifts. We expect that the latter also induces homotopy equivalence in $\Wgp$. We will address this elsewhere.

In \cref{t: geometric Fukaya A_infty-pre-category is good}\,(3) a strict equivalence $\F\colon \W_0 \to \W_1$ is a bijection on objects $\F^0$ that establishes a bijection on transverse sequences together with isomorphisms of graded $\k$-modules $\F^1\colon \W_0(-, -) \to \W_1(\F^0 -, \F^0 -)$ on transverse pairs such that $\F^1 \circ \bmu_0^n = \bmu_1^n \circ (\F^1)^{\otimes n}$.
An isotopy of $\eta$ is any path in $C^\infty(S, \PP TS)$.
An isotopy of $\C$ is induced by a path starting at $\id_S$ in the space of diffeomorphisms of $S$ fixing $Q$ with the $C^\infty$-topology.

\subsection*{Geometric Fukaya strictly unital \texorpdfstring{$A_\infty$}{A∞}-category}
Kontsevich and Soibelman~\cite{KS01}*{Conjecture 4} conjectured that the quasi-equivalence classes of $\k$-linear $A_\infty$-pre-categories with enough quasi-isomorphisms are in bijection with quasi-equivalence classes of strictly unital $\k$-linear $A_\infty$-categories.
This was proven when $\k$ is a field by Efimov~\cite{E11}, but remains open for arbitrary commutative rings, which is the setting we are in.
We therefore resort to a hands-on way to construct an $A_\infty$-category spanned by a collection of objects in an $A_\infty$-pre-category that works well for our geometric setup.

We sketch the construction.
Call a collection $\bL^0$ of Lagrangian branes \emph{tame} if the branes in it are supported on a no more than countable set of transverse curves.
Call a tame collection $\bL^0$ \emph{full} if any Lagrangian brane in $(\S; \C)$ is homotopy equivalent in $\Wgp$ to a Lagrangian brane from $\bL^0$.
For any tame collection $\bL^0 = \{\LL^0_\alpha\}_{\alpha \in A}$  we construct a collection
$\bL^\bullet = \{\LL^j_\alpha\}_{\alpha \in A}^{j = 0, 1, \ldots}$
of Lagrangian branes and homotopy equivalences
${e^{j_0,j_1}_\alpha \in \Wgp(\LL^{j_0}_\alpha, \LL^{j_1}_\alpha)}$
for each $j_0 \ne j_1$ that behave like strict units, see \cref{subsubs: system of strict units}.
We call such $\bL^\bullet$ a \emph{system of strict units} and associate to it a strictly unital $A_\infty$-category $\W_{\bL^\bullet}$ with $\Ob \W_{\bL^\bullet} = \bL^0$, see \cref{subsubs: A_infty-category spanned by a collection of objects}.
This system of strict units is obtained by taking small perturbations of the underlying curve of $\LL^0_\alpha$, see \cref{subsubs: systems of strict units in Fukaya category}.
The result of the construction is the following.

\begin{theorem} \label{t: geometric Fukaya strictly unital A_infty-category}
	Let $(\S; \C)$ be a weighted Liouville surface.
	Then
	\begin{enumerate}
		\item For any system of strict units $\bL^\bullet$ the strictly unital $A_\infty$-category $\W_{\bL^\bullet}$ extends the full subcategory of $\Wgp$ on $\bL^0$, i.e.\ $\Ob \W_{\bL^\bullet} = \bL^0$, $\W_{\bL^\bullet}(-, -) = \Wgp(-, -)$ on transverse pairs, and $\bmu_{\W_{\bL^\bullet}} = \bmu_{\Wgp}$ on transverse sequences.
		Moreover, two objects of $\W_{\bL^\bullet}$ are homotopy equivalent in $\W_{\bL^\bullet}$ if and only if they are homotopy equivalent in $\Wgp$.
		\item For any two systems of strict units $\bL^\bullet_0 \subset \bL^\bullet_1$ we have a natural strict full embedding ${\W_{\bL_0^\bullet} \hookrightarrow \W_{\bL_1^\bullet}}$.
		\item For any two systems of strict units $\bL^\bullet_0, \bL^\bullet_1$ there exists a strict autoequivalence $\I$ of $\Wgp$ such that $\bL^\bullet_0 \sqcup \I(\bL^\bullet_1)$ is a system of strict units.
		\item Any tame collection can be extended to a full tame collection.
		\item If $\bL^0$ is a full tame collection, then $\W_{\bL^\bullet}$ is triangulated.
	\end{enumerate}
\end{theorem}

A strict full embedding above is a strict $A_\infty$-functor which is injective on objects and an isomorphism of graded modules on morphisms.
\cref{t: geometric Fukaya strictly unital A_infty-category} implies that the homotopy equivalence class of $\W_{\bL^\bullet}$ for a system of strict units $\bL^\bullet$ on a full tame collection $\bL^0$ does not depend on these choices.
We call the resulting category the \emph{geometric Fukaya strictly unital $A_\infty$-category} $\Wgsu(\S; \C)$ of $\S$ with coefficients in $\C$.
Note that $\Wgsu(\S; \C)$ is a small strictly unital triangulated $\ZZ$-graded $\k$-linear $A_\infty$-category.

We now tie back to the topological Fukaya category.
Any arc brane can be interpreted as a Lagrangian brane and an arc brane system $\bGamma$ gives a tame collection of Lagrangian branes $\bGamma^0$.

\begin{theorem} \label{t: topological Fukaya category inside geometric}
	Let $(\S; \C)$ be a weighted Liouville surface.
	Then
	\begin{enumerate}
		\item For any arc brane system $\bGamma$ the $A_\infty$-category $\V_\bGamma$ is strictly equivalent to $\W_{\bGamma^\bullet}$ for a system of strict units $\bGamma^\bullet$ on the corresponding tame collection $\bGamma^0$.
		\item For a full tame collection $\bL^0$ the subcategory of $\W_{\bL^\bullet}$ generated by $\bGamma^0 \subset \bL^0$ corresponding to a full arc brane system $\bGamma$ does not depend on such $\bGamma$.
	\end{enumerate}
\end{theorem}

The subcategory generated by $B \subset \Ob \W$ is the smallest full subcategory of $\W$ containing $B$ and closed under taking shifts, cones, and homotopy equivalent objects.
In the proof of \cref{t: topological Fukaya category inside geometric}\,(1) we do not use that $\V_\bGamma$ satisfies the $A_\infty$-relations hence \cref{t: topological Fukaya category inside geometric} implies \cref{t: A_infty-categories of arc systems} and moreover that $\Wt$ is homotopy equivalent to a subcategory of $\Wgsu$, well-defined up to homotopy autoequivalence of $\Wgsu$.
We expect that $\Wt$ is in fact homotopy equivalent to $\Wgsu$. We will address this elsewhere.

\subsection*{Structure of the paper and other results}
\label{subs: Structure of the paper}
In \cref{s: Combinatorial definition} we define weighted surfaces and their topological Fukaya categories. In \cref{p: formal full arc system} we give a sufficient condition for $\Wt(\S; \C)$ to have a formal generator.
We describe the relation of $\Wt(\S; \C)$ to classical topological Fukaya categories in \cref{subsubs: Classical topological Fukaya category} and to relative wrapped Fukaya categories of surfaces in \cref{subsubs: Relative wrapped Fukaya category of a surface}.

In \cref{s: Mirror symmetry examples} we prove \cref{t: HMS}.

In \cref{s: Geometric preparations} we study Liouville structures on a given surface, Lagrangian structures on a given curve and deformations of such structures.
We describe the moduli space of graded Liouville marked structures on a given surface in \cref{p: graded marked surfaces up to deformation} and
prove that diffeomorphisms of marked surfaces deform to symplectomorphisms in \cref{p: Diff retracts on Symp}.
We characterise Lagrangian structures on a given curve in terms of linear inequalities on actions and areas in \cref{p: criterion for existence of exact immersion in a given component,p: Lagrangian supported on exact curve with prescribed signs}.
In \cref{p: characterisation of isotopy} we characterise isotopies of Lagrangians topologically.
We then define the relevant polygons bounded by curves and discuss the grading.

In \cref{s: Obstruction algebra} we define Lagrangian pre-branes and their obstruction algebras.
In \cref{p: isotopy induces isomorphism of obstruction algebras} we prove that the obstruction algebra does not change under isotopies of Lagrangians and another move that allows to exchange Reeb chords for double points.
In \cref{p: obstruction algebra is a nilpotent A infinity algebra} we prove that the obstruction algebra is a curved $A_\infty$-algebra.

In \cref{s: A_infty-machinery} we study $A_\infty$-pre-categories with enough homotopy equivalences.
After establishing basic facts about homotopy equivalences in \cref{l: 2-out-of-6 property for homotopy equivalences,l: homotopy equivalent transverse objects,l: exact morphisms and contractible objects in A_infty-pre-category,l: homotopy equivalence can be tested pointwise}, we explain how to get a strictly unital $A_\infty$-category from a certain $A_\infty$-pre-category in \cref{p: A_infty category spanned by a collection of objects is well-defined} and give a sufficient condition for this construction to be canonical in \cref{p: W_SU is independent of choices}.
We define the triangulated envelope of an $A_\infty$-pre-category in \cref{subsubs: twisted complexes} and prove that homotopy equivalences behave well for strongly triangulated $A_\infty$-categories in \cref{p: homotopy equivalence in triangulated A_infty-pre-category}.

In \cref{s: Geometric Fukaya category} we define the geometric Fukaya $A_\infty$-pre-category in \cref{subs: Fukaya A_infty-pre-category} and prove \cref{t: A_infty-categories of arc systems,t: geometric Fukaya A_infty-pre-category is good,t: geometric Fukaya strictly unital A_infty-category,t: topological Fukaya category inside geometric}.

\subsection*{Conventions}
We adopt the following conventions throughout.

For geometric considerations we work in the smooth category unless otherwise stated.
Manifold means (smooth) manifold with (possibly empty) boundary.
A smooth map on a subset $A$ of a manifold is the restriction of a smooth map on a neighbourhood of $A$.

Given a smooth map $\ph\colon M \to N$ the pullback on $n$-forms is
${\ph^*\colon \Omega^n(N) \to \Omega^n(M)}$,
for any other fibre bundle $E$ on $N$ the pullback on sections is
${\ph^*\colon C^\infty(N, E) \to C^\infty(M, \ph^*E)}$.
The pullback of $a$ along an embedding is denoted by $\res{a}{M}$.

A component of a topological space $T$ is its path-connected component, the set of all components is denoted by $\pi_0 T$.

For algebraic considerations we work over a unital commutative ring $\k$.
A linear map is a morphism of (left) $\k$-modules.
The tensor product over $\k$ is denoted by $\otimes$.

All (dg-, $A_\infty$-, 2-)categories except $\dgcat_\k$, $Ch(\k)$, and $\Dp X$ are small.
All 2-categories and 2-functors are strict.

Compositions of morphisms and concatenations of paths are written from right to left, e.g.\ the concatenation of a path $a$ from $p$ to $q$ and a path $b$ from $q$ to $r$ is the path $ba$ from $p$ to $r$.

\subsection*{Acknowledgements}
This work started during my PhD at the London School of Geometry and Number Theory and incorporates part of my thesis.
It is a pleasure to thank my supervisor Yank{\i} Lekili for countless discussions and generous sharing of ideas.

I am also grateful to Ilaria Di Dedda, Ezra Getzler, Jeff Hicks, Thomas Kragh, and Michael Wemyss for discussing my work and sharing their insights.

This work was supported by the Swedish Research Council Centre of Excellence grant VR 2022-06593: The Centre of Excellence in Geometry and Physics at Uppsala University and the Engineering and Physical Sciences Research Council [EP/S021590/1], The EPSRC Centre for Doctoral Training in Geometry and Number Theory (The London School of Geometry and Number Theory), University College London.

\subsection*{Declaration on the use of AI}
AI tools (Anthropic's Claude and OpenAI's ChatGPT) were used while preparing this manuscript: for checking proofs and computations, detecting errors and inconsistencies, suggesting corrections to formulas and arguments of a routine nature, literature search, and implementing notational and editorial changes.
All main ideas, constructions and results are the author's.
The author has verified all AI-assisted content and takes full responsibility for the paper.
The wording of this declaration was suggested by AI.

\addtocontents{toc}{\protect\setcounter{tocdepth}{2}}
\section{Topological Fukaya category}
\label{s: Combinatorial definition}
In \cref{subs: Weighted system of dg-categories} we define weighted systems of dg-categories on a marked surface.
In \cref{subs: Fukaya category from arc systems} we define the topological Fukaya category of a weighted surface, give a sufficient condition for existence of a formal generator in it, and explain how to specialise it to earlier constructions.

\subsection{Weighted system of dg-categories}
\label{subs: Weighted system of dg-categories}
Let $S$ be an oriented smooth surface and $Q \subset S$ a finite set, called the \emph{marking}.
In \cref{subsubs: Directed relative path 2-category} we define the \emph{directed relative path 2-category} $\P_2^+(S, Q)$  by modifying the definition of the path 2-groupoid of a smooth manifold by Schreiber--Waldorf~\cite{SW11}*{\Sec 2--3}.
In \cref{subsubs: Weighted system of dg-categories on a graded marked surface} we define a weighted system of dg-categories on a marked surface $(S, Q)$ as a 2-functor $\C \colon \P_2^+(S, Q) \to \dgcat_\k$ for a commutative ring $\k$.
In \cref{subsubs: Composition along a 2-morphism} we define the composition along a 2-morphism in $\P_2^+(S, Q)$ with respect to $\C$, which is the key construction to define all $A_\infty$-operations in the paper.
In \cref{subsubs: Fundamental and homological 2-category} we describe two more manageable quotients of $\P_2^+(S, Q)$, the directed relative fundamental and homological 2-categories $\Pi_2^+(S, Q)$ and $\H_2^+(S, Q)$.

\subsubsection{Directed relative path 2-category}
\label{subsubs: Directed relative path 2-category}
First recall general notions on lazy paths, see Schreiber--Waldorf~\cite{SW11} for details.
Let $M$ be a smooth manifold.
A \emph{lazy 1-path} (a.k.a.\ a path with sitting instants) from $p_0 \in M$ to $p_1 \in M$ is a smooth map $\gamma\colon [0, 1] \to M$ such that
\begin{enumerate}
	\item $\gamma(s) = p_0$ for $s$ in a neighbourhood of $0$ and
	\item $\gamma(s) = p_1$ for $s$ in a neighbourhood of $1$.
\end{enumerate}
We use notation $p_0 \xrightarrow{\gamma} p_1$ to denote such a lazy 1-path.

A \emph{lazy 2-path} from a lazy 1-path $p_0 \xrightarrow{\gamma_0} p_1$ to a lazy 1-path $p_0 \xrightarrow{\gamma_1} p_1$ is a smooth map $\sigma\colon [0, 1] \times [-1, 0] \to M$ such that
\begin{enumerate}
	\item $\sigma(s, t) = p_i$ for fixed $i = 0, 1$ and $(s, t)$ in a neighbourhood of $\{i\} \times [-1, 0]$ and
	\item $\sigma(s, t) = \gamma_i(s)$ for fixed $i = 0, 1$ and $(s, t)$ in a neighbourhood of $[0, 1] \times \{-i\}$.
\end{enumerate}
We use notation $\gamma_0 \xRightarrow{\sigma} \gamma_1$ to denote such a lazy 2-path.
\begin{remark}
	The choice of the domain for a lazy 2-path is made to satisfy the following conditions.
	\begin{enumerate}
		\item The first coordinate is horizontal and the second coordinate is vertical (so the horizontal and vertical compositions are the respective stackings of the domains).
		\item The orientation of the domain $[0, 1] \times \{-1\}$ of the \emph{target} 1-path coincides with its orientation as boundary of $[0, 1] \times [-1, 0]$.
	\end{enumerate}
\end{remark}

A lazy 2-path $\sigma$ is called a \emph{thin homotopy} if $\rk d\sigma \le 1$ everywhere.
Observe that thin homotopy is an equivalence relation on lazy 1-paths.

Specialising to $Q \subset S$ as above, we say that a lazy 2-path $\gamma_0 \xRightarrow{\sigma} \gamma_1$ in $S$ is \emph{positive} if
\begin{enumerate}
	\item $\gamma_0$ and $\gamma_1$ are lazy 1-paths in $S \setminus Q$ and
	\item
	$\left(\frac{\pd \sigma}{\pd s}, \frac{\pd \sigma}{\pd t}\right)$
	is a positive basis of $T_{\sigma(s, t)} S$ for each $(s, t) \in [0, 1] \times [-1, 0]$ with $\sigma(s, t) \in Q$.
\end{enumerate}

Note that a positive $2$-path cannot have $Q \cap \pd S$ in its image.
We therefore divide $Q$ into two parts, the \emph{divisor} $D = Q \cap \Int S$ and the \emph{stop} $\Lambda = Q \cap \pd S$.

An \emph{isotopy} from a positive 2-path $\gamma_0^0 \xRightarrow{\sigma^0} \gamma_1^0$ to a positive 2-path $\gamma_0^1 \xRightarrow{\sigma^1} \gamma_1^1$ is a smooth map $h\colon [0, 1] \times [-1, 0] \times [0, 1] \to S$ such that
\begin{enumerate}
	\item $\res{h}{[0, 1] \times \{-i\} \times [0, 1]}$ is a thin homotopy from $\gamma_i^0$ to $\gamma_i^1$ in $S \setminus Q$ for $i = 0, 1$;
	\item $h(s, t, r) = \sigma^i(s, t)$ for fixed $i = 0, 1$ and $(s, t, r)$ in a neighbourhood of $[0, 1] \times [-1, 0] \times \{i\}$
	\item $\res{h}{[0, 1] \times [-1, 0] \times \{r\}}$ is a positive 2-path for any $r \in [0, 1]$.
\end{enumerate}
Observe that isotopy induces an equivalence relation on positive 2-paths.
Note that there is no notion of thinness for isotopies since $\dim S = 2$.

\begin{defn} \label{d: directed relative path 2-category}
Let $S$ be an oriented surface and $Q \subset S$ a finite set.
The \emph{directed relative path 2-category} $\P_2^+(S, Q)$ is the following (strict) 2-category.
\begin{enumerate}
	\item $\Ob \P_2^+(S, Q) = S \setminus Q$.
	\item For $p_0, p_1 \in S \setminus Q$ the set of $1$-morphisms is
	$$
		\P_2^+(S, Q)(p_0, p_1) = \frac{\left\{\text{lazy 1-paths from $p_0$ to $p_1$ in $S \setminus Q$}\right\}}{\text{thin homotopy}}.
	$$
	\item For $[\gamma_0], [\gamma_1] \in \P_2^+(S, Q)(p_0, p_1)$ the set of 2-morphisms is
	$$
		\P_2^+(S, Q)(p_0, p_1)([\gamma_0], [\gamma_1]) = \frac{\left\{\text{positive 2-paths $\gamma_0' \xRightarrow{\sigma} \gamma_1'$ in $S$ with $\gamma_i' \in [\gamma_i]$ for $i = 0, 1$}\right\}}{\text{isotopy}}.
	$$
\end{enumerate}
The composition $\circ$ of $1$-morphisms is induced by the concatenation of paths; the vertical and horizontal compositions (we denote both by $\cup$ if it does not cause confusion) of $2$-morphisms are induced by vertical and horizontal stackings of the domains $[0, 1]^2$ respectively.
\end{defn}

We refer the reader to Kelly--Street~\cite{KS74} for the 2-categorical notation used.
That \cref{d: directed relative path 2-category} indeed gives a 2-category is a straightforward check analogous to the case of the path 2-groupoid, see \cite{SW11}*{\Sec 3.1} for details of the latter.

Observe that any piece-wise smooth path $\gamma$ in $S \setminus Q$ defines a 1-morphism $[\gamma]$ in $\P_2^+(S, Q)$.
Moreover, if we have a smooth map $u \colon \Pi \to M$ where $\Pi$ is a compact region in $\RR^2$ bounded by a piecewise-smooth Jordan curve, then for any two points $v_0, v_1 \in \pd \Pi$ we get a 2-morphism $[u]\colon [\gamma^0] \to [\gamma^1]$ where $\gamma^0$ (resp.\ $\gamma^1$) is the restriction of $u$ to the part of $\pd \Pi$ from $v_0$ to $v_1$ that has $\Pi$ on the right (resp.\ left).

Note that any 1-morphism in $\P_2^+(S, Q)$ is invertible with the inverse for a lazy 1-path $\gamma$ given by $\gamma^{-1}(s) = \gamma(1 - s)$.
Also observe that the \emph{degree function}
$$
\deg_\sigma \colon D \to \ZZ_{\ge 0}, \, \deg_\sigma(q) = |\sigma^{-1}(q)|
$$
of a positive 2-path $\gamma_0 \xRightarrow{\sigma} \gamma_1$ is constant in isotopy classes and additive with respect to vertical and horizontal composition.
In particular, if $\deg_\sigma \ne 0$, then $\sigma$ is not invertible.
On the other hand, if $\deg_\sigma = 0$, then $\sigma(s, -1 - t)$ is the inverse of $\sigma(s, t)$ (with respect to vertical composition).

\subsubsection{Weighted system of dg-categories on a graded marked surface}
\label{subsubs: Weighted system of dg-categories on a graded marked surface}
We define the surfaces we consider in this paper.
\begin{defn} \label{d: graded marked surface}
	A \emph{graded marked surface} is the data of
	\begin{enumerate}
		\item A compact smooth oriented surface $S$ such that each component of $S$ has non-empty boundary.
		\item A finite set $Q \subset S$ called the \emph{marking}.
		\item A line field $\eta \in C^{\infty}(S, \PP TS)$ called the \emph{grading}.
	\end{enumerate}
\end{defn}

Let $\S = (S, Q, \eta)$ be a graded marked surface.
We associate to it the 2-category $\P^+_2(\S)$ from \cref{d: directed relative path 2-category}.
We recall the na{\"i}ve 2-category structure on dg-categories.

Fix a commutative ring $\k$ and recall that the category $Ch(\k)$ of dg-$\k$-modules is symmetric monoidal.
A \emph{dg-category} $\C$ \emph{(over $\k$)} is a category enriched over $Ch(\k)$, a \emph{dg-functor} is a functor enriched over $Ch(\k)$, see Keller~\cite{K06}*{\Sec 2} for details.
$Z^0(\C)$ and $H^0(\C)$ are the categories with the same objects as $\C$ and morphisms given by
$$
Z^0(\C)(M, N) = Z^0(\C(M, N)) \text{ and } H^0(\C)(M, N) = H^0(\C(M, N))
$$
where $Z^0 \text{ and } H^0\colon Ch(\k) \to \k\text{-Mod}$ are the $0$-cycle and $0$-cohomology functors respectively.
We will only consider small dg-categories.
They form a closed symmetric monoidal category with the internal hom-functor $\Fun_{dg}$.
By \emph{morphisms} of dg-functors from $\C$ to $\D$ we mean that in the category $Z^0(\Fun_{dg}(\C, \D))$.
A dg-functor $F\colon \C \to \D$ is a \emph{(strict) dg-equivalence} if there exists a dg-functor $G\colon \D \to \C$ such that $F \circ G$ is isomorphic to $\id_{\D}$ (resp.\ $F \circ G = \id_{\D}$) and $G \circ F$ is isomorphic to $\id_{\C}$ (resp.\ $G \circ F = \id_{\C}$).
We denote by $\dgcat_\k$ the 2-category of dg-categories where an object is a small dg-category over $\k$, a 1-morphism is a dg-functor, and a 2-morphism is a morphism of dg-functors.

\begin{defn} \label{d: weighted system of dg-categories}
	A \emph{weighted system of dg-categories} on the graded marked surface $\S$ is the data of
	\begin{enumerate}
		\item a commutative ring $\k$;
		\item a strict 2-functor $\C\colon \P_2^+(\S) \to \dgcat_{\k}$.
	\end{enumerate}
	We call the pair $(\S; \C)$ a \emph{weighted surface}.

	For $p \in \Ob \P_2^+(\S) = S \setminus Q$ we denote the corresponding dg-category, called the \emph{fibre} at $p$, by $\C_p$.
	For a 1-morphism in $\P_2^+(\S)$ represented by a lazy 1-path $\gamma$ we denote the corresponding dg-functor, called the \emph{transition functor} along $\gamma$, by $\C(\gamma)$.
	For a 2-morphism in $\P_2^+(\S)$ represented by a positive 2-path $\sigma$ we denote the corresponding morphism of dg-functors, called the \emph{weight} of $\sigma$, by $\C(\sigma)$.
\end{defn}

Note that all transition functors in a weighted system of dg-categories are strict dg-equivalences and all weights of positive 2-paths that do not intersect $Q$ are isomorphisms of dg-functors.

\begin{remark}
	As Fukaya categories are $A_\infty$-categories, it might seem more natural to consider $A_\infty$-category coefficients instead.
	We do not do this for three reasons.
	First, similarly to the case of tensor product, the formulae in the case of $A_\infty$-categories are complicated and unenlightening, see Ornaghi~\cite{O25} for an overview of the tensor product constructions.
	Secondly, to obtain an object of the Fukaya category, one chooses a local system of modules on a Lagrangian, and modules over an $A_\infty$-category form a dg-category, so we can simply consider that dg-category and local system of objects over it.
	Thirdly, $A_\infty$-categories do not form a 2-category at the chain level since the vertical composition of $A_\infty$-natural transformations is not strictly associative.
\end{remark}

\subsubsection{Composition along a 2-morphism}
\label{subsubs: Composition along a 2-morphism}
Fix a weighted system of dg-categories
${\C\colon \P^+_2(S, Q) \to \dgcat_\k}$.
Let $p^0, \ldots, p^n \in S \setminus Q$, $\gamma^0, \ldots, \gamma^n$ be lazy 1-paths in $S \setminus Q$ such that $\gamma^0$ goes from $p^0$ to $p^n$ and $\gamma^i$ goes from $p^{i - 1}$ to $p^i$ for $i = 1, \ldots, n$, and $\M(p^i)$ be an object in $\C_{p^i}$ for $i = 0, \ldots, n$.
Let $\gamma$ be a lazy parametrisation of $\gamma^n\ldots\gamma^1$ and consider a positive 2-path $\sigma\colon \gamma^0 \Rightarrow \gamma$.
The \emph{composition along} $\sigma$ is the morphism
$$
\xymatrix{
\C_{p^n}(\C(\gamma^n)\M(p^{n - 1}), \M(p^n)) \otimes \ldots \otimes \C_{p^1}(\C(\gamma^1)\M(p^0), \M(p^1)) \ar[d]_{\C(\gamma^0; \gamma^n, \ldots, \gamma^1; \sigma)}\\
\C_{p^n}(\C(\gamma^0)\M(p^0), \M(p^n))
}
$$
given by
\begin{equation} \label{eq: composition along 2-morphism}
	\C(\gamma^0; \gamma^n, \ldots, \gamma^1; \sigma)(\ph^n, \ldots, \ph^1) = \ph^{\gamma^n, \ldots, \gamma^1} \cdot \C(\sigma)_{\M(p^0)}
\end{equation}
where $\ph^{\gamma^n, \ldots, \gamma^1} = \ph^{\gamma^n, \ldots, \gamma^1}(\ph^n, \ldots, \ph^1)$ is the composition
\begin{equation} \label{eq: phi^gamma^n...gamma^1}
	\C(\gamma^n\ldots\gamma^1)\M(p^0) \xrightarrow{\C(\gamma^n\ldots\gamma^2)\ph^1} \C(\gamma^n\ldots\gamma^2)\M(p^1) \xrightarrow{\C(\gamma^n\ldots\gamma^3)\ph^2} \ldots \xrightarrow{\ph^n} \M(p^n)
\end{equation}
in $\C_{p^n}$.

We list elementary properties of composition along a positive 2-path for future reference.

\begin{lemma}~ \label{l: properties of composition along a 2-path}
	Let $\gamma^0, \ldots, \gamma^n$ and $\sigma$ be as above. Then
	\begin{enumerate}
		\item $\C(\gamma^0; \gamma^n, \ldots, \gamma^1; \sigma)$ only depends on the thin homotopy classes of $\gamma^0, \ldots, \gamma^n$ and the isotopy class of $\sigma$.
		\item $\C(\gamma^0; \gamma^0; \id) = \id$.
		\item If $\C(\gamma^i; \gamma^i_{m^i}, \ldots, \gamma^i_1; \sigma^i)$ is defined for $i = 1, \ldots, n$, then $\C(\gamma^0; \gamma^n_{m^n}, \ldots, \gamma^n_1, \ldots, \gamma^1_{m^1}, \ldots, \gamma^1_1; \sigma \cup \sigma^n \cup \ldots \cup \sigma^1)$ is defined and equals $\C(\gamma^0; \gamma^n, \ldots, \gamma^1; \sigma) \circ (\C(\gamma^n; \gamma^n_{m^n}, \ldots, \gamma^n_1; \sigma^n), \ldots, \C(\gamma^1; \gamma^1_{m^1}, \ldots, \gamma^1_1; \sigma^1))$.
	\end{enumerate}
\end{lemma}
\begin{proof}
	(1) and (2) are immediate from construction.
	(3) is the computation
	$$
	\begin{aligned}
		&\C(\gamma^0; \gamma^n, \ldots, \gamma^1; \sigma) \circ (\C(\gamma^n; \gamma^n_{m^n}, \ldots, \gamma^n_1; \sigma^n), \ldots, \C(\gamma^1; \gamma^1_{m^1}, \ldots, \gamma^1_1; \sigma^1))(\ph^n_{m^n}, \ldots, \ph^1_1) =\\
		&\C(\gamma^0; \gamma^n, \ldots, \gamma^1; \sigma)\left(\ph^{\gamma^i_{m^i}, \ldots, \gamma^i_1}(\ph^i_{m^i}, \ldots, \ph^i_1) \cdot \C(\sigma^i)_{\M(p^{i - 1})}\right)_{i = 1}^n =\\
		&\left(\prod_{i = 1}^n \C(\gamma^n\ldots\gamma^{i + 1}) \left(\ph^{\gamma^i_{m^i}, \ldots, \gamma^i_1} \cdot \C(\sigma^i)_{\M(p^{i - 1})}\right)\right) \cdot \C(\sigma)_{\M(p^0)} =\\
		&\left(\prod_{i = 1}^n \C(\gamma^n\ldots\gamma^{i + 1}) \left(\ph^{\gamma^i_{m^i}, \ldots, \gamma^i_1}\right) \cdot \C(\id_{\gamma^n\ldots\gamma^{i + 1}} \cup \sigma^i)_{\M(p^{i - 1})}\right) \cdot \C(\sigma)_{\M(p^0)} =^*\\
		&\left(\prod_{i = 1}^n \C(\gamma^n_{m^n}\ldots\gamma^{i + 1}_1) \left(\ph^{\gamma^i_{m^i}, \ldots, \gamma^i_1}\right)\right) \cdot
		\left(\prod_{i = 1}^n \C(\id_{\gamma^n\ldots\gamma^{i + 1}} \cup \sigma^i)_{\C(\gamma^{i - 1}_{m^{i - 1}}\ldots\gamma^1_1)\M(p^0)}\right) \cdot \C(\sigma)_{\M(p^0)} = \\
		&\ph^{\gamma^n_{m^n}, \ldots, \gamma^1_1}(\ph^n_{m^n}, \ldots, \ph^1_1) \cdot
		\left(\prod_{i = 1}^n \C(\id_{\gamma^n\ldots\gamma^{i + 1}} \cup \sigma^i \cup \id_{\gamma^{i - 1}_{m^{i - 1}}\ldots\gamma^1_1})_{\M(p^0)}\right) \cdot \C(\sigma)_{\M(p^0)} =\\
		&\C(\gamma^0; \gamma^n_{m^n}, \ldots, \gamma^1_1; \sigma \cup \sigma^n \cup \ldots \cup \sigma^1)(\ph^n_{m^n}, \ldots, \ph^1_1)
	\end{aligned}
	$$
	where $\prod_{i = 1}^n x_i = x_n\ldots x_1$ and for $=^*$ we repeatedly use that
	$$
	\C(\tau)_{M_1} \cdot \C(\delta_0)\ph = \C(\delta_1)\ph \cdot \C(\tau)_{M_0} \in \C_{q_1}(\C(\delta_0)M_0, \C(\delta_1)M_1)
	$$
	for any $q_0, q_1 \in S \setminus Q$, lazy 1-paths $\delta_0, \delta_1$ from $q_0$ to $q_1$, positive 2-path $\tau\colon \delta_0 \Rightarrow \delta_1$, $M_0, M_1 \in \C_{q_0}$, and $\ph \in \C_{q_0}(M_0, M_1)$
\end{proof}

\subsubsection{Fundamental and homological 2-category}
\label{subsubs: Fundamental and homological 2-category}
The equivalence relation on paths in \cref{d: directed relative path 2-category} is the finest that allows us to define a Fukaya category.
There are two coarser and more manageable 2-categories which one can associate to $Q \subset S$ as above.

\begin{defn} \label{d: directed relative fundamental 2-category}
Let $S$ be an oriented surface and $Q \subset S$ a finite set.
The \emph{directed relative fundamental 2-category} $\Pi_2^+(S, Q)$ is the following 2-category.
\begin{enumerate}
	\item $\Ob \Pi_2^+(S, Q) = S \setminus Q$.
	\item For $p_0, p_1 \in S \setminus Q$ the set of $1$-morphisms is
	$$
		\Pi_2^+(S, Q)(p_0, p_1) = \frac{\left\{\text{lazy 1-paths from $p_0$ to $p_1$ in $S \setminus Q$}\right\}}{\text{homotopy rel endpoints in } S \setminus Q}.
	$$
	\item For $[\gamma_0], [\gamma_1] \in \Pi_2^+(S, Q)(p_0, p_1)$ the set of 2-morphisms is
	$$
		\Pi_2^+(S, Q)(p_0, p_1)([\gamma_0], [\gamma_1]) = \frac{\left\{\text{positive 2-paths $\gamma_0' \xRightarrow{\sigma} \gamma_1'$ in $S$ with $\gamma_i' \in [\gamma_i]$ for $i = 0, 1$}\right\}}{\text{isotopy not necessarily thin on $[0, 1] \times \{-i\} \times [0, 1]$}}.
	$$
\end{enumerate}
\end{defn}

Note that the 1-category truncation of $\Pi_2^+(S, Q)$ is the fundamental groupoid $\Pi_1(S \setminus Q)$.
Also observe that the identity on objects, lazy 1-paths, and positive 2-paths induce a 2-functor $\P_2^+(S, Q) \to \Pi_2^+(S, Q)$.

\begin{defn} \label{d: directed relative homology 2-category}
Let $S$ be an oriented surface and $Q \subset S$ a finite set.
The \emph{directed relative homology 2-category} $\H_2^+(S, Q)$ is the following 2-category.
\begin{enumerate}
	\item $\Ob \H_2^+(S, Q) = S \setminus Q$.
	\item For $p_0, p_1 \in S \setminus Q$ the set of $1$-morphisms is
	$$
		\H_2^+(S, Q)(p_0, p_1) = \frac{\left\{x \in C_1(S \setminus Q) \colon \pd x = p_1 - p_0\right\}}{\pd C_2(S \setminus Q)}
	$$
	where $(C_*(S \setminus Q), \partial)$ is the chain complex of singular smooth simplices in $S \setminus Q$ with coefficients in $\ZZ$.
	\item For $[x_0], [x_1] \in \H_2^+(S, Q)(p_0, p_1)$ the set of 2-morphisms is
	$$
		\H_2^+(S, Q)(p_0, p_1)([x_0], [x_1]) = \left\{v \in \left(\ZZ_{\ge 0}\right)^{D} \colon x_1 - x_0 - \sum_{q \in D} v(q) \gamma_q \in \pd C_2(S \setminus Q)\right\}
	$$
	where $\gamma_q$ is the boundary of a small disc in $S$ around $q \in D$.
\end{enumerate}
The composition of $1$-morphisms is induced by addition in $C_1(S \setminus Q)$ and the vertical and horizontal compositions are given by addition in $\left(\ZZ_{\ge 0}\right)^{D}$.
\end{defn}

Note that the set of 2-morphisms between two fixed 1-morphisms has at most one element since $\{[\gamma_q]\}_{q \in D}$ is a linearly independent set in $H_1(S \setminus Q)$. Also observe that the identity on $S \setminus Q$, sending a 1-path to the corresponding 1-simplex, and sending a positive 2-path $\sigma$ to its degree $\deg_\sigma$ induce a 2-functor $\Pi_2^+(S, Q) \to \H_2^+(S, Q)$.

\subsection{Fukaya category from arc systems}
\label{subs: Fukaya category from arc systems}
Fix a weighted surface $(\S = (S, Q, \eta); \C)$.

In \cref{subsubs: Arc systems and boundary paths} we define arc brane systems, and degrees of boundary paths.
In \cref{subsubs: A_infty-category associated to an arc system} we construct an $A_\infty$-category from an arc brane system.
In \cref{subsubs: Basic properties} we define topological Fukaya categories and give a sufficient condition for existence of a formal generator.
In \cref{subsubs: Classical topological Fukaya category,subsubs: Relative wrapped Fukaya category of a surface} we recall the classical topological Fukaya category and the relative wrapped Fukaya category of a surface respectively, and discuss their relation to our topological Fukaya category.

\subsubsection{Arc systems and boundary paths}
\label{subsubs: Arc systems and boundary paths}
An \emph{arc} in $\S$ is an embedded curve $\gamma\colon [0, 1] \to S \setminus Q$ such that $\gamma\{0, 1\} \subset \pd S$, $\gamma(0, 1) \subset \Int S$, and $\gamma$ is transverse to $\pd S$.
An \emph{arc system} is a no more than countable set of arcs such that any two arcs from the set are  disjoint.
An arc system is called \emph{full} if its arcs cut $S$ into topological discs and each such disc has at most one point from $Q$.
Note that full arc systems exist for any marked surface $(S, Q)$ since we can first cut $S$ into a single disc by an arc system, and then for each $q \in Q$ add an arc that cuts off $q$ from the disc.

Let $\gamma_0$ and $\gamma_1$ be two (not necessarily distinct) arcs from an arc system $\Gamma$.
A \emph{boundary path} from $\gamma_0$ to $\gamma_1$ is a  homotopy class of lazy 1-paths $y\colon [0, 1] \to \pd S \setminus \Lambda$ such that $y(0) = \gamma_0(y_+)$, $y(1) = \gamma_1(y_-)$ for $y_+, y_- \in \{0, 1\}$, and $dy/dt$ is always non-negative with respect to the natural orientation of $\pd S$.

A \emph{grading} on an arc $\gamma$ is a homotopy class of paths from $\gamma^*\eta$ to the tangent line $\dot\gamma$, i.e.\ an element
$\gr \in \Pi_1C^{\infty}([0, 1], \gamma^*\PP T S)(\gamma^*\eta, \dot\gamma)$.
If the arcs $\gamma_0$ and $\gamma_1$ are graded by $\gr_0$ and $\gr_1$ respectively, then we define the \emph{degree} $|y| \in \ZZ$ of a boundary path $y$ from $\gamma_0$ to $\gamma_1$ as follows.
First, obtain a path from $\eta(y(1))$ to $T_{y(1)}\pd S$ in $\PP(T_{y(1)}S)$ as a concatenation of $\gr_1(y_-)$ and the smallest counter-clockwise rotation from $\dot \gamma(y_-)$ to $T_{y(1)}\pd S$.
Secondly, transport it along $y$ to a path from $\eta(y(0))$ to $T_{y(0)}\pd S$ in $\PP(T_{y(0)}S)$.
Thirdly, compose it with the smallest clockwise rotation from $T_{y(0)}\pd S$ to $\dot \gamma(y_+)$ and then with the inverse of $\gr_0(y_+)$ to obtain an element in $\pi_1(\PP T_{y(0)}, \eta(y(0)))$.
Finally, identify $\pi_1(\PP T_{y(0)}, \eta(y(0)))$ with $\ZZ$ by sending the counter-clockwise rotation by $\pi$ to $1$ to obtain the degree of $y$.
Note that the degree is additive with respect to concatenation of paths.

An \emph{arc brane} is a graded arc $(\gamma, \gr)$ together with an object $M^0$ of $\C_{\gamma(0)}$.
Note that the parallel transport along $\gamma$ defines an object $M^t$ of $\C_{\gamma(t)}$ for any $t \in [0, 1]$.
An \emph{arc brane system} is a set of arc branes $\bGamma$ such that
$$
\Gamma_\bGamma = \{\gamma \mid (\gamma, \gr, M) \in \bGamma\}
$$
is an arc system.
For $\gamma \in \Gamma_\bGamma$ set
$$
\M_\gamma = \{M \mid (\gamma, \gr, M) \in \bGamma\} \subset \Ob \C_{\gamma(0)}.
$$

We say that $\M \subset \Ob \C$ \emph{generates} the dg-category (or $A_\infty$-category) $\C$ if the smallest full subcategory of the triangulated envelope $\Tw \C$ containing $\M$ and closed under taking shifts, cones of closed degree zero morphisms, and homotopy equivalent objects, contains $\C$.
An arc brane system $\bGamma$ is called \emph{full} if $\Gamma_\bGamma$ is a full arc system and for any $\gamma \in \Gamma_\bGamma$ the set $\M_\gamma$ generates $\C_{\gamma(0)}$.
Our conventions for the triangulated envelope $\Tw$ are given in \cref{subsubs: twisted complexes}.
To define $\Tw$ for dg-categories we view them as $A_\infty$-categories via $\mu^1(a) = (-1)^{|a|}da$ and $\mu^2(a, b) = (-1)^{|b|}a \cdot b$.
This coincides with the usual definition of twisted complexes over a dg-category (see e.g.\ Drinfeld~\cite{D04}*{\Sec 2.4}) when the shift by $m$ is realised as $\k[m] \otimes -$.

Note that for any full arc system $\Gamma$ we can take $\M_\gamma = \Ob \C_{\gamma(0)}$ for each $\gamma \in \Gamma$ and obtain a full arc brane system.

\subsubsection{\texorpdfstring{$A_\infty$}{A∞}-category associated to an arc system}
\label{subsubs: A_infty-category associated to an arc system}
We generalise the construction of the $A_\infty$-category of an arc system from \cite{HKK17}*{\Sec3.3} to our setting.

Given an arc brane system $\bGamma$ in $(\S; \C)$ we define the strictly unital $\k$-linear $A_\infty$-category $\V_\bGamma$ as follows.

\paragraph{Objects} $\Ob \V_\bGamma = \bGamma$.

\paragraph{Morphisms}
Fix two arc branes $\LL_i = (\gamma_i, \gr_i, M_i) \in \bGamma, \; i = 0, 1$.

For a boundary path $y$ from $\gamma_0$ to $\gamma_1$ define a graded $\k$-module
$$
\V_{\bGamma}(y) = \C_{y(1)}\left(\C(y)(M_0^{y_+}), M_1^{y_-}\right).
$$
The morphism space between arc branes is the graded $\k$-module
\begin{equation} \label{eq: A_Gamma(L_0 L_1)}
	\V_{\bGamma}\left(\LL_0, \LL_1\right) = \left. \bigoplus_{y} \V_{\bGamma}(y)[-|y|] \right/ \sim
\end{equation}
where $y$ runs through boundary paths from $\gamma_0$ to $\gamma_1$ up to reparametrisation, and $\sim$ identifies the summands for the two constant paths via the transition functor if $\gamma_0 = \gamma_1$.

\paragraph{$A_\infty$-operations}
There are three sources of $A_\infty$-operations: operations in the fibres of $\C$ (i.e.\ the differential $d$ and the composition $\cdot$), concatenation of boundary paths, and polygons bounded by the arcs.

If $\gamma_1, \ldots, \gamma_n$ are arcs from $\Gamma_\bGamma$ and $y^i, \, i = 1, \ldots, n$ are boundary paths from $\gamma_{i - 1}$ to $\gamma_i$ (the index is mod $n$) such that $y^i_+ \ne y^{i - 1}_-$, then the concatenation $\gamma_ny^n\gamma_{n - 1}y^{n - 1}\ldots\gamma_1y^1$ (whenever we write $\gamma_i$ in a concatenation we tacitly assume that it can be travelled in either direction) is a closed curve.
This closed curve can bound at most one immersed disc, and if it does, then $(y^n, \ldots, y^1)$ is called a \emph{disc sequence}.

For $n \ge 1$ and
$\LL_0 = (\gamma_0, \gr_0, M_0), \ldots, \LL_n = (\gamma_n, \gr_n, M_n) \in \bGamma$
the $A_\infty$-operation
$$
\bmu^n\colon \V_\bGamma(\LL_{n - 1}, \LL_n) \otimes \ldots \otimes \V_\bGamma(\LL_0, \LL_1) \to \V_\bGamma(\LL_0, \LL_n)
$$
is determined by its entries
$$
\bmu^{y^0; y^n, \ldots, y^1}\colon \V_\bGamma(y^n) \otimes \ldots \otimes \V_\bGamma(y^1) \to \V_\bGamma(y^0)
$$
where $y^i$ is a boundary path from $\gamma_{i - 1}$ to $\gamma_i$ for $i = 1, \ldots, n$, and $y^0$ is a boundary path from $\gamma_0$ to $\gamma_n$.

For
$\ph^i \in \C_{y^i(1)}\left(\C(y^i)(M_{i - 1}^{y^i_+}), M_i^{y^i_-}\right), \, i = 1, \ldots, n$
the element
$$
\ph^0 = \bmu^{y^0; y^n, \ldots, y^1}(\ph^n, \ldots, \ph^1) \in \C_{y^0(1)}\left(\C(y^0)(M_0^{y^0_+}), M_n^{y^0_-}\right)
$$
is possibly non-zero in the following cases.
\begin{enumerate}
    \item $n = 1$ and $y^0 = y^1$.
    Then
	\begin{equation} \label{eq: mu for y^0 = y^1}
			\ph^0 = (-1)^{|\ph^1| + |y^0|}d\ph^1.
	\end{equation}

    \item $n = 2$ and $y^0 = y^2y^1$.
    If $y^1_- = 1$ and $y^0_- = 0$, then
	\begin{equation} \label{eq: mu for y^0 = y^2y^1}
		\ph^0 = (-1)^{|\ph^1|(|y^2| - 1)}\ph^2 \cdot \C(y^2)(\ph^1).
	\end{equation}

    \item There exists a boundary path $z$ from $\gamma_{n - 1}$ to $\gamma_0$ such that $y^n = y^0z$ and $(z, y^{n - 1}, \ldots, y^1)$ is a disc sequence bounding an immersed polygon $\bu$.
    If $\gamma_1, \ldots, \gamma_{n - 1}$ are oriented as the boundary of $\bu$ and $y^0_- = 0$, then
	\begin{equation} \label{eq: mu for a polygon with y^n = y^0z}
	\begin{aligned}
		\ph^0 = (-1)^{\sum_{i < j} |\ph^i|(|y^j| - 1)}
		\C\big(y^0; &y^n, &\gamma_{n-1}, &y^{n-1}, &\ldots, &\gamma_1, &y^1, &\gamma_0; \sigma_{\bu}\big)\\
		\big(&\ph^n,& \id, &\ph^{n-1}, &\ldots, &\id, &\ph^1, &\id\big)
	\end{aligned}
	\end{equation}
	where the 2-morphism
	$\sigma_{\bu}\colon [y^0] \Rightarrow [y^n\gamma_{n-1}y^{n-1}\ldots\gamma_1 y^1\gamma_0]$
	is obtained from $\bu$ (which can be viewed as a 2-morphism
	$[y^0(0)] \Rightarrow [z\gamma_{n-1}y^{n-1}\ldots\gamma_1 y^1\gamma_0]$)
	by post-composing with $[y^0]$.

    \item There exists a boundary path $z$ from $\gamma_{n}$ to $\gamma_1$ such that $y^1 = zy^0$ and $(z, y^n, \ldots, y^2)$ is a disc sequence bounding an immersed polygon $\bu$.
    If $\gamma_1, \ldots, \gamma_n$ are oriented as the boundary of $\bu$, then
	\begin{equation} \label{eq: mu for a polygon with y^1 = zy^0}
	\begin{aligned}
		\ph^0 = (-1)^{\sum_{i < j} |\ph^i|(|y^j| - 1)}
		\C\big(y^0; &\gamma_n, &y^n, &\gamma_{n-1}, &y^{n-1}, &\ldots, &\gamma_1, &y^1; \sigma_{\bu}\big)\\
		\big(&\id, &\ph^n, &\id, &\ph^{n-1}, &\ldots, &\id, &\ph^1\big)
	\end{aligned}
	\end{equation}
	where the 2-morphism
	$\sigma_{\bu}\colon [y^0] \Rightarrow [\gamma_ny^n\gamma_{n-1}y^{n-1}\ldots\gamma_1 y^1]$
	is obtained from $\bu$ (which can be viewed as a 2-morphism
	$[y^0(1)] \Rightarrow [\gamma_ny^n\gamma_{n-1}y^{n-1}\ldots\gamma_1 z]$)
	by pre-composing with $[y^0]$.

	\item In the cases (2)--(4) above the change of orientation of $\gamma_i$ for $i = 1, \ldots, n - 1$ (resp.\ $\gamma_n$) multiplies the right-hand side in \cref{eq: mu for y^0 = y^2y^1,eq: mu for a polygon with y^n = y^0z,eq: mu for a polygon with y^1 = zy^0} by $(-1)^{|y^i|}$ (resp.\ $(-1)^{|y^n| + |y^0|}$).
\end{enumerate}

It is not hard to see that the cases above are mutually exclusive except when $y^0$ is a constant path and $(y^n, \ldots, y^1)$ is a disc sequence.
The formulae in \cref{eq: mu for a polygon with y^1 = zy^0,eq: mu for a polygon with y^n = y^0z} are then related by the transition functor along $\gamma_0 = \gamma_n$ and hence get identified under the equivalence relation in \cref{eq: A_Gamma(L_0 L_1)}.

The strict units in $\V_\bGamma$ are given by
$\id \in \C_{y(1)}(M^{y_+}, M^{y_-}) = \V_\bGamma(y) \subset \V_\bGamma(\LL, \LL)$
for a constant path $y$ (for which $y_+ = y_-$).
That $\V_\bGamma$ is a strictly unital $A_\infty$-category is the content of \cref{t: A_infty-categories of arc systems}\,(1).
It can be proven similarly to \cite{HKK17}*{Proposition 3.1}, we get it as a corollary of \cref{t: topological Fukaya category inside geometric}.

\subsubsection{Basic properties}
\label{subsubs: Basic properties}
Recall from the introduction that \cref{t: A_infty-categories of arc systems} implies that we have a well-defined up to homotopy equivalence \emph{topological Fukaya category} ${\Wt(\S; \C)}$, which is given by $\Tw \V_\bGamma$ for an arbitrary full arc brane system $\bGamma$ in $(\S; \C)$.
We prove \cref{t: A_infty-categories of arc systems} in \cref{subs: Proofs of the main theorems}.
We now give conditions under which $\Wt(\S; \C)$ has a dg- or formal generator.

\begin{prop} \label{p: formal full arc system}
Let $(\S = (S, Q, \eta); \C)$ be a weighted surface.
\begin{enumerate}
	\item If each component of $S$ containing no points from $\Lambda = Q \cap \pd S$ is a disc or an annulus, then there exists a full arc brane system $\bGamma$ such that $\V_\bGamma$ is a dg-category (i.e.\ $\bmu^n = 0$ for $n \ge 3$).
	\item If in addition to (1) each disc component of $S$ has at least two points from $\Lambda$ and the fibres of $\C$ have zero differential, then $\bGamma$ can be chosen so that $\V_\bGamma$ is a graded $\k$-linear category (i.e.\ $\bmu^n = 0$ for $n \ne 2$).
\end{enumerate}

\end{prop}

\begin{figure}[!htbp]
	\centering
	\begin{subfigure}[t]{0.5\linewidth}
		\centering
		\includegraphics[width=\linewidth]{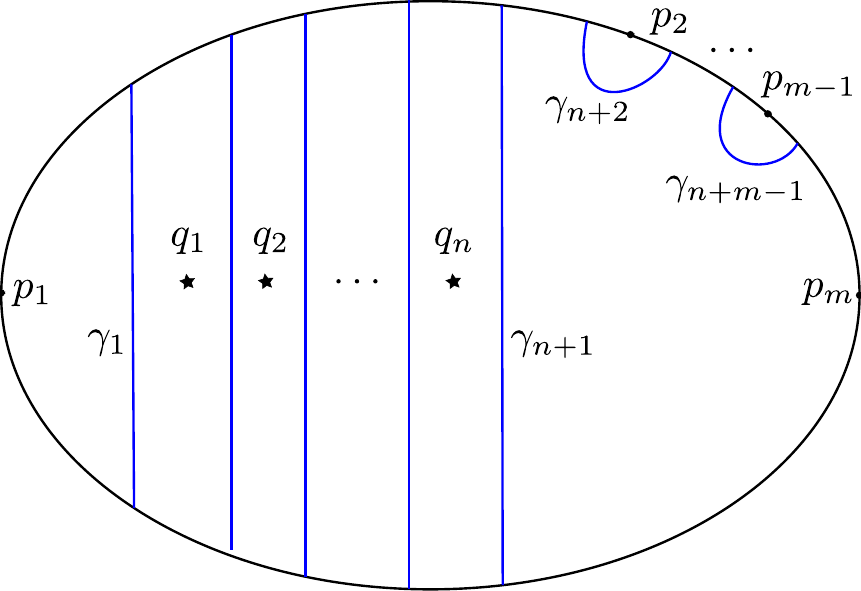}
		\caption{$S$ is a disc.
		$\Gamma = \{\gamma_1, \ldots, \gamma_{n + m - 1}\}$.}
		\label{subfig: formal arc system on a disc}
	\end{subfigure}%
	\hfill%
	\begin{subfigure}[t]{0.5\linewidth}
		\centering
		\includegraphics[width=0.9\linewidth]{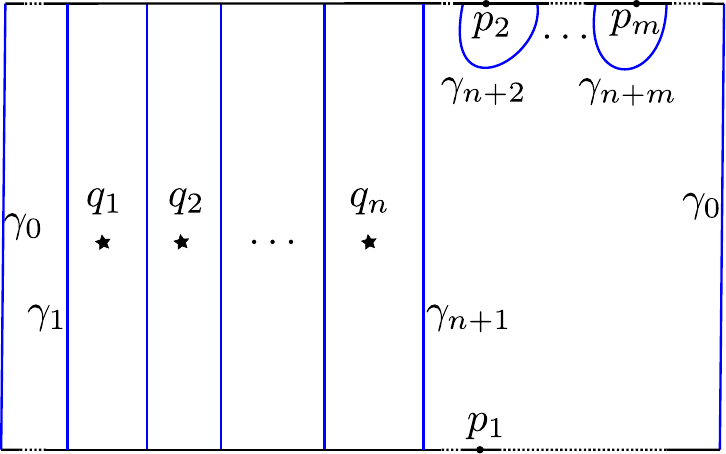}
		\caption{$S$ is not a disc.
		$\Gamma = \Gamma_0 \cup \{\gamma_1, \ldots, \gamma_{n + m}\}$.}
		{\small The continuous black line are $\pd S$, the dotted lines can belong to $\pd S$ or $\Gamma_0$.}
		\label{subfig: formal arc system}
	\end{subfigure}
	\caption{Full arc systems giving a dg-category.
	$\Lambda = \{p_1, \ldots, p_m\}$, $D = \{q_1, \ldots, q_n\}$, the arcs are in \textcolor{blue}{blue}.}
	\label{fig: formal arc systems}
\end{figure}

\begin{proof}
	We prove (1).
	We can deal with each component of $S$ separately, so we assume $S$ is connected.
	If $S$ is a disc, then we choose the full arc system as in \cref{subfig: formal arc system on a disc}.
	Now suppose $S$ is not a disc and $\Lambda \ne \varnothing$.
	Then we can cut $S$ by an arc system $\Gamma_0$ to get a single disc and then glue back along $\gamma_0 \in \Gamma_0$ to obtain an annulus.
	We then add to $\Gamma_0$ the arcs as in \cref{subfig: formal arc system} to obtain a full arc system $\Gamma$.
	If $S$ is an annulus and $\Lambda = \varnothing$, then $\Gamma = \{\gamma_1, \ldots, \gamma_{n}\}$ in \cref{subfig: formal arc system} is a full arc system.

	In each case we obtained a full arc system $\Gamma$ in $S$ such that each region cut out by the arcs either has a point from $\Lambda$ or is a quadrilateral or is a bigon (the bigons exist only when $S$ is a disc with at most one point from $\Lambda$).
	Let $\bGamma$ be a full arc brane system supported on $\Gamma$.
	As quadrilaterals and bigons can only glue to quadrilaterals and bigons, the only contributions to $\bmu_{\V_\bGamma}$ from \cref{eq: mu for a polygon with y^1 = zy^0,eq: mu for a polygon with y^n = y^0z} are to $\bmu^2$ and $\bmu^1$.
	We conclude that $\bmu^n = 0$ for $n \ge 3$ in $\V_\bGamma$.

	To prove (2) observe that $\bmu^1$ in $\V_\bGamma$ only comes from \cref{eq: mu for a polygon with y^1 = zy^0,eq: mu for a polygon with y^n = y^0z} for bigons or \cref{eq: mu for y^0 = y^1} hence vanishes by the assumption.
\end{proof}

\subsubsection{Classical topological Fukaya category}
\label{subsubs: Classical topological Fukaya category}
Our definition of the topological Fukaya category $\Wt(\S; \C)$ is a direct generalisation of the classical topological Fukaya category introduced by Haiden--Katzarkov--Kontsevich~\cite{HKK17}.
To specialise back to \cite{HKK17} take $Q = \Lambda \subset \pd S$, $\k$ a field and $\C$ a trivial constant system (i.e.\ all fibres are the same dg-category, all transition functors are the identity dg-functor, and all weights are the identity morphism of dg-functors) with fibre $\k$ viewed as a dg-category with one object.

Note that the signs in our definition and in \cite{HKK17}*{\Sec3.3} are different.
The discrepancy comes from the fact that in \cite{HKK17} the arcs are not oriented and the orientation at an end of an arc is taken to point out of $S$, this convention does not extend naturally to the case when arcs are allowed to intersect.
Since in our setup we eventually extend the category to include all immersed Lagrangian branes, we choose a different convention, where the orientation is induced by the orientation of the domain.

Classical topological wrapped Fukaya categories are equivalent to other categories coming from algebraic geometry and representation theory.
These include
\begin{enumerate}
	\item Categories of matrix factorisations, see Abouzaid--Auroux--Efimov--Katzarkov--Orlov~\cite{AAEFO13} and Pascaleff--Sibilla~\cite{PS19} for Landau--Ginzburg models, and Bocklandt~\cite{B16} for a non-commutative version.
	\item Derived categories of (stacky) curves, see Kontsevich~\cite{K09}*{Pictures}, Lekili--Polishchuk~\cite{LP18}, and Habermann~\cite{H25}.
	\item Perfect derived categories of dg-modules over graded gentle algebras, see \cite{HKK17}*{\Sec3.4}, Lekili--Polishchuk~\cite{LP20}, and Opper--Plamondon--Schroll~\cite{OPS25}.
	\item Global sections of constructible cosheaves or of parametrised perverse schobers on ribbon graphs, see \cite{HKK17}*{\Sec3.6} and Christ--Haiden--Qiu~\cite{CHQ26}*{\Sec7.2}.
\end{enumerate}

The first examples of implementing ideas from (2) in our setting are given in \cref{t: HMS},~(1)~and~(2).
(3) extends to our setting: it is the content of \cref{p: formal full arc system}.
We do not address (1) and (4) in this paper.

\subsubsection{Relative wrapped Fukaya category of a surface}
\label{subsubs: Relative wrapped Fukaya category of a surface}
Relative Fukaya categories were introduced by Seidel \cite{S02} to relate a Fukaya category of a closed symplectic manifold to a Fukaya category of an open submanifold in order to establish HMS, e.g.\ for the quartic surface~\cite{S15}.
A similar idea of weighting disc counts according to how many times a disc passes through a point goes back to the work of Ozsv\'{a}th--Szab\'{o} on Heegaard Floer homology \cite{OS04}.
Relative Fukaya categories are usually defined in a formal neighbourhood of a point, i.e.\ $\k$ is a complete local ring, see e.g.\ Perutz--Sheridan~\cite{PS23}.

Relative wrapped Fukaya categories of surfaces were introduced by Lekili--Perutz~\cite{LP12} and then developed by Lekili--Polishchuk~\cite{LP17} to prove HMS between the `$n$-Tate curve' and the $n$-punctured torus over $\ZZ[[t_1, \ldots, t_n]]$.
Bocklandt--van de Kreeke~\cites{BvdK24,BvdK25} gave an algebraic definition of relative wrapped Fukaya category linear over $\CC[[q_1, \ldots, q_n]]$ where $q_1,\ldots q_n$ are punctures.
Evans--Lekili~\cite{EL26} gave a combinatorial definition of the relative wrapped Fukaya category over an arbitrary commutative ring, but were only able to prove that the resulting category does not depend on choices after completion.
\cref{t: A_infty-categories of arc systems} resolves this issue.

To specialise to the setup of \cite{EL26} from ours take $Q = D \subset \Int S$ and $\C\colon \P^+_2(\S) \to \dgcat_\k$ such that all fibres of $\C$ are $\k$ viewed as a dg-category with one object, all transition functors are the identity, and $\C$ factors through $\H^+_2(\S)$ from \cref{d: directed relative homology 2-category}.
Indeed, unwrapping the definitions, we see that the data of such $\C$ is equivalent to specifying a weight
$\w(q) \in Z^0(\Fun_{dg}(\k, \k))(\id_{\k}, \id_{\k}) = \k$
for each $q \in Q$.

\section{Mirror symmetry examples}
\label{s: Mirror symmetry examples}
In this section we prove \cref{t: HMS}.
We prove part (1) in \cref{subs: Total spaces of projective completions of line bundles}, part (2) in \cref{subs: Bamboo forest blow-ups of affine plane}, and part (3) in \cref{subs: Cotangent bundles of spheres}.

\subsection{Total spaces of projective completions of line bundles}
\label{subs: Total spaces of projective completions of line bundles}
To prove \cref{t: HMS}\,(1) we use the fact that for a 2-dimensional vector bundle $\E$ on $X$ the category $\Dpdg \PP(\E)$ is the gluing of $\Dpdg X$ with itself along a bimodule.
When $\E = \O \oplus \L$ we can realise this gluing as the topological Fukaya category of
${\S = (S = \{z \in \CC \colon 1 \le |z| \le 2\}, Q = \{1, -2\}, \eta \text{ radial})}$
with coefficients in a system of dg-categories with fibre $\Dpdg X$.
To this end we recall a particular model for $\Dpdg X$ in \cref{subsubs: Dg-category of perfect complexes}, recall the gluing of dg-categories and identify the relevant bimodule on the B-side in \cref{subsubs: Gluing of dg-categories}, and in \cref{subsubs: Proof of HMS (1)} identify the bimodule on the A-side and prove \cref{t: HMS}\,(1).

\subsubsection{Dg-category of perfect complexes}
\label{subsubs: Dg-category of perfect complexes}
We recall a \v{C}ech-style dg-enhancement of the derived category of perfect complexes on a separated Noetherian scheme $X$ by Wei~\cite{W16}.
To use the enhancement as a fibre of a weighted system of dg-categories we need it to be small and to enhance the tensor product with a line bundle to a strict dg-autoequivalence.

\begin{remark}
	One can probably do the same with a dg-enhancement via h-injective complexes by carefully choosing a small subcategory of that large dg-enhancement as in Kuznetsov--Lunts~\cite{KL15}*{\Sec3.11} and then apply a dg-version of strictification of autoequivalences as in Shinder~\cite{Sh18}*{Theorem 5.4}.
	That should work in larger generality, possibly for any algebraic stack $X$.
	We do not pursue this direction as it looks less geometric and enlightening.
\end{remark}

Fix a separated Noetherian scheme $X$ over $\k$ and a finite affine open cover ${\U = \{U_i\}_{i \in I}}$ of $X$.
As usual, denote $U_{i_0\ldots i_p} = U_{i_0} \cap \ldots \cap U_{i_p}$.

Recall the small dg-category $\Uproj^\bullet$ from \cite{W16}*{\Sec2.2} (we describe a simplified version sufficient for our setting, in particular, all products become direct sums).
An object of $\Uproj^\bullet$ is a formal direct sum $\E = \bigoplus_{i \in I} \E_i$ where $\E_i$ is a graded $\O_{U_i}$-module that is a direct summand subsheaf of $\bigoplus_{\alpha = -n}^{n} \O_{U_i}^{\oplus n}[-\alpha]$ for some $n \ge 0$.
The morphism space is
\begin{equation} \label{eq: Uproj(E F)}
	\Uproj^\bullet(\E, \F) = \bigoplus_{p \ge 0} \bigoplus_{(i_0, \ldots, i_p) \in I^{p + 1}} \O_{U_{i_0\ldots i_p}}\text{-}\Mod^\bullet\left(\res{\E_{i_p}}{U_{i_0\ldots i_p}}, \res{\F_{i_0}}{U_{i_0\ldots i_p}}\right).
\end{equation}
For $u \in \Uproj^\bullet(\E, \F)$ denote by $u^{p, q}_{i_0\ldots i_p}$ the corresponding component.
The \emph{total degree} is given by $|u^{p, q}_{i_0\ldots i_p}| = p + q$.

The differential is given by
\begin{equation} \label{eq: check d on Uproj(E F)}
	(\check{d} u)^{p + 1, q}_{i_0\ldots i_{p + 1}} = \sum_{k = 1}^{p} (-1)^k \res{u^{p, q}_{i_0\ldots \hat i_k \ldots i_{p + 1}}}{U_{i_0\ldots i_{p + 1}}}.
\end{equation}
Note that the index $k$ starts from $1$ and ends at $p$ in the sum above.
The composition is given by
\begin{equation} \label{eq: composition in Uproj}
	(u \cdot v)^{p + r, q + s}_{i_0\ldots i_{p + r}} = (-1)^{qr} \res{u^{p, q}_{i_0\ldots i_p}}{U_{i_0\ldots i_{p + r}}} \cdot \res{v^{r, s}_{i_p \ldots i_{p + r}}}{U_{i_0\ldots i_{p + r}}}.
\end{equation}
Note that $\Uproj^\bullet$ is the additive enlargement of its full subcategory on objects
$\bigsqcup_{i \in I} \Ob \left(\O_{U_i}\text{-}\mathrm{proj}\right)$.

We now describe the triangulated dg-category $\Twp (\U)$ of \emph{twisted perfect complexes} from \cite{W16}*{\Sec2.3}.
Note that the construction is different from twisted complexes we consider elsewhere in the paper in that the differential of a twisted perfect complex is not upper-triangular and its diagonal entries satisfy a certain non-degeneracy condition.

Namely, an object of $\Twp (\U)$ is a pair $(\E, a)$ where $\E \in \Ob \Uproj^\bullet$ and $a \in \Uproj^1(\E, \E)$ such that the Maurer--Cartan equation
$
\check da + a \cdot a = 0
$
is satisfied and $a^{1, 0}_{ii}$ is chain homotopic to $\id_{(\E_i, a^{0, 1}_i)}$ for each $i \in I$.

The morphism space is the dg-$\k$-module
$$
\Twp (\U)((\E, a), (\F, b)) = \Uproj^\bullet(\E, \F) \text{ with } du = \check{d}u + b \cdot u - (-1)^{|u|} u \cdot a.
$$

Denote by $\Shp(X)$ the dg-category of perfect complexes of $\O_X$-modules.
Recall the \emph{sheafification dg-functor} $\SS \colon \Twp(\U) \to \Shp(X)$ from \cite{W16}*{\Sec3.1}.
First, for $\E \in \Ob \Uproj^\bullet$ and an open $V$ in $X$ define
$$
\check C(\U, \E; V) = \bigoplus_{p \ge 0} \bigoplus_{q \in \ZZ} \check C^p(\U, \E^q; V) \text{ where } \check C^p(\U, \E^q; V) = \bigoplus_{(i_0, \ldots, i_p) \in I^{p + 1}} \Gamma(U_{i_0\ldots i_p} \cap V, \E^q_{i_0}).
$$
For $c\in \check C(\U, \E; V)$ denote by $c^{p, q}_{i_0\ldots i_p}$ the corresponding component.
We also have $\Uproj^\bullet(\E, \F; V)$ defined as in \cref{eq: Uproj(E F)} but with $U_{i_0\ldots i_p} \cap V$ in place of $U_{i_0\ldots i_p}$.
The same modification of \cref{eq: check d on Uproj(E F),eq: composition in Uproj} gives the dg-category $\Uproj^\bullet(-; V)$.

Moreover, $\check C(\U, \E; V)$ has a differential given by
$$
(\check{d} c)^{p + 1, q}_{i_0\ldots i_{p + 1}} = \sum_{k = 1}^{p + 1} (-1)^k \res{c^{p, q}_{i_0\ldots \hat i_k \ldots i_{p + 1}}}{U_{i_0\ldots i_{p + 1}} \cap V}
$$
and the map
\begin{equation} \label{eq: Uproj otimes C to C}
	\Uproj^\bullet(\E, \F; V) \otimes \check C(\U, \E; V) \to \check C(\U, \F; V)
\end{equation}
given by \cref{eq: composition in Uproj} with $c$ in place of $v$ and $U_{i_0\ldots i_p} \cap V$ in place of $U_{i_0\ldots i_p}$ satisfies the graded Leibniz rule.

Now, for a twisted perfect complex $(\E, a)$ set
$$
\SS(\E, a) = \left(\bigoplus_{n \in \ZZ} \SS^n(\E), \check d + a \cdot\right) \text{ where } \SS^n(\E)(V) = \bigoplus_{p + q = n} \check C^p(\U, \E^q; V)
$$
and define $\SS$ on morphisms via restriction and \cref{eq: Uproj otimes C to C}.
By \cite{W16}*{Theorem 3.32, Corollary 6.6, and Appendix 2} the induced exact functor
$$
\SS \colon H^0(\Twp(\U)) \to \Dp(X)
$$
is an equivalence for any finite affine open cover $\U$.
That is, $\Twp(\U)$ is a dg-enhancement of $\Dp(X)$.
We write $\Dpdg X$ for $\Twp(\U)$ when we do not need to specify the cover.

Now let $\L$ be a line bundle on $X$.
Take $\U$ trivialising $\L$.
Let $s \in \check C^0(\U, \L)$ be a \v{C}ech cochain of non-vanishing sections and $l \in \check Z^1(\U, \O_X^\times)$ the corresponding \v{C}ech cocycle with $l_{i_0i_1} = s_{i_1}/s_{i_0}$.
Define $F_l\colon \Uproj^\bullet \to \Uproj^\bullet$ by identity on objects and
\begin{equation} \label{eq: F_l}
	F_l(u)^{p, q}_{i_0\ldots i_p} = l_{i_0i_p} u^{p, q}_{i_0\ldots i_p}
\end{equation}
on morphisms.

\begin{lemma}~ \label{l: automorphisms induced by line bundle}
	\begin{enumerate}
		\item $F_l$ is a strict dg-autoequivalence of $\Uproj^\bullet$ with $F_l^{-1} = F_{l^{-1}}$.
		\item The induced strict dg-autoequivalence $F_l$ of $\Twp(\U)$ enhances $- \otimes \L$.
		Namely,
		$$
			\ph_{(\E, a)}\colon \SS(\E, a) \otimes \L \to \SS(\E, F_l(a)), \; c^{p, q}_{i_0\ldots i_p} \otimes \sigma_V \mapsto \frac{\sigma_V}{s_{i_0}} c^{p, q}_{i_0\ldots i_p}
		$$
		defines an isomorphism of dg-functors $\SS(-)\otimes \L \to \SS \circ F_l$ from $\Twp(\U)$ to $\Shp X$.
		\item For $\sigma \in \Gamma(X, \L)$ the \v{C}ech cochain $\sigma/s \in \check C^0(\U, \O_X)$ defines a morphism of dg-functors ${\sigma/s\colon \id_{\Twp(\U)} \to F_l}$ which $\ph$ identifies with the multiplication by $\sigma\colon \SS \to \SS(-) \otimes \L$.

	\end{enumerate}
\end{lemma}

\begin{proof}
	(1) First prove $F_l$ is a dg-functor.
	Compatibility with composition \cref{eq: composition in Uproj} follows from the cocycle condition $l_{i_0i_p} l_{i_pi_{p + r}} = l_{i_0i_{p + r}}$.
	Compatibility with $\check d$ holds since \cref{eq: check d on Uproj(E F)} does not change the outer indices.
	Next, $F_{l^{-1}}$ is a \v{C}ech cocycle corresponding to $\L^{-1}$, so is also a dg-functor.
	Clearly, $F_l \circ F_{l^{-1}} = \id_\B = F_{l^{-1}} \circ F_l$.

	(2) Since $F_l$ is a strict dg-autoequivalence and $F_l(a)^{1, 0}_{ii} = a^{1, 0}_{ii}$, it preserves $\Twp(\U)$.
	Note that $\SS$ and $\SS \circ F_l$ assign the same graded sheaf $\SS^\bullet(\E)$ to the twisted perfect complex $(\E, a)$.
	The map $\ph_{\E, a}$ commutes with $\check d$ since $\check d$ keeps $i_0$, and it intertwines $a \cdot$ with $F_l(a) \cdot$ since
	$\frac{s_{i_p}}{s_{i_0}} \cdot \frac{\sigma_V}{s_{i_p}} = \frac{\sigma_V}{s_{i_0}}$,
	so it is a map of complexes.
	The same computation with $u$ in place of $a$ gives naturality of $\ph$.
	Each $\ph_{\E, a}$ is a local isomorphism as it multiplies each component by a unit.

	(3) The cochain $\sigma/s$ is $\check d$-closed since \cref{eq: check d on Uproj(E F)} vanishes in \v{C}ech degree $0$, and $F_l(u) \cdot \sigma/s = \sigma/s \cdot u$ since
	$\frac{s_{i_p}}{s_{i_0}} \cdot \frac{\sigma}{s_{i_p}} = \frac{\sigma}{s_{i_0}}$.
	Taking $u = a$ shows that $\sigma/s$ is closed in $\Twp(\U)$, and arbitrary $u$ gives naturality.
	Finally, $\ph(c \otimes \sigma) = \sigma/s \cdot c$ by definition.
\end{proof}

\subsubsection{Gluing of dg-categories}
\label{subsubs: Gluing of dg-categories}

We recall the gluing of dg-categories by Tabuada~\cite{T07}*{Definition 3.57} and Orlov~\cite{O16}*{Definitions 3.4, 3.5}.
Let $\A$ and $\B$ be two small dg-categories over a commutative ring $\k$ and $\T$ an $\A$-$\B$-bimodule (i.e.\ a dg-functor $\A^{op} \otimes \B \to Ch(\k)$).
The \emph{upper triangular dg-category} $\A \with{\T} \B$ is given by
$\Ob \A \with{\T} \B = \Ob \A \sqcup \Ob \B$ and
$$
\A \with{\T} \B(M, N) =\begin{cases}
	\A(M, N) &\text{if } M, N \in \Ob \A;\\
	\B(M, N) &\text{if } M, N \in \Ob \B;\\
	\T(M, N) &\text{if } M \in \Ob \A \text{ and } N \in \Ob \B;\\
	0 &\text{otherwise}
\end{cases}
$$
with composition given by that in $\A$ and $\B$, and the bimodule structure on $\T$.
The \emph{gluing of $\A$ and $\B$ along $\T$} is the triangulated envelope $\A \oright{\T}\B = \Tw \left(\A \with{\T} \B\right)$.

We recast \cite{O16}*{Example 3.9} for our setting.

\begin{lemma} \label{l: gluing for projectivisation of line bundle}
	Let $\L$ be a line bundle on a separated Noetherian scheme $X$ over $\k$ and $\U^X = \{U_i\}_{i \in I}$ a finite affine open cover of $X$ trivialising $\L$.
	Then $\Dpdg \PP(\O_X \oplus \L)$ is quasi-equivalent to $\B \oright{\T} \B$ where ${\B = \Twp(\U^X)}$ and the $\B$-$\B$-bimodule $\T$ is given by
	$$
	\T(M, N) = \B(M, N) \oplus \B(F_l(M), N) \text{ for } M, N \in \Ob \B.
	$$
\end{lemma}
\begin{proof}
	Let $\pi\colon Y =\PP(\O_X \oplus \L) \to X$
	be the projection and
	$$
	x_\O \in \Gamma(Y, \O_Y(1)), \, x_\L \in \Gamma(Y, \O_Y(1) \otimes \pi^*\L)
	$$
	the coordinates of the tautological inclusion
	$\O_Y(-1) \hookrightarrow \O_Y \oplus \pi^*\L$.
	By Orlov's semi-orthogonal decomposition~\cite{O92}*{Theorem 2.6}, which holds in our setting by Bergh--Schn\"urer~\cite{BS20}*{Corollary 6.8}, we have
	$$
		\Dp(Y) = \langle \pi^* \Dp(X), \pi^* \Dp(X) \otimes \O_Y(1)\rangle.
	$$
	By \cite{O16}*{Proposition 3.8} (it is stated over a field, but the proof works verbatim over any commutative ring $\k$)
	it is enough to find a finite affine cover $\U^Y$ of $Y$, two dg-functors $\pi^*_0, \pi^*_1\colon \B \to \Twp(\U^Y)$ enhancing $\pi^*$ and $\pi^*(-) \otimes \O_Y(1)$ respectively, and a quasi-isomorphism of $\B$-$\B$-bimodules
	\begin{equation} \label{eq: morphism of B-B-bimodules}
		\B(-_0, -_1) \oplus \B(F_l(-_0), -_1) \to \Twp(\U^Y)(\pi^*_0 (-_0), \pi^*_1 (-_1)).
	\end{equation}

	For $i \in I$ define $U^Y_{i, \O} = \pi^{-1}(U_i) \cap \{x_\L \ne 0\}$ and ${U^Y_{i, \L} = \pi^{-1}(U_i) \cap \{x_\O \ne 0\}}$ and set
	${\U^Y = \{U^Y_{i, \O}, U^Y_{i, \L}\}_{i \in I}}$.
	Let $s \in \check C^0(\U, \L)$ be a \v{C}ech cochain of non-vanishing sections.
	Then $s^Y \in \check C^0(\U^Y, \O_Y(1))$
	given by $s^Y_{i, \O} = x_\L \cdot \pi^*s_i^{-1}$ and $s^Y_{i, \L} = x_\O$ also consists of non-vanishing sections.
	Let $l \in \check Z^1(\U, \O_X^\times)$ and $l^Y \in \check Z^1(\U, \O_Y^\times)$ be the corresponding \v{C}ech cocycles, and $F_l$, $F_{l^Y}$ the dg-functors from \cref{eq: F_l}.

	Define $\pi_0^*$ by pulling back the twisted perfect complexes along $U^Y_{i, \O}, U^Y_{i, \L} \to U_i$.
	Set $\pi_1^* = F_{l^Y} \circ \pi_0^*$.
	\cref{l: automorphisms induced by line bundle} shows that $\pi_0^*$ and $\pi_1^*$ enhance $\pi^*$ and $\pi^*(-) \otimes \O_Y(1)$ respectively.

	The morphism of $\B$-$\B$-bimodules in \cref{eq: morphism of B-B-bimodules} is given by
	\begin{equation} \label{eq: morphism of B-B-modules on (u v)}
		(u, v) \mapsto \pi^*_1(u) \cdot \frac{x_\O}{s^Y}\pi^*_0 + \pi^*_1(v) \cdot \frac{x_\L}{s^Y \cdot \pi^*s}\pi^*_0
	\end{equation}
	where
	$$
	\frac{x_\O}{s^Y}\pi^*_0\colon \pi^*_0 \to F_{l^Y} \circ \pi_0^* = \pi^*_1 \text{ and } \frac{x_\L}{s^Y \cdot \pi^*s} \pi_0^* \colon \pi^*_0 \to F_{l^Y} \circ F_{\pi^*l} \circ \pi_0^* = \pi^*_1 \circ F_l
	$$
	are morphisms of dg-functors given by \cref{l: automorphisms induced by line bundle}\,(3).
	\cref{eq: morphism of B-B-modules on (u v)} defines a quasi-isomorphism since on cohomology it is the adjunction isomorphism composed with the isomorphism ${\O_X \oplus \L^{-1} \to \pi_*\O_Y(1)}$ given by $x_\O$ and $x_\L$.
\end{proof}

\subsubsection{Proof of \cref{t: HMS}\,(1)}
\label{subsubs: Proof of HMS (1)}
Before proving \cref{t: HMS}\,(1) we give an abstract dg-category construction realising the gluing as in \cref{l: gluing for projectivisation of line bundle} as a topological Fukaya category.
Let $\B$ be a small dg-category over $\k$ and $F\colon \B \to \B$ a strict dg-autoequivalence.
Define the $\B$-$\B$-bimodule $\T = \T({\id_\B} \oplus F)$ by $\T(M, N) = \B(M, N) \oplus \B(F(M), N)$ and the corresponding upper-triangular dg-category $\D = \B \with{\T} \B$.
We describe the weighted system $\C \colon \P^+_2(\S) \to \dgcat_\k$ on
$$
\S = (S = \{z \in \CC \colon 1 \le |z| \le 2\}, Q = \{1, -2\}, \eta \text{ radial})
$$
such that $\Wt(\S; \C)$ is homotopy equivalent to $\Tw \D$.

The fibres of $\C$ are equal to $\B$.
For a lazy 1-path $\gamma\colon [0, 1] \to S \setminus Q$ lift the polar coordinate on $S \subset \CC^\times$ along $\gamma$ to $\theta_\gamma \colon [0, 1] \to \RR$ and set
$$
n(\gamma) = \left\lfloor \frac{\theta_\gamma(1)}{2\pi}\right\rfloor - \left\lfloor \frac{\theta_\gamma(0)}{2\pi}\right\rfloor \text{ and } \C(\gamma) = F^{n(\gamma)}.
$$
Note that $n(\gamma)$ is independent of the lift, invariant under homotopy-rel-endpoints, and additive with respect to concatenation of paths.
Hence, by weighting all 2-morphisms in $\P^+_2(\S)$ by the identity, we obtain a weighted system of dg-categories $\C$.
Note that $\C$ factors through the directed relative fundamental 2-category $\Pi^+_2(\S)$ from \cref{d: directed relative fundamental 2-category}.

\begin{lemma} \label{l: gluing along bimodule is W(S C)}
	In the above notation $\Tw \D$ is homotopy equivalent to $\Wt(\S; \C)$.
\end{lemma}

\begin{figure}[!htbp]
	\includegraphics[width=0.4\linewidth]{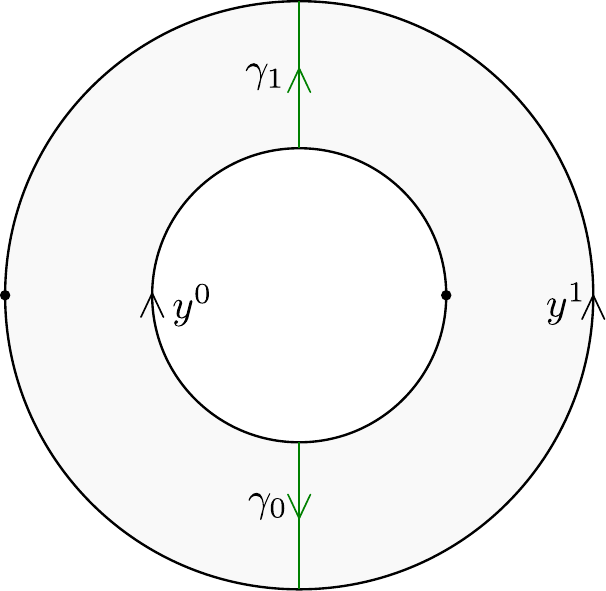}
	\caption{A full arc system $\{\gamma_0, \gamma_1\}$ in $\S$.}
	\label{fig: full arc system in annulus with two stops}
\end{figure}

\begin{proof}
	Let $\Gamma = \{\gamma_0, \gamma_1\}$ be the full arc system in $\S$ from \cref{fig: full arc system in annulus with two stops} and consider the full arc brane system
	$$
		\bGamma = \{\LL_i(M) = (\gamma_i, const, M) \mid i \in \{0, 1\}, \, M \in \Ob \B\}
	$$
	($const$ above means the constant path as $\gamma_i^*\eta = \dot \gamma$).
	A straightforward computation shows that for any $M, N \in \Ob \B$ we have
	$$
	\begin{aligned}
		&\V_\bGamma(\LL_i(M), \LL_i(N)) = \B(M, N) \text{ for } i = 0, 1;\\
		&\V_\bGamma(\LL_0(M), \LL_1(N)) = \V_\bGamma(y^0) \oplus \V_\bGamma(y^1) = \B(M, N) \oplus \B(F(M), N);\\
		&\V_\bGamma(\LL_1(M), \LL_0(N)) = 0
	\end{aligned}
	$$
	The $A_\infty$-operations in $\V_\bGamma$ only come from \cref{eq: mu for y^0 = y^1,eq: mu for y^0 = y^2y^1} and, using dg-category notation for $\V_\bGamma$, i.e.\ $da = (-1)^{|a|}\bmu^1(a)$ and $a \cdot b = (-1)^{|b|}\bmu^2(a, b)$, we identify $\V_\bGamma$ with $\D$.
	Thus $\Tw \D = \Tw \V_\bGamma$ is homotopy equivalent to $\Wt(\S; \C)$.
\end{proof}

\begin{proof}[Proof of \cref{t: HMS}\,(1)]
	In the notation of \cref{l: gluing for projectivisation of line bundle} $\Dpdg \PP(\O_X \oplus \L)$ is quasi-equivalent to $\B \oright{\T} \B$ and in the notation of \cref{l: gluing along bimodule is W(S C)} $\B \oright{\T} \B = \Tw \D$ is quasi-equivalent to $\Wt(\S; \C)$.
\end{proof}

\subsection{Bamboo forest blow-ups of affine plane}
\label{subs: Bamboo forest blow-ups of affine plane}
Suppose we are in the setup of \cref{t: HMS}\,(2), i.e.\ $X$ is a bamboo forest blow-up of $\AA^2_{\k}$ for an infinite field $\k$.

First, we introduce some notation for $X$.
Let $\pi\colon X \to \AA^2$ be the projection, $\{p_1, \ldots, p_n\} \subset \AA^2(\k)$ the image of the exceptional locus, and $\pi^{-1}(p_i) = C_i^1\cup \ldots \cup C_{i}^{k_i}$ with
$$
C_i^{j_0} \cdot C_i^{j_1} = \begin{cases}
	-2& \text{if } j_0 = j_1 < k_i;\\
	-1& \text{if } j_0 = j_1 = k_i;\\
	1& \text{if } |j_0 - j_1| = 1;\\
	0& \text{otherwise}
\end{cases}
$$
and $\iota_i^j\colon \PP^1 \hookrightarrow X$ parametrising $C_i^j$ in such a way that $\iota_i^1(\infty) = C^1_i \cap C^2_i = \iota_i^2(0)$, $\ldots$,
${\iota_i^{k_i - 1}(\infty) = C^{k_i - 1}_i \cap C^{k_i}_i = \iota_i^{k_i}(0)}$.

\begin{figure}[!htbp]
	\begin{subfigure}[t]{0.5\linewidth}
		\includegraphics[width=0.9\linewidth]{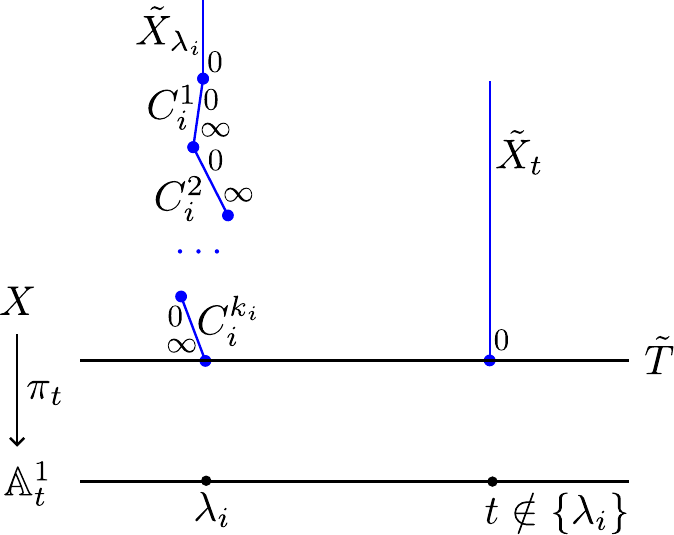}
		\caption{Bamboo forest blow-up as a relative curve.}
		\label{subfig: Plane blow-up}
	\end{subfigure}%
	\hfill%
	\begin{subfigure}[t]{0.5\linewidth}
		\includegraphics[width=0.98\linewidth]{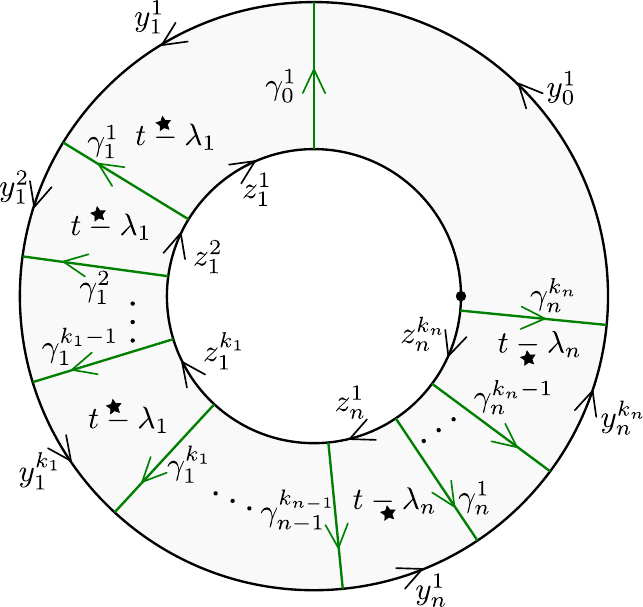}
		\caption{$(\S; \C)$ with a full arc system
		$\Gamma = \{\gamma_0^1, \gamma_1^1, \ldots, \gamma_1^{k_1}, \ldots, \gamma_n^1, \ldots, \gamma_n^{k_n}\}$.}
		\label{subfig: Annulus for plane blow-up}
	\end{subfigure}
	\caption{}
	\label{fig: Plane blow-up and its mirror}
\end{figure}

Since $\k$ is infinite we can find linear coordinates $(t, x)$ on $\AA^2$ such that the $t$-coordinate $\lambda_i$ of $p_i$ is different for different $i$ and $\iota_i^1(\infty)$ does not lie on the strict transform $\tilde X_{\lambda_i}$ of $\{t = \lambda_i\}$.
Composing with an automorphism $(t, x) \mapsto (t, x - \ph(t))$, where $\ph \in \k[t]$ has a prescribed $(k_i - 1)$-jet at each $\lambda_i$, we can assume that $p_i = (\lambda_i, 0)$, $\tilde X_{\lambda_i}$ meets $\pi^{-1}(p_i)$ only at $\iota_i^1(0)$, the strict transform $\tilde T$ of $\{x = 0\}$ meets $\pi^{-1}(p_i)$ only at $\iota_i^{k_i}(\infty)$, and all $\iota_i^j$ are the natural parametrisations coming from the successive blow-up structure.
See \cref{subfig: Plane blow-up} for the picture of the projection $\pi_t = \pr_1 \circ \pi \colon X \to \AA^1_t$ that makes $X$ into a relative curve (i.e.\ a flat morphism with 1-dimensional fibres) over $\AA^1$.

Secondly, we specify the graded marked surface $\S$ and the weighted system of dg-categories $\C\colon \P^+_2(S, Q) \to \dgcat_{\k}$ from \cref{t: HMS}\,(2) corresponding to $X$.
Namely,
$$
\S = (S = \{z \in \CC \colon 1 \le |z| \le 2\}, Q, \eta \text{ radial})
$$
where $\Lambda = Q \cap \pd S = \{1\}$ and
$D = Q \cap \Int S = \{q_1^1, \ldots, q^{k_1}_1, \ldots, q_n^1, \ldots, q_n^{k_n}\}$.
The weighted system of dg-categories $\C\colon \P^+_2(S, Q) \to \dgcat_{\k}$ factors through $\H^+_2(S, Q)$ from \cref{d: directed relative homology 2-category}, has fibres $\k[t]$ viewed as a dg-category with one object, all transition functors are the identity, and the weight of $v \colon D \to \ZZ_{\ge 0}$ is
$$
\prod_{i = 1}^n \prod_{j = 1}^{k_i}(t - \lambda_i)^{v(q_i^j)} \in \k[t] = Z^0(\Fun_{dg}(\k[t], \k[t]))(\id, \id).
$$
See \cref{subfig: Annulus for plane blow-up} for the picture of $(\S; \C)$.

Thirdly, we describe a $\k[t]$-linear category $\B$ such that $\Tw \B$ is quasi-equivalent to $\Dp X$ and $\Wt(\S; \C)$.
$\B$ is the quotient of the $\k[t]$-linear path category of the quiver
\begin{equation} \label{eq: quiver for annulus with one stop}
\xy
	(0.0,28.0)*+{\binom{1}{0}}="v0";
	(-13.0,24.8)*+{\binom{1}{1}}="v1";
	(-23.0,15.9)*+{\binom{2}{1}}="v2";
	(-27.8,3.4)*+{\vdots}="d1";
	(-26.2,-9.9)*+{\binom{k_1 - 1}{1}}="v3";
	(-18.6,-21.0)*+{\binom{k_1}{1}}="v4";
	(-6.7,-27.2)*+{\cdots}="d2";
	(6.7,-27.2)*+{\binom{k_{n - 1}}{n - 1}}="v5";
	(18.6,-21.0)*+{\binom{1}{n}}="v6";
	(26.2,-9.9)*+{\vdots}="d3";
	(27.8,3.4)*+{\binom{k_n - 1}{n}}="v7";
	(23.0,15.9)*+{\binom{k_n}{n}}="v8";
	{\ar@/_/_{y^1_1} "v0";"v1"};
	{\ar@/_/_{z^1_1} "v1";"v0"};
	{\ar@/_/_{y^2_1} "v1";"v2"};
	{\ar@/_/_{z^2_1} "v2";"v1"};
	{\ar@/_/ "v2";"d1"};
	{\ar@/_/ "d1";"v2"};
	{\ar@/_/ "d1";"v3"};
	{\ar@/_/ "v3";"d1"};
	{\ar@/_/_{y^{k_1}_1} "v3";"v4"};
	{\ar@/_/_{z^{k_1}_1} "v4";"v3"};
	{\ar@/_/ "v4";"d2"};
	{\ar@/_/ "d2";"v4"};
	{\ar@/_/ "d2";"v5"};
	{\ar@/_/ "v5";"d2"};
	{\ar@/_/_{y^1_n} "v5";"v6"};
	{\ar@/_/_{z^1_n} "v6";"v5"};
	{\ar@/_/ "v6";"d3"};
	{\ar@/_/ "d3";"v6"};
	{\ar@/_/ "d3";"v7"};
	{\ar@/_/ "v7";"d3"};
	{\ar@/_/_{y^{k_n}_n} "v7";"v8"};
	{\ar@/_/_{z^{k_n}_n} "v8";"v7"};
	{\ar@/_/_{y_0^1} "v8";"v0"};
\endxy
\end{equation}
by the ideal generated by
\begin{equation} \label{eq: relations for the annulus quiver}
	y_i^j \cdot z_i^j = (t - \lambda_i)\id \text{ and } z_i^j \cdot y_i^j = (t - \lambda_i)\id \text{ for } i = 1, \ldots, n, \, j = 1, \ldots, k_i.
\end{equation}
Note that the paths that contain only $y$-arrows or only $z$-arrows give a $\k[t]$-basis of morphisms in $\B$.

We are now ready to prove the theorem.

\begin{proof}[Proof of \cref{t: HMS}\,(2)]
	We first show that $\Tw \B$ is homotopy equivalent to $\Wt(\S; \C)$.
	Consider the full arc system
	$\Gamma = \{\gamma_0^1, \gamma_1^1, \ldots, \gamma_1^{k_1}, \ldots, \gamma_n^1, \ldots, \gamma_n^{k_n}\}$
	in $\S$ from \cref{subfig: Annulus for plane blow-up}.
	Adding the constant grading and $\k[t]$ to each $\gamma \in \Gamma$ we obtain a full arc brane system
	$\bGamma = \{\LL_0^1, \LL_1^1, \ldots, \LL_1^{k_1}, \ldots, \LL_n^1, \ldots, \LL_n^{k_n}\}$
	in $(\S; \C)$.
	Note that all boundary paths have degree $0$ and the only non-zero $\bmu$ in $\V_\bGamma$ is $\bmu^2$.
	We therefore view $\V_\bGamma$ as the ordinary $\k[t]$-linear category with $a\cdot b = \bmu^2(a, b)$.
	The evident bijections between $\Ob \B$ and $\bGamma$, and between the basis of paths in the quiver \cref{eq: quiver for annulus with one stop} and boundary paths between arcs from $\Gamma$ then establish an isomorphism $\B \cong \V_\bGamma$ as path concatenation in the quiver corresponds to \cref{eq: mu for y^0 = y^2y^1} and \cref{eq: relations for the annulus quiver} corresponds to \cref{eq: mu for a polygon with y^1 = zy^0,eq: mu for a polygon with y^n = y^0z}.

	We now prove that $\Tw \B$ is quasi-equivalent to $\Dpdg X$.
	Construct $X$ from $X_0^1 = \AA^2$ by successive blow-ups at one point that create
	$C^1_1, \ldots, C^{k_1}_1, \ldots, C^{1}_n, \ldots, C^{k_n}_n$
	in this order.
	Denote the varieties we get in the process $X^j_i$, note that they are all bamboo forest blow-ups of $\AA^2$.
	Denote by $\tilde C_i^j = C_i^j + \ldots + C_i^{k_i}$ the total transform of $C_i^j \subset X_i^j$ in $X$.
	Define line bundles $\L_0^1 = \O_X$ and $\L_i^j = \O_X(\tilde C^1_1 + \ldots + \tilde C_i^j)$ for $i = 1, \ldots, n$, $j = 1, \ldots, k_i$.

	We have the bijection
	$$
		F\colon \Ob \B \to \{\L_0^1, \ldots, \L_n^{k_n}\}, \; \binom{j}{i} \mapsto \L_i^j.
	$$
	Observe that the principal divisor of $(t - \lambda_i) \in \Gamma(X, \O_X)$ is $\tilde X_{\lambda_i} + \tilde C_i^1$ and of $x \in \Gamma(X, \O_X)$ is $\tilde T + \tilde C^1_1 + \ldots + \tilde C_n^{k_n}$.
	We therefore can define $F$ on the arrows of the quiver \cref{eq: quiver for annulus with one stop} by
	$$
		F(y_0^1) = x \in \Gamma(X, \O_X(-\tilde C^1_1 - \ldots - \tilde C_n^{k_n})) = \Hom_X(\L_n^{k_n}, \L_0^1)
	$$
	and for $i = 1, \ldots, n$, $j = 1, \ldots, k_i$ by
	$$
	\begin{aligned}
		&F(y_i^j) = s_{\tilde C_i^j} \in \Gamma(X, \O_X(\tilde C_i^j)) = \Hom_X(\L_i^{j - 1}, \L_i^j) \text{ where } s_{\tilde C_i^j} \text{ is the canonical section};\\
		&F(z_i^j) = (t - \lambda_i) \in \Gamma(X, \O_X(-\tilde C_i^j)) = \Hom_X(\L_i^j, \L_i^{j - 1}).
	\end{aligned}
	$$
	where we use the notation $\L_i^0 = \L_{i - 1}^{k_{i - 1}}$ and $k_0 = 1$.
	Note that
	$$
		s_{\tilde C_i^j} \cdot (t - \lambda_i) = (t - \lambda_i) \cdot s_{\tilde C_i^j} = (t - \lambda_i) \in \Gamma(X, \O_X),
	$$
	so $F$ respects the relations \cref{eq: relations for the annulus quiver} hence defines a functor from $\B$ to the $\Hom_X$-category on objects $\L_0^1, \ldots, \L_n^{k_n}$.
	The following claim finishes the proof since $\B$ is intrinsically formal.
	\begin{claim}~
		\begin{enumerate}
			\item $F$ is an isomorphism from $\B$ to the $\Ext^*$-category on objects $\L_0^1, \ldots, \L_n^{k_n}$.
			\item The smallest full subcategory of $\Dp X$ containing $\L_0^1, \ldots, \L_n^{k_n}$ and closed under shifts, cones, and isomorphisms is $\Dp X$.
		\end{enumerate}
	\end{claim}
\end{proof}

\begin{proof}[Proof of Claim]
	We argue by induction on the number of blow-ups.
	All $X_i^j$ are bamboo forest blow-ups in normal form with respect to the same coordinates $(t, x)$, so the construction above gives $\k[t]$-linear categories $\B_i^j$, line bundles on $X_i^j$, and functors $F_i^j$.
	The pullback along $X \to X_i^j$ is fully faithful on $\Dp$ and sends these line bundles to the corresponding $\L$'s, so it is enough to prove the claim for every $X_i^j$.
	For $X_0^1 = \AA^2$ we have
	$\B_0^1 = \k[t][y_0^1] \cong \k[t, x] = \Ext^*_{\AA^2}(\O, \O)$,
	and every perfect complex on $\AA^2$ is quasi-isomorphic to a bounded complex of free modules by the Quillen--Suslin theorem.

	Let $\rho\colon Y' = X_i^j \to Y = X_i^{j - 1}$ (where $X_i^0 = X_{i - 1}^{k_{i - 1}}$) be the blow-up at $p$ with exceptional curve $E = C_i^j$, let $\L$ be the last line bundle on $Y$, so that $\L' = \rho^*\L(E)$ is the last line bundle on $Y'$, and assume the claim for $Y$.

	We first prove (2).
	By Orlov's blow-up formula~\cite{O92}*{Theorem 4.3}, which holds in our setting by \cite{BS20}*{Corollary 6.10}, we have $\Dp Y' = \langle \O_E(-1), \rho^*\Dp Y\rangle$, so every object of $\Dp Y'$ is a cone of a morphism between an object of $\rho^*\Dp Y$ and a finite direct sum of shifts of $\O_E(-1)$.
	As $(\rho^*\L)|_E$ is trivial, $\O_E(-1) = \L'|_E$ is the cone of $s_E\colon \rho^*\L \to \L'$, and we are done by the induction hypothesis.

	Now prove (1).
	Since $\B_i^{j - 1}$ and $\B_i^j$ have bases of $y$- and $z$-paths, sending $y_0^1$ to $y_0^1 \cdot y_i^j$ identifies $\B_i^{j - 1}$ with the full subcategory of $\B_i^j$ on all objects except $\binom{j}{i}$, compatibly with $F_i^{j - 1}$, $F_i^j$, and $\rho^*$.
	It remains to compare morphisms to and from $\L'$.
	The short exact sequence $0 \to \O \to \O(E) \to \O_E(-1) \to 0$ gives $\rho_*\O(E) = \O$ and $\rho_*\O(-E) = \J_p$, the ideal sheaf of $p$ (all functors are derived).
	Hence, by the projection formula, twisting by $\O(E)$ induces an isomorphism
	\begin{equation} \label{eq: Ext(L' L') = Ext(L L)}
		\Ext^*_{Y'}(\O_{Y'}, \O_{Y'}) = \Ext^*_{Y'}(\L', \L') \cong \Ext^*_Y(\L, \L) = \Ext^*_{Y}(\O_{Y}, \O_{Y}),
	\end{equation}
	and for every line bundle $\M$ on $Y$ from our collection the multiplication by $F(y_i^j) = s_E$ induces isomorphisms
	\begin{equation} \label{eq: Ext (rho^*M L')}
		\Ext^*_{Y'}(\rho^*\M, \L') \cong \Ext^*_Y(\M, \L)
	\end{equation}
	and
	\begin{equation} \label{eq: Ext (L' rho^*M)}
		\Ext^*_{Y'}(\L', \rho^*\M) \cong \Ext^*_Y(\L, \M \otimes \J_p).
	\end{equation}

	\cref{eq: Ext(L' L') = Ext(L L)} matches $\B_i^j$ since both sides are $\k[t, x]$ with the $y$-loop at $\binom{j}{i}$ mapping to $x$.
	\cref{eq: Ext (rho^*M L')} matches $\B_i^j$ since every path to $\binom{j}{i}$ from another vertex ends with $y_i^j$.
	For \cref{eq: Ext (L' rho^*M)}, every path from $\binom{j}{i}$ starts with $y_0^1$ or $z_i^j$, so $- \cdot y_i^j$ identifies $\B_i^j(\binom{j}{i}, -)$ with the $\k[t]$-span of the $y$-paths and of $(t - \lambda_i)$ times the $z$-path in $\B_i^{j - 1}(\binom{j - 1}{i}, -)$.
	On the other hand, $\Ext^*_Y(\L, \M \otimes \J_p)$ is the fibre of the evaluation map
	$\mathrm{ev}_p\colon \Hom(\L, \M) \to (\M \otimes \L^{-1})|_p \cong \k$.
	It thus suffices to show that $\mathrm{ev}_p$ vanishes on the images of the $y$-paths but not on the image of the $z$-path.
	The images of the $y$-paths are divisible by $x$ and hence vanish on $p \in \tilde T$.
	The image of a $z$-path is a product of the $(t - \lambda_l)$, so its divisor is supported on the $\tilde X_{\lambda_l}$ and on $C_i^1, \ldots, C_i^{j - 2}$, none of which contains $p$ by the choice of coordinates.
\end{proof}

\subsection{Cotangent bundles of spheres}
\label{subs: Cotangent bundles of spheres}
Abouzaid~\cite{A11} proved that the wrapped Fukaya category $\W(T^*S^n)$ is quasi-equivalent to the triangulated envelope $\Tw C_{-*}\Omega S^n$ of the chains on the loop space of $S^n$.
The Adams--Hilton model for $C_{-*}\Omega S^n$ for $n \ge 2$ is the dg-algebra $\B_{m} = (\k\langle z\rangle, \, dz = 0, \, |z| = m)$ with $m = 1 - n$.
Therefore, to prove \cref{t: HMS}\,(3) we need to find a weighted system of dg-categories $\C \colon \P^+_2(\S) \to \dgcat_\k$ on
$$
\S = (S = \{z \in \CC \colon |z| \le 2\}, Q = \{1, -1\}, \eta \text{ vertical})
$$
such that $\Tw \B_{1 - n}$ is quasi-equivalent to $\Wt(\S; \C)$.
We first do it for any $\B_m$ with $m \in \ZZ$ in \cref{subsubs: From B_m to B_m - 1} and then give a geometric heuristic explaining why the quasi-equivalence is not a coincidence in \cref{subsubs: Lefschetz fibration and symplectic reduction}.

\subsubsection{From $\B_m$ to $\B_{m - 1}$}
\label{subsubs: From B_m to B_m - 1}
\begin{figure}[!htbp]
	\centering
	\includegraphics[width=0.4\linewidth]{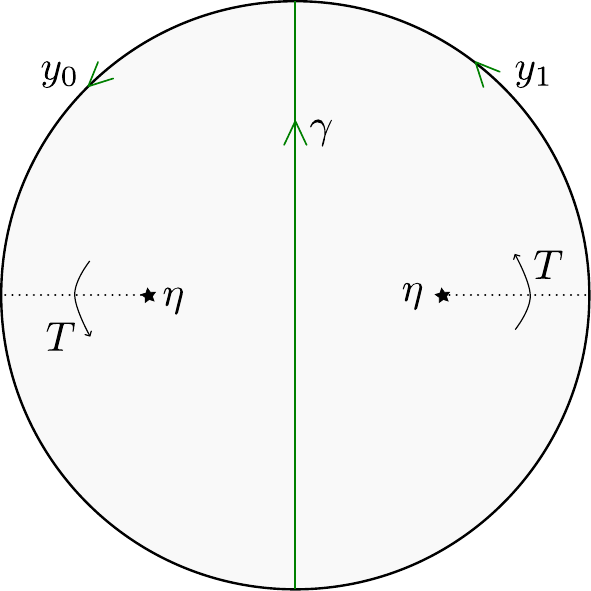}
	\caption{$(\S; \C_m)$ and the full arc system $\{\gamma\}$.}
	{\small Crossing a dotted cut in the direction of the arrow gives the transition functor $T$, and the weights at $\pm 1$ are $\eta$.}
	\label{fig: weighted system on disc with two marked points}
\end{figure}

We describe $\C_{m}\colon \P^+_2(\S) \to \dgcat_\k$ that gives $\Wt(\S; \C_m)$ quasi-equivalent to $\Tw \B_{m - 1}$.
All fibres of $\C_m$ are $\Tw \B_m$.
Define a strict dg-autoequivalence $T = [m] \circ \Tw (z \mapsto (-1)^m z)$ of $\Tw \B_m$
and a morphism of dg-functors $\eta\colon \id \Rightarrow T$ whose component from $\B_m[n]$ to $\B_m[n + m]$ is $(-1)^{mn}z[m]$.
A straightforward computation shows that $\eta T = T \eta$.
$\C_m$ is obtained from the trivial system of dg-categories on $S \setminus ([-2, -1] \cup [1, 2])$ by gluing with $T$ along $[-2, -1) \cup (1, 2]$ and weighting $\pm 1$ by $\eta$, see \cref{fig: weighted system on disc with two marked points}.

Namely, for a lazy 1-path $\gamma\colon [0, 1] \to S \setminus Q$ let $\theta^{\pm 1}_\gamma\colon [0, 1] \to \RR$ be the lift of the polar coordinate of $\gamma \mp 1$ along $\gamma$.
Set
$$
n(\gamma) = \left\lfloor \frac{\theta^{+1}_\gamma(1)}{2\pi}\right\rfloor - \left\lfloor \frac{\theta^{+1}_\gamma(0)}{2\pi}\right\rfloor + \left\lfloor \frac{\theta^{-1}_\gamma(1) - \pi}{2\pi}\right\rfloor - \left\lfloor \frac{\theta^{-1}_\gamma(0) - \pi}{2\pi}\right\rfloor
$$
and note that $n(\gamma)$ is independent of the lift, invariant under homotopy-rel-endpoints, and additive with respect to concatenation of paths.
The transition functor is $\C_m(\gamma) = T^{n(\gamma)}$.

For a positive $2$-path $\sigma$ from $\gamma_0$ to $\gamma_1$ we have
$n(\gamma_1) - n(\gamma_0) = \deg_\sigma(-1) + \deg_\sigma(1) =: \deg(\sigma)$.
The weight $\C_m(\sigma)$ is
$$
\eta^{\deg(\sigma)}T^{n(\gamma_0)}\colon T^{n(\gamma_0)} \Rightarrow T^{n(\gamma_1)}
$$
where the power of $\eta$ is with respect to the horizontal composition.

\begin{lemma} \label{l: B_m as topological Fukaya category}
	$\Wt(\S; \C_m)$ is quasi-equivalent to $\Tw \B_{m - 1}$.
\end{lemma}
\begin{proof}
	$\{\gamma\}$ from \cref{fig: weighted system on disc with two marked points} is a full arc system in $\S$ and hence $\bGamma = \{\LL = (\gamma, const, \B_m)\}$ is a full arc brane system in $(\S; \C_m)$.
	The only non-zero $\bmu$ in $\V_\bGamma$ is $\bmu^1$ and $\bmu^2$ and $\V_\bGamma$ has one object, so we describe $\V_\bGamma$ as a dg-algebra with $da = (-1)^{|a|}\bmu^1(a)$ and $ab = (-1)^{|b|}\bmu^2(a, b)$.
	Denote by $b \in \B_m = \V_\bGamma(c)$ the generator of $\B_m$, where $c$ is a constant boundary path, and by $a_i$ the element
	$1[-|y_i|] \in \Tw \B_m(T\B_m, \B_m)[-|y_i|] = \V_\bGamma(y_i)[-|y_i|]$
	for $i = 0, 1$.
	Then $|b| = m$, $|a_0| = |a_1| = m - 1$ and a straightforward computation using \cref{eq: mu for y^0 = y^2y^1,eq: mu for a polygon with y^1 = zy^0} shows
	$$
		\V_\bGamma = \frac{\k\langle b, a_0, a_1\rangle}{\left(a_0^2, \, a_1^2, \, ba_0 - (-1)^ma_0b, \, ba_1 - (-1)^ma_1b\right)}, \quad db = 0, \quad da_0 = -da_1 = (-1)^{m - 1}b,
	$$
	where the sign of $da_1$ comes from $\gamma$ not being oriented as the boundary of the region to its right, and the sign in $ba_i = (-1)^ma_ib$ from the sign in $T$.

	Observe that $\V_\bGamma$ has a basis given by the prefixes of the words $b^na_ia_{1 - i}a_i\ldots$ for $n \ge 0$ and $i = 0, 1$.
	We have $db^n = 0$, $d(b^na_i) = \pm b^{n + 1}$, and on a basis word with $k \ge 2$ letters $a$ we have
	$d(b^{n}a_i\ldots a_j) = \pm b^{n + 1}(a_i\ldots + \ldots a_j)$.
	Therefore,
	$$
		\B_{m - 1} = (\k\langle z\rangle, dz = 0, |z| = m - 1)\to \V_\bGamma, \; z \mapsto a_0 + a_1
	$$
	is a quasi-isomorphism and we are done.
\end{proof}

\begin{proof}[Proof of \cref{t: HMS}\,(3)]
	By \cref{l: B_m as topological Fukaya category} with $m = 2 - n$ we have
	$$
	\Wt(\S; \C_{2 - n}) \simeq \Tw \B_{1 - n} \simeq \Tw C_{-*}\Omega S^n,
	$$
	 which is quasi-equivalent to $\W(T^*S^n; \k)$ by \cite{A11}.
	The fibres of $\C_{2 - n}$ are $\Tw \B_{2 - n} \simeq \W(T^*S^{n - 1}; \k)$ for $n \ge 3$ by the same argument.
\end{proof}

\subsubsection{Lefschetz fibration and symplectic reduction}
\label{subsubs: Lefschetz fibration and symplectic reduction}
We now give the geometric heuristic.
Let $n \ge 3$ (the case $n = 2$ is covered in \cite{LS26} and does not quite fit in the discussion below since $\W(T^*S^1) \cong \Tw \k[z, z^{-1}] \not \cong \Tw \B_{0}$).
Identify $T^*S^n$ with the quadric
${V = \{x_0^2 + \ldots + x_n^2 = 1\} \subset \CC^{n + 1}}$.
The projection onto the first coordinate $\pi\colon V \to \CC_t$ is a Lefschetz fibration with regular fibres $V_t \cong T^*S^{n - 1}$ for $t \ne \pm 1$ and singular fibres $V_{\pm 1}$ which are $T^*S^{n - 1}$ with the zero section collapsed to a point.
It is not hard to see that
\begin{enumerate}
	\item the zero section $S^n_0$ of $T^*S^n$ is the matching sphere $V \cap \RR^{n + 1}$ above $[-1, 1]$;
	\item the set
	$$
		S^{n - 1}_t = \{t\} \times \sqrt{1 - t^2}S^{n - 1}
	$$
	is the zero section of $T^*S^{n - 1}$ under the identification $T^*S^{n - 1} \cong V_t$ for $t \ne \pm 1$ and for $t = \pm 1$ it is the collapsed zero section;
	\item for $r \in \RR$ and $y \in S^{n - 1} \subset \RR^n$ the cotangent fibre $T^*_yS^{n - 1}$ under the identification $T^*S^{n - 1} \cong V_{ir}$ is
	$$
		V_{ir, y} = \{(ir, x) \in V \colon \Re x \in \RR_{> 0}\cdot y\};
	$$
	\item for $y \in S^{n - 1}$ the cotangent fibre $T^*_{(0, y)}S^n$ under the identification $T^*S^n \cong V$ is
	$
		\bigcup_{r \in \RR} V_{ir, y}.
	$
\end{enumerate}
Fix the regular fibre $V_0$ and vanishing paths $[0, \pm 1]$.
The monodromy around $+1$ and $-1$ is the Dehn twist along the zero section $S_0^{n - 1}$.
On $\W(T^*S^{n - 1})$ the Dehn twist gives the spherical twist along the zero section.
Under $\W(T^*S^{n - 1})\cong \Tw \B_{2 - n}$ the spherical twist becomes $T$, which is the spherical twist along the twisted complex $\bS = (\B_{2 - n}[n - 1] \oplus \B_{2 - n}, z[1 - n])$ and $\eta$ fits into the triangle
$$
\Tw \B_{2 - n} (\bS, X) \otimes \bS \to X \xrightarrow{\eta_X} T(X).
$$

Also observe that the standard action of $\GL(n + 1, \CC)$ on $\CC^{n + 1}$ restricts to a Hamiltonian action of $\SO(n)$ on $V$ (the action is on the last $n$ coordinates) with moment map
$$
\mu\colon V \to \so(n)^* \cong \Lambda^2 \RR^n, \; \mu(x_0, \ldots, x_n) = \left(\Im \bar x_j x_k\right)_{1 \le j < k \le n}.
$$
This action preserves the fibres $V_t$ and
$$
\mu^{-1}(0) = V \cap \left(\CC \times (\CC \cdot S^{n - 1})\right) = \bigcup_{t \in \CC} S^{n - 1}_t.
$$
The action is transitive on $\mu^{-1}(0) \cap V_t$ with stabilisers conjugate to $\SO(n - 1)$ and fixes the node for $t = \pm 1$.
Hence set-theoretically we have $(\mu^{-1}(0) \cap \pi^{-1}(S))/\SO(n) = S$.

Moreover, a cotangent fibre $T^*_{(0, y)} S^n$ from (4) above intersects $\mu^{-1}(0)$ transversally and we get a bijection
$$
\pi\colon T^*_{(0, y)} S^n \cap \mu^{-1}(0) \cap \pi^{-1}(S) \to \Im \gamma
$$
where $\gamma$ is as in \cref{fig: weighted system on disc with two marked points}.
That is, $(\S; \C)$ can be thought of as a `singular symplectic reduction' of $(T^*S^n, \k)$ with respect to the action of $\SO(n)$ under which the cotangent fibre generating $\W(T^*S^n; \k)$ reduces to the arc brane $\LL$ generating $\Wt(\S; \C)$.

\section{Geometric preparations}
\label{s: Geometric preparations}
In this section we set up the geometric underpinning for Fukaya categories of surfaces with coefficients.

Let $S$, $Q = \Lambda \sqcup D$ be as in \cref{d: graded marked surface}, i.e.\ $S$ is a compact oriented surface with $\pd S' \ne \varnothing$ for each $S' \in \pi_0 S$, and $\Lambda \subset \pd S$ and $D \subset \Int S$ are finite sets.
We call such data $(S, Q)$ a \emph{marked surface}.

In \cref{subs: Surface} we define a compatible Liouville form on a marked surface; the choice of such form is contractible (see \cref{p: graded marked surfaces up to deformation}) and the diffeomorphism group of the marked surfaces retracts onto the exact symplectomorphism group (see \cref{p: Diff retracts on Symp}).
In \cref{subs: Curves} we first define Lagrangians in a Liouville marked surface, which are certain immersed curves with additional data, these will give rise to objects in the Fukaya category; then we describe the self-intersection information of a Lagrangian, consisting of its double points and Reeb chords between its ends, these will give morphisms in the Fukaya category.
We characterise Lagrangian structures on a given curve in terms of linear inequalities on actions and areas in \cref{p: criterion for existence of exact immersion in a given component,p: Lagrangian supported on exact curve with prescribed signs}.
In \cref{p: characterisation of isotopy} we characterise isotopies of Lagrangians topologically.
In \cref{subs: Polygons} we describe the polygons bounded by a Lagrangian, their counts will give rise to $A_\infty$-operations.
In \cref{subs: Grading} we add the grading information to surfaces and curves to make the Fukaya category $\ZZ$-graded.

\subsection{Surface}
\label{subs: Surface}
In \cref{subsubs: Liouville structure on a marked surface} we define Liouville structures on a marked surface $(S, Q)$, i.e.\ Liouville forms on $S$ compatible with $Q$.
In \cref{subsubs: Liouville marked structures on a given surface} we prove that the space of markings of a given surface with prescribed numbers of points on each boundary component and in the interior is aspherical, and that the space of Liouville structures on a given marked surface is non-empty and convex, see \cref{p: graded marked surfaces up to deformation}.
In \cref{subsubs: Hamiltonian diffeomorphisms} we define Hamiltonian diffeomorphisms of a Liouville marked surface and prove that they form the identity component of the group of exact symplectomorphisms in \cref{l: Symp_0 = Ham}.
In \cref{subsubs: Deforming diffeomorphisms to exact symplectomorphisms} we prove that exact symplectomorphisms of Liouville marked surfaces form a strong deformation retract of diffeomorphisms of the underlying marked surfaces in \cref{p: Diff retracts on Symp}; in particular, any two Liouville structures on the same marked surface are related by an exact symplectomorphism isotopic to the identity.
In \cref{subsubs: Contactisation} we describe the two contact manifolds associated to a Liouville marked surface $(S, Q, \lambda)$: the boundary $\pd S \setminus \Lambda$ and the contactisation $(S \setminus Q) \times \RR$.

\subsubsection{Liouville structure on a marked surface}
\label{subsubs: Liouville structure on a marked surface}
We describe a Liouville structure on $S$ compatible with the marking $Q$.

Recall that a $1$-form $\lambda$ on a $2n$-dimensional manifold $M$ is called a \emph{Liouville form} if
$(d \lambda)^{\wedge n}$
is a volume (i.e.\ nowhere vanishing) form on $M$.
The \emph{Liouville vector field} $Z = Z_\lambda$ is defined by
$d\lambda(Z, -) = \lambda$, and its flow is called the \emph{Liouville flow}.
Note that if $M$ admits a Liouville form, then it is orientable (by the volume form) and $M$ cannot be closed (by the Stokes theorem).

Let $\lambda$ be a Liouville form on the surface $S$.
We impose the following conditions on $\lambda$ along $\pd S$, $\Lambda$, and $D$.

The first set of conditions ensures that the Liouville flow $\Phi_Z^t$ is defined for all $t \le 0$ and that $Q$ does not move under the flow.
\begin{enumerate}
	\item\label{item: lambda vanishes on D} $\lambda$ (equivalently, $Z$) vanishes on $D$.
	\item\label{item: lambda vanishes on Lambda} $\Lambda = \left\{p \in \pd S \colon \lambda(p) = 0\right\}$.
	\item\label{item: Z points outwards} $Z$ points outwards along $\pd S \setminus \Lambda$.
\end{enumerate}
Note that there can be zeroes of $\lambda$ in $\Int S$ not lying in $D$, but all zeroes of $\lambda$ on $\pd S$ lie in $\Lambda$.
Observe that condition \cref{item: Z points outwards} is equivalent to $\res{\lambda}{\pd S \setminus \Lambda}$ being a volume form that induces the orientation of $\pd S \setminus \Lambda$ agreeing with the orientation of $\pd S$ as the boundary of $S$.

We impose one more condition to ensure that introducing $\lambda$ does not affect the mapping class group.
On one hand, the natural notion of automorphism of the marked surface $(S, Q)$ is an orientation-preserving diffeomorphism $\Phi$ of $S$ such that $\Phi(Q) = Q$.
Let $\Diff(S, Q)$ be the group of such maps $\Phi$.
On the other hand, an automorphism of the exact symplectic manifold $(S, \lambda)$ is a diffeomorphism $\Phi$ of $S$ such that $\Phi^*\lambda - \lambda$ is exact.
Let $\Symp_{\ex}(S, \lambda)$ be the group of such maps $\Phi$.
It is therefore natural to define
\begin{equation} \label{eq: def of Aut(Liouville marked surface)}
	\Symp_{\ex}(S, Q, \lambda) = \Diff(S, Q) \cap \Symp_{\ex}(S, \lambda).
\end{equation}
We would like the mapping class group $\pi_0 \Diff(S, Q) = \MCG(S, Q)$ to be equal to the mapping class group $\pi_0 \Symp_{\ex}(S, Q, \lambda) = \MCG(S, Q, \lambda)$, but there is an obstruction to that.
Namely, if $C_0$ and $C_1$ are two components of $\pd S$ with $|C_0 \cap Q| = |C_1 \cap Q|$, then there exists $\Phi \in \Diff(S, Q)$ mapping $C_0$ to $C_1$, but if $\Phi \in \Symp_{\ex}(S, \lambda)$, then $\int_{C_0} \lambda = \int_{C_1} \lambda$.
We require this obstruction to vanish.
\begin{enumerate}[resume]
	\item\label{item: integral of lambda over C} $\int_C \lambda = 1$ for each connected component $C$ of $\pd S$.
\end{enumerate}
This is enough for $\Symp_{\ex}(S, Q, \lambda)$ to be a strong deformation retract of $\Diff(S, Q)$, see \cref{p: Diff retracts on Symp}.

We summarise the discussion above into the following definition.
\begin{defn} \label{d: Liouville structure on marked surface}
	A \emph{Liouville structure} on the marked surface
	$(S, Q)$
	is a Liouville form $\lambda$ on $S$ satisfying conditions \cref{item: lambda vanishes on D,item: lambda vanishes on Lambda,item: Z points outwards,item: integral of lambda over C} above.
	$\S = (S, Q, \lambda)$ is then called a \emph{Liouville marked surface}.
\end{defn}

\begin{figure}[!htbp]
	\centering
	\begin{subfigure}[t]{0.36\linewidth}
		\centering
		\includegraphics[width=0.8\linewidth]{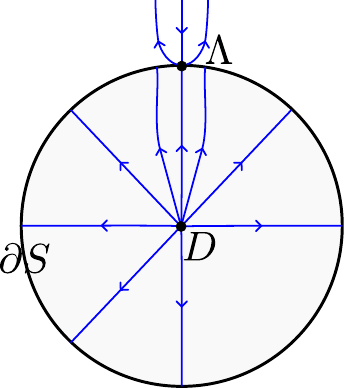}
		\caption{$S$ is a disc and $|\Lambda| = |D| = 1$.}
		\label{subfig: Liouville marked disc}
	\end{subfigure}%
	\hspace{0.04\linewidth}%
	\begin{subfigure}[t]{0.60\linewidth}
		\centering
		\includegraphics[width=\linewidth]{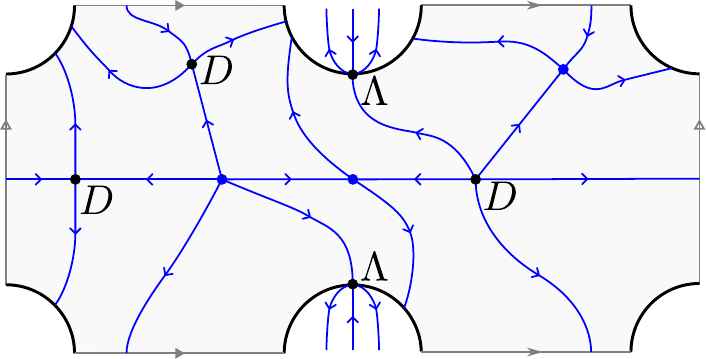}
		\caption{$S$ is a surface of genus $1$ glued along the \textcolor{Gray}{grey} lines marked with the same arrow, $|\pi_0 \pd S| = 2$, $|\Lambda| = 2$, and $|D| = 3$.}
		\label{subfig: Liouville marked torus}
	\end{subfigure}
	\caption{Examples of Liouville marked surfaces.}
	{\small The dots are the zeroes and the \textcolor{blue}{blue} lines are the integral lines of the Liouville vector field.
	The zeroes that belong to $\Lambda$ or $D$ are marked by the corresponding letter.}
	\label{fig: Liouville marked surfaces}
\end{figure}

In \cref{fig: Liouville marked surfaces} the zeroes and integral lines of the Liouville vector field are sketched for two Liouville marked surfaces.

If $S$ is orientable and $\pd S \ne \varnothing$ (which is necessary for the existence of a Liouville form on $S$), then the marked surface $(S, Q)$ admits a Liouville structure.
Moreover, the space of Liouville structures on $(S, Q)$ such that $d\lambda$ induces a given orientation of $S$ is convex in $\Omega^1(S)$, see \cref{p: graded marked surfaces up to deformation}.

\subsubsection{Liouville marked structures on a given surface}
\label{subsubs: Liouville marked structures on a given surface}
We study the space of Liouville marked surface structures on a connected surface $S$, the generalisation to disconnected case is straightforward.

For a non-negative integer $n$ and a function $\sigma$ from the set of connected components of $\pd S$ to $\ZZ_{\ge 0}$ let
\begin{equation} \label{eq: M_sigma n(S)}
	\M(S; \sigma, n) = \left\{
		\begin{aligned}
			\Lambda &\subset \pd S\\
			D &\subset \Int S
		\end{aligned}
		\;\colon
		\begin{aligned}
			&|\Lambda \cap C| = \sigma (C) \text{ for all } C \in \pi_0 \pd S\\
			&\text{and } |D| = n
		\end{aligned}
	\right\}
\end{equation}
be the set parametrising the choice of a marked surface structure $(S, Q)$ with $Q = \Lambda \sqcup D$ on $S$.

For a finite set $Q \subset S$ define
\begin{equation}
	\LS(S, Q) = \left\{
		\begin{aligned}
			\lambda \in \Omega^1(S)
		\end{aligned}
		\;\colon
		\begin{aligned}
			&\lambda \text{ is a Liouville structure on } (S, Q)\\
			&\text{and } d\lambda \text{ induces the orientation of } S
		\end{aligned}
	\right\}.
\end{equation}

\begin{prop} \label{p: graded marked surfaces up to deformation}
	Let $S$ be a compact connected oriented surface with non-empty boundary.
	Then
	\begin{enumerate}
		\item $\M(S; \sigma, n)$ is a smooth aspherical manifold of dimension $2n + \sum_{C \in \pi_0 \pd S} \sigma(C)$ for any $\sigma\colon \pi_0 \pd S \to \ZZ_{\ge 0}$ and $n \in \ZZ_{\ge 0}$.
		\item $\LS(S, Q)$ is a non-empty convex subset of $\Omega^1(S)$ for any finite set $Q \subset S$.
	\end{enumerate}

\end{prop}
\begin{proof}
	We prove (1).
	$\M = \M(S; \sigma, n)$ is the orbit space of the trivial-stabiliser locus of the action of the product of symmetric groups
	$$
	S_n \times \prod_{C \in \pi_0 \pd S} S_{\sigma(C)} \quad \text{on} \quad (\Int S) ^ n \times \prod_{C \in \pi_0 \pd S} C^{\sigma(C)}
	$$
	where each symmetric group acts on the corresponding factor by permuting entries.
	Since the trivial-stabiliser locus (i.e.\ the complement of the big diagonal) is open, $\M$ is a smooth manifold of dimension $2n + \sum_{C \in \pi_0 \pd S} \sigma(C)$.
	$\M$ is connected since both the complement of finitely many embedded closed intervals in $S$ and the configuration space of a fixed number of (unordered) points in $S^1$ are connected.
	That $\pi_k(\M) = 0$ for $k \ge 2$ follows by induction on $n + \sum_{C \in \pi_0 \pd S} \sigma(C)$ from the homotopy sequence of fibration and the fact that the complement of finitely many points in $S$ (resp.\ in a component of $\pd S$) is aspherical (resp.\ a disjoint union of contractible spaces).

	To prove (2) we first show $\LS = \LS(S, Q)$ is non-empty.
	Recall that $S$ admits a Stein structure, i.e.\ a complex structure $J$ together with a Morse subharmonic function $\phi\colon S \to \RR$ such that $\pd S$ is a regular level set of $\phi$ where $\phi$ attains its maximum.
	Moreover, up to reparametrisation of $\RR$ we can choose $\phi$ to be arbitrary as long as it does not have critical points of index $2$.
	In particular, we can ensure that $\phi$ attains a non-degenerate local minimum at each point of $D$, see Cieliebak--Eliashberg~\cite{CE12}*{Chapter 8} for details.
	The $1$-form $\lambda' = -d\phi \circ J$ is almost what we need, except the behaviour of $\lambda'$ along the boundary.
	By pushing $S$ inside its Liouville completion and rescaling we can ensure that $\int_C \lambda' = 1$ for all components $C$ of $\pd S$.
	Next, for each $C \in \pi_0 \pd S$ fix an orientation-preserving diffeomorphism
	$$
	\tau_C\colon C \to \RR/\ZZ \text{ such that } \res{\lambda'}{C} = \tau_C^*(dt)
	$$
	and pick a function $g_C\colon \RR/\ZZ \to \RR$ such that $g_C'(t) \ge -1$ for all $t \in \RR/\ZZ$ and the solution set of $g_C'(t) = -1$ is $\tau_C(C \cap \Lambda)$ (this can be done since we can ensure $\int_{\RR/\ZZ} g'_C(t)dt = 0$).
	Extend the collection of functions $\{g_C \circ \tau_C\}_{C \in \pi_0 \pd S}$ to a function $g\colon S \to \RR$ supported near $\pd S$ such that $(\lambda' + dg)(p) = 0$ for each $p \in \Lambda$ (this can be done since $\res{(\lambda' + dg)}{\pd S}(\Lambda) = 0$ and we can prescribe the derivative of $g$ along $\pd S$ in the normal direction).
	We obtain the form $\lambda = \lambda' + dg$ that lies in $\LS$.

	Now show that $\LS$ is convex.
	Observe that $\lambda \in \Omega(S)$ belongs to $\LS$ if and only if
	\begin{enumerate}
		\item $\lambda$ vanishes on $Q$;
		\item $\int_C \lambda = 1$ for each $C \in \pi_0 \pd S$;
		\item $d\lambda$ is a volume form on $S$ inducing the orientation of $S$;
		\item $\res{\lambda}{\pd S \setminus \Lambda}$ is a volume form on $\pd S \setminus \Lambda$ inducing the orientation of $\pd S$.
	\end{enumerate}
	(1) is linear in $\lambda$, (2) is affine in $\lambda$, and (3), (4) are convex in $\lambda$, hence the result.
\end{proof}

\subsubsection{Hamiltonian diffeomorphisms}
\label{subsubs: Hamiltonian diffeomorphisms}
Fix a Liouville marked surface $\S = (S, Q, \lambda)$.
Recall from \cref{eq: def of Aut(Liouville marked surface)} that $\Diff(\S)$ is the group of diffeomorphisms of $S$ that preserve $Q$ setwise and $\Symp_{\ex}(\S)$ is its subgroup consisting of exact symplectomorphisms.
We describe the identity component $\Symp_{\ex}(\S)_0$ of $\Symp_{\ex}(\S)$.

Let $H = \{H_s\colon S \to \RR\}_{s \in [0, 1]}$ be a smooth family such that $dH_s$ vanishes on $Q$ and $\res{H_s}{\pd S}$ is locally constant.
Then, since $d\lambda$ is non-degenerate, the equation $d\lambda(X_s, -) = dH_s$ defines a family of vector fields
$$
[0, 1] \to C^\infty(S, TS), s \mapsto X_s^H.
$$
As $X_s^H$ vanishes on $Q$ and is tangent to $\pd S$, it integrates to a family
$$
[0, 1] \to \Symp_{\ex}(\S), s \mapsto \Phi_s^H.
$$
We therefore obtain a map
$$
C^\infty_{[0, 1]}(\S, \RR) \to \Symp_{\ex}(\S), \; H = \{H_s\}_{s \in [0, 1]} \mapsto \Phi_1^H
$$
where $C^\infty_{[0, 1]}(\S, \RR)$ denotes the space of smooth families of functions $\{H_s\}_{s \in [0, 1]}$ described above.
The image of this map is a connected subgroup $\Ham(\S)$ of $\Symp_{\ex}(\S)$ consisting of \emph{Hamiltonian diffeomorphisms} of $\S$.

\begin{lemma} \label{l: Symp_0 = Ham}
	$\Ham(\S) = \Symp_{\ex}(\S)_0$.
\end{lemma}

\begin{proof}
	$\Ham(\S)$ is connected hence we only need to prove that the right hand side is contained in the left hand side.
	Let $\Phi_s \in \Symp_{\ex}(\S), s \in [0, 1]$ be a family of symplectomorphisms generated by a family of vector fields
	${X_s \in C^\infty(S, TS), s \in [0, 1]}$.
	Since a function on $S$ is determined up to a constant (on each component of $S$) by its differential, we get a family $g_s \in C^\infty(S, \RR), s \in [0, 1]$ such that
	$$
		\Phi_s^*\lambda - \lambda = dg_s \text{ for } s \in [0, 1].
	$$
	Differentiating with respect to $s$ and using Cartan's magic formula we obtain
	$$
		\Phi_s^*\left(d(\lambda(X_s)) + d\lambda(X_s, -)\right) = d\left(\frac{dg_s}{ds}\right).
	$$
	Therefore,
	$$
		d\lambda(X_s, -) = dH_s \text{ where } H_s = \frac{dg_s}{ds} \circ \Phi_s^{-1} - \lambda(X_s).
	$$
	Since $\Phi_s$ is a diffeomorphism of $S$ fixing $Q$ pointwise, $X_s$ vanishes on $Q$ and is tangent to $\pd S$.
	Therefore, $dH_s$ vanishes on $Q$ and $H_s$ is locally constant on $\pd S$.
	That is, $H = \{H_s\}_{s \in [0, 1]} \in C^\infty_{[0, 1]}(\S, \RR)$ and hence $\Phi_s$ is a path in $\Ham(\S)$.
\end{proof}

\subsubsection{Deforming diffeomorphisms to exact symplectomorphisms}
\label{subsubs: Deforming diffeomorphisms to exact symplectomorphisms}
The purpose of proving \cref{l: Symp_0 = Ham} was mostly to warm up before the following statement.

\begin{prop} \label{p: Diff retracts on Symp}
	Let $\S_i = (S_i, Q_i, \lambda_i)$, $i = 0, 1$ be two Liouville marked surfaces.
	Then
	$\Symp_{\ex}(\S_0, \S_1)$ is a strong deformation retract of $\Diff(\S_0, \S_1)$.
\end{prop}

$\Diff(\S_0, \S_1)$ in \cref{p: Diff retracts on Symp} is the space of orientation-preserving diffeomorphisms $\Phi$ from $S_0$ to $S_1$ such that $\Phi(Q_0) = Q_1$, with the $C^\infty$-topology.
$\Symp_{\ex}(\S_0, \S_1)$ is the subspace of $\Diff(\S_0, \S_1)$ consisting of exact symplectomorphisms, i.e.\ such that $\Phi^*\lambda_1 - \lambda_0 = dg$ for some $g \in C^\infty(S_0; \RR)$.

Applying \cref{p: Diff retracts on Symp} to $\S_0 = \S_1$ we get
\begin{cor} \label{cor: Diff retracts on Symp}
	If $\S$ is a Liouville marked surface, then $\Symp_{\ex}(\S)$ is a strong deformation retract of $\Diff(\S)$.
\end{cor}

Applying \cref{p: Diff retracts on Symp} to Liouville structures on the same marked surface we get
\begin{cor} \label{cor: Liouville structures on the same marked surface are symplectomorphic}
If $\lambda_0$ and $\lambda_1$ are Liouville structures on a marked surface $(S, Q)$, then there exists a smooth family $\Phi_s \in \Diff(S, Q)$ such that $\Phi_0 = \id_S$ and
$\Phi_1 \in \Symp_{\ex}((S, Q, \lambda_0), (S, Q, \lambda_1))$.
\end{cor}

\begin{remark}
	A related result by Bourgeois--Ekholm--Eliashberg~\cite{BEE12}*{Lemma 1.1} is that a symplectomorphism between Liouville manifolds with cylindrical ends can be deformed to an exact symplectomorphism.

	Note that for a Liouville domain $M$ the natural map $\pi_0 \Symp_{\ex}(M) \to \pi_0 \Symp(M)$ is not necessarily surjective.
	For instance, on the cylinder $M = [-1, 3]_r \times (\RR/\ZZ)_{\theta}$ with the Liouville form $\lambda = rd\theta$ the symplectomorphism $(r, \theta) \mapsto (2 - r, -\theta)$ swaps the components of the boundary, but there is no exact symplectomorphism $\Phi$ swapping the components since we would have
	$$
		0 = \int_{\{3\} \times \RR/\ZZ} (\Phi^*\lambda - \lambda) = \int_{-\{-1\} \times \RR/\ZZ}-d\theta - \int_{\{3\} \times \RR/\ZZ} 3d\theta = -2.
	$$

	On the other hand, a slight adaptation of our proof yields that for any Liouville domain $M$ the group $\Symp_{\ex}(M)$ is a strong deformation retract of the group
	$$
		\{\Phi \in \Symp(M) \colon \res{(\Phi^*\lambda - \lambda)}{\pd M} \text{ is exact}\}
	$$
	of symplectomorphisms exact on the boundary.
\end{remark}

To prove \cref{p: Diff retracts on Symp} we use a Moser argument and the following lemma.

\begin{lemma} \label{l: mu -> g_mu}
	Let $(S, Q)$ be a marked surface.
	Then there exists a continuous (in the $C^\infty$-topology) map $\mu \mapsto g_\mu$ from
	\begin{equation} \label{eq: tilde Omega^1}
		\tilde \Omega^1(S) = \left\{\mu \in \Omega^1(S) \colon
		\begin{aligned}
			&\int_C \mu = 0 \text{ for all } C \in \pi_0 \pd S
		\end{aligned}\right\}
	\end{equation}
	to $C^\infty(S, \RR)$
	such that $\res{\mu}{\pd S} = \res{dg_\mu}{\pd S}$ for any $\mu \in \tilde \Omega^1(S)$, $\mu(p) = d_pg_{\mu}$ for any $p \in Q$, and $dg_\mu = \mu$ whenever $\mu$ is exact.
\end{lemma}

\begin{proof}[Proof of \cref{p: Diff retracts on Symp}]
	For any $\Phi \in \Diff(\S_0, \S_1)$ we have ${\Phi^*\lambda_1 - \lambda_0 \in \tilde \Omega^1(S)}$ and hence a function $g_{\Phi^*\lambda_1 - \lambda_0}$ provided by \cref{l: mu -> g_mu} that depends on $\Phi$ continuously.
	The equation
	$$
		\Phi^*\lambda_1 - \lambda_0 + d\tilde \lambda_s(X_s, -) = dg_{\Phi^*\lambda_1 - \lambda_0} \text{ where } \tilde \lambda_s = s \Phi^*\lambda_1 + (1 - s)\lambda_0
	$$
	defines a family of vector fields
	$$
		[0, 1] \to C^{\infty}(S, TS), \; s \mapsto X^\Phi_s = X_s
	$$
	since $d\tilde\lambda_s$ is non-degenerate.
	As $d\tilde\lambda_s$ is skew-symmetric and $\res{d\tilde \lambda_s(X_s, -)}{\pd S}$ vanishes, $X_s$ is tangent to $\pd S$.
	Since $\Phi^*\lambda_1 - \lambda_0$ and hence $dg_{\Phi^*\lambda_1 - \lambda_0}$ vanish on $Q$, so does $X_s$.
	Therefore, $X_s$ integrates to a family of diffeomorphisms
	$$
		{[0, 1] \to \Diff(\S_0), \; s \mapsto \Psi^\Phi_s = \Psi_s}; \; \Psi_0 = \id_{S_0}.
	$$
	Integrating
	$$
		\frac{d}{ds} \Psi_s^*\tilde\lambda_s = \Psi_s^*\left(\frac{d}{ds}\tilde \lambda_s + d(\tilde \lambda_s(X_s)) + d\tilde \lambda_s(X_s, -)\right) = \Psi_s^*\big(d(\tilde \lambda_s(X_s)) + dg_{\Phi^*\lambda_1 - \lambda_0}\big)
	$$
	we obtain
	$$
		(\Phi \circ \Psi_1)^*\lambda_1 - \lambda_0 = \Psi_1^*\tilde \lambda_1 - \lambda_0 = d h \text{ where } h = \int_0^1 (\tilde \lambda_s(X_s) + g_{\Phi^*\lambda_1 - \lambda_0}) \circ \Psi_s ds.
	$$
	Therefore, $\Phi \circ \Psi_1 \in \Symp_{\ex}(\S_0, \S_1)$.
	If ${\Phi \in \Symp_{\ex}(\S_0, \S_1)}$, then $\Psi^\Phi_s = \id_{S_0}$ since in that case
	$dg_{\Phi^*\lambda_1 - \lambda_0} = \Phi^*\lambda_1 - \lambda_0$.
	Thus the map
	$$
		\Diff(\S_0, \S_1) \times [0, 1] \to \Diff(\S_0, \S_1), \quad (\Phi, s) \mapsto \Phi \circ \Psi^\Phi_s
	$$
	defines the desired strong deformation retraction.
\end{proof}

\begin{proof}[Proof of \cref{l: mu -> g_mu}]
	We perform the construction in three steps.
	The first step takes care of exactness using the Helmholtz decomposition, the second of the behaviour on $Q$ and $\pd S$ using a partition of unity argument, and the third combines the two.

	\emph{Step 1.}
	Fix a Riemannian metric on $S$ and recall that the Helmholtz decomposition on $S$ yields for any $\mu \in \Omega^1(S)$ an $L^2$-orthogonal splitting
	$
	\mu = \mu_{ex} + \mu_{co}
	$
	where $\mu_{ex}$ is exact and $\mu_{co}$ is co-closed with vanishing normal component, see e.g.\ Schwarz~\cite{S95}*{p. 87}.
	Since $g \in C^\infty(S, \RR)$ is determined by its differential up to a constant we obtain a continuous linear map
	$$
		\Omega^1(S) \to C^\infty(S, \RR), \ \mu \mapsto g^0_{\mu}
	$$
	such that $\mu = dg^0_{\mu}$ if and only if $\mu$ is exact.

	\emph{Step 2.}
	Note that the integral condition in \cref{eq: tilde Omega^1} is equivalent to $\res{\mu}{\pd S}$ being exact.
	Therefore, the behaviour of $\mu \in \tilde \Omega^1(S)$ on $Q$ and $\pd S$ can be recorded by a continuous linear map
	$$
		\mu \mapsto \left((\mu(q))_{q \in D}, (\mu(p))_{p \in \Lambda}, f_\mu\right)
	$$
	from $\tilde \Omega^1(S)$ to
	$$
		\T = \left\{
		\begin{aligned}
		&\left((\eta_q)_{q \in D}, (\eta_p)_{p \in \Lambda}, f\right) \colon\\
		&d_pf = \res{\eta_p}{T_p\pd S} \text{ for all } p \in \Lambda
		\end{aligned} \right\}
		\subset \prod_{q \in D} T^*_q S \times \prod_{p \in \Lambda} T^*_p S \times C^{\infty}(\pd S, \RR)
	$$
	satisfying $\res{\mu}{\pd S} = df_{\mu}$.

	A standard partition of unity argument gives a continuous linear map
	$$
		\T \to C^\infty(S, \RR), \ \bt \mapsto g_{\bt}
	$$
	such that $\res{g_{\bt}}{\pd S} = f$ and $d_pg_{\bt} = \eta_p$ for all $p \in Q$ for any
	$\bt = \left((\eta_q)_{q \in D}, (\eta_p)_{p \in \Lambda}, f\right)$.

	Composing the two maps above we get
	$$
		\tilde \Omega^1(S) \to C^{\infty}(S; \RR), \ \mu \mapsto g^1_{\mu}
	$$
	satisfying $\res{\mu}{\pd S} = \res{dg^1_{\mu}}{\pd S}$ and $\mu(p) = d_pg^1_{\mu}$ for any $p \in Q$.

	\emph{Step 3.} Define the continuous linear map
	$$
		\tilde \Omega^1(S) \to C^\infty(S, \RR), \ \mu \mapsto g_{\mu} = g^0_{\mu} + g^1_{\mu - dg^0_{\mu}}.
	$$
	Then
	$$
		d_p g_{\mu} = d_pg^0_{\mu} + (\mu - dg^0_{\mu}) = \mu \text{ for any } p \in Q
	$$
	and
	$\res{dg_\mu}{\pd S} = \res{dg^0_\mu}{\pd S} + \res{(\mu - dg^0_{\mu})}{\pd S} = \res{\mu}{\pd S}$.
	Moreover, if $\mu$ is exact, then $dg^0_\mu = \mu$ and hence $dg_\mu = \mu$.
	Thus $\mu \mapsto g_\mu$ satisfies all the requirements and we are done.
\end{proof}

\subsubsection{Contactisation}
\label{subsubs: Contactisation}
Liouville structures are closely connected to contact structures.
Recall that a $1$-form $\alpha$ on a $(2n + 1)$-dimensional manifold $N$ is called a \emph{contact form} if
$\alpha \wedge (d \alpha)^{\wedge n}$
is a volume form on $N$.
The pair $(N, \alpha)$ is then called a \emph{(co-oriented) contact manifold}.
The \emph{Reeb vector field} $R$ is defined by
$d\alpha(R, -) = 0$ and $\alpha(R) = 1$,
its integral curves are called \emph{Reeb trajectories}.

An $n$-dimensional submanifold of a $(2n + 1)$-dimensional contact manifold $(N, \alpha)$ is called a \emph{Legendrian} if it is everywhere tangent to the \emph{contact distribution} $\ker \alpha$.
A Reeb trajectory connecting two points of a Legendrian is called a \emph{Reeb chord}.
The \emph{length} of a Reeb chord $y\colon [0, l] \to N$ is $l = \int_{0}^l y^*\alpha$.
Unless stated otherwise, we only consider non-constant Reeb chords.

Given a Liouville marked surface
$\S = (S, Q, \lambda)$
we obtain two contact manifolds from it.
The first one is
$(\pd S \setminus \Lambda, \res{\lambda}{\pd S \setminus \Lambda})$,
which is a disjoint union of open intervals and circles.
The second is the contact $3$-manifold
\begin{equation} \label{eq: def of S^+}
	(S^+, \lambda^+) = \left((S \setminus Q) \times \RR_z, \lambda + dz\right),
\end{equation}
called the \emph{contactisation} of $\S$.

\subsection{Curves}
\label{subs: Curves}
Fix a Liouville marked surface $\S = (S, Q, \lambda)$.
We describe the curves in $S$ that will give rise to objects in the Fukaya category.

In \cref{subsubs: Lagrangians} we define Lagrangians in $\S$, which are immersed curves together with a potential, i.e.\ a Legendrian lift to the contactisation.
In \cref{subsubs: Double points and Reeb chords} we describe the self-intersection information of a Lagrangian, consisting of its double points and Reeb chords between its ends together with their actions, and introduce slightly degenerate Lagrangians, whose ends are allowed to collide.
In \cref{subsubs: Exact curves isotopic to a given curve,subsubs: Lagrangian structures on a given exact curve} we characterise in terms of linear inequalities on actions and areas when an isotopy class of immersed curves contains an exact curve and which Reeb signs and Reeb chords of positive action a Lagrangian supported on a given exact curve can have, see \cref{p: criterion for existence of exact immersion in a given component,p: Lagrangian supported on exact curve with prescribed signs}.
In \cref{subsubs: Lagrangian isotopy} we parametrise a neighbourhood of a Lagrangian in the space of Lagrangians and characterise isotopies of Lagrangians topologically.
In \cref{subsubs: Isotopy of slightly degenerate Lagrangians} we prove that any Lagrangian is isotopic through slightly degenerate Lagrangians to one without Reeb chords of positive action in \cref{l: sufficiently self-wrapped model for Lagrangian} and describe what happens to a Lagrangian when an isotopy passes through a slightly degenerate one; this will allow us to exchange Reeb chords for double points in \cref{p: isotopy induces isomorphism of obstruction algebras}.

\subsubsection{Lagrangians}
\label{subsubs: Lagrangians}
Let $I$ be a $1$-dimensional compact manifold, i.e.\ a finite disjoint union of circles and closed segments.
To avoid set-theoretic issues we tacitly only work with
\begin{equation} \label{eq: I m n}
	I = I_{m, n} = \bigsqcup_{i = 1}^m \RR/\ZZ \sqcup \bigsqcup_{i = 1}^n [0, 1] \text{ for } m, n \in \ZZ_{\ge 0}.
\end{equation}
As a byproduct, we obtain a nowhere vanishing $1$-form $dt$ on $I$, which, in particular, determines the orientation of $I$.

Consider an immersion $\gamma\colon I \to S$.
The first set of conditions we impose is that $\gamma$ is a generic properly immersed curve in $S \setminus Q$:
\begin{enumerate}
	\item\label{item: gamma avoid D and Lambda} $\gamma(I) \cap Q = \varnothing$.
	\item\label{item: boundary of gamma} $\gamma^{-1}(\pd S) = \pd I$ and $\gamma$ is transverse to $\pd S$;
	\item\label{item: ends of gamma are distinct} $\res{\gamma}{\pd I}$ is an embedding into $\pd S$;
	\item\label{item: gamma has at worst transcerse double points} all singularities of $\gamma$ are transverse double points.
\end{enumerate}

The curve $\gamma$ is called an \emph{exact immersion} if there exists a function ${P\colon I \to \RR}$ such that
\begin{enumerate}[resume]
	\item\label{item: exact Lagrangian equation} $\gamma^*\lambda + dP = 0$.
\end{enumerate}
Note that this implies that the \emph{Legendrian lift} $\gamma^+ = (\gamma, P)\colon I \to S^+$ is indeed Legendrian in the contactisation $(S^+, \lambda^+)$, i.e.\ $\gamma^+$ is tangent to the contact distribution $\ker \lambda^+$.

We impose further general position requirements, now on $P$:
\begin{enumerate}[resume]
	\item\label{item: P at double point} for any double point of $\gamma$ the values of $P$ at its two preimages are different;
	\item\label{item: P for a Reeb chord} for any Reeb chord
	$y\colon [0, l] \to \pd S \setminus \Lambda$
	connecting $\gamma(t_0) = y(0)$ and $\gamma(t_1) = y(l)$, the value $P(t_1)$ is different from $P(t_0) - l$.
\end{enumerate}
Note that \cref{item: P at double point} ensures the Legendrian lift $\gamma^+$ is embedded and that \cref{item: P for a Reeb chord} is in a sense a variant of \cref{item: P at double point}: for $t_0 \in \pd I$ if we attach the maximally extended Reeb trajectory $y\colon [0, l_0) \to \pd S \setminus \Lambda$ with $y(0) = t_0$ to $\gamma$ and extend the potential to satisfy (5), then the potential at a point $l \in [0, l_0)$ will be equal to $P(t_0) - \int_{0}^l y^*\lambda = P(t_0) - l$.

Combining all conditions together we obtain the following definition.

\begin{defn} \label{d: Lagrangian}
	A (exact immersed) \emph{Lagrangian} $L$ (with a fixed Legendrian lift) in the Liouville marked surface $\S$ is the data of
	\begin{itemize}
		\item a $1$-dimensional oriented compact manifold $I$, called the \emph{domain};
		\item an immersion $\gamma\colon I \to S$, called the \emph{underlying curve};
		\item a function $P\colon I \to \RR$, called the \emph{potential}
	\end{itemize}
	satisfying conditions \cref{item: boundary of gamma,item: gamma avoid D and Lambda,item: gamma has at worst transcerse double points,item: ends of gamma are distinct,item: exact Lagrangian equation,item: P at double point,item: P for a Reeb chord} above.
\end{defn}

Note that the potential is part of the data.
For a generic exact immersion $\gamma\colon I \to S$ the choice of potential $P$ is parametrised by the complement of a hyperplane arrangement in $\RR^{\pi_0 I}$, see \cref{subsubs: Transverse Lagrangians}.
Choices of the potential from different connected components of this complement can give non-equivalent obstruction algebras and hence different objects in the Fukaya category.

Observe that the orientation of $I$ was not used in \cref{d: Lagrangian}.
We only choose it for the sake of convenience, different orientations of $I$ will give isomorphic obstruction algebras and equivalent objects in the Fukaya category.

\subsubsection{Double points and Reeb chords}
\label{subsubs: Double points and Reeb chords}
Let
$L = (I, \gamma, P)$
be a Lagrangian in a Liouville marked surface $\S$.
We describe the self-intersection information of $L$.
It consists of double points of $\gamma$ and certain Reeb chords connecting the ends of $L$.
They will give rise to the obstruction algebra and to morphism spaces in the Fukaya category.

Let $X = X_L$ be the set of double points of $\gamma$.
Fix $x \in X$ and let
\begin{equation} \label{eq: x_+ and x_-}
	x_+, x_- \in I \text{ be such that } \gamma(x_+) = x = \gamma(x_-) \text{ and } P(x_+) > P(x_-).
\end{equation}
This makes sense by condition~\cref{item: P at double point}.
The \emph{action} of $x$ is
\begin{equation} \label{eq: def of A(x) for double point}
	A(x) = P(x_+) - P(x_-) > 0.
\end{equation}

Note that double points of $L$ are in bijection with Reeb chords of its Legendrian lift $L^+$ and that the action of a double point is equal to the length of the corresponding Reeb chord.

\begin{figure}[!htbp]
	\centering
	\begin{subfigure}[t]{0.48\linewidth}
		\centering
		\includegraphics[width=0.9\linewidth]{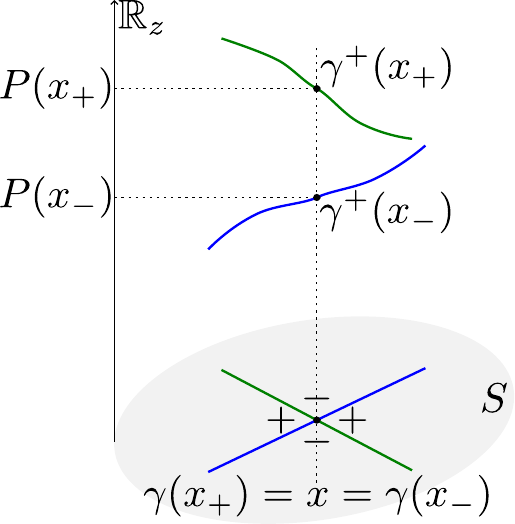}
		\caption{Double point $x$ of $L$ with Reeb signs indicated.}
		\label{subfig: double point}
	\end{subfigure}%
	\hfill%
	\begin{subfigure}[t]{0.48\linewidth}
		\centering
		\includegraphics[width=\linewidth]{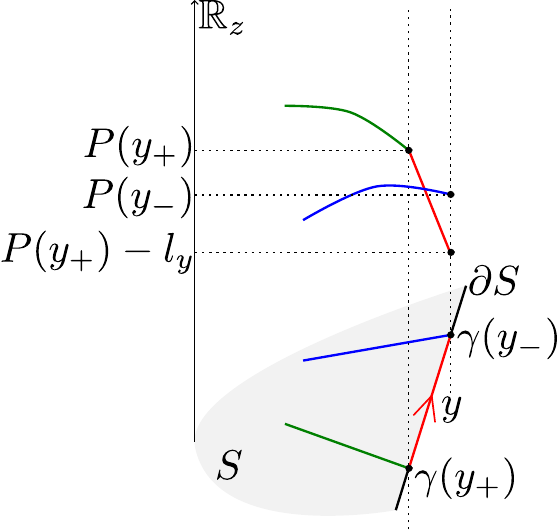}
		\caption{Reeb chord $y$ between two ends of $L$ and its Legendrian lift.}
		\label{subfig: Reeb chord}
	\end{subfigure}
	\caption{Lifts of a double point and a Reeb chord.}
	\label{fig: double point and Reeb chord}
	{\small Lagrangian $L = (I, \gamma, P)$ in $S$ and its Legendrian lift $\gamma^+ = (\gamma, P)$ in the contactisation $S^+$ around a double point and a Reeb chord.}
\end{figure}

The image of a small neighbourhood of $x_+$ (resp.\ $x_-$) under $\gamma$ is called the \emph{overstrand (resp.\ understrand)} at $x$.
To arrange the counts of discs later it is convenient to encode the over- and understrand information at $x$ by \emph{Reeb signs}.
These are the assignments of $+$ or $-$ to each of the four quadrants cut out by the tangent lines $\dot \gamma(x_+)$ and $\dot \gamma(x_-)$ from $T_x S$.
Namely, we assign $+$ (resp.\ $-$) to the two quadrants that lie to the left of the curve that first follows the understrand (resp.\ overstrand) and then jumps to the overstrand (resp.\ understrand) at $x$, see \cref{subfig: double point}.
The quadrants with Reeb sign $+$ (resp.\ $-$) are called \emph{positive} (resp.\ \emph{negative}).

Note that $\res{\gamma}{\pd I}$ gives a Legendrian submanifold in the contact manifold
$(\pd S \setminus \Lambda, \res{\lambda}{\pd S \setminus \Lambda})$.
Let $Y = Y_L$ be the set of Reeb chords $y\colon [0, l_y] \to \pd S \setminus \Lambda$ that connect two ends of $L$.
Consider a Reeb chord $y \in Y$.
Let
\begin{equation} \label{eq: y_+ and y_-}
	y_+, y_- \in \pd I \text{ be such that } \gamma(y_+) = y(0) \text{ and } \gamma(y_-) = y(l_y).
\end{equation}
The \emph{action} of $y$ is
\begin{equation} \label{eq: def of A(y)}
	A(y) = P(y_+) - l_y - P(y_-).
\end{equation}
Note that $A(y)$ is always non-zero by condition~\cref{item: P for a Reeb chord}.

The definitions of action for double points and Reeb chords are compatible in the following sense.
Given a Reeb chord $y \in Y$, consider the Reeb trajectory $\tilde y\colon [0, l) \to \pd S \setminus \Lambda$ that maximally extends $y$.
Extend $\gamma\colon I \to S$ to a piece-wise smooth curve
$\tilde \gamma = \gamma \cup \tilde y \colon I \underset{y_+ = 0}{\cup} [0, l) \to S$
and $P\colon I \to \RR$ to a continuous function $\tilde P = I \underset{y_+ = 0}{\cup} [0, l) \to \RR$ satisfying \cref{item: exact Lagrangian equation} everywhere except $y_+ = 0$.
Then $\tilde \gamma(y_-) = \tilde \gamma(l_y)$ is a self-intersection point of $\tilde \gamma$ and $A(y) = \tilde P(l_y) - \tilde P(y_-)$, see \cref{subfig: Reeb chord}.

Allow for a moment $\res{\gamma}{\pd I}$ not to be an embedding and let $z_+, z_- \in \pd I$ be such that $\gamma(z_+) = \gamma(z_-) = z \in \pd S$ and $P(z_+) > P(z_-)$.
Then $z$ can be considered both as a self-intersection point of $L$ and as a constant Reeb chord between the ends of $L$.
The two definitions of the action clearly agree.
Motivated by this observation, we give the following definition.

\begin{defn} \label{d: slightly degenerate Lagrangian}
	A \emph{slightly degenerate Lagrangian} $L$ in the Liouville marked surface $\S$ is the same thing as a Lagrangian from \cref{d: Lagrangian} except that $\res{\gamma}{\pd I}$ is not required to be an embedding and at each double point $z \in \pd S$ the unique quadrant of $z$ contained in $S$ is negative.
\end{defn}

Given a slightly degenerate Lagrangian $L$ we denote by $X_L$ the set of double points of $L$ in $\Int S$, by $Y_L$ the set of Reeb chords of positive length, and by $XY_{L}$ the set of double points of $L$ in $\pd S$, which we call \emph{point-chords} of $L$.

Let
$y^i\colon [0, l_i] \to \pd S \setminus \Lambda, \, i = 0, 1$
be two Reeb chords connecting the ends of $L$ such that $y^0(l_0) = y^1(0)$.
Their \emph{concatenation} is the Reeb chord
\begin{equation}
	y^1y^0\colon [0, l_0 + l_1] \to \pd S \setminus \Lambda, \;
	s \mapsto \begin{cases}
		y^0(s), & \text{if } s \in [0, l_0];\\
		y^1(s - l_0), & \text{if } s \in [l_0, l_0 + l_1].
	\end{cases}
\end{equation}
The action is additive with respect to concatenation.

We make the following simple observation
\begin{lemma} \label{l: Y^>C is finite}
	Let $B$ be a real number.
	Then
	\begin{enumerate}
		\item The set
		$$
			Y^{>B} = \{y \in Y \colon A(y) > B\}
		$$
		is finite.
		\item The set
		$$
			Y^{<B} = \{y \in Y \colon A(y) < B\}
		$$
		is finite if and only if for any $t \in \pd I$ there is a stop on the component $C$ of $\pd S$ with $\gamma(t) \in C$.
	\end{enumerate}
\end{lemma}
\begin{proof}
	$\pd I$ is finite and for fixed
	$y_+, y_- \in \pd I$
	the set of Reeb chords from $\gamma(y_+)$ to $\gamma(y_-)$ either contains at most one element or is equal to $\{yz^k \mid k \in \ZZ_{\ge 0}\}$ where $y$ is the shortest Reeb chord from $\gamma(y_+)$ to $\gamma(y_-)$ and $z$ is the Reeb chord from $\gamma(y_+)$ to itself that goes around the component of $\pd S$ once.
	As $A(yz^{k+1}) - A(yz^k) = -l_z < 0$,
	only finitely many of these chords will have action larger than $B$.
	On the other hand, if $\gamma(t) \in C$ for some $t \in \pd I$ and a component $C \in \pi_0 \pd S$ with $C \cap \Lambda = \varnothing$, then $A(z^k) < B$ for $k \gg 0$ where $z$ is the Reeb chord from $\gamma(t)$ to itself going around $C$ once.
\end{proof}

The sets $X$ and $XY$ are finite by compactness and so is the set $Y^{>0}$ of Reeb chords of positive action by \cref{l: Y^>C is finite}.
We denote $Z = Z_L = X \cup Y^{>0} \cup XY$.

\begin{figure}[!htbp]
	\includegraphics[width=\linewidth]{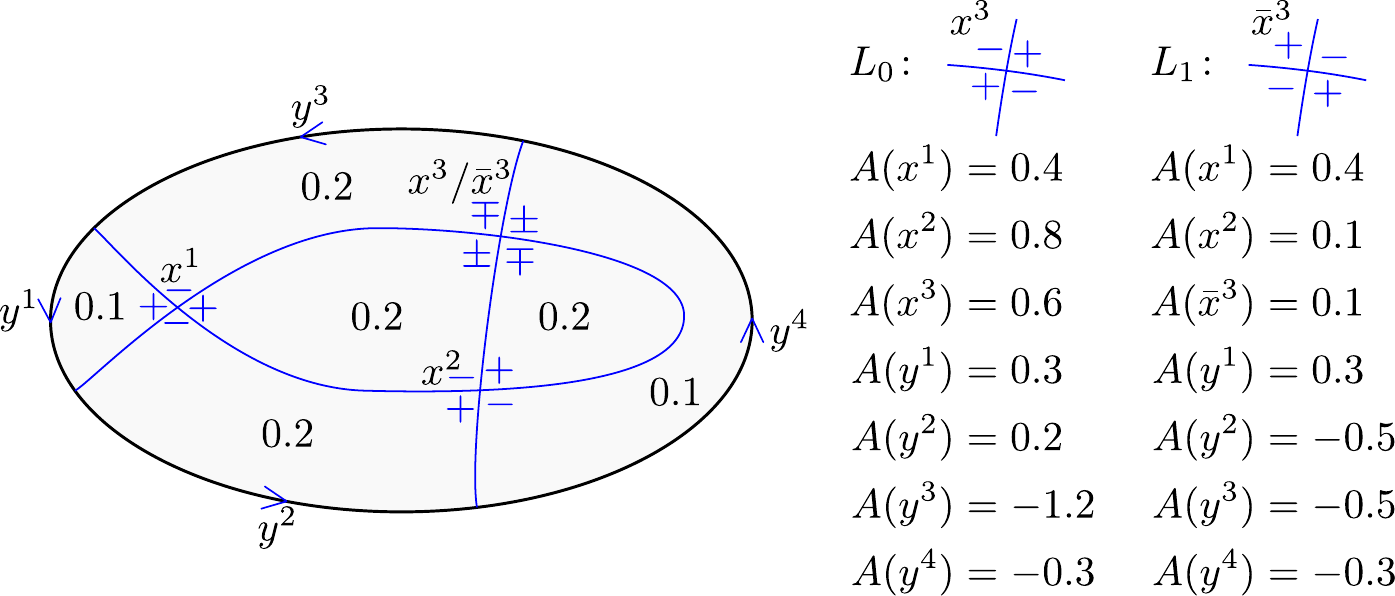}
	\caption{Two Lagrangians with the same underlying curve.}
	{\small The \textcolor{blue}{blue} curve is the underlying curve of both $L_0$ and $L_1$, the \textcolor{blue}{blue} signs at each double point are Reeb signs, the numbers inside the regions cut out from the disc by the curve are their areas.}
	\label{fig: Lagrangian with actions}
\end{figure}

\begin{example} \label{ex: two Lagrangians with the same underlying curve}
In \cref{fig: Lagrangian with actions} two Lagrangians $L_0$ and $L_1$ in a disc with some Liouville form $\lambda$ are depicted.
The domain of both Lagrangians is the disjoint union of two segments and their underlying curves are the same.
$L_1$ is obtained from $L_0$ by changing the potential on one of the segments by a constant.
We have $Z_{L_0} = \{x^1, x^2, x^3, y^1, y^2, y^2y^1, y^4y^2y^1\}$ and $Z_{L_1} = \{x^1, x^2, \bar x^3, y^1\}$.

This example is continued in \cref{ex: polygons bounded by two Lagrangians with the same underlying curve}.
\end{example}

The Liouville form in \cref{fig: Lagrangian with actions} is only represented by areas (i.e.\ integrals of $d\lambda$) and actions.
This numerical information is enough to determine whether a Lagrangian with these areas and actions exists as we explain in \cref{p: criterion for existence of exact immersion in a given component,p: Lagrangian supported on exact curve with prescribed signs} below.
Moreover, the Reeb signs and the signs of the actions of Reeb chords between the ends of a Lagrangian $L$ determine it uniquely up to isotopy, see \cref{p: characterisation of isotopy}.

\subsubsection{Exact curves isotopic to a given curve}
\label{subsubs: Exact curves isotopic to a given curve}

Fix a marked surface $(S, Q)$.
Let $I$ be a compact $1$-dimensional manifold and $\gamma\colon I \to S \setminus Q$ a generic enough properly immersed curve, i.e.\ $\gamma$ satisfies conditions \cref{item: boundary of gamma,item: gamma avoid D and Lambda,item: gamma has at worst transcerse double points,item: ends of gamma are distinct} from \cref{d: Lagrangian}.
Let $\Imm_I(S, Q)$ be the space of all such immersions with $C^\infty$-topology.
We say that two curves in $\Imm_I(\S)$ are \emph{isotopic} if they lie in the same component of $\Imm_I(\S)$.

Note that $\Diff(\S)$ acts on $\Imm_I(\S)$ on the left by post-composition and $\Diff(I)$ acts on the right by pre-composition.
Make the following simple observation.

\begin{lemma} \label{l: components of the space of immersions are orbits}
	$$
		\pi_0 \Imm_I(\S) = \Diff(\S)_0\backslash \Imm_I(\S) / \Diff(I)_0.
	$$
\end{lemma}
\begin{proof}
	On one hand, each orbit is contained in one component since $\Diff(\S)_0$ and $\Diff(I)_0$ are connected.
	On the other hand, if we have an isotopy $\gamma_s \in \Imm_I(S)$, $s \in [0, 1]$, then we can first make the double points fixed throughout the isotopy by an isotopy $\Phi_s \in \Diff(\S)$, $s \in [0, 1]$, $\Phi_0 = \id_S$, then make the preimage of each double point in $I$ fixed throughout the isotopy by an isotopy ${\psi_s \in \Diff(I)}$, ${s \in [0, 1]}$, ${\psi_0 = \id_I}$, and then extend the time-dependent vector field
	${\res{\frac{d\gamma_s(t)}{ds}}{s = s_0} \in C^\infty(I, \gamma^*TS)}$,
	$s_0 \in [0, 1]$ to a time-dependent vector field on $S$ that integrates to an isotopy $\tilde \Phi_s \in \Diff(\S)$, $s \in [0, 1]$, $\tilde \Phi_0 = \id_S$.
\end{proof}

Now let $\lambda$ be a Liouville structure on $\S$.
We give a criterion for the existence of an exact (i.e.\ such that ${\gamma'}^*\lambda$ is exact) immersion $\gamma'$ isotopic to a given $\gamma \in \Imm_I(\S)$, cf.\ Chekanov \cite{C02}*{Lemma 11.4}.

Note that if $\gamma \in \Imm_I(\S)$ is exact, then the restriction of $\gamma$ to any union of connected components of $I$ is exact and that we can add a number of embedded segments to $\gamma$ and it will remain exact.
It therefore suffices to study the case when ${I = I_0 \sqcup I_1 \sqcup I_2}$ where $I_0$ is a disjoint union of circles, $I_1$ and $I_2$ are disjoint unions of segments, and $\res{\gamma}{I_2}$ is an embedded curve such that $S \setminus \gamma(I_2)$ is contractible.
We call such $\gamma$ \emph{totally separating}.
Observe that in a given connected component $\C$ of $\Imm_I(\S)$ either all curves are totally separating or none are.
Note that $I_2$ has $l = 2g + b - 1$ components where $g$ is the genus and $b$ is the number of boundary components of $S$.
In particular, $I_2 = \varnothing$ if and only if $S$ is a disc.

Let $\R$ be the set of regions cut out by $\gamma$ from $S$.
Let $\Sigma\colon \RR^{\R} \to \RR$ be the linear map $(\alpha_R)_{R \in \R} \mapsto \sum_{R \in \R} \alpha_R$.
Note that if $(I, \gamma, P)$ is a Lagrangian, then
$$
\Sigma\left(\int_R d\lambda \right)_{R \in \R} = b
$$
by the Stokes theorem and condition \cref{item: integral of lambda over C} in \cref{d: Liouville structure on marked surface}.

For each $t_+ \in \pd I_2$ let $y_{t_+}\colon [0, 1] \to \pd S$ be a path on the boundary that starts at $\gamma(t_+)$ and ends at $\gamma(t_-)$ where $t_- \in \pd I_2$ is the point following $t_+$ in the order of the orientation of $\pd S$ ($y_{t_+}$ can pass through $\Lambda$ and $t_-$ can be equal to $t_+$).
Let $Y_2$ be the set of such paths, it has order $2l$.
For $y \in Y_2$ denote by $y_\pm \in \pd I$ the points satisfying $\gamma(y_+) = y(0)$ and $\gamma(y_-) = y(1)$ as in \cref{subsubs: Double points and Reeb chords}.
If $(I_2, \res{\gamma}{I_2}, P_2)$ is a Lagrangian, then we can define the actions $A(y) = P_2(y_+) - P_2(y_-) - \int_0^1 y^*\lambda$ for $y \in Y_2$.

Let $B\colon \RR^{Y_2} \to \RR^{\pi_0 \pd S}$ be the linear map
\begin{equation} \label{eq: def of B}
	(A_y)_{y \in Y_2} \mapsto \left(\sum_{\substack{y \in Y_2\colon \Im y \subset C}} A_y \right)_{C \in \pi_0 \pd S}.
\end{equation}
Note that if $(I_2, \res{\gamma}{I_2}, P_2)$ is a Lagrangian, then
${B(A(y))_{y \in Y_2} = (-1, \ldots, -1)}$
by condition \cref{item: integral of lambda over C} in \cref{d: Liouville structure on marked surface}.
If $Y_2 = \varnothing$ we define $B\colon \{0\} \to \RR$ by $B(0) = -1$.

Fix an orientation of $I$.
Then we have a linear map $E\colon \RR^{\R} \oplus \RR^{Y_2} \to \RR^{\pi_0 I_0}$ that, given the areas $(\alpha_R)_{R \in \R}$ of the regions and the actions $(A_y)_{y \in Y_2}$ of the paths computes the integrals $(\int_C \gamma^*\lambda)_{C \in \pi_0 I_0}$ by the Stokes theorem.
Note that the definitions of $\Sigma$, $B$, and $E$ only depends on the component of $\gamma$ in $\Imm_I(\S)$.

\begin{prop} \label{p: criterion for existence of exact immersion in a given component}
	Let $\C$ be a component of $\Imm_I(\S)$ consisting of totally separating immersions.
	Then there exists an exact immersion in $\C$ if and only if there exist $\alpha_R > 0, \, R \in \R$ and $A_y \in \RR, \, y \in Y_2$ such that
	$$
		\Sigma(\alpha_R)_{R \in \R} = b,\ B(A_y)_{y \in Y_2} = (-1, \ldots, -1), \text{ and } E\left((\alpha_R)_{R \in \R}, (A_y)_{y \in Y_2}\right) = 0.
	$$
\end{prop}
\begin{proof}
	The `only if' part is clear from the discussion above.
	Assume we have $(\alpha_R)_{R \in \R}$ and $(A_y)_{y \in Y_2}$ as in the statement.
	We want to find an immersion $\gamma \in \C$ and a function $P\colon I \to \RR$ such that $\gamma^*\lambda + dP = 0$.

	Pick an immersion $\gamma_0 \in \C$.
	By the condition on $B$ we can pick $2l$ numbers $P(t), \, t \in \pd I_2$ such that $A(y) = P(y_+) - P(y_-) - \int_{0}^1 y^*\lambda = A_y$ for each $y \in Y_2$.
	Recall that $\res{\gamma_0}{I_2}$ cuts $S$ to form one polygon.
	By deforming the interiors of two sides of this polygon corresponding to a component $C$ of $I$ we can achieve any value of $\int_C \gamma^*\lambda$.
	Since this deformation can be realised by an isotopy in $\Diff(\S)$, we can obtain a new immersion $\gamma_1 \in \C$ such that $P(t_1) - P(t_0) + \int_{t_0}^{t_1} \gamma_1^*\lambda = 0$ for each component $C \in \pi_0 I_2$ where $\pd C = \{t_0, t_1\}$.
	Hence we extend $P$ to $I_2$ to satisfy $\res{(\gamma^*\lambda + dP)}{I_2} = 0$.

	We now proceed by removing the regions $R$ from $\R$ one by one after ensuring that $\int_R d\lambda = \alpha_R$.
	We start with the polygon $\T_1 = (T, D)$ cut out by $\res{\gamma_1}{I_2}$ from $S$, the curve $\tilde \gamma_1 = \res{\gamma_1}{I \setminus I_2}$, and $\R_1 = \R$.
	At each step we have a polygon $\T_k = (T_k \subset T, D_k = D \cap T_k)$ and a curve $\tilde \gamma_k$ cutting $T_k$ into regions $\R_k \subset \R$ such that
	$$
	\sum_{R \in \R_k} \int_R d\lambda = \sum_{R \in \R_k} \alpha_R.
	$$
	We take a region $R \in \R_k$ adjacent to $\pd T_k$ such that $T_k \setminus R$ is connected.
	If $\pd R \subset \pd T_k$, then $R = T_k$, $\int_R d\lambda = \alpha_R$ and we are done.
	Otherwise, we can find an isotopy
	$$
	\Phi_s \in \Diff(\T) = \{\Phi \in \Diff(T) \colon \res{\Phi}{\pd S} = \id,\ \Phi(D_k) = D_k\},\ s \in [0, 1],\ \Phi_0 = \id_T
	$$
	such that $\int_{\Phi_1(R)} d\lambda = \alpha_R$.
	We set $T_{k + 1}$ to be the closure (in $T_k$) of $T_k \setminus \Phi_1(R)$ and $\tilde \gamma_{k + 1}$ to be the part of $\Phi_1 \circ \tilde \gamma_k$ that lies inside $T_{k + 1}$.

	After the procedure above we obtain an immersion $\gamma \in \C$ such that $\int_R d\lambda = \alpha_R$ for each $R \in \R$ and $A(y) = A_y$ for each $y \in Y_2$.
	The vanishing of $E$ now implies that $\gamma$ is exact and we are done.
\end{proof}

\subsubsection{Lagrangian structures on a given exact curve}
\label{subsubs: Lagrangian structures on a given exact curve}
Let $\gamma \in \Imm_I(\S)$ be an exact immersed curve.
We study Lagrangians in $\S$ with underlying curve $\gamma$.

Let $L = (I, \gamma, P)$ be such Lagrangian.

Given $L$ we obtain Reeb signs at each double point of $\gamma$ (recall that this is the same information as specifying under- and over-crossings) and the set $Y^{>0}_L$ of Reeb chords of positive action.

Let $\gamma \in \Imm_I(\S)$ be a totally separating exact curve with Reeb signs at each double point $x \in X$ and a finite set $Y^{>0} \subset Y$ of Reeb chords between the ends of $\gamma$ specified.
We give a criterion for when there exists a Lagrangian supported on $\gamma$ inducing this information, cf.\ Chekanov \cite{C02}*{Proposition 11.2}.

Let $Y_{short}$ be the set of paths $y\colon [0, 1] \to \pd S$ that follow the orientation of $\pd S$, start at $\gamma(y_+)$, end at $\gamma(y_-)$ for some $y_{\pm} \in \pd I$, and do not have the ends of $\gamma$ in the interior.
The number of paths in $Y_{short}$ is twice the number of segment connected components of $I$.

Recall that $\R$ is the set of regions cut out by $\gamma$ from $S$.
We have a linear map $\Area\colon \RR^{X \cup Y_{short}} \to \RR^{\R}$ that computes the area of the region $R$ by the Stokes theorem \cref{eq: area formula for an immersed polygon}.
We also have the map $F\colon \RR^{Y_{short}} \to \RR^{Y}$ computing the action of each $y \in Y$ from the actions on $Y_{short}$ by the additivity with respect to concatenation and the map
$\tilde B\colon \RR^{Y_{short}} \to H^1(\gamma(I_2) \cup \pd S; \RR)$
that, given a closed curve $\sigma$ in $\gamma(I_2) \cup \pd S$ thought of as a $1$-cycle, computes $\sum_{y \in Y_{short}} \deg_\sigma(y)$ where $\deg_\sigma(y)$ is the signed number of times $\sigma$ covers $y$.

\begin{prop} \label{p: Lagrangian supported on exact curve with prescribed signs}
	Let $\gamma \in \Imm_I(\S)$ be a totally separating exact curve with Reeb signs at each double point $x \in X$ and a finite set $Y^{>0} \subset Y$ of Reeb chords between the ends of $\gamma$ specified.
	Then there exists a potential $P\colon I \to \RR$ such that the Lagrangian $L = (I, \gamma, P)$ induces the same Reeb signs and $Y^{>0}_L = Y^{>0}$ if and only if there exists $A_x > 0, \, x \in X$ and $A_y \in \RR, \, y \in Y_{short}$ such that
	\begin{enumerate}
		\item $\left(F(A_y)_{y \in Y_{short}}\right)_y > 0$ for $y \in Y^{>0}$;
		\item $\left(F(A_y)_{y \in Y_{short}}\right)_y < 0$ for $y \in Y \setminus Y^{>0}$;
		\item
		$\Area\left((A_x)_{x \in X}, (A_y)_{y \in Y_{short}}\right) = \left(\int_R d\lambda\right)_{R \in \R}$;
		\item
		$\tilde B(A_y)_{y \in Y_{short}} = \res{\lambda}{\gamma(I_2) \cup \pd S}$
		in $H^1(\gamma(I_2) \cup \pd S; \RR)$.
	\end{enumerate}
\end{prop}

Note that it is enough to check condition (2) on the shortest Reeb chord $yz^k \in Y \setminus Y^{>0}$ from the infinite set $\{yz^k \mid k \ge 0\}$ where $z$ goes around a component of $\pd S$ once.
Hence there are finitely many conditions to check in \cref{p: Lagrangian supported on exact curve with prescribed signs}.

\begin{proof}
	The `only if' part is clear.
	To prove the `if' part we construct a function $P\colon I \to \RR$ such that $\gamma^*\lambda + dP = 0$, $A(x) = P(x_+) - P(x_-) = A_x$ for each $x \in X$, and $A(y) = P(y_+) - P(y_-) - \int_0^1 y^*\lambda = A_y$ for each $y \in Y_{short}$.
	Then $L = (I, \gamma, P)$ is a Lagrangian since $A(z) \ne 0$ for each $z \in X \cup Y$, it induces the correct Reeb signs since $A_x > 0$ for $x \in X$ and $Y^{>0}_L = Y^{>0}$ by (1) and (2).

	Define $P(t) = 0$ for some $t \in I$ and spread using $\gamma^*\lambda + dP = 0$ and $A_z$ for $z \in X \cup Y_{short}$ to all $I$ (if $\Im \gamma \cup \pd S$ is not connected, then do this separately for each component).
	No conflict arises by conditions (3) and (4) since any closed curve in $\Im \gamma \cup \pd S$ is a linear combination (as a $1$-chain) of boundaries of the regions from $Y$ and closed curves in $\gamma(I_2) \cup \pd S$.
\end{proof}

\subsubsection{Isotopy of Lagrangians}
\label{subsubs: Lagrangian isotopy}
For a compact $1$-dimensional manifold $I$ denote by $\Lag_I(\S)$ the space of Lagrangians (in the sense of \cref{d: Lagrangian}) in $\S$ with the $C^\infty$-topology.
Two Lagrangians with the same domain $I$ are \emph{isotopic} if they lie in the same component of $\Lag_I(\S)$.

The group $\Symp_{\ex}(\S)$ acts on $\Lag_I(\S)$ by
$$
\Phi_*(I, \gamma, P) = (I, \Phi \circ \gamma, P - \gamma^*g)
$$
where $g \in C^{\infty}(S, \RR)$ is such that $\Phi^*\lambda - \lambda = dg$ and $g(q_0) = 0$ for some fixed basepoint $q_0 \in S$.
Note that $\Phi_*L$ is isotopic to $L$ for $\Phi \in \Ham(\S)$, but the converse is not true (even after allowing vertical shifts) since the actions of double points and Reeb chords are the same for the corresponding points or chords of $L$ and $\Phi_*L$, but can change under an isotopy.

Fix $L_0 = (I, \gamma_0, P_0)$.
We describe a neighbourhood of $L_0$ in $\Lag_I(\S)$.
By the Weinstein neighbourhood theorem there exists a local diffeomorphism $\Phi\colon I \times (-\eps, \eps) \to S \setminus Q$ such that $\Phi(t, 0) = \gamma_0(t)$, $\Phi^*d\lambda = dt \wedge ds$ and $\Phi^{-1}(\pd S) = \pd I \times (-\eps, \eps)$.
As $\Phi^*\lambda + sdt$ is closed and restricts to an exact form $-dP_0$ on $I \times \{0\}$, there is a function $\tilde P\colon I \times (-\eps, \eps)$ such that $\tilde P(t, 0) = P_0(t)$ and $\Phi^*\lambda + sdt + d\tilde P = 0$.

Therefore, for $F \in C^\infty(I, \RR)$ with $F'$ close to $0$ we obtain a Lagrangian $(I, \gamma_F, P_F)$ where
$$
\gamma_F(t) = \Phi(t, F'(t)) \text{ and } P_F(t) = F(t) + \tilde P(t, F'(t)) \text{ for } t \in I.
$$

\begin{lemma} \label{l: neighbourhood of L in Lag_I(S)}
	There exists a neighbourhood $\U$ of $\{0\} \times \{\id\}$ in $C^\infty(I, \RR) \times \Diff_0(I)$ such that the map
	$$
		\U \to \Lag_I(\S), \; (F, \ph) \mapsto L(F, \ph) = (I, \gamma_F \circ \ph, P_F \circ \ph)
	$$
	is a homeomorphism onto a neighbourhood of $L_0$ in $\Lag_I(\S)$.
\end{lemma}
\begin{proof}
	We construct the inverse map.
	Let $L = (I, \gamma, P)$ be a Lagrangian sufficiently $C^\infty$-close to $L_0$ for $\gamma$ to be contained in the image of $\Phi$ and $\pr_1 \circ \Phi^{-1} \circ \gamma$ to be a diffeomorphism.
	As $\Phi$ is a local diffeomorphism, there exist a unique diffeomorphism $\ph_L\colon I \to I$ and a map $h_L\colon I \to (-\eps, \eps)$ such that $\Phi(\ph_L(t), h_L(\ph_L(t))) = \gamma(t)$ for any $t \in I$.
	Set $F_L(t) = P (\ph^{-1}(t)) - \tilde P(t, h(t))$ for $t \in I$.
	Since $\gamma^*\lambda + dP = 0$ and $\Phi^*\lambda + sdt + d\tilde P = 0$, we have
	$$
	\ph_L^* dF_L= dP - \ph_L^* (t \mapsto (t, h_L(t)))^* d\tilde P = \ph_L^* (t \mapsto (t, h_L(t)))^* \left(-\Phi^*\lambda - d\tilde P\right) = \ph^*_L (h_Ldt).
	$$
	Therefore, $F_L' = h_L$ hence $L(F_L, \ph_L) = L$, $\ph_{L(F, \ph)} = \ph$, and $F_{L(F, \ph)} = F$.
\end{proof}

We give a topological characterisation of isotopy.
The result is folklore in $\RR^2$ (see e.g.\ K\'alm\'an~\cite{K05}*{p. 2023}) and not hard in our setting, but we could not find a suitable reference.

\begin{prop} \label{p: characterisation of isotopy}
	Let $L_i = (I, \gamma_i, P_i)$, $i = 0, 1$ be two Lagrangians in the Liouville marked surface $\S$.
	Then the following are equivalent.

	\begin{enumerate}
		\item $L_0$ and $L_1$ are isotopic.
		\item $\gamma_0$ and $\gamma_1$ lie in the same component of $\Imm_I(S)$, the Reeb signs at the corresponding double points of $L_0$ and $L_1$ agree, and the actions of the corresponding Reeb chords between the ends of $L_0$ and between the ends of $L_1$ have the same sign.
	\end{enumerate}
\end{prop}
\begin{proof}
	We prove that (1) follows from (2), the converse is clear.
	By \cref{l: components of the space of immersions are orbits} and the fact that the isotopy of the domain induces an isotopy, we can assume that there exists an isotopy $\Phi_s \in \Diff(\S)$, $s \in [0, 1]$ such that $\Phi_0 = \id$ and $\Phi_1 \circ \gamma_0 = \gamma_1$.
	To prove the proposition it suffices to find a path $\Phi_s, \, s \in [1, 2]$ in $\Diff(\S)$ such that $\Phi_2 \in \Symp_{\ex}(\S)_0$ and the Lagrangian $L_2 = (\Phi_2)_* L_0$ is isotopic to $L_1$.
	Indeed, \cref{cor: Diff retracts on Symp} then implies that there exists a path $\Theta_s, \, s \in [0, 1]$ in $\Symp_{\ex}(\S)$ such that $\Theta_0 = \id_S$ and $\Theta_1 = \Phi_2$ hence an isotopy from $L_0$ to $L_1$ can be obtained as the concatenation of $s \mapsto (\Theta_s)_*L_0$ and the isotopy from $L_2$ to $L_1$.

	We now find the desired path $\Phi_s, \, s \in [1, 2]$ in $\Diff(\S)$.
	Consider the equation
	$$
		(\Phi_1^{-1})^*\lambda - \lambda + d\lambda_s(X_s, -) = dg_{(\Phi_1^{-1})^*\lambda - \lambda}, \text{ where } s \in [0, 1].
	$$
	where
	$\lambda_s = s (\Phi_1^{-1})^*\lambda + (1 - s)\lambda$
	and $g_{(\Phi_1^{-1})^*\lambda - \lambda}$ is provided by \cref{l: mu -> g_mu}.
	As in the proof of \cref{p: Diff retracts on Symp}, this equation determines a time-dependent vector field $X_s$ on $S$ which integrates to a family of diffeomorphisms
	$$
		[0, 1] \to \Diff(\S),\ s \mapsto \Psi_s,\ \Psi_0 = \id_S.
	$$
	Cartan's magic formula yields
	$$
		s(\Phi_1^{-1})^*\lambda + (1 - s)\lambda = (\Psi_s^{-1})^*\lambda + dh_s \text{ where } h_s = sg_{(\Phi_1^{-1})^*\lambda - \lambda} + \int_{0}^s \lambda_t(X_t) dt
	$$
	We define the isotopy $L_u = (I, \gamma_u, P_u), u \in [1, 2]$ by
	$$
		\gamma_{s + 1} = \Psi_s^{-1} \circ \gamma_1, \, P_{s + 1} = sP_0 + (1 - s)P_1 + \gamma_1^* h_s.
	$$
	A straightforward calculation shows that $\gamma_{s + 1}^*\lambda + dP_{s + 1} = 0$.
	For any self-intersection point $x$ of $\gamma_{s + 1}$ we have $\gamma_{s + 1}(x_+) = x = \gamma_{s + 1}(x_-)$ for some $x_+, x_- \in I$ such that $P_i(x_+) > P_i(x_-)$ for $i = 0, 1$ by assumption of the proposition.
	Therefore,
	\begin{align*}
		&P_{s + 1}(x_+) =
		sP_0(x_+) + (1 - s)P_1(x_+) + h_s(x) > \\
		&sP_0(x_-) + (1 - s)P_1(x_-) + h_s(x) =
		P_{s + 1}(x_-).
	\end{align*}
	For any Reeb chord $y\colon [0, l] \to \pd S \setminus \Lambda$ connecting $\gamma_{s + 1}(y_+) = y(0)$ and $\gamma_{s + 1}(y_-) = y(l)$ we have that $P_0(y_+) - l - P_0(y_-)$ and $P_1(y_+) - l - P_1(y_-)$ are non-zero and have the same sign by assumption, therefore, a computation similar to the above shows that $P_{s + 1}(y_+) - l - P_{s + 1}(y_-)$ is non-zero.
	Hence $L_u, u \in [1, 2]$ is indeed an isotopy.

	We set $\Phi_{s + 1} = \Psi_s^{-1} \circ \Phi_1, s \in [0, 1]$.
	Then
	$$
		\Phi_2^*\lambda - \lambda = \Phi_1^*\left((\Psi_1^{-1})^*\lambda - (\Phi_1^{-1})^*\lambda\right) = d(-\Phi_1^*h_1) \text{ hence } \Phi_2 \in \Symp_{\ex}(\S)_0
	$$
	and
	$$
		(\Phi_2)_*L_0 = (I, \Psi_1^{-1} \circ \Phi_1 \circ \gamma_0, P_0 - \gamma_0^*(-\Phi_1^*h_1)) =
		(I, \Psi_1^{-1} \circ \gamma_1, P_0 + \gamma_1^*h_1) = L_2,
	$$
	as required.
\end{proof}

\subsubsection{Isotopy of slightly degenerate Lagrangians}
\label{subsubs: Isotopy of slightly degenerate Lagrangians}

For a compact 1-dimensional manifold $I$ denote by $\Lag_I^{sg}(\S)$ the space of slightly degenerate Lagrangians (see \cref{d: slightly degenerate Lagrangian}) in $\S$ with the $C^\infty$-topology.
Two slightly degenerate Lagrangians with the same domain $I$ are \emph{isotopic} if they lie in the same component of $\Lag_I^{sg}(\S)$.

The following lemma is the reason we introduce slightly degenerate Lagrangians.
We call a Lagrangian $L$ \emph{sufficiently self-wrapped} if there are no Reeb chords of positive action with ends on $L$, i.e.\ $Y^{>0}_L = \varnothing$.

\begin{lemma} \label{l: sufficiently self-wrapped model for Lagrangian}
	Any slightly degenerate Lagrangian $L$ is isotopic to a sufficiently self-wrapped Lagrangian.
\end{lemma}
\begin{figure}[!htbp]
		\centering
		\includegraphics[width=0.46\linewidth]{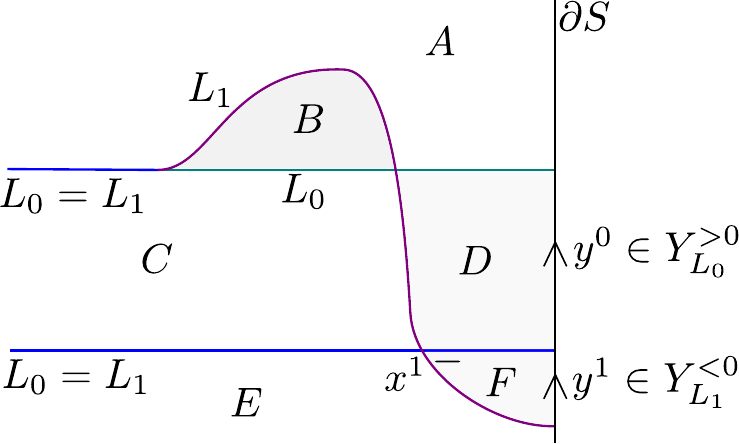}
		\caption{Replacing a Reeb chord by a double point.}
		\label{fig: Move II>0}
\end{figure}
\begin{proof}
	This follows by induction on $|Y_L^{>0} \cup XY_{L}|$ since we can always exchange a chord by a double point.
	Indeed, if $Y_L^{>0} \cup XY_{L} \ne \varnothing$ and $y^0$ is the shortest Reeb chord in $Y_L^{>0} \cup XY_{L}$ then either it is a point-chord, or all its proper sub-chords must have negative action, but as the action is additive with respect to concatenation, $y^0$ does not have proper sub-chords hence no ends of $L$ in the interior.
	In any case, $L$ in a neighbourhood of $y^0$ looks like $L_0 = (I, \gamma_0, P_0)$ in \cref{fig: Move II>0} (if $y^0$ is a point-chord, then we have a slightly different non-diffeomorphic picture, but the reasoning stays the same).
	Now replace $L_0$ with a Lagrangian $L_1 = (I, \gamma_1, P_1)$ in \cref{fig: Move II>0} in such a way that the move happens on a subsegment of $I$ (i.e.\ the \textcolor{teal}{teal} segment is changed to the \textcolor{violet}{purple} segment), the darker shaded area is equal to the lighter shaded area (i.e.\ $\int_B d\lambda = \int_{D \cup F} d\lambda$), and $\int_{F} d\lambda < A(y^0)$.

	Note that for double points we have $X_{L_1} = X_{L_0} \cup \{x^1\}$.
	For the unique $p \in \pd I$ such that $\gamma_0(p) \ne \gamma_1(p)$ let $r_p = y^0y^1$, otherwise let $r_p$ be the constant path at $\gamma_0(p) = \gamma_1(p)$.
	Define the map $\psi\colon X_{L_0} \cup Y_{L_0} \to X_{L_1} \cup Y_{L_1}$ by
	\begin{equation} \label{eq: def of psi}
		\psi(z) = \begin{cases}
			z, & \text{ if } z \in X_{L_0};\\
			x^0, & \text{ if } z = y^0;\\
			z', & \text{ such that } zr_{z_-} = r_{z_+}z' \text{ if } z \in Y_{L_0} \setminus \{y^0\}.
		\end{cases}
	\end{equation}
	In other words, on $Y_{L_0} \setminus \{y^0\}$ the map $\psi$ is defined by continuously deforming the Reeb chord by the isotopy from $L_0$ to $L_1$.

	It is clear that $\psi$ is a bijection from $X_{L_0} \cup Y_{L_0}$ to $X_{L_1} \cup Y_{L_1} \setminus \{y^1\}$.
	If $z \in X_{L_0}$, then $\psi(z) = z$, so $A(z) = A(\psi(z))$.
	If $z \in Y_{L_0} \setminus \{y^0\}$, then
	$$
		A(\psi(z)) = (P_{0}(z_+) - l_{r_{z_+}}) - (P_{0}(z_-) - l_{r_{z_-}}) - (l_z + l_{r_{z_-}} - l_{r_{z_+}}) = A(z).
	$$
	Finally, applying the Stokes theorem to region $D$ and using ${\int_B d\lambda = \int_{D} d\lambda + \int_{F} d\lambda}$, we get
	$$
		A(x^0) = A(y^0) - \int_F d\lambda > 0.
	$$
	In particular, $\psi$ restrict to a bijection $Z_{L_0} \to Z_{L_1}$ and
	$|Y_{L_0}^{>0} \cup XY_{L_0}| - 1 = |Y_{L_1}^{>0} \cup XY_{L_1}|$,
	which finishes the proof of the inductive step.
\end{proof}

We will also need the following observation.

\begin{lemma} \label{l: family of maps from isotopy of slighly degenerate Legendrians}
	Let $L_s = (I, \gamma_s, P_s), s \in [0, 1]$ be a path of slightly degenerate Lagrangians such that $L_s$ is a Lagrangian for $s \ne 1$ and $L_1$ has a unique point-chord $z^1 \in XY_{L_1}$ and $L_s$ is transverse to the slightly degenerate locus.
	Then there exist smooth families of maps
	$\Phi^s \colon S \to S,\,  \psi^s\colon I \to I,\, s \in [0, 1]$
	such that $\Phi^0 = \id$, $\psi^0 = \id$, $\Phi^s \circ \gamma_0 = \gamma_s \circ \psi^s$ for $s \in [0, 1]$, $\Phi^s \in \Diff(\S), \, \psi^s \in \Diff(I)$ for $s \in [0, 1)$, and
	$\res{\Phi^1}{S \setminus (\Phi^1)^{-1}(z^1)}, \, \res{\psi^1}{I \setminus (\psi^1)^{-1}\{z^1_+, z^1_-\}}$
	are diffeomorphisms onto their images.
\end{lemma}
\begin{figure}[!htbp]
	\centering
	\begin{subfigure}[t]{0.48\linewidth}
		\centering
		\includegraphics[width=\linewidth]{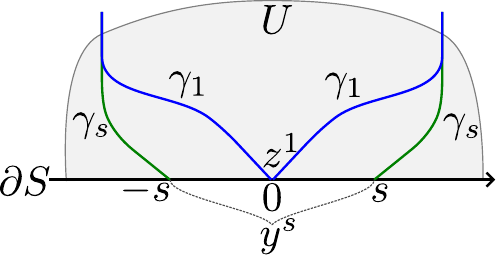}
		\caption{Local model for $y^s \in Y_{L_s}^{>0}$.}
		\label{subfig: contracting y>0}
	\end{subfigure}%
	\hfill%
	\begin{subfigure}[t]{0.48\linewidth}
		\centering
		\includegraphics[width=\linewidth]{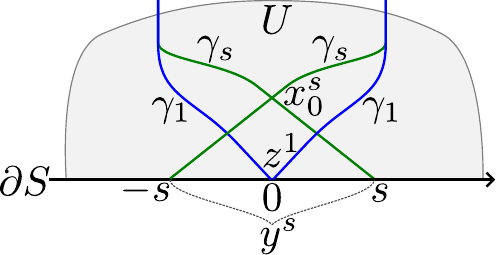}
		\caption{Local model for $y^s \in Y_{L_s}^{<0}$}
		\label{subfig: contracting y<0}
	\end{subfigure}
	\caption{}
	\label{fig: contracting y}
\end{figure}
\begin{proof}
	Let $y^s, s \in [0, 1]$ be the Reeb chord in $L_s$ that collapses to the point-chord $z^1$.
	There are two cases.

	\emph{Case 1.} $A(y^s) > 0$ for all $s \in [0, 1]$.
	In this case we can trace each double point to obtain a family $x^s \in X_{L_s}$.
	Correct $\gamma_s$ by a family in $\Diff(I)$ in such a way that the points $x^s_{\pm} \in I$ do not depend on $s$ for all families of double points $x^s$.
	Then by a family of diffeomorphisms $\Phi^s_0 \in \Diff(\S), \, s \in [0, 1]$ we can identify neighbourhoods of $y^s$ with a neighbourhood $U$ of $z^1$ and make $\Phi^s_0 \circ \gamma_0 = \gamma_s$ outside of $U$, see \cref{subfig: contracting y>0}.
	Applying a family of compactly supported diffeomorphisms of $U$ we can choose local coordinates $(x, y)$ in $U$ centred at $z^1$ with $U$ given by $\{y \ge 0\}$ such that the two branches of $\gamma_s$ are parametrised by $[0, 1)$ with $t \mapsto (s + t, t)$ and $t \mapsto (-s - t, t)$ near $\pd S$.
	There is a family of smooth functions $\ph^s\colon \RR \to \RR, s \in [0, 1]$ such that $\ph^s(x) = sx$ for $x \in [-1, 1]$, $\ph^s(x) = x$ for large $|x|$, $\ph^s$ is a diffeomorphism for $s < 1$, and $\res{\ph^1}{\RR \setminus [-1, 1]}$ is a diffeomorphism.
	We extend this to a family of smooth functions $\Phi^s_1\colon \RR^2 \to \RR^2, s \in [0, 1]$ such that $\Phi^s_1 \circ \gamma_0 = \gamma_s$ near $\pd S$, $\Phi^s_1(x, y) = (x, y)$ for large $|(x, y)|$, $\Phi^s_1$ is a diffeomorphism for $s < 1$ and $\res{\Phi^1_1}{\RR^2 \setminus [-1, 1] \times \{0\}}$ is a diffeomorphism.
	Finally, we correct by a family of diffeomorphisms compactly supported on $U \setminus \pd S$ to make $\gamma_s$ agree with $\Phi^s \circ \gamma_0$ everywhere.

	\emph{Case 2.} $A(y^s) < 0$ for all $s \in [0, 1]$.
	In this cases there is also a double point $x^s_0 \in X_{L_s}$ that collapses to the point-chord $z^1$.
	Any other double points we can track to obtain a family $x^s \in X_{L_s} \setminus \{x^s_0\}$.
	Choose a family of smooth maps $\psi^s\colon I \to I, \, s \in [0, 1]$ such that $\psi^s(x^0_{\pm}) = x^s_{\pm}$ for any $x^s \in X_{L_s} \cup XY_{L_s}$, $\psi^s$ is a diffeomorphism for $s < 1$, and $\res{\psi^1}{(\psi^1)^{-1}\{z^1_+, z^1_-\}}$ is a diffeomorphism.
	Then by a family of diffeomorphisms $\Phi^s_0 \in \Diff(\S), \, s \in [0, 1]$ we can identify neighbourhoods of $y^s$ with a neighbourhood $U$ of $z^1$ and make $\Phi^s_0 \circ \gamma_0 = \gamma_s \circ \psi^s$ outside of $U$, see \cref{subfig: contracting y<0}.
	Applying a family of compactly supported diffeomorphisms of $U$ we can choose local coordinates $(x, y)$ in $U$ centred at $z^1$ with $U$ given by $\{y \ge 0\}$ such that the two branches of $\gamma_s$ are parametrised by $[0, 1)$ with $t \mapsto (s - t, t)$ and $t \mapsto (-s + t, t)$ near $\pd S$, and $x_0^s = (0, s)$.
	Let $T_s$ be the triangle with vertices $(s, 0)$, $(0, s)$, $(-s, 0)$, i.e.\ the triangle shrinking to $z^1$.
	There is a family of smooth functions $\Phi^s_1\colon \RR^2 \to \RR^2, s \in [0, 1]$ such that $\Phi^s_1 \circ \gamma_0 = \gamma_s \circ \psi^s$ near $T_s$, $\Phi^s_1(x, y) = (x, y)$ for large $|(x, y)|$, $\Phi^s_1$ is a diffeomorphism for $s < 1$ and $\res{\Phi^1_1}{\RR^2 \setminus T_1}$ is a diffeomorphism.
	Finally, we correct by a family of diffeomorphisms compactly supported on $U \setminus T_s$ to make $\gamma_s \circ \psi^s$ agree with $\Phi^s \circ \gamma_0$ everywhere.
\end{proof}

\subsection{Polygons}
\label{subs: Polygons}
Fix a Lagrangian
$L = (I, \gamma, P)$
in a Liouville marked surface
${\S = (S, Q, \lambda)}$.
We describe immersed polygons bounded by $L$.
Their counts will give $A_\infty$-operations.

In \cref{subsubs: Polygons bounded by a Lagrangian} we define immersed polygons bounded by $L$, their virtual dimension, and their positive and negative corners, and relate the area of a polygon to the actions of its corners.
In \cref{subsubs: Polygons with one positive corner} we define the sets $\Delta(L; z^0; z^n, \ldots, z^1)$ of polygons with a positive corner at $z^0$ and negative corners at $z^1, \ldots, z^n$, where we also allow Reeb chords as corners, and prove that only finitely many of them are non-empty in \cref{l: moduli spaces are finite}.
In \cref{subsubs: 2-morphism induced by a polygon} we associate to each such polygon a 2-morphism in $\P_2^+(S, Q)$, which determines the contribution of the polygon to the $A_\infty$-operations via composition along a 2-morphism from \cref{subsubs: Composition along a 2-morphism}.

\subsubsection{Polygons bounded by a Lagrangian}
\label{subsubs: Polygons bounded by a Lagrangian}

\begin{defn} \label{d: polygon bounded by L}
	A \emph{polygon bounded by $L$} is a tuple $(\Pi, u, v)$ such that
	\begin{itemize}
		\item $\Pi$ is the compact region in $\RR^2$ bounded by a piecewise-smooth Jordan curve;
		\item
		$u\colon (\Pi, \pd \Pi) \to (S, \Im \gamma \cup (\pd S \setminus \Lambda))$
		is an orientation-preserving immersion;
		\item $v$ is a vertex (i.e.\ non-smooth point) of $\Pi$.
	\end{itemize}

	Two such polygons $(\Pi_i, u_i, v_i), \, i = 0, 1$ are \emph{equivalent} if there exists a diffeomorphism $\phi\colon \Pi_0 \to \Pi_1$ such that $u_1 \circ \phi = u_0$ and $\phi(v_0) = v_1$.
\end{defn}

Let $X$ be the set of self-intersection points of $L$ and $Y$ the set of Reeb chords between the ends of $L$, as in \cref{subsubs: Double points and Reeb chords}.
Fix a polygon $\bu = (\Pi, u, v^0)$ bounded by $L$.
Note that the set of vertices $V$ of $\Pi$ is mapped by $u$ to $X \cup \gamma(\pd I)$.
Let $W(\bu) = \{w^0, \ldots w^n\}$ be the set of connected components of $V \cup u^{-1}(\pd S)$, numbered in the counter-clockwise order on $\pd \Pi$, with $v^0 \in w^0$.
For each $w \in W(\bu)$ either $u(w)$ is a double point of $L$ or $u(w)$ is a Reeb chord between the ends of $L$, i.e.\ $u(w) \in X \cup Y$.

Let $v$ be a vertex of $\Pi$.
We say that $v$ is \emph{convex} (resp.\ \emph{concave}) if the angle of the tangent cone of $\Pi$ at $v$ is in $(0, \pi)$ (resp.\ in $(\pi, 2\pi)$).
Note that concave vertices of $\Pi$ are necessarily mapped to $X$ by $u$.
Call the number of concave vertices of $\Pi$ the \emph{virtual dimension of $\bu$} and denote it by $\vdim \bu$.
To define the $A_\infty$-operations we will only need polygons of virtual dimension zero.
To prove the $A_\infty$-relations we will also need polygons of virtual dimension one.

Fix $w \in W(\bu)$ such that $u(w)$ is a double point of $L$ (i.e.\ $u(w)\in X$).
Then the image under $u$ of a small neighbourhood of $w$ in $\Pi$ covers one quadrant of $u(w)$ if $w$ is convex and three quadrants of $u(w)$ if $w$ is concave.
If the majority of these quadrants is positive (negative), then we say that $\bu$ has a \emph{positive (negative) corner} at $w$.
We obtain the partition
\begin{equation}
	\{w \in W(\bu) \colon u(w) \in X\} = W^+(\bu) \sqcup W^-(\bu)
\end{equation}
into positive and negative corners.
Our definition of positive and negative corners follows the convention in Legendrian contact homology, the convention in Lagrangian Floer theory is, unfortunately, the opposite.

The Stokes theorem gives
\begin{equation} \label{eq: area formula for an immersed polygon}
\sum_{w \in W^+(\bu)} A(u(w)) -
\sum_{\substack{w \in W(\bu) \colon\\u(w) \in Y}} A(u(w)) -
\sum_{w \in W^-(\bu)} A(u(w)) =
\int\limits_{\Pi} u^*d\lambda > 0.
\end{equation}

\subsubsection{Polygons with one positive corner}
\label{subsubs: Polygons with one positive corner}

Formula \cref{eq: area formula for an immersed polygon} suggests that the Reeb chords on the boundary of an immersed polygon can be thought of as negative corners.
In fact, certain Reeb chords not on the boundary of an immersed polygon can be thought of as its positive corners.
In the definition below we make this idea precise by grouping immersed polygons bounded by $L$ based on where their negative and positive `corners' are.
\begin{defn} \label{d: moduli spaces of polygons}
Fix $n \ge 0$, and $z^0, z^1, \ldots, z^n \in Z = X \cup Y^{>0}$.
We define the set $\Delta(L; z^0; z^n, \ldots, z^1)$ of polygons bounded by $L$ with a `positive corner' at $z^0$ and `negative corners' at $z^1, \ldots, z^n$.

\begin{enumerate}
	\item If $z^0$ is a double point of $L$ (i.e.\ $z^0 \in X$), let $\Delta(L; z^0; z^n, \ldots, z^1)$ be the set of equivalence classes of polygons $\bu = (\Pi, u, v^0)$ bounded by $L$ such that
	$W(\bu) = \{w^0, w^1, \ldots, w^n\}$,
	$W^+(\bu) = w^0$, and
	$u(w^i) = z^i$ for $i = 0, \ldots, n$.
	\item If $z^0$ is a Reeb chord between the ends of $L$ (i.e.\ $z^0 \in Y$), let $\Delta(L; z^0; z^n, \ldots, z^1)$ be the union of the following:
	\begin{enumerate}
		\item if $z^1 \in Y^{>0}$ and there exists (necessarily unique) $y \in Y$ such that the concatenation $yz^0$ is equal to $z^1$, then the set of equivalence classes of polygons $\bu = (\Pi, u, v^0)$ bounded by $L$ such that
		$W(\bu) = \{w^0, w^1, \ldots, w^{n - 1}\}$,
		$W^+(\bu) = \varnothing$,
		$u(w^0) = y$, and
		$u(w^i) = z^{i + 1}$ for $i = 1, \ldots, n - 1$;
		\item if $z^n \in Y^{>0}$ and there exists (necessarily unique) $y \in Y$ such that the concatenation $z^0y$ is equal to $z^n$, then the set of equivalence classes of polygons $\bu = (\Pi, u, v^0)$ bounded by $L$ such that
		$W(\bu) = \{w^0, w^1, \ldots, w^{n - 1}\}$,
		$W^+(\bu) = \varnothing$,
		$u(w^i) = z^{i + 1}$ for $i = 0, \ldots, n - 2$, and
		$u(w^{n - 1}) = y$;
		\item if $n = 2$, $z^1, z^2 \in Y^{>0}$, and $z^0 = z^2z^1$, then a singleton $\{*\}$, we call its unique element a \emph{degenerate polygon}.
		Set $\vdim(*) = 0$.
	\end{enumerate}
\end{enumerate}

For $k \ge 0$ denote by $\Delta^k(L; z^0; z^n, \ldots, z^1)$ the subset of $\Delta(L; z^0; z^n, \ldots, z^1)$ consisting of the elements of virtual dimension $k$.
\end{defn}

\begin{figure}[!htbp]
	\centering
	\begin{subfigure}[t]{0.48\linewidth}
		\centering
		\includegraphics[width=0.3\linewidth]{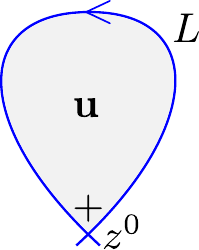}
		\caption{$\Delta^0(L; z^0;)$.}
		\label{subfig: Delta^0(L z^0)}
	\end{subfigure}%
	\hspace{0.04\linewidth}%
	\begin{subfigure}[t]{0.48\linewidth}
		\centering
		\includegraphics[width=0.3\linewidth]{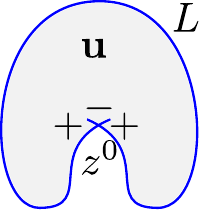}
		\caption{$\Delta^1(L; z^0;)$.}
		\label{subfig: Delta^1(L z^0)}
	\end{subfigure}

	\bigskip%
	\begin{subfigure}[t]{\linewidth}
		\includegraphics[width=\linewidth]{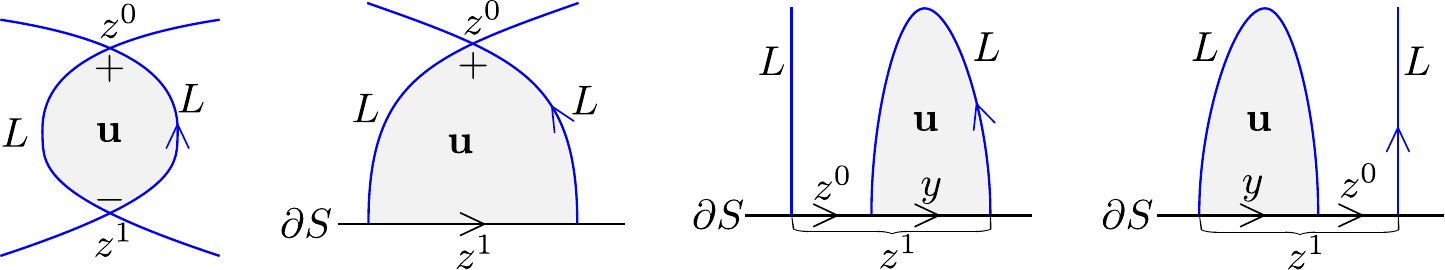}
		\caption{$\Delta^0(L; z^0; z^1)$.}
		\label{subfig: Delta^0(L z^0 z^1)}
	\end{subfigure}

	\bigskip%
	\begin{subfigure}[t]{\linewidth}
		\centering
		\includegraphics[width=0.7\linewidth]{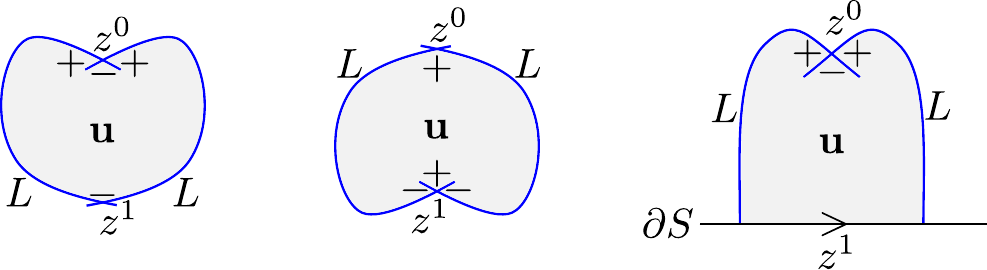}
		\caption{$\Delta^1(L; z^0; z^1)$.}
		\label{subfig: Delta^1(L z^0 z^1)}
	\end{subfigure}

	\bigskip%
	\begin{subfigure}[t]{\linewidth}
		\includegraphics[width=\linewidth]{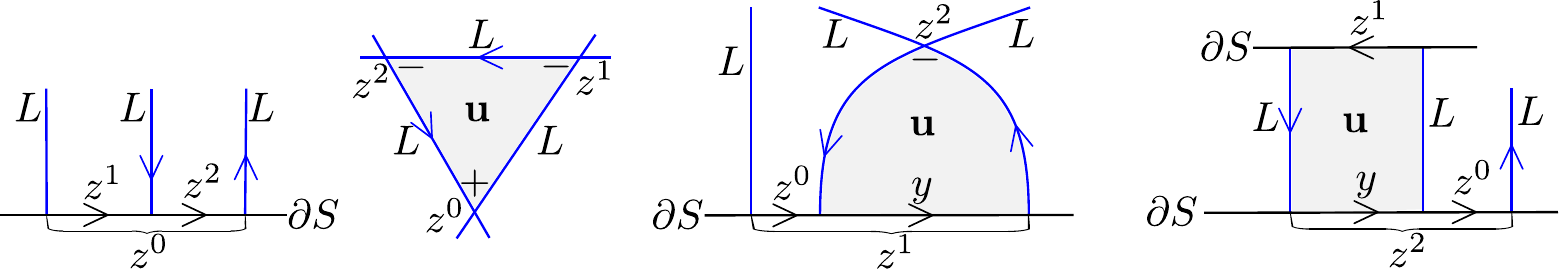}
		\caption{$\Delta^0(L; z^0; z^2, z^1)$.
		Four out of nine types are listed.}
		\label{subfig: Delta^0(L z^0 z^2 z^1)}
	\end{subfigure}
	\caption{Examples of polygons in $\Delta(L; z^0; z^n, \ldots, z^1)$.}
	\label{fig: Examples of polygons bounded by L}
\end{figure}

See \cref{fig: Examples of polygons bounded by L} for examples of polygons from \cref{d: moduli spaces of polygons}.
The pictures of polygons bounded by $L$ that we draw are of the pull-back of a neighbourhood of the polygon $\bu$ and Reeb chords involved from $S$ to $\RR^2$ via an extension of $u$.
For a polygon $\bu = (\Pi, u, v^0) \in \Delta(L; z^0; z^n, \ldots, z^1)$ the vertex $v^0$ is uniquely determined by the rest of the data, so we suppress it from the notation and write $\bu = (\Pi, u)$.

\begin{example} \label{ex: polygons bounded by two Lagrangians with the same underlying curve}
We continue \cref{ex: two Lagrangians with the same underlying curve}.
For the Lagrangian $L_0$ in \cref{fig: Lagrangian with actions} only the following $\Delta^0(L_0; \ldots)$ are non-empty:
$$
\begin{aligned}
&\Delta^0(L_0; x^1;), \Delta^0(L_0; x^2; x^3),\\
&\Delta^0(L_0; x^1; y^1), \Delta^0(L_0; x^2; y^2y^1), \Delta^0(L_0; x^3; y^2), \Delta^0(L_0; x^3; y^4y^2y^1),\\
&\Delta^0(L_0; y^2y^1; y^4y^2y^1), \Delta^0(L_0; y^2y^1; y^2, y^1), \Delta^0(L_0; x^2; y^2, x^1).
\end{aligned}
$$
For the Lagrangian $L_1$ in \cref{fig: Lagrangian with actions} only the following $\Delta^0(L_1; \ldots)$ are non-empty:
$$
\begin{aligned}
&\Delta^0(L_1; x^1;), \Delta^0(L_1; x^1; y^1), \Delta^0(L_1; x^1; \bar x^3, x^2).
\end{aligned}
$$
These non-empty sets are all singletons.

This example is continued in \cref{ex: degree for two Lagrangians with the same underlying curve}.
\end{example}

\begin{remark}
	If we were to consider moduli spaces of pseudoholomorphic discs, then the set $\Delta^k(L; z^0; z^n, \ldots, z^1)$ would match the set of connected components of the corresponding moduli space for $k = 0, 1$.
\end{remark}

Below we list elementary properties of $\Delta(L; z^0; z^n, \ldots, z^1)$.

\begin{lemma} \label{l: moduli spaces are finite}~
	Let $L$ be a Lagrangian in a Liouville marked surface $\S$.
	\begin{enumerate}
		\item Suppose that $\Delta(L; z^0; z^n, \ldots, z^1) \ne \varnothing$.
		Then
		$$
			0 \le \sum_{i = 1}^n A(z^i) \le A(z^0).
		$$
		Moreover, the first inequality is strict unless $n = 0$ and the second inequality is strict unless $n = 2$ and $\Delta(L; z^0; z^2, z^1)$ contains a degenerate polygon.
		In particular, if $\Delta(L; z^0; z^2, z^1)$ contains a degenerate polygon, then $\Delta(L; z^0; z^2, z^1) = \{*\}$.
		\item The set
		$$
			\bigcup_{\substack{n \ge 0,\\z^0, z^1, \ldots, z^n \in Z_L}} \Delta(L; z^0; z^n, \ldots, z^1)
		$$
		is finite.
	\end{enumerate}
\end{lemma}
\begin{proof}
	Suppose the case (c) in \cref{d: moduli spaces of polygons} holds, i.e.\ $z^0 = z^2z^1$ for $z^0, z^1, z^2 \in Y^{>0}$.
	Then $0 < A(z^1) + A(z^2) = A(z^0)$.

	We prove (1) in the remaining cases.
	Pick a polygon $\bu = (\Pi, u) \in \Delta(L; z^0; z^n, \ldots, z^1)$.
	From \cref{eq: area formula for an immersed polygon} and additivity of action with respect to the concatenation of Reeb chords we get
	\begin{equation} \label{eq: upper bound on the area of immersed polygon}
		0 < \int\limits_{\Pi} u^*d\lambda = A(z^0) - \sum_{i = 1}^n A(z^i) \le A(z^0).
	\end{equation}
	(1) follows since $A(z^i) > 0$ for all $i$.

	To show (2), as $Z_L$ is finite, \cref{eq: upper bound on the area of immersed polygon} implies that it suffices to prove that for any $C > 0$ there are only finitely many (equivalence classes of) polygons bounded by $L$ of area not exceeding $C$.
	This is indeed the case since each polygon bounded by $L$ is obtained by gluing copies of the regions into which $L$ divides $S$ along the sides of these regions (a side is a smooth segment of the boundary of the region between two double points of $L$ or between a double point of $L$ and an end of $L$).
	There are finitely many such regions, each of them has positive area and finitely many sides, so there are only finitely many ways to glue a number of their copies to get the total area of no more than $C$.
\end{proof}

\subsubsection{2-morphism induced by a polygon}
\label{subsubs: 2-morphism induced by a polygon}
We continue studying polygons with one positive corner bounded by the Lagrangian $L = (I, \gamma, P)$.

Note that each $z \in Z$ gives a path from $\gamma(z_+)$ to $\gamma(z_-)$, which we again denote by $z$.
Namely, if $z$ is a Reeb chord, then the path is $z$ itself and if $z$ is a double point, then the path is the constant path at $z$.
We denote by $\gamma(z_+) \xrightarrow{[z]} \gamma(z_-)$ the corresponding 1-morphism in the path 2-category $\P_2^+(S, Q)$ (see \cref{d: directed relative path 2-category}).

Fix $\bu \in \Delta(L; z^0; z^n, \ldots, z^1)$.
We construct a 2-morphism $[z^0] \xRightarrow{\sigma_\bu}[\gamma_\bu]$.

If $\bu = (\Pi, u)$ is a polygon bounded by $L$, then the curve $\gamma_\bu$ is a piecewise-smooth curve from $\gamma(z^0_+)$ to $\gamma(z^0_-)$ given by the concatenation
$$
	\gamma(a^n)z^n\ldots\gamma(a^1)z^1\gamma(a^0)
$$
where $a^i = a^i_{\bu}$ for $i = 0, \ldots, n$ are paths in $I$ defined as follows.
\begin{itemize}
	\item If $\gamma(z^0_+) = \gamma(z^1_+)$, then $a^0$ is the constant path at $z^0_+$, otherwise there is a side of $\bu$ going from $\gamma(z^0_+)$ to $\gamma(z^1_+)$, and $a^0$ is the path in $I$ parametrising this side by $\gamma$.
	\item For $i = 1, \ldots, n - 1$ if $\gamma(z^i_-) = \gamma(z^{i + 1}_+)$, then $a^i$ is the constant path at $z^i_-$, otherwise there is a side of $\bu$ going from $\gamma(z^i_-)$ to $\gamma(z^{i + 1}_+)$, and $a^i$ is the path in $I$ parametrising this side by $\gamma$.
	\item If $\gamma(z^n_-) = \gamma(z^0_-)$, then $a^n$ is the constant path at $z^n_-$, otherwise there is a side of $\bu$ going from $\gamma(z^n_-)$ to $\gamma(z^0_-)$ and $a^n$ is the path in $I$ parametrising this side by $\gamma$.
\end{itemize}
The 2-morphism $\sigma_{\bu} \colon [z^0] \Rightarrow [\gamma_\bu]$ is induced by $\bu$.
Namely,
\begin{itemize}
	\item in case (1) of \cref{d: moduli spaces of polygons} $u \colon \Pi \to S$ itself is (after a lazy reparametrisation) a 2-morphism from $[z^0]$ to $[\gamma_\bu]$;
	\item in case (2a) of \cref{d: moduli spaces of polygons} $u$ gives a 2-morphism from $[\gamma(z^0_-)]$ to $[\gamma(a^n)z^n\ldots\gamma(a^1)y]$ where $z^1 = yz^0$, and $\sigma_{\bu}$ is obtained from that by pre-composition with $[z^0]$;
	\item in case (2b) of \cref{d: moduli spaces of polygons} $u$ gives a 2-morphism from $[\gamma(z^0_+)]$ to $[y\gamma(a^{n - 1})\ldots\gamma(a^1)z^1\gamma(a^0)]$ where $z^n = z^0y$, and $\sigma_{\bu}$ is obtained from that by post-composition with $[z^0]$.
	\item in case (2c) of \cref{d: moduli spaces of polygons} we have $[\gamma_{\bu}] = [z^0] = [z^2z^1]$ and $\sigma_{\bu} = \id_{[z^0]}$.
\end{itemize}

\subsection{Grading}
\label{subs: Grading}
We describe the last piece of geometric data, which is needed to make obstruction algebras and Fukaya categories $\ZZ$-graded.

\begin{defn} \label{d: graded surface and Lagrangian}
A \emph{grading} on a marked surface $(S; Q)$ is a \emph{line field} $\eta$ on $S$, i.e.\ an element of $C^\infty(S, \PP TS)$.

A \emph{grading} on an immersed curve $\gamma\colon I \to S \setminus Q$ is a homotopy class of paths from $\gamma^*\eta$ to the tangent line $\dot\gamma$, i.e.\ an element
$\gr \in \Pi_1C^{\infty}(I, \gamma^*\PP T S)(\gamma^*\eta, \dot\gamma)$.
\end{defn}

Following Lekili--Polishchuk~\cite{LP20}*{\Sec 1.1}, we observe that ${\pi_0 C^{\infty}(S, \PP TS)}$ is a torsor over $H^1(S; \ZZ)$ and each component of $C^{\infty}(S, \PP TS)$ is a $K(\ZZ, 1)$ space.
For a line field $\eta$ on $S$ the generator of
${\pi_1(C^{\infty}(S, \PP TS), \eta) = \ZZ}$
is given by a counter-clockwise rotation of $\eta$ by $\pi$ in each tangent space $T_pS, \, p \in S$.

In \cref{subsubs: Gradings on a given curve} we characterise the immersed curves admitting a grading in terms of winding numbers and prove that the gradings on such a curve form a torsor over $\ZZ^{\pi_0 I}$ in \cref{l: gradings on an immersed curve}.
In \cref{subsubs: Degrees of double points and Reeb chords} we define the degrees of double points and Reeb chords of a graded Lagrangian, check that the two definitions agree on point-chords of slightly degenerate Lagrangians, and relate the degrees of the corners of a polygon to its virtual dimension.
In \cref{subsubs: Graded diffeomorphisms} we recall graded diffeomorphisms and deduce from \cref{p: Diff retracts on Symp} that graded exact symplectomorphisms form a strong deformation retract of graded diffeomorphisms in \cref{cor: graded Diff retracts on Symp}.
\subsubsection{Gradings on a given curve}
\label{subsubs: Gradings on a given curve}

Consider an immersed curve $\gamma\colon \RR/\ZZ \to S \setminus Q$.
The tangent line $\dot \gamma$ is a section of $\gamma^*\PP TS$ and $\eta$ gives a trivialisation of $\PP TS$, hence we obtain the \emph{winding number} $w_\eta(\gamma) \in \ZZ$ as the degree of $\dot \gamma$.
Note that the winding number is an invariant of the isotopy class of an immersed curve and of the homotopy class of the line field.

We have the following easy lemma characterising gradings on a given curve.

\begin{lemma} \label{l: gradings on an immersed curve}
	Let $I$ be a compact 1-dimensional manifold and $\gamma\colon I \to S \setminus Q$ be an immersed curve in a graded surface $(S, Q, \eta)$.
	Then
	\begin{enumerate}
		\item $\gamma$ admits a grading
		$\gr \in \Pi_1C^{\infty}(I, \gamma^*\PP T S)(\gamma^*\eta, \dot\gamma)$
		if and only if the winding number $w_{\eta}(\res{\gamma}{C})$ vanishes for each circle component $C$ of $I$.
		\item If $\gamma$ admits a grading, then the set of gradings on $\gamma$ is a torsor over $\ZZ^{\pi_0 I}$.
	\end{enumerate}
\end{lemma}

\begin{proof}
	(1) is clear since given a path $\gr_t$ from $\gamma^*\eta(t)$ to $\dot \gamma(t)$ for some $t \in I$ there is a homotopically unique way to extend it to a path $i^*\gr$ from $i^*\gamma^*\eta$ to $i^*\dot \gamma$ for any immersion $i\colon [0, 1] \to C$ with $i(0) = t$.
	If $i\colon [0, 1] \to C$ is a parametrisation of a circle component $C$ of $I$ with $i(0) = i(1) = t$, then the winding number $w_{\eta}(\res{\gamma}{C})$ is equal to
	$$
		i^*\gr(0)^{-1} \cdot i^*\gr(1) \in \pi_1(\PP T_{\gamma(t)}S, \eta(\gamma(t))) = \ZZ.
	$$

	(2) follows from the above since
	$$
		\pi_1\left(C^{\infty}(C, \gamma^*\PP T S), \gamma^*\eta\right) = \pi_1(\PP T_{\gamma(t)}S, \eta(\gamma(t))) = \ZZ
	$$
	for any connected component $C$ of $I$ and $t \in C$.
\end{proof}

\cref{l: gradings on an immersed curve} implies that grading is discrete data defined for any immersed curve in $S \setminus Q$ hence it is transported along isotopies of such curves.

\subsubsection{Degrees of double points and Reeb chords}
\label{subsubs: Degrees of double points and Reeb chords}
Let $L = (I, \gamma, P; \gr)$ be a graded Lagrangian in a graded Liouville marked surface $\S = (S, Q, \eta, \lambda)$.
We define the degree of double points of $L$ and Reeb chords between the ends of $L$.

The \emph{degree} of a double point $x \in X_L$ is an integer $|x|$ determined as follows.
First, obtain an element in $\pi_1(\PP T_x \S, \eta(x))$ as a concatenation of $\gr(x_-)$ (which rotates $\eta(x)$ to $\dot \gamma(x_-)$), followed by the smallest counter-clockwise rotation from $\dot \gamma(x_-)$ to $\dot \gamma(x_+)$, and then followed by the inverse of $\gr(x_+)$.
Now identify $\pi_1(\PP T_x S, \eta(x))$ with $\ZZ$ by sending the counter-clockwise rotation by $\pi$ to $1$ to obtain the degree of $x$.

The \emph{degree} of a Reeb chord $y \in Y_L$ is an integer $|y|$ determined as follows.
First, obtain a path from $\eta(\gamma(y_-))$ to $\dot y (l_y)$ in $\PP(T_{\gamma(y_-)}S)$ as a concatenation of $\gr(y_-)$ and the smallest counter-clockwise rotation from $\dot \gamma(y_-)$ to $\dot y (l_y)$.
Secondly, transport it along $[0, l_y]$ to a path from $\eta(\gamma(y_+))$ to $\dot y (0)$ in $\PP(T_{\gamma(y_+)}S)$.
Thirdly, compose it with the smallest clockwise rotation from $\dot y(0)$ to $\dot \gamma(y_+)$ and then with the inverse of $\gr(y_+)$ to obtain an element in $\pi_1(\PP T_{\gamma(y_+)}, \eta(\gamma(y_+)))$.
Finally, identify $\pi_1(\PP T_{\gamma(y_+)}, \eta(\gamma(y_+)))$ with $\ZZ$ by sending the counter-clockwise rotation by $\pi$ to $1$ to obtain the degree of $y$.
Note that the degree is additive with respect to concatenation of Reeb chords.

\begin{lemma} \label{l: degree is well-defined for slightly degenerate Lagrangian}
	Let $L$ be a slightly degenerate Lagrangian (see \cref{d: slightly degenerate Lagrangian}).
	Then for a point-chord $z \in XY_{L}$ the degrees of $z$ as a double point and as a Reeb chord are equal.
\end{lemma}
\begin{proof}
	The two definitions of degree agree if and only if the composition of the smallest counter-clockwise rotation from $\dot \gamma(z^0_-)$ to $T_{z}\pd S$ with the smallest clockwise rotation from $T_{z}\pd S$ to $\dot \gamma(z^0_+)$ is homotopic to the smallest counter-clockwise rotation from $\dot \gamma(z^0_-)$ to $\dot \gamma(z^0_+)$.
This is equivalent to the unique quadrant of $z$ contained in $S$ being negative.
\end{proof}

Let $\bu = (\Pi, u)$ be a polygon bounded by $L$ as in \cref{d: polygon bounded by L}.
The Poincar\'e--Hopf index theorem gives
\begin{equation} \label{eq: degree formula for an immersed polygon}
	\sum_{w \in W^+(\bu)} |u(w)| +
	\sum_{w \in W^-(\bu)} (1 - |u(w)|) +
	\sum_{\substack{w \in W(\bu) \colon\\u(w) \in Y}} (1 - |u(w)|) =
	2 + \vdim \bu.
\end{equation}
Regrouping the terms for $\bu \in \Delta^k(L; z^0; z^n, \ldots, z^1)$ we obtain

\begin{equation} \label{eq: |z^0| from |z^1| + ... + |z^n|}
|z^0| = 2 - n + |z^1| + \ldots + |z^n| + k.
\end{equation}

\begin{figure}[!htbp]
	\includegraphics[width=0.6\linewidth]{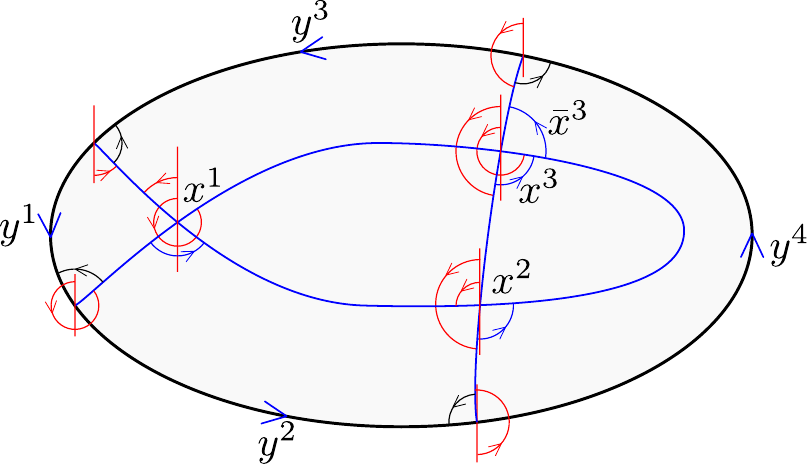}
	\caption{Graded Lagrangian in a graded surface.}
	{\small The vertical \textcolor{red}{red} lines are the line field $\eta$, the \textcolor{red}{red} arcs are the paths from $\eta$ to the tangent line, the \textcolor{blue}{blue} arcs are the smallest counter-clockwise rotations of the understrand to the overstrand, the black arcs are the smallest counter-clockwise rotations from the curve to the boundary.}
	\label{fig: Lagrangian with degrees}
\end{figure}

\begin{example} \label{ex: degree for two Lagrangians with the same underlying curve}
	We continue \cref{ex: polygons bounded by two Lagrangians with the same underlying curve}.

	We fix a line field on the disc to be the vertical line at every point and pick the same grading for $L_0$ and $L_1$, as indicated by the red arcs in \cref{fig: Lagrangian with degrees}.
	The only difference between $L_0$ and $L_1$ that affects the degree is that the understrand and the overstrand for $x^3$ and $\bar x^3$ are swapped, the corresponding smallest counter-clockwise rotations are marked with $x^3$ for $L_0$ and $\bar x^3$ for $L_1$.

	The degrees of the elements of $Z_{L_0}$ are
	$$
		|x^1| = 2,\ |x^2| = 1,\ |x^3| = 0,\ |y^1| = 1,\ |y^2| = -1,\ |y^2y^1| = 0,\ |y^4y^2y^1| = -1.
	$$
	The degrees of the elements of $Z_{L_1}$ are
	$$
		|x^1| = 2,\ |x^2| = 1,\ |\bar x^3| = 1,\ |y^1| = 1.
	$$
\end{example}

\subsubsection{Graded diffeomorphisms}
\label{subsubs: Graded diffeomorphisms}
We recall the well-known way to incorporate gradings into diffeomorphisms and symplectomorphisms, see Seidel~\cite{S00}*{\Sec 2.b} for a more general discussion.

Recall that a diffeomorphism of surfaces $\Phi\colon S_0 \to S_1$ induces a bijection
$$
\Phi_*\colon C^\infty(S_0, \PP TS_0) \to C^\infty(S_1, \PP TS_1), \; (\Phi_*\eta)(\Phi(p)) = d_{p}\Phi (\eta(p)).
$$

Let $\eta_i$ be a line field on $S_i$ for $i = 0, 1$.
Then $\Phi$
is called \emph{gradable} if $\Phi_*\eta_0$ and $\eta_1$ lie in the same component of $C^{\infty}(S_1, \PP TS_1)$.
Let $\Phi^s, \, s \in [0, 1]$ be a family in $\Diff(S_0, S_1)$ and $\tilde \Phi^0$ a path from $(\Phi^0)_*\eta_0$ to $\eta_1$ in $C^{\infty}(S_1, \PP TS_1)$.
Then we can get a path $\tilde \Phi^1$ from $(\Phi^1)_*\eta_0$ to $\eta_1$ by concatenating $s \mapsto \Phi^{1 - s}_*\eta_0$ and $\tilde \Phi^0$.
This implies that in a given component of $\Diff(S_0, S_1)$ either all diffeomorphisms are gradable or none is.

Let $\Diff^{\gr}((S_0, \eta_0), (S_1, \eta_1))$ be the space of pairs $(\Phi, \tilde \Phi)$ where $\Phi \in \Diff(S_0, S_1)$ is gradable and $\tilde \Phi$ is a homotopy class of paths from $\Phi_*\eta_0$ to $\eta_1$, i.e.\ an element of $\Pi_1 C^\infty(S_1, \PP TS_1)(\Phi_*\eta_0, \eta_1)$.
The composition is given by
$$
(\Phi, \tilde \Phi) \circ (\Psi, \tilde \Psi) = (\Phi \circ \Psi, \tilde \Phi \cdot \Phi_*\tilde \Psi)
$$
where $\cdot$ denotes concatenation.

Consider a graded marked surface $(S, Q, \eta)$.
Each component of $S$ can be treated separately, so we assume $S$ is connected.
The observation after \cref{d: graded surface and Lagrangian} gives an exact sequence
$$
\xymatrix{
0 \ar[r] & \ZZ \ar[r] & \Diff^{\gr}(S, Q, \eta) \ar[r] & \Diff(S, Q) \ar[r] & H^1(S; \ZZ)}
$$
where all maps except the one to $H^1(S; \ZZ)$ are group homomorphisms and $H^1(S; \ZZ)$ is considered as a pointed set.
Moreover, the projection map $(\Phi, \tilde \Phi) \mapsto \Phi$ sends each component of $\Diff^{\gr}(\S)$ homeomorphically onto a component of $\Diff(\S)$.
In particular, the identity component $\Diff(\S)_0$ can be identified with the identity component $\Diff^{\gr}(\S)_0$.

Now let $\S_i = (S, Q, \eta, \lambda_i)$ be two Liouville surfaces supported on the same $(S, Q, \eta)$.
Denote by $\Symp_{\ex}^{\gr}(\S_0, \S_1)$ the subset of $\Diff^{\gr}(S, Q, \eta)$ consisting of the pairs $(\Phi, \tilde \Phi)$ with $\Phi \in \Symp_{\ex}(\S_0, \S_1)$.
\cref{p: Diff retracts on Symp} and the above discussion yield a strong deformation retraction for graded diffeomorphisms.
\begin{cor} \label{cor: graded Diff retracts on Symp}
	$\Symp_{\ex}^{\gr}(\S_0, \S_1)$ is a strong deformation retract of $\Diff^{\gr}(S, Q, \eta)$.
\end{cor}

\section{Obstruction algebra}
\label{s: Obstruction algebra}
In this section we define Lagrangian branes in a weighted Liouville surface, which will be the objects of the geometric Fukaya $A_\infty$-pre-category.

In \cref{subs: Lagrangian pre-branes and branes} we define Lagrangian pre-branes and their obstruction algebras, which generalise the Chekanov--Eliashberg algebra to dg-coefficients and to Legendrians with ends on the boundary.
In \cref{subs: Invariance of the obstruction algebra} we prove that the obstruction algebra does not change under isotopies of slightly degenerate Lagrangians, which in particular allow us to exchange Reeb chords for double points, see \cref{p: isotopy induces isomorphism of obstruction algebras}.
In \cref{subs: Obstruction algebra is A_infty} we prove that the obstruction algebra is a structured curved $A_\infty$-algebra in \cref{p: obstruction algebra is a nilpotent A infinity algebra}.
In \cref{subs: Lagrangian branes} we define Lagrangian branes as Lagrangian pre-branes together with a bounding cochain in the obstruction algebra and relate bounding cochains to augmentations of the Chekanov--Eliashberg algebra.

\subsection{Obstruction algebra of a Lagrangian pre-brane}
\label{subs: Lagrangian pre-branes and branes}
We now combine the algebraic data in \cref{subs: Weighted system of dg-categories} with the geometric setup described in \cref{s: Geometric preparations}.
Namely, we fix a \emph{graded Liouville marked surface} $\S = (S, Q, \eta, \lambda)$ (see \cref{d: Liouville structure on marked surface,d: graded marked surface}) and a weighted system of dg-categories $\C\colon \P^+_2(\S) \to \dgcat_\k$ (see \cref{d: weighted system of dg-categories}). Recall from \nameref{s: Introduction} that we call $(\S; \C)$ a \emph{weighted Liouville surface}.

In \cref{subsubs: Lagrangian pre-branes} we define Lagrangian pre-branes, which are graded Lagrangians together with a local system of objects over the pullback of $\C$.
In \cref{subsubs: curved A_infty-algebra} we recall curved $A_\infty$-algebras and introduce structured $A_\infty$-algebras, in which sufficiently nested expressions in $\mu$ vanish, so that no completion is needed over an arbitrary commutative ring $\k$.
In \cref{subsubs: Obstruction algebra} we define the obstruction algebra $\A_\L$ of a Lagrangian pre-brane $\L$, whose entries are given by double points and Reeb chords of positive action and whose structure maps are given by compositions along the 2-morphisms induced by polygons with one positive corner, and relate it to the Chekanov--Eliashberg algebra.

\subsubsection{Lagrangian pre-branes}
\label{subsubs: Lagrangian pre-branes}
Let $\C$ be a local system of dg-categories on a manifold $I$, i.e.\ a functor from the fundamental groupoid $\Pi_1(I)$ to $\dgcat_\k$.
A \emph{local system of objects} over $\C$ on $I$ is the data of
\begin{itemize}
	\item $\M(t) \in \Ob \C_t$ for each $t \in I$ and
	\item $\M(a) \in Z^0\C_{a(1)}(\C(a)\M(a(0)), \M(a(1)))$ for each path $a\colon [0, 1] \to I$
\end{itemize}
such that
\begin{itemize}
	\item $\M(a)$ only depends on the homotopy-rel-endpoints class of $a$;
	\item $\M(\text{constant path})$ is the identity;
	\item $\M(ba) = \M(b) \circ \C(b)\M(a)$ for any two paths $a, b\colon [0, 1] \to I$ with $a(1) = b(0)$.
\end{itemize}

Now let $L = (I, \gamma, P)$ be a Lagrangian in $\S$.
Then $\C$ induces a local system $\gamma^*\C$ of dg-categories on $I$, namely, $\gamma^*\C_t = \C_{\gamma(t)}$ for $t \in I$ and the transition functor
${\gamma^*\C(a)\colon \gamma^*\C_{a(0)} \to \gamma^*\C_{a(1)}}$
associated to a path $a\colon [0, 1] \to I$ only depends on the homotopy-rel-endpoints class of $a$ since all lazy 2-paths in $I$ are thin homotopies.

\begin{defn} \label{d: Lagrangian pre-brane}
	A \emph{Lagrangian pre-brane} in $(\S; \C)$ is the data of
	\begin{itemize}
		\item a graded Lagrangian $L = (I, \gamma, P; \gr)$ in $\S$ (see \cref{d: Lagrangian,d: graded surface and Lagrangian}) and
		\item a local system $\M$ of objects over $\gamma^*\C$ on $I$.
	\end{itemize}
\end{defn}

Note that our Lagrangian pre-branes do not carry a  $\Pin$ structure since the signs coming from the non-triviality of the $\Pin$ structure are incorporated into the transition maps of the local system.

\subsubsection{Curved \texorpdfstring{$A_\infty$}{A∞}-algebras}
\label{subsubs: curved A_infty-algebra}
Our primary reference for curved $A_\infty$-algebras is Fukaya--Oh--Ohta--Ono~\cite{FOOO}*{\Sec\Sec 3.2, 3.6, 4.2, 4.3}.
The difference with \cite{FOOO} is that there the curved $A_\infty$-algebras are filtered and supported on a free graded module over the Novikov ring, while our curved $A_\infty$-algebras are structured and supported on an arbitrary graded module over an arbitrary commutative ring $\k$.

A \emph{($\ZZ$-graded) curved $A_\infty$-algebra (over $\k$)} is a $\ZZ$-graded $\k$-module
$$
\A = \bigoplus_{d \in \ZZ} \A^d
$$
together with a collection of linear maps
$
	\mu^n\colon \A^{\otimes n} \to \A
$
of degree $2 - n$ for each $n \ge 0$ such that the $A_n$-relation
\begin{equation} \label{eq: A-infinity relations}
	\sum_{n \ge j \ge i \ge 0} (-1)^{\dagger}\mu^{n - j + 1 + i}\left(a^n, \ldots, a^{j + 1}, \mu^{j - i}(a^{j}, \ldots, a^{i + 1}), a^i, \ldots, a^1\right) = 0
\end{equation}
where $\dagger = |a^i| + \ldots + |a^1| - i$ holds for any $n \ge 0$ and homogeneous elements $a^1, \ldots, a^n \in \A$.

$\mu^0()$ in \cref{eq: A-infinity relations} is interpreted as $\mu^0(1_\k) \in \A^2$.
The $A_0$-relation \cref{eq: A-infinity relations} reads $\mu^1(\mu ^ 0()) = 0$ and the non-vanishing of $\mu^0$ is the obstruction to $\mu^1$ being a differential as
$\mu^1(\mu^1(a)) = (-1)^{|a|} \mu^2(\mu^0(), a) - \mu^2(a, \mu^0())$
by the $A_1$-relation.

An \emph{$A_\infty$-algebra} is a curved $A_\infty$-algebra with $\mu^0 = 0$.
We sometimes use the term \emph{uncurved $A_\infty$-algebra} for emphasis.
In an uncurved $A_\infty$-algebra the $A_0$-relation is vacuous and the first three non-trivial relations are

\begin{equation*} \label{eq: A_1-relation} \tag{$A_1$}
		\mu^1(\mu^1(a_1)) = 0;
\end{equation*}
\begin{equation*} \label{eq: A_2-relation} \tag{$A_2$}
	\mu^1(\mu^2(a_2, a_1)) + \mu^2\left(a_2, \mu^1(a_1)\right) + (-1)^{|a_1| - 1}\mu^2\left(\mu^1(a_2), a_1\right) = 0;
\end{equation*}
\begin{equation*} \label{eq: A_3-relation} \tag{$A_3$}
\begin{aligned}
		&\mu^2\left(a_3, \mu^2(a_2, a_1)\right) +
		(-1)^{|a_1| - 1}\mu^2\left(\mu^2(a_3, a_2), a_1\right) +
		\mu^1\left(\mu^3(a_3, a_2, a_1)\right)\\
		&+\mu^3\left(a_3, a_2, \mu^1(a_1)\right) +
		(-1)^{|a_1| - 1}\mu^3\left(a_3, \mu^1(a_2), a_1\right) +\\
		&(-1)^{|a_1| + |a_2| - 2}\mu^3\left(\mu^1(a_3), a_2, a_1\right) =
		0.
\end{aligned}
\end{equation*}

These relations imply that $a \cdot b = (-1)^{|b|} \mu^2(a, b)$
satisfies the graded Leibniz rule with respect to the differential
$da = (-1)^{|a|}\mu^1(a)$ and descends to an associative product on the cohomology of $(\A, d)$.

In curved $A_\infty$-algebras without any additional structure many natural notions (such as bounding cochains, which we define in \cref{subsubs: bounding cochains}) do not make sense since they involve infinite sums.
The additional structure we use is the following.

\begin{defn} \label{d: structured A_infty-algebra}
A \emph{structured module} is a graded $\k$-module $\A$ together with the data of
\begin{itemize}
	\item a finite set $Z$ and a graded $\k$-module $\A(z)$ for each $z \in Z$;
	\item a \emph{degree} function $|\cdot|\colon Z \to \ZZ$;
	\item a decomposition
	$
		\A = \bigoplus_{z \in Z} \A(z)[-|z|].
	$
\end{itemize}

The data of a structured module $\A$ and a curved $A_\infty$-algebra $(\A, \mu)$ is called a \emph{structured $A_\infty$-algebra} if there exists an \emph{action} function $A\colon Z \to \RR_{> 0}$ such that if the entry
$$
	\mu^{z^0; z^n, \ldots, z^1}\colon \A(z^n) \otimes \ldots \otimes \A(z^1) \to \A(z^0)
$$
of $\mu^n$ is non-zero, then $A(z^0) \ge A(z^n) + \ldots + A(z^1)$ and the inequality is strict except possibly in the following two cases:
\begin{enumerate}
	\item $n = 1$ and $z^0 = z^1$;
	\item $n = 2$.
\end{enumerate}
\end{defn}

Note that the action function is not part of the data in \cref{d: structured A_infty-algebra}, we explain the reason for that in \cref{r: isomorphism does not preserve action}.
Observe that in a structured $A_\infty$-algebra $\mu^n = 0$ for $n \gg 0$ and moreover any nested expression involving $\mu$ with sufficiently many occurrences of $\mu$ or sufficiently many inputs vanishes.

Recall that a curved $A_\infty$-algebra $(\A, \mu)$ is called \emph{nilpotent} if it possesses a filtration $\A = F^1\A \supset F^2\A \supset \ldots \supset F^N\A = 0$ such that
$$
	\mu^n(F^{i_1}\A \otimes \ldots \otimes F^{i_n}\A) \subset F^{i_1 + \ldots + i_n}\A.
$$
The following example shows that a nilpotent $A_\infty$-algebra supported on a structured module is not necessarily structured.

\begin{example}
	Let $\k = \ZZ[x]/x^2$, $Z = \{a, b\}$, $\A(a) = \A(b) = \k$, $|b| = |a| + 1 = 0$ and the only non-zero entries of $\mu$ be
	$$
		\mu^{b; a} = x \text{ and } \mu^{a; b, b, b} = x.
	$$
	The $A_\infty$-relations are satisfied since they are quadratic and $x^2 = 0$.
	$\A$ is nilpotent as witnessed by the filtration
	$$
		\A(a)[1] \oplus \A(b) \supset x\A(a)[1] \supset x\A(a)[1] \supset 0.
	$$
	However, $\A$ is not structured since we would have $A(b) > A(a) > 3A(b) > 0$ if it was.
\end{example}

\subsubsection{Obstruction algebra}
\label{subsubs: Obstruction algebra}
Let $\L = (L=(I, \gamma, P, \gr), \M)$ be a Lagrangian pre-brane in the weighted Liouville surface $(\S; \C)$.
We define the \emph{obstruction algebra} $\A_\L$, which is a structured $A_\infty$-algebra over $\k$.

Let $Z$ be the union of the set $X$ of double points of $L$ and the set $Y^{>0}$ of Reeb chords of positive action between the ends of $L$, as in \cref{subsubs: Double points and Reeb chords}.
Recall that $Z$ is a finite set by \cref{l: Y^>C is finite}.
We defined the action function $A\colon Z \to \RR_{>0}$ in \cref{subsubs: Double points and Reeb chords} and the degree function $|\cdot|\colon Z \to \ZZ$ in \cref{subs: Grading}.
As structured $\k$-module we define
\begin{equation} \label{eq: def of A_L}
	\A_\L = \bigoplus_{z \in Z} \A_\L(z)[-|z|]
\end{equation}
where
\begin{equation} \label{eq: def of A_L(z)}
	\A_\L(z) = \begin{cases}
		\C_{z}(\M(z_+), \M(z_-)) & \text{if } z \text{ is a double point};\\
		\C_{\gamma(z_-)}(\C(z)\M(z_+), \M(z_-)) & \text{if } z \text{ is a Reeb chord}.
	\end{cases}
\end{equation}
To avoid chasing cases in \cref{eq: def of A_L(z)} we will write $\C(z)$ for the identity functor when $z$ is a double point, so $\A_\L(z) = \C_{\gamma(z_-)}(\C(z)\M(z_+), \M(z_-))$ for any $z \in Z$.

The structure maps $\mu^n\colon \A_\L^{\otimes n} \to \A_\L[2 - n]$ are defined entry-wise.
For $n \ge 0$ and $z^0, \ldots, z^n \in Z_L$
the entry
\begin{equation} \label{eq: mu^z^0 z^n ... z^1}
	\mu^{z^0; z^n, \ldots, z^1}\colon \A_{\L}(z^n) \otimes \ldots \otimes \A_{\L}(z^1) \to \A_{\L}(z^0)
\end{equation}
is given by
\begin{equation} \label{eq: mu^z^0 z^0}
	\mu^{z^0; z^0}(\ph) = (-1)^{|\ph| + |z^0|}d_{\A_{\L}(z^0)}(\ph)
\end{equation}
if $n = 1$ and $z^1 = z^0$, and in all other cases by
\begin{equation} \label{eq: mu from polygons}
	\mu^{z^0; z^n, \ldots, z^1}(\ph^n, \ldots, \ph^1) = (-1)^{\sum_{j > k} (|z^j| - 1)|\ph^k|} \sum_{\bu \in \Delta^0(L; z^0; z^n, \ldots, z^1)}  \mu^{\bu}(\ph^n, \ldots, \ph^1) \cdot \prod_{i = 0}^n (\eps^{\bu}_i)^{|z^i|}
\end{equation}
where for $\bu \in \Delta^0(L; z^0; z^n, \ldots, z^1)$ the map $\mu^{\bu}$ and the signs $\eps^{\bu}_i$ are defined as explained below.

Note that $\Delta^0(L; z^0; z^n, \ldots, z^1)$ is non-empty only if $z^0_+ \in C_0, \, z^0_- \in C_n$ and $z^i_+ \in C_{i - 1}, \, z^i_- \in C_i$ for $i = 1, \ldots, n$ for some components $C_0, \ldots, C_n$ of $I$.
We first describe a choice of orientation of $C_1, \ldots, C_n$ induced by $\bu$ that extends the one from Seidel~\cite{S08}*{Figure 13.1}.
The sign $\eps^\bu_i$ is equal to $1$ (resp.\ $(-1)$) if the orientation of the component $C$ with $z^i_- \in C$ as a component of $I$ coincides with (resp.\ is opposite to) the orientation of $C$ induced by $\bu$.

\begin{defn}~\label{d: orientation induced by u}
	\begin{itemize}
		\item Suppose $\bu$ is a degenerate polygon, i.e.\ $n = 2$, $z^0, z^1, z^2 \in Y^{>0}$ and $z^0 =  z^2z^1$.
		The \emph{orientation induced by $\bu$} of $C_2$ (resp.\ $C_1$) is orientation induced by a tangent vector to $\gamma(C_2)$ at $\gamma(z^0_-) \in \pd S$ pointing inside $S$ (resp.\ to $\gamma(C_1)$ at $\gamma(z^1_-) \in \pd S$ pointing outside $S$).
		\item Suppose $\bu = (\Pi, u)$ is a polygon bounded by $L$.
		For $i = 1, \ldots, n$ let the orientation of $C_i$ \emph{induced by $\bu$} be the orientation of $\pd \Pi$ if there is a vertex of $\pd \Pi$ corresponding to $z^i_-$, or if this is not the case, i.e.\ we need to determine the orientation of $C_n$ in the case (2b) of \cref{d: moduli spaces of polygons}, then the orientation of $C_n$ induced by a tangent vector to $\gamma(C_n)$ at $\gamma(z^n_-) \in \pd S$ pointing inside $S$.
	\end{itemize}

	See the blue arrows in \cref{fig: Examples of polygons bounded by L} for examples of orientations defined above.
\end{defn}

Finally, given $\ph^i \in \C_{\gamma(z^i_-)}(\C(z^i)\M(z^i_+), \M(z^i_-))$
for $i = 1, \ldots, n$, we define $\mu^{\bu}(\ph^n, \ldots, \ph^1)$ as the composition along a 2-morphism (see \cref{subsubs: Composition along a 2-morphism})
\begin{equation} \label{eq: def of mu^u}
	\begin{aligned}
\C(z^0; &\gamma(a^n), &z^n, &\ldots, &\gamma(a^1), &z^1, &\gamma(a^0)&; \sigma_{\bu})\\
(&\M(a^n), &\ph^n, &\ldots, &\M(a^1), &\ph^1, &\M(a^0)&)
\end{aligned}
\end{equation}
where the curves $a^i = a^i_{\bu}$ for $i = 0, \ldots, n$, the curve $\gamma_{\bu}$, and the 2-morphism $\sigma_\bu \colon [z^0] \Rightarrow [\gamma_\bu]$ in $\P^+_2(S, Q)$ are defined in \cref{subsubs: 2-morphism induced by a polygon}.

We prove that $\A_{\L}$ is indeed a structured $A_\infty$-algebra in \cref{p: obstruction algebra is a nilpotent A infinity algebra}.

\begin{remark} \label{r: obstruction algebra and Chekanov-Eliashberg algebra}
	Assume that $\A_\L(z)$ is a finitely generated projective graded $\k$-module for each $z \in Z$.
	Applying the component-wise $\k$-linear dual to the bar construction on the obstruction algebra $(\A_{\L}, \mu)$ we get a differential graded algebra structure on
	$$
		B\A_{\L}^{\lor} = \bigoplus_{d = 0}^\infty \left(\A_{\L}^{\lor}[-1]\right)^{\otimes d}.
	$$
	More concretely, the algebra structure comes from the tensor product and the differential is given by
	$$
	\pd (az^0) = \pm (\mu^{z^0; z^n, \ldots, z^1})^{\lor}(a) \cdot z^1\ldots z^n \text{ for } z^0, \ldots, z^n \in Z_L.
	$$
	and the graded Leibniz rule.
	The $A_\infty$-relations translate to $\pd^2 = 0$.

	What we get is a version of the Chekanov--Eliashberg algebra.
	Namely, if $S$ is a disc, $\k = \FF_2$, $\C$ is the trivial 2-functor with fibre $\k$ (thought of as dg-category with one object), then
	\begin{itemize}
		\item when $Q = \varnothing$ and the domain of $L$ is closed, we obtain the DGA of a Legendrian link;
		\item when $D = \varnothing$, $|\Lambda| = 2$, the domain of $L$ is a segment, and the two ends of $L$ lie in different components of $\pd S \setminus \Lambda$,  we obtain the DGA of a long Legendrian knot
	\end{itemize}
	both defined by Chekanov~\cite{C02}.
	More generally, modulo sign conventions, when $Q = \varnothing$, $\k$ is any commutative ring, $\C$ is the trivial 2-functor with fibre $\k$, and the domain of $L$ is closed, we obtain the DGA defined by Ekholm--Etnyre--Sullivan~\cites{EES05,EES07}.
\end{remark}

\subsection{Invariance of the obstruction algebra}
\label{subs: Invariance of the obstruction algebra}
We prove that the obstruction algebra is invariant under isotopies of slightly degenerate Lagrangian pre-branes (see \cref{d: slightly degenerate Lagrangian}).
For this first observe that in \cref{subsubs: Obstruction algebra} the definition goes through if we allow Reeb chords of zero length, so we define the obstruction algebra of a slightly degenerate Lagrangian pre-brane by treating the point-chords as Reeb chords.
Throughout the subsection we work with a fixed weighted Liouville surface ${(\S = (S, Q, \eta, \lambda); \C)}$.

In \cref{subsubs: Naive isomorphisms} we define na\"ive isomorphisms, the sense in which the obstruction algebra is invariant.
In \cref{subsubs: isotopy induces isomorphism of obstruction algebras} we transport the grading and the local system of a Lagrangian pre-brane along an isotopy of slightly degenerate Lagrangians and prove that the obstruction algebra is invariant in \cref{p: isotopy induces isomorphism of obstruction algebras}.

\begin{remark}
	We expect that, similarly to the classical Chekanov--Eliashberg algebra, the obstruction algebra is invariant with respect to more general isotopies of the Legendrian lift of a Lagrangian pre-brane.
	The invariance in that case should be understood as up to an $A_\infty$-isomorphism.
	We will not treat it in this paper.
\end{remark}

\subsubsection{Na\"ive isomorphisms}
\label{subsubs: Naive isomorphisms}
By invariance of the obstruction algebra we mean that all data defining it is preserved in the na\"ive sense.
Namely, a \emph{na\"ive morphism} of two structured modules with higher products (which are not required to satisfy the $A_\infty$-relations) $(\A_i, Z_i, \mu_i^0, \mu_i^1, \ldots), \, i = 0, 1$ is
\begin{itemize}
	\item a map $\psi\colon Z_0 \to Z_1$;
	\item a morphism of graded $\k$-modules $\F^z\colon \A_0(z) \to \A_1(\psi(z))$ for each $z \in Z_0$
\end{itemize}
such that
\begin{itemize}
	\item $|\psi(z)| = |z|$ for each $z \in Z_0$;
	\item
	$\mu_1^{\psi(z^0); \psi(z^n), \ldots, \psi(z^1)} (\F^{z^n}(-), \ldots, \F^{z^1}(-)) = \F^{z^0} \circ \mu_0^{z^0; z^n, \ldots, z^1}(-, \ldots, -)$
	for all $z^0, \ldots, z^n \in Z_0$.
\end{itemize}
Structured modules with higher products together with na\"ive morphisms form a category under the component-wise composition and a \emph{na\"ive isomorphism} is an isomorphism in this category.
We have the following standard observation.

\begin{lemma} \label{l: characterisation of naive isomorphism}
	A na\"ive morphism $(\psi, (\F^z)_{z \in Z_0})$ is an isomorphism if and only if $\psi$ is a bijection and $\F^z$ is an isomorphism for each $z \in Z_0$.
\end{lemma}

\subsubsection{Invariance under isotopies}
\label{subsubs: isotopy induces isomorphism of obstruction algebras}
Let $\L_0 = (I, \gamma_0, P_0, \gr_0, \M_0)$ be a Lagrangian pre-brane and $L_1 = (I, \gamma_1, P_1)$ a Lagrangian isotopic to $L_0 = (I, \gamma_0, P_0)$ as a slightly degenerate Lagrangian, see \cref{subsubs: Isotopy of slightly degenerate Lagrangians}.

We want to transport $\gr_0$ and $\M_0$ to get a Lagrangian pre-brane $\L_1 = (L_1, \gr_1, \M_1)$ and prove that obstruction algebras $(\A_{\L_0}, \mu_0)$ and $(\A_{\L_1}, \mu_1)$ are isomorphic.
Since a generic isotopy intersects the slightly degenerate locus transversely, it is enough to consider the case when $L_s = (I, \gamma_s, P_s), s \in [0, 1]$ is a path of slightly degenerate Lagrangians such that $L_s$ is a Lagrangian for $s \ne 1$ and $L_1$ is either a Lagrangian or has a unique point-chord $z^1 \in XY_{L_1}$, i.e.\ we are in the setting of \cref{l: family of maps from isotopy of slighly degenerate Legendrians}.

By \cref{l: gradings on an immersed curve} the grading $\gr_0$ transports to a grading $\gr_1$ on $L_1$.
We define the local system $\M_1$ by $\M_1(p) = \C[s \mapsto \gamma_s(p)]\M_0(p)$ for $p \in I$ and for a path $a\colon [0, 1] \to I$, the transition map
$$
\M_1(a) \in \C_{\gamma_1a(1)}(\C(\gamma_1a)\M_1(a(0)), \M_1(a(1)))
$$
is the composition along a 2-morphism (see \cref{subsubs: Composition along a 2-morphism})
\begin{equation} \label{eq: def of M_1}
	\C(\gamma_1 a; s \mapsto \gamma_s a(1), \gamma_0 a, s \mapsto \gamma_{1 - s} a(0); \sigma_a)(\id_{\M_1(a(1))}, \M_0(a), \id_{\M_0(a(0))})
\end{equation}
where $\sigma_a$ is the 2-morphism in $\P^+_2(S, Q)$ given by $\gamma_s a (t)$.

A na\"ive isomorphism $(\psi, \F)\colon \A_{\L_0} \to \A_{\L_1}$ is constructed as follows.
Recall from \cref{l: family of maps from isotopy of slighly degenerate Legendrians} that we have the Reeb chords $y^s_0, s \in [0, 1]$ of $L_s$ collapsing to the point-chord $z^1$ and that there are two cases.

\emph{Case 1.} $A(y_0^s) > 0$ for all $s \in [0, 1]$.
In this case we can track each double point $x^0 \in X_{L_0}$ to obtain a family $x^s \in X_{L_s}, s \in [0, 1]$ and can track each Reeb chord $y^0 \in Y_{L_0} \setminus \{y^0_0\}$ to get a family $y^s \in Y_{L_s}, s \in [0, 1]$, and moreover $A(y^s)$ does not change sign in the family.
We then define $\psi(x^0) = x^1$ for each $x^0 \in X_{L_0}$, $\psi(y^0) = y^1$ for each $y^0 \in Y_{L_0}^{>0} \setminus \{y^0_0\}$ and $\psi(y_0^0) = z^1$.

\emph{Case 2.} $A(y_0^s) < 0$ for all $s \in [0, 1]$.
In this case there is a family of double points $x^s_0 \in X_{L_s}$ that collapses to the point-chord $z^1$.
We can track each double point $x^0 \in X_{L_0} \setminus \{x_0^0\}$ to obtain a family $x^s \in X_{L_s}, s \in [0, 1]$ and can track each Reeb chord $y^0 \in Y_{L_0} \setminus \{y^0_0\}$ to get a family $y^s \in Y_{L_s}, s \in [0, 1]$, and moreover $A(y^s)$ does not change sign in the family.
We then define $\psi(x^0) = x^1$ for each $x^0 \in X_{L_0} \setminus \{x_0^0\}$, $\psi(y^0) = y^1$ for each $y^0 \in Y_{L_0}^{>0}$ and $\psi(x_0^0) = z^1$.

Fix $z^0 \in Z_0$.
From the reasoning above we can always track $z^0$ to obtain the family ${z^s \in Z_s, \, s \in [0, 1]}$.
To define
$$
\F^{z^0}\colon \C_{\gamma_0(z^0_-)}(\C(z^0)\M_0(z^0_+), \M_0(z^0_-)) \to \C_{\gamma_1(z^1_-)}(\C(z^1)\M_1(z^1_+), \M_1(z^1_-))
$$
first take into account the moving of $z^s_{\pm}$ in $I$.
Namely, consider the morphism
$$
\begin{aligned}
\M(s \mapsto \gamma_s(z^s_{-})) =
\C\big(s \mapsto \gamma_s(z^s_{-})&; s \mapsto \gamma_s(z^1_-)&,& s \mapsto \gamma_0(z^s_{-})&;& [\gamma_s(z^t_{-}), s \le t]\big)\\
&\big(\id&,& \M_0(s \mapsto z^s_{-})&\big)&
\end{aligned}
$$
in
$\C_{\gamma_1(z^1_{-})}(\C(s \mapsto \gamma_s(z^s_{-}))\M_0(z^0_-), \M_1(z^1_-))$
and the morphism
$$
\begin{aligned}
\M(s \mapsto \gamma_{1 - s}(z^{1 - s}_{+})) =
\C\big(s \mapsto \gamma_{1 - s}(z^{1 - s}_{+})&; s \mapsto \gamma_{0}(z^{1 - s}_{+})&,& s \mapsto \gamma_{1 - s}(z^{1}_{+})&;& [\gamma_{1 - t}(z^{1 - s}_{+}), s \le t]\big)\\
&\big(\M_0(s \mapsto z^{1 - s}_{+})&,& \id&\big)&
\end{aligned}
$$
in
$\C_{\gamma_0(z^0_{+})}(\C(s \mapsto \gamma_{1 - s}(z^{1 - s}_{+}))\M_1(z^1_+),\M_0(z^0_+))$.
Note that \cref{l: properties of composition along a 2-path} implies that they are both isomorphisms since the $2$-morphisms $[\gamma_s(z^t_{-}), s \le t]$ and $[\gamma_{1 - t}(z^{1 - s}_{+}), s \le t]$ are invertible in $\P^+_2(\S)$.

Now define $\F^{z^0}(\ph)$ to be
$$
\begin{aligned}
\C\big(z^1&; s \mapsto \gamma_s(z^s_{-})&,& z^0&,& s \mapsto \gamma_{1 - s}(z^{1 - s}_{+})&;& [(s, t) \mapsto z^s (t)]\big)\\
&\big(\M(s \mapsto \gamma_s(z^s_{-}))&,& \ph&,& \M(s \mapsto \gamma_{1 - s}(z^{1 - s}_{+}))&\big)&.
\end{aligned}
$$
Note that $\sigma_{z^0} = [(s, t) \mapsto z^s (t)]$ is thin (because either it does not depend on $t$ or it is contained in $\pd S$), so in fact
$[(s \mapsto \gamma_s(z^s_-))z^0(s \mapsto \gamma_{1 - s}(z^{1 - s}_+))] = [z^1]$
and
\begin{equation} \label{eq: F^z^0(ph)}
	\F^{z^0}(\ph) = \M(s \mapsto \gamma_s(z^s_{-})) \cdot \C(s \mapsto \gamma_s(z^s_{-}))\ph \cdot \C((s \mapsto \gamma_s(z^s_{-}))z^0)\M(s \mapsto \gamma_{1 - s}(z^{1 - s}_{+})).
\end{equation}

\begin{prop} \label{p: isotopy induces isomorphism of obstruction algebras}
	$(\psi, \F)\colon \A_{\L_0} \to \A_{\L_1}$ defined above is a na\"ive isomorphism.
\end{prop}
\begin{remark} \label{r: isomorphism does not preserve action}
	The isomorphism in \cref{p: isotopy induces isomorphism of obstruction algebras} does not necessarily preserve the action, i.e.\ if the double points or Reeb chords $x^0, y^0 \in Z_0$ correspond to $x^1, y^1 \in Z_1$, then it is possible to have $A(x^0) < A(y^0)$, but $A(x^1) > A(y^1)$.
\end{remark}

\begin{proof}
$\psi$ is a bijection by construction, and it preserves degree by \cref{l: gradings on an immersed curve,l: degree is well-defined for slightly degenerate Lagrangian}.
As $\M(s \mapsto \gamma_s(z^s_{-}))$ and $\M(s \mapsto \gamma_{1 - s}(z^{1 - s}_{+}))$ are isomorphisms, \cref{eq: F^z^0(ph)} implies that $\F^z$ is an isomorphism of dg-modules.
Therefore, by \cref{l: characterisation of naive isomorphism} it remains to prove that the higher products given by \cref{eq: mu from polygons} are preserved.

Fix $z_0^0, \ldots, z_n^0 \in Z_0$ and note that the sign before the summation in \cref{eq: mu from polygons} is the same since $\psi$ and $\F^z$ preserve the degree.

Recall that by \cref{l: family of maps from isotopy of slighly degenerate Legendrians} we have smooth families of maps
$\Phi^s \colon S \to S,\,  \psi^s\colon I \to I,\, s \in [0, 1]$
such that $\Phi^0 = \id$, $\psi^0 = \id$, $\Phi^s \circ \gamma_0 = \gamma_s \circ \psi^s$ for $s \in [0, 1]$, $\Phi^s \in \Diff(\S), \, \psi^s \in \Diff(I)$ for $s \in [0, 1)$, and
$\res{\Phi^1}{S \setminus (\Phi^1)^{-1}(z^1)}, \, \res{\psi^1}{I \setminus (\psi^1)^{-1}\{z^1_+, z^1_-\}}$
are diffeomorphisms onto their images.

\begin{claim}
	$\Phi^1$ establishes a bijection
	$$
		\Delta^0(L_0; z^0_0; z^0_n, \ldots, z^0_1) \to \Delta^0(L_1; z^1_0; z^1_n, \ldots, z^1_1), \, \bu^0 \mapsto \bu^1
	$$
\end{claim}
\begin{proof}[Proof of Claim]
	If none of $z^1_0, \ldots, z^1_n$ is the point-chord $z^1$, then there is nothing to prove.
	If we are in Case 1 above, i.e.\ the family of Reeb chord $y^s_0$ collapsing to $z^1$ have positive action, then there is also nothing to prove since we treat the point-chord as a Reeb chord.
	Therefore suppose we are in Case 2, and let $x^s_0$ be the family of double points collapsing to $z^1$.

	If $\bu^1$ is a degenerate polygon, i.e.\ $z^1_0 = z^1_2z^1_1$, then, since $z^1$ has no non-trivial subchords, either $z^1_1 = z^1$ or $z^1_2 = z^1$.
	In either case, $\bu^0$ is the triangle with vertices $x^0_0$, $\gamma_0(y^0_{0+})$, $\gamma_0(y^0_{0-})$, see \cref{fig: Delta(L_1 z_0 z_2z^1),fig: Delta(L_1 z_0 z^1z_1)}.

	\begin{figure}[!htbp]
	\includegraphics[width=0.8\linewidth]{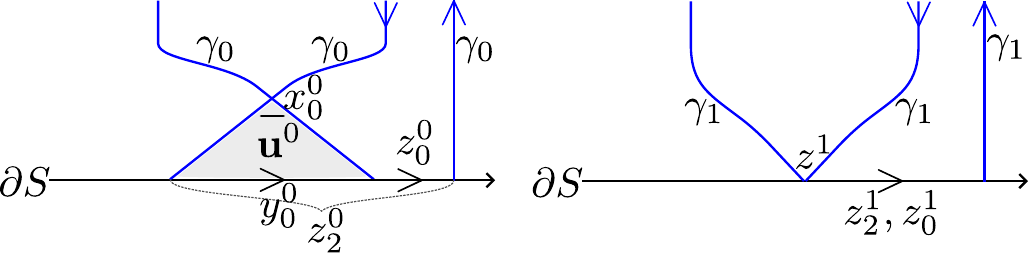}
	\caption{$\Delta(L_0; z^0_0; z_2^0, x^0_0) \to \Delta(L_1; z^1_0; z^1_2, z^1)$.}
	\label{fig: Delta(L_1 z_0 z_2z^1)}
	\end{figure}

	\begin{figure}[!htbp]
	\includegraphics[width=0.8\linewidth]{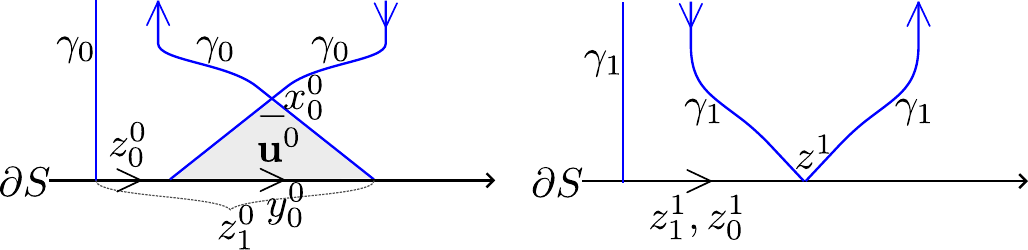}
	\caption{$\Delta(L_0; z^0_0; x^0_0, z_1^0) \to \Delta(L_1; z^1_0; z^1, z^1_1)$.}
	\label{fig: Delta(L_1 z_0 z^1z_1)}
	\end{figure}

	Now suppose $\bu^1$ is a non-degenerate polygon.
	If $z^1 = z^1_i$, then $z^0_i = x^0_0$ and the polygons $\bu^0$ and $\bu^1$ look like in \cref{fig: Delta(L_1 z_0 ...z^1...)}.

	\begin{figure}[!htbp]
	\includegraphics[width=0.8\linewidth]{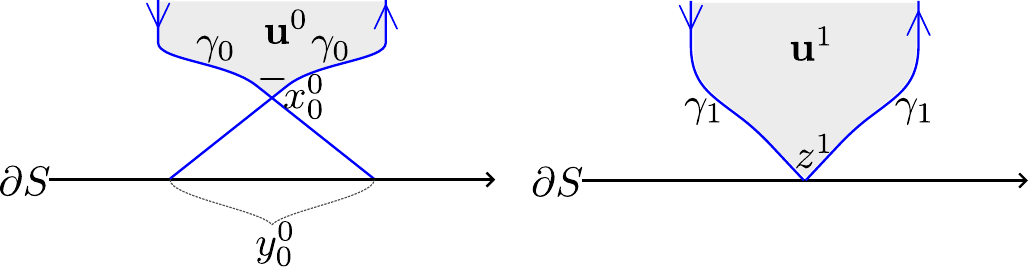}
	\caption{$\Delta(L_0; z_0^0; \ldots, x^0_0, \ldots) \to \Delta(L_1; z_0^1; \ldots, z^1, \ldots)$.}
	\label{fig: Delta(L_1 z_0 ...z^1...)}
	\end{figure}

	Finally, if $z^1 = z^1_0$, then, depending on whether we are in the case (2a) or case (2b) of \cref{d: moduli spaces of polygons}, the polygons $\bu^0$ and $\bu^1$ look like in \cref{fig: Delta(L_1 z^1 ...)right,fig: Delta(L_1 z^1 ...)left}.

	\begin{figure}[!htbp]
	\includegraphics[width=0.8\linewidth]{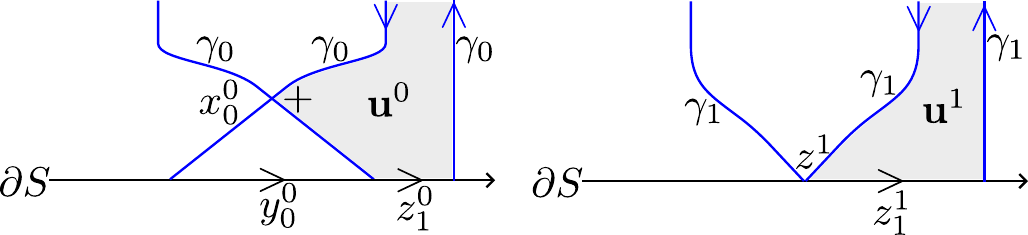}
	\caption{$\Delta(L_0; x^0_0; \ldots, z^0_1) \to \Delta(L_1; z^1; \ldots, z^1_1)$.}
	\label{fig: Delta(L_1 z^1 ...)right}
	\end{figure}

	\begin{figure}[!htbp]
	\includegraphics[width=0.8\linewidth]{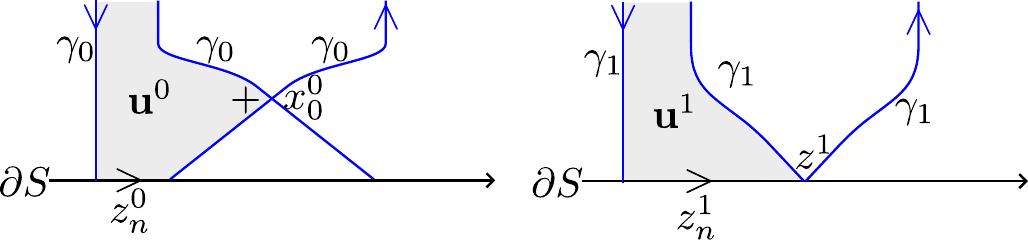}
	\caption{$\Delta(L_0; x^0_0; z^0_n, \ldots) \to \Delta(L_1; z^1; z^1_n, \ldots)$.}
	\label{fig: Delta(L_1 z^1 ...)left}
	\end{figure}

	We organised the cases above depending on what happens at $s = 1$.
	Observe that at $s = 0$ the cases depend on what corner of $x_0^0$ is covered by $\bu^0$, and we covered all possibilities.
	This observation finishes the proof of the claim.
\end{proof}

We therefore have a bijection between the summands in \cref{eq: mu from polygons} for $\L_0$ and $\L_1$.
Fix $\bu^0 \in \Delta^0(L_0; z^0_0; z^0_n, \ldots, z^0_1)$.
For the signs we have
$$
\prod_{i = 0}^n \left(\eps_i^{\bu^0}\right)^{|z^0_i|} = \prod_{i = 0}^n \left(\eps_i^{\bu^1}\right)^{|z^1_i|}
$$
 since $\Phi^1$ is orientation-preserving and in \cref{fig: Delta(L_1 z_0 z_2z^1),fig: Delta(L_1 z_0 z^1z_1),fig: Delta(L_1 z_0 ...z^1...),fig: Delta(L_1 z^1 ...)right,fig: Delta(L_1 z^1 ...)left} the orientations from \cref{d: orientation induced by u} are the same except for \cref{fig: Delta(L_1 z_0 z^1z_1)}, for which we have
$$
\prod_{i = 0}^2 \left(\frac{\eps_i^{\bu^1}}{\eps_i^{\bu^0}}\right)^{|z^0_i|} = (-1)^{|z^0_0|}(-1)^{|z^0_1|}(-1)^{|z^0_2|} = (-1)^{|z^0_1| + |z^0_2|}(-1)^{|z^0_1|}(-1)^{|z^0_2|} = 1.
$$
It therefore remains to check that
\begin{equation} \label{eq: mu^u for isotopy}
	\mu^{\bu^1}(\F^{z_n^0}(\ph^n), \ldots, \F^{z_1^0}(\ph^1)) = \F^{z_0^0} \circ \mu^{\bu^0}(\ph^n, \ldots, \ph^1).
\end{equation}

This follows by \cref{l: properties of composition along a 2-path} from the isotopy given by $\Phi^s$ in \cref{fig: Isotopy of 2-morphism from isotopy}.

\begin{figure}[!htbp]
	\includegraphics[width=0.4\linewidth]{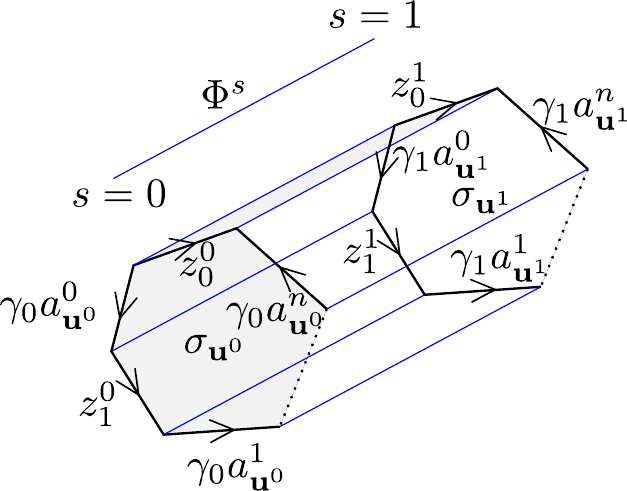}
	\caption{Isotopy of 2-morphism from $\Phi^s$.}
	\label{fig: Isotopy of 2-morphism from isotopy}
\end{figure}

Namely, the right-hand side of \cref{eq: mu^u for isotopy} is given by
\begin{equation} \label{eq: shaded part of the prism}
\begin{aligned}
	\C\Big(z^1_0; &[s \mapsto \gamma_s(z^s_{0-})], &\gamma_0(a^n_{\bu^0}), &\ldots, &\gamma_0(a^1_{\bu^0}), &z^0_1, &\gamma_0(a^0_{\bu^0}), &[s \mapsto \gamma_{1 - s}(z^{1 - s}_{0+})]; \sigma_{z_0^0} \cup \sigma_{\bu^0}\Big)\\
	\Big(&\M(s \mapsto \gamma_s(z^s_{0-})), &\M_0(a^n_{\bu^0}), &\ldots, &\M_0(a^1_{\bu^0}), &\ph^1, &\M_0(a^0_{\bu^0}), &\M(s \mapsto \gamma_{1 - s}(z^{1 - s}_{0+}))\Big),
\end{aligned}
\end{equation}

and the left-hand side of \cref{eq: mu^u for isotopy} is given by
\begin{equation} \label{eq: unshaded part of the prism}
\begin{aligned}
	\C\Big(z^1_0; &\gamma_1(a^n_{\bu^1}), &\ldots, &[s \mapsto \gamma_s(z^s_{1-})], &z^0_1, &[s \mapsto \gamma_{1 - s}(z^{1 - s}_{1+})], &\gamma_1(a^0_{\bu^1}); \sigma\Big)\\
	\Big(&\M_1(a^n_{\bu^1}), &\ldots, &\M(s \mapsto \gamma_s(z^s_{1-})), &\ph^1, &\M(s \mapsto \gamma_{1 - s}(z^{1 - s}_{1+})), &\M_1(a^0_{\bu^1})\Big)
\end{aligned}
\end{equation}

where
$\sigma = \sigma_{\bu^1} \cup \sigma_{z^0_1} \cup \ldots \cup \sigma_{z^0_n}$.
Expanding $\M$ and $\M_1$ by definition we see that \cref{eq: shaded part of the prism} is equal to the composition along the shaded part of the boundary of the prism in \cref{fig: Isotopy of 2-morphism from isotopy} and \cref{eq: unshaded part of the prism} is equal to the composition along the unshaded part of the boundary of the prism.
Thus the solid prism $\Phi^s \circ \sigma_{\bu^0}$ provides the required isotopy.
\end{proof}

\subsection{\texorpdfstring{$A_\infty$}{A∞}-relations}
\label{subs: Obstruction algebra is A_infty}

Let $\L = (L=(I, \gamma, P), \gr, \M)$ be a Lagrangian pre-brane in $(\S; \C)$.
Recall that in \cref{subsubs: Obstruction algebra} we defined the obstruction algebra $\A_\L$ of $\L$.

\begin{prop} \label{p: obstruction algebra is a nilpotent A infinity algebra}
The obstruction algebra $\A_\L$ is a structured $A_\infty$-algebra over $\k$.
\end{prop}

The argument is based on the decomposition of polygons with one concave corner by Chekanov~\cite{C02}*{Theorem 3.3}.
There are versions of this argument that take into account signs or local systems, see e.g.\ Abouzaid~\cite{A08}*{\Sec 3} for the case of immersed curves in closed surfaces that do not bound any teardrops and Etnyre--Ng~\cite{EN22}*{\Sec 3} for an overview of various flavours of Chekanov--Eliashberg DGAs.

\begin{proof}
	The set $Z = Z_L$ of double points of $L$ and Reeb chords between the ends of $L$ of positive action is finite by \cref{l: Y^>C is finite}.
	The sums defining the structure maps $\mu^n$ in \cref{eq: mu from polygons} are finite by \cref{l: moduli spaces are finite}(2).
	The degree of any $\mu^n, \, n \ge 0$ is equal to $2 - n$ by \cref{eq: |z^0| from |z^1| + ... + |z^n|}.
	The action property from \cref{d: structured A_infty-algebra} is satisfied by \cref{l: moduli spaces are finite}(1).
	It remains to check the $A_\infty$-relations \cref{eq: A-infinity relations}.
	By \cref{l: sufficiently self-wrapped model for Lagrangian,p: isotopy induces isomorphism of obstruction algebras} we can assume $L$ is sufficiently self-wrapped.

	Let $C_0, \ldots, C_n$ be (not necessarily distinct) connected components of $I$ and $x^0, \ldots, x^n \in X_L = Z_L$ be such that
	$$
		x^0_+ \in C_0, \, x^0_- \in C_n \text{ and } x^i_+ \in C_{i - 1}, \, x^i_- \in C_i \text{ for } i = 1, \ldots, n.
	$$
	We show that the entry
	\begin{equation} \label{eq: entry of A_infty relation}
	\A_{\L}(x^n) \otimes \ldots \otimes \A_{\L}(x^1) \to \A_{\L}(x^0)
	\end{equation}
	of the map in the left hand side of the $A_\infty$-relation \cref{eq: A-infinity relations} vanishes.

	We first handle the part of the $A_\infty$-relations coming from the dg-coefficients, i.e.\ given by \cref{eq: mu^z^0 z^0}.
	$\mu^{x^0; x^0}$ squares to $0$ since $d_{\A_\L(x^0)}$ is a differential.
	Other summands of \cref{eq: entry of A_infty relation} involving \cref{eq: mu^z^0 z^0} group into
	$$
	\begin{aligned}
		&\mu^{x^0; x^0} \circ \mu^{\bu}(\ph^n, \ldots, \ph^1) + \sum_{i = 1}^n (-1)^{\sum_{j < i}(|x^j| + |\ph^j| - 1) + \sum_{j > i}(|x^j| - 1)}\mu^{\bu}(\ph^n, \ldots, \mu^{x^i; x^i} (\ph^i), \ldots, \ph^1) = \\
		&(-1)^{\sum_{i = 1}^n (|\ph^i| + |x^i| - 1)} d(\C(x^0; \gamma^n_{\bu}, x^n, \ldots, \gamma^1_{\bu}, x^1, \gamma^0_{\bu}; \sigma_{\bu})(\M(\gamma^n_{\bu}), \ph^n, \ldots, \M(\gamma^1_{\bu}), \ph^1, \M(\gamma^0_{\bu}))) +\\
		&\sum_{i = 1}^n (-1)^{\sum_{j \le i} |\ph^j| + \sum_{j = 1}^n (|x^j| - 1) + 1} \C(x^0; \ldots, \gamma^i_{\bu}, x^i, \gamma^{i - 1}_{\bu}, \ldots; \sigma_{\bu})(\ldots, \M(\gamma^i_{\bu}), d\ph^i, \M(\gamma^{i - 1}_{\bu}), \ldots)
	\end{aligned}
	$$
	where $\bu \in \Delta^0(L; x^0; x^n, \ldots, x^1)$ and we omitted the common sign.
	Recall from \cref{eq: composition along 2-morphism} that the composition along a 2-morphism is a certain composition in $\C_{x^0}$.
	Since all $\M(\gamma^i_{\bu})$ are closed morphisms of degree $0$ and $\sigma_{\bu}$ is a closed morphism of dg-functors, the above expression vanishes by the Leibniz rule.

	Each other summand of the entry \cref{eq: entry of A_infty relation} is given by
	\begin{equation} \label{eq: summand of A-infinity relation}
		\pm \mu^{\bu_1}(\ph^n, \ldots, \ph^{k + j + 1}, \mu^{\bu_0}(\ph^{k + j}, \ldots, \ph^{k + 1}), \ph^k, \ldots, \ph^1)
	\end{equation}
	where
	\begin{itemize}
		\item $i, j, k \in \ZZ_{\ge 0}$ and $i + j + k = n$;
		\item $x$ is a double point of $L$;
		\item
		$\bu_0 = (\Pi_0, u_0) \in \Delta^0(L; x; x^{k + j}, \ldots, x^{k + 1})$;
		\item
		$\bu_1 = (\Pi_1, u_1) \in \Delta^0(L; x^0; x^n, \ldots, x^{k + j + 1}, x, x^{k}, \ldots, x^1)$.
	\end{itemize}
	Denote by $T$ the set of the tuples $\bt = ((i, j, k), x, \bu_0, \bu_1)$ as above.

	By \cite{C02}*{\Sec 7} there is an involution $\iota$ on $T$ and a $\iota$-invariant double cover $\theta\colon T \to \Delta^1(L; x^0; x^n, \ldots, x^1)$ such that
	\begin{itemize}
		\item if $\theta((i, j, k), x, \bu_0, \bu_1) = \bu = (\Pi, u)$, then the pullback of $\bu_0$ and $\bu_1$ to $\Pi$ via $u$ is diffeomorphic to one of the pictures in \cref{fig: Concave Polygon 1,fig: Concave Polygon 2,fig: Concave Polygon 3,fig: Concave Polygon 4,fig: Concave Polygon 5} where $\Pi$ has vertices at each point marked by $x^i$ for $i = 0, \ldots, n$, the corners which are not depicted are all negative;
		\item the involution $\iota$ swaps the pictures in \cref{fig: Concave Polygon 1} with that in \cref{fig: Concave Polygon 2}, the two pictures in \cref{fig: Concave Polygon 3}, and the pictures in \cref{fig: Concave Polygon 4} with that in \cref{fig: Concave Polygon 5}.
	\end{itemize}

	\begin{figure}[!htbp]
	\centering
	\begin{subfigure}[t]{0.24\linewidth}
		\centering
		\includegraphics[width=\linewidth]{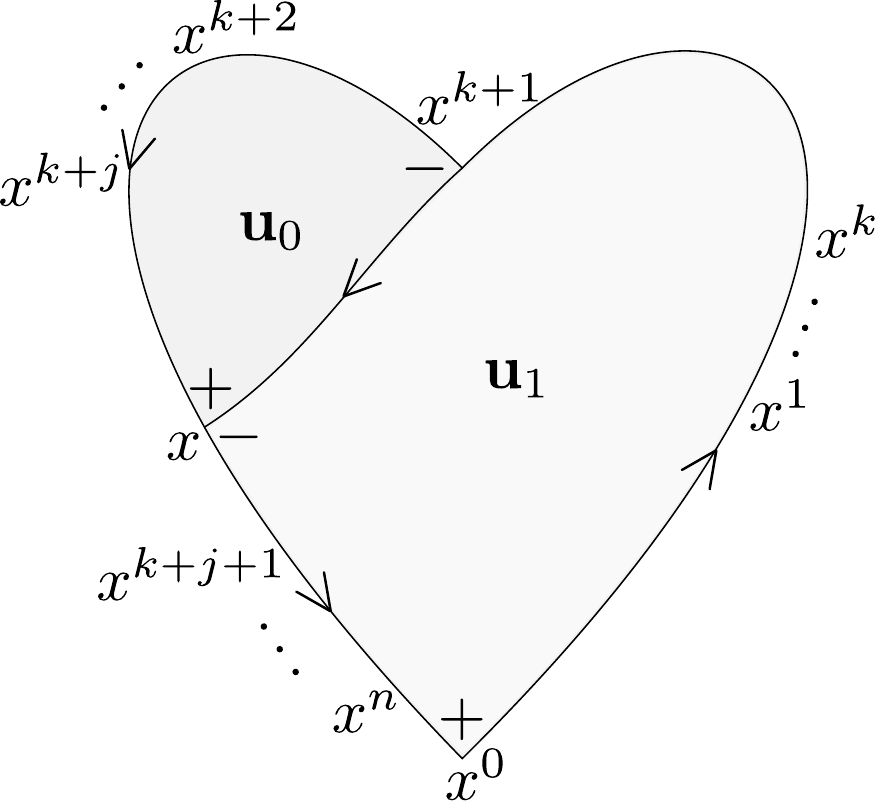}
		\caption{}
		\label{subfig: ConcavePolygon_1a}
	\end{subfigure}%
	\hfill%
	\begin{subfigure}[t]{0.25\linewidth}
		\centering
		\includegraphics[width=\linewidth]{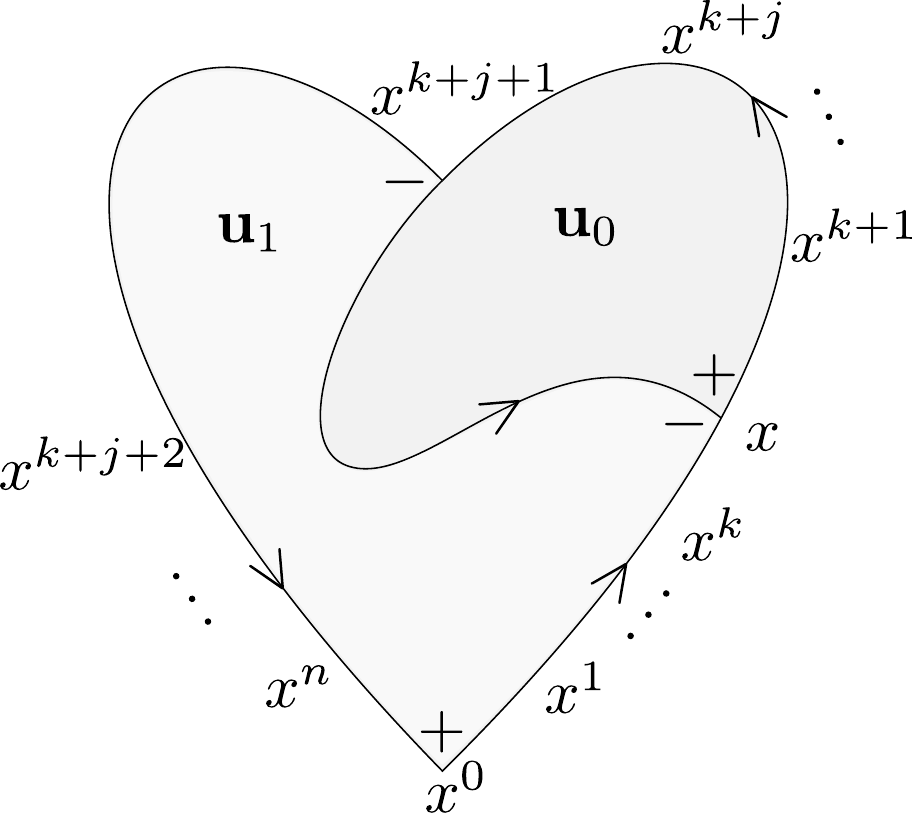}
		\caption{}
		\label{subfig: ConcavePolygon_1b}
	\end{subfigure}%
	\hfill%
	\begin{subfigure}[t]{0.24\linewidth}
		\centering
		\includegraphics[width=\linewidth]{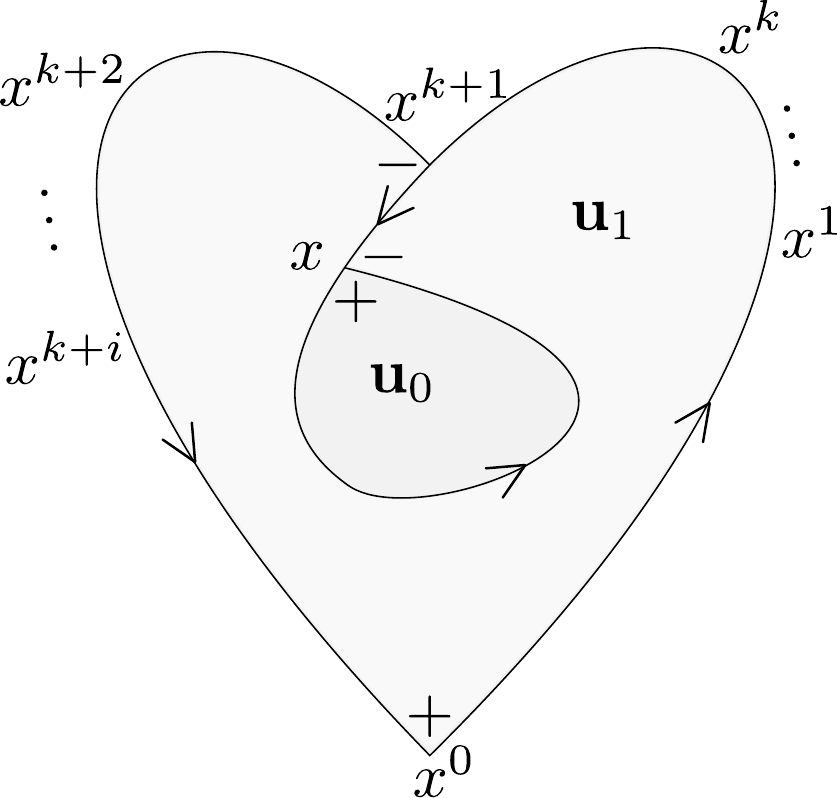}
		\caption{$j = 0$}
		\label{subfig: ConcavePolygon_1c}
	\end{subfigure}%
	\hfill%
	\begin{subfigure}[t]{0.24\linewidth}
		\centering
		\includegraphics[width=\linewidth]{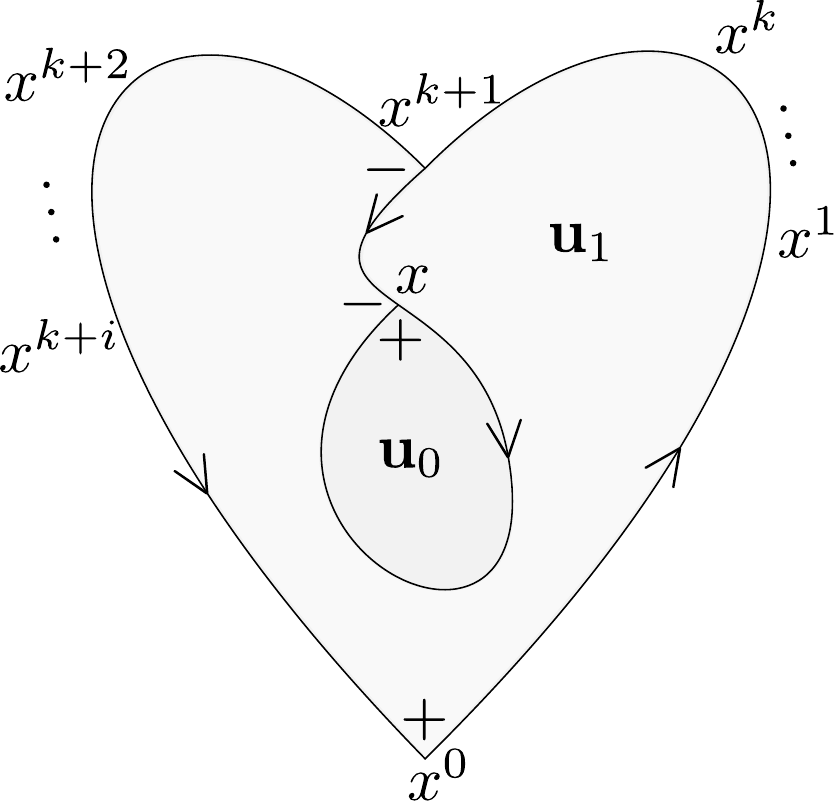}
		\caption{$j = 0$}
		\label{subfig: ConcavePolygon_1d}
	\end{subfigure}
	\caption{Extension of the incoming side of a negative concave corner.}
	\label{fig: Concave Polygon 1}
	\end{figure}

	\begin{figure}[!htbp]
	\centering
	\begin{subfigure}[t]{0.24\linewidth}
		\centering
		\includegraphics[width=\linewidth]{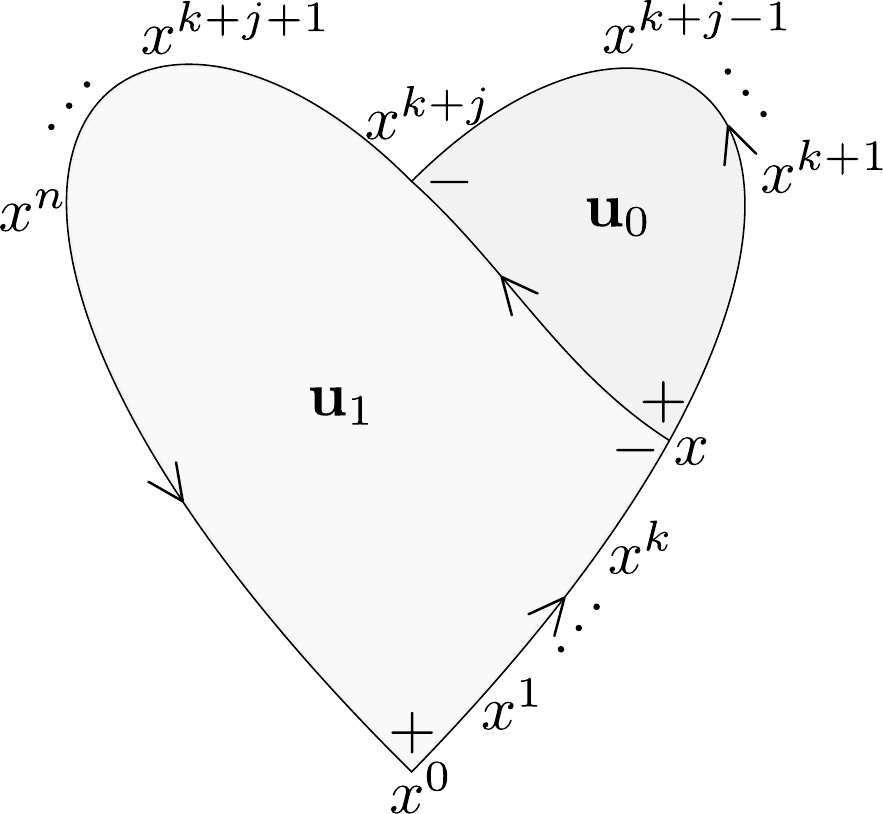}
		\caption{}
		\label{subfig: ConcavePolygon_2a}
	\end{subfigure}%
	\hfill%
	\begin{subfigure}[t]{0.24\linewidth}
		\centering
		\includegraphics[width=\linewidth]{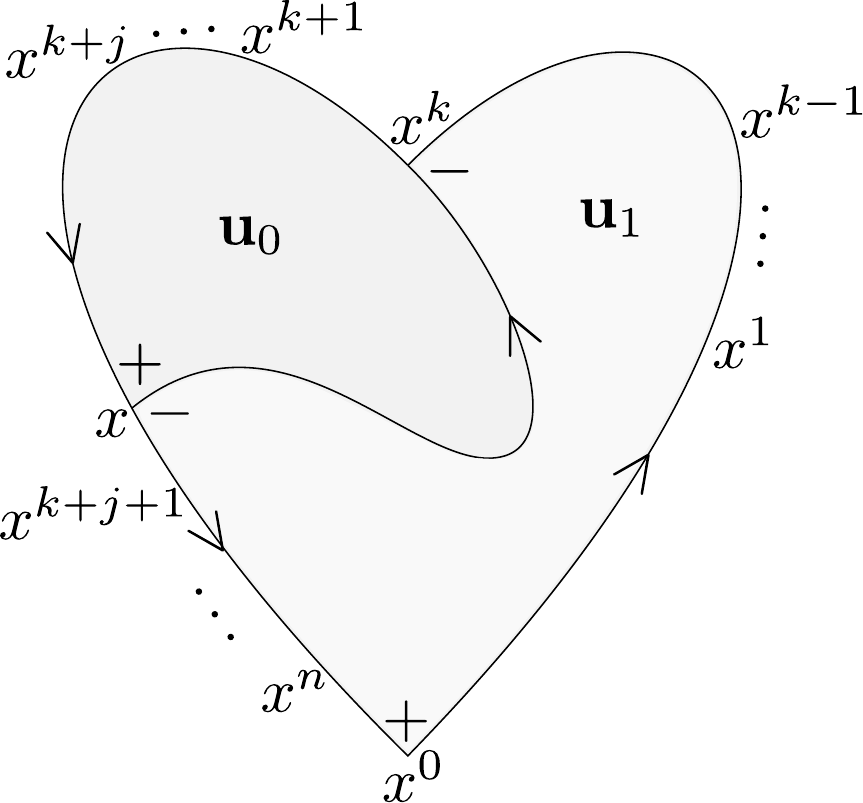}
		\caption{}
		\label{subfig: ConcavePolygon_2b}
	\end{subfigure}%
	\hfill%
	\begin{subfigure}[t]{0.24\linewidth}
		\centering
		\includegraphics[width=\linewidth]{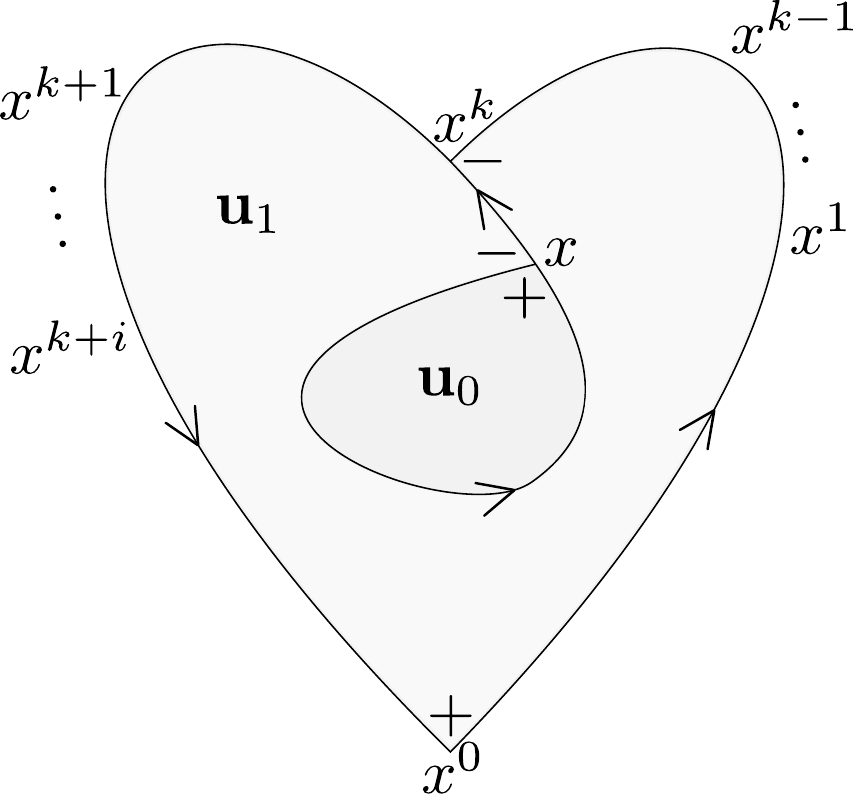}
		\caption{$j = 0$}
		\label{subfig: ConcavePolygon_2c}
	\end{subfigure}%
	\hfill%
	\begin{subfigure}[t]{0.25\linewidth}
		\centering
		\includegraphics[width=\linewidth]{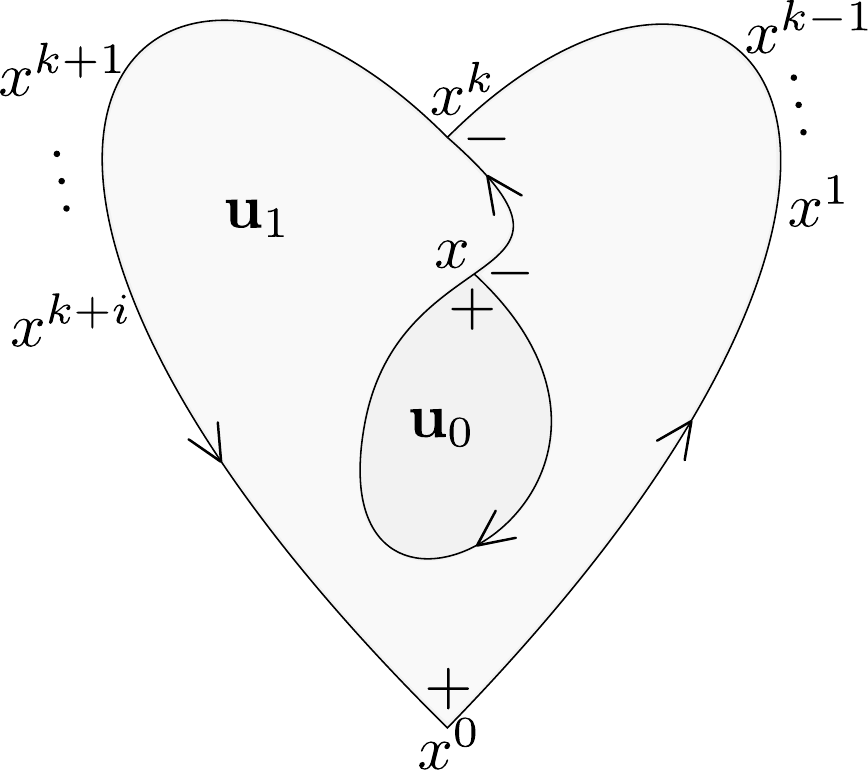}
		\caption{$j = 0$}
		\label{subfig: ConcavePolygon_2d}
	\end{subfigure}
	\caption{Extension of the outgoing side of a negative concave corner.}
	\label{fig: Concave Polygon 2}
	\end{figure}

	\begin{figure}[!htbp]
	\centering
	\begin{subfigure}[t]{0.25\linewidth}
		\centering
		\includegraphics[width=\linewidth]{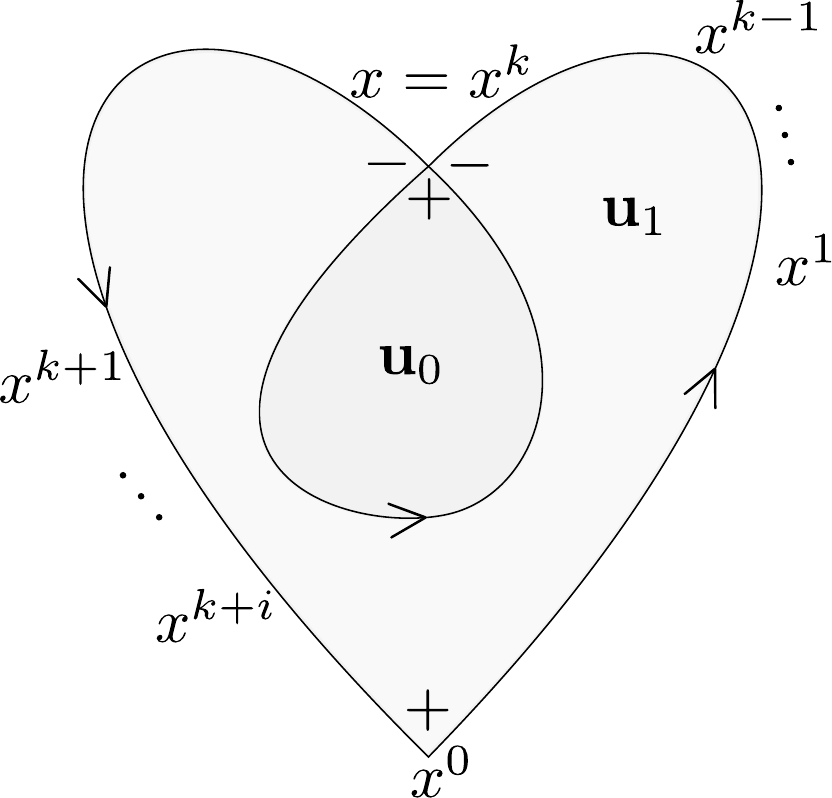}
		\caption{$j = 0$}
		\label{subfig: ConcavePolygon_3a}
	\end{subfigure}%
	\hspace{0.10\linewidth}%
	\begin{subfigure}[t]{0.25\linewidth}
		\centering
		\includegraphics[width=\linewidth]{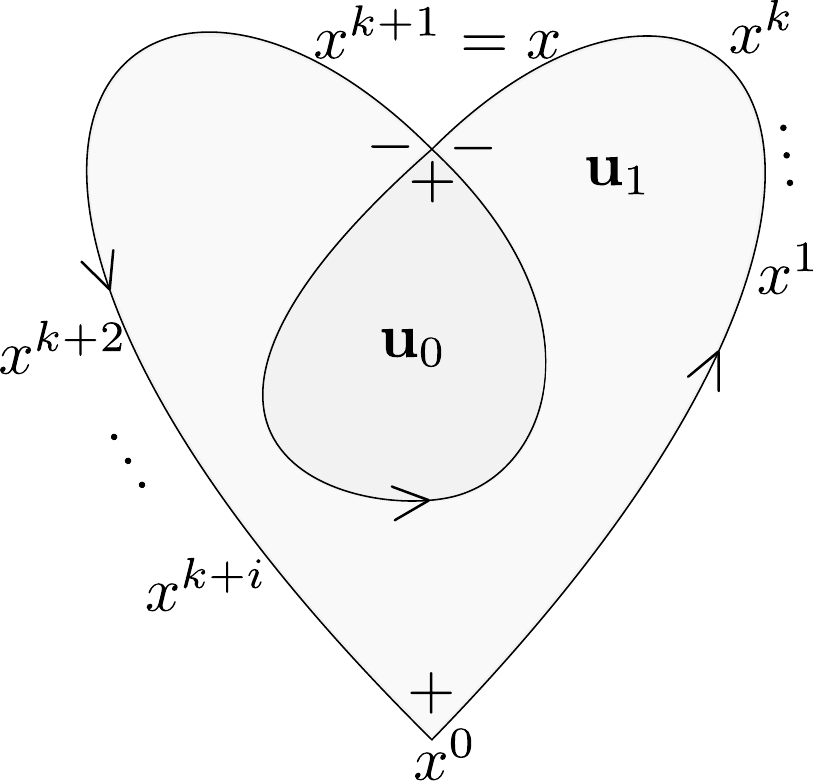}
		\caption{$j = 0$}
		\label{subfig: ConcavePolygon_3b}
	\end{subfigure}%
	\caption{Decompositions of a polygon with one negative concave corner arising when extensions of the incoming and outgoing sides coincide.}
	\label{fig: Concave Polygon 3}
	\end{figure}

	\begin{figure}[!htbp]
	\centering
	\begin{subfigure}[t]{0.25\linewidth}
		\centering
		\includegraphics[width=\linewidth]{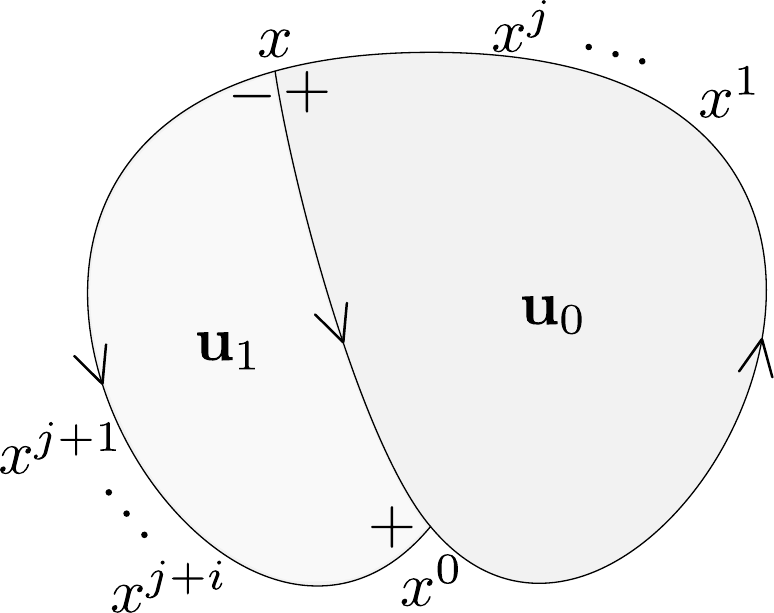}
		\caption{$k = 0$}
		\label{subfig: ConcavePolygon_4a}
	\end{subfigure}%
	\hspace{0.10\linewidth}%
	\begin{subfigure}[t]{0.23\linewidth}
		\centering
		\includegraphics[width=\linewidth]{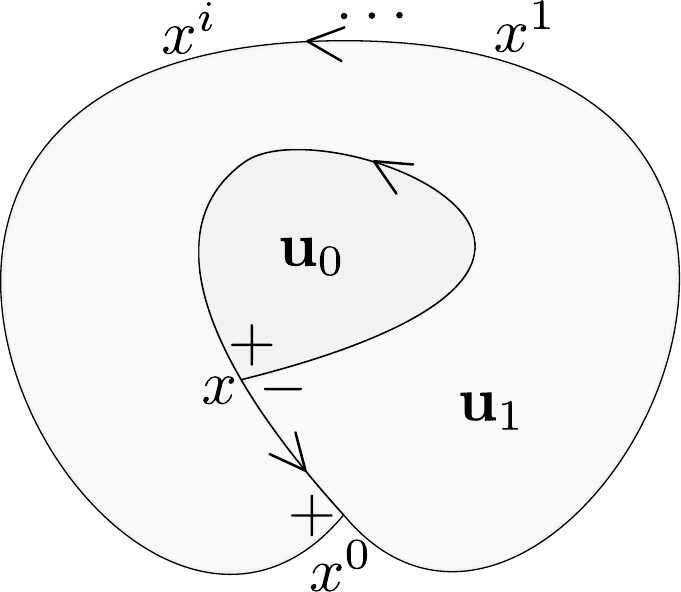}
		\caption{$j = k = 0$}
		\label{subfig: ConcavePolygon_4b}
	\end{subfigure}%
	\hspace{0.10\linewidth}%
	\begin{subfigure}[t]{0.23\linewidth}
		\centering
		\includegraphics[width=\linewidth]{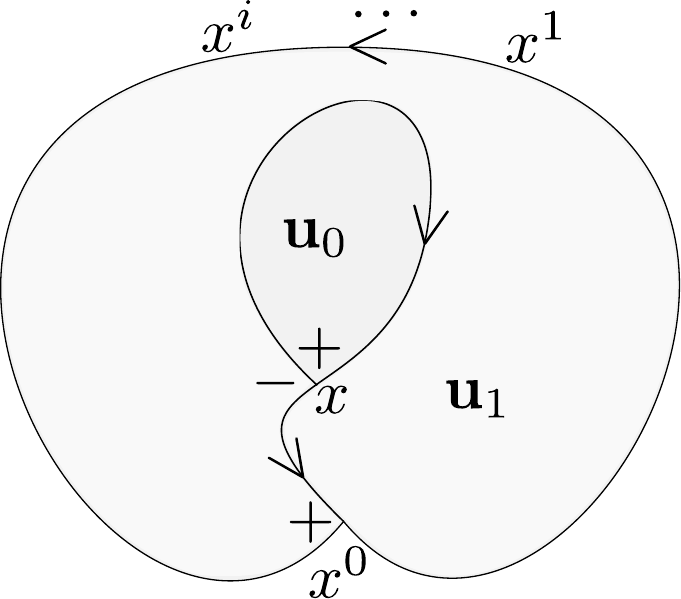}
		\caption{$j = k = 0$}
		\label{subfig: ConcavePolygon_4c}
	\end{subfigure}%
	\caption{Extension of the outgoing side of a positive concave corner.}
	\label{fig: Concave Polygon 4}
	\end{figure}

	\begin{figure}[!htbp]
	\centering
	\begin{subfigure}[t]{0.25\linewidth}
		\centering
		\includegraphics[width=\linewidth]{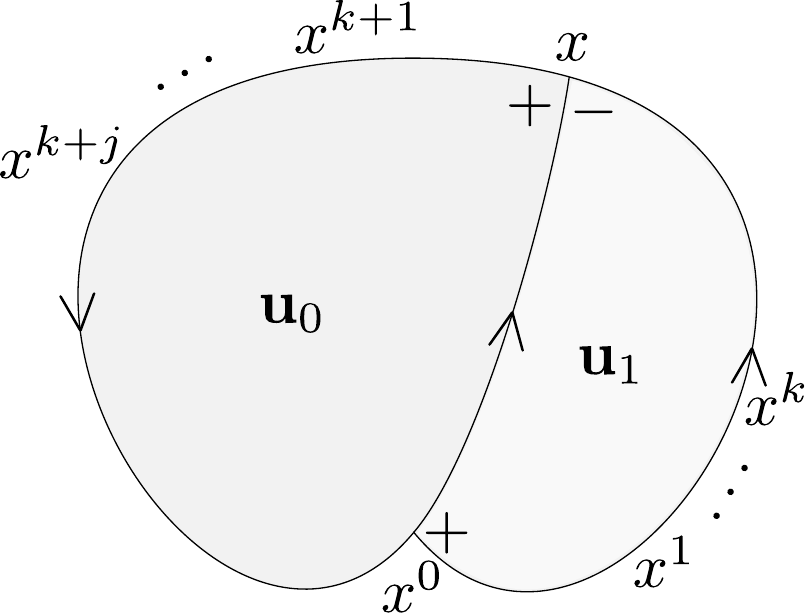}
		\caption{$i = 0$}
		\label{subfig: ConcavePolygon_5a}
	\end{subfigure}%
	\hspace{0.10\linewidth}%
	\begin{subfigure}[t]{0.23\linewidth}
		\centering
		\includegraphics[width=\linewidth]{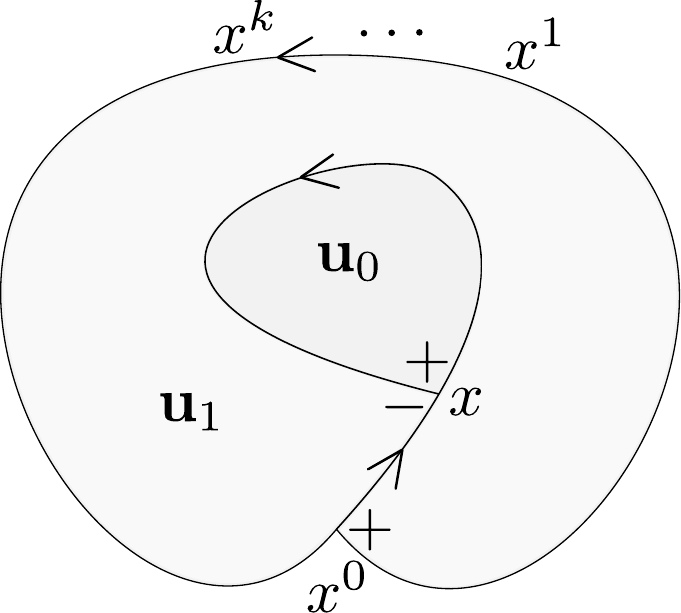}
		\caption{$i = j = 0$}
		\label{subfig: ConcavePolygon_5b}
	\end{subfigure}%
	\hspace{0.10\linewidth}%
	\begin{subfigure}[t]{0.23\linewidth}
		\centering
		\includegraphics[width=\linewidth]{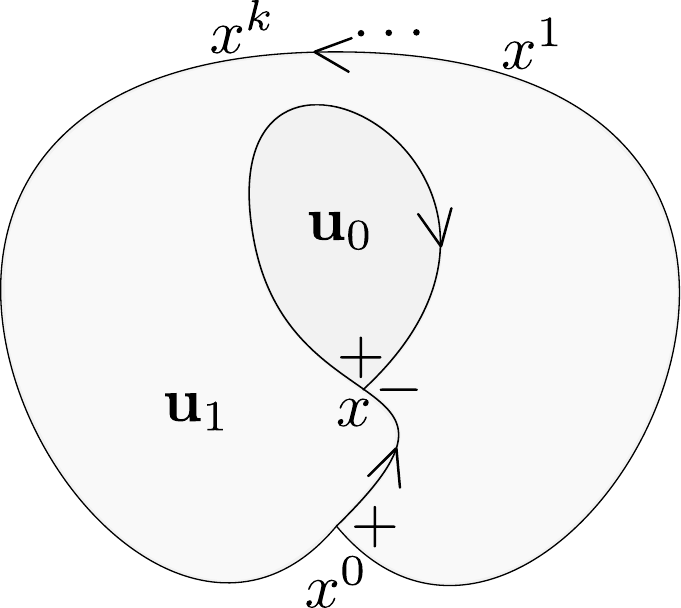}
		\caption{$i = j = 0$}
		\label{subfig: ConcavePolygon_5c}
	\end{subfigure}%
	\caption{Extension of the incoming side of a positive concave corner.}
	\label{fig: Concave Polygon 5}
	\end{figure}

	It remains to prove that for each $\bt = ((i, j, k), x, \bu_0, \bu_1) \in T$ the summands given by $\bt$ and $\iota(\bt)$ sum to zero.
	Let
	$\bu = (\Pi, u) = \theta(\bt) \in \Delta^1(L; x^0; x^n, \ldots, x^1)$
	and $\ph = \mu^{\bu_0}(\ph^{k + j}, \ldots, \ph^{k + 1})$.
	By \cref{l: properties of composition along a 2-path} we have
	$$
	\begin{aligned}
		&\mu^{\bu_1}(\ph^n, \ldots, \ph^{k + j + 1}, \mu^{\bu_0}(\ph^{k + j}, \ldots, \ph^{k + 1}), \ph^k, \ldots, \ph^1) = \\
		&\C(x^0; \gamma^{n}_{\bu}, x^n, \ldots, \gamma^1_{\bu}, x^1, \gamma^0_{\bu}; \sigma_{\bu})(\M(\gamma^{n}_{\bu}), \ph^n, \ldots, \M(\gamma^1_{\bu}), \ph^1, \M(\gamma^0_{\bu}))
	\end{aligned}
	$$
	where $\gamma^i_{\bu}$ is the part of the boundary of $\bu$ going from $x^i$ to $x^{i + 1}$ (with $x^{n + 1} = x^0$).
	That is, up to sign, the contributions of $\bt$ and $\iota(\bt)$ are the same.
	It remains to compute the sign in \cref{eq: summand of A-infinity relation}.

	The sign comes from the following.
	\begin{enumerate}
		\item From $A_\infty$-relation \cref{eq: A-infinity relations}:
		$$
			(-1)^{\sum_{0 < \alpha \le k} (|\ph^\alpha| + |x^\alpha| - 1)}.
		$$
		\item From the sign before the sum in \cref{eq: mu from polygons} for $\bu_0$:
		$$
			(-1)^{\sum\limits_{k < \alpha < \beta \le k + j} |\ph^\alpha|(|x^\beta| - 1)}.
		$$
		\item From the sign before the sum in \cref{eq: mu from polygons} for $\bu_1$:
		$$
			(-1)^{\sum\limits_{0 < \alpha \le k} |\ph^\alpha| \left(\sum\limits_{\alpha < \beta \le k} (|x^\beta| - 1) + |x| - 1 + \sum\limits_{k + j < \beta \le n} (|x^\beta| - 1)\right) +
			|\ph|\sum\limits_{k + j < \beta \le n} (|x^\beta| - 1) +
			\sum\limits_{k + j < \alpha < \beta \le n} |\ph^\alpha|(|x^\beta| - 1)}.
		$$
		\item From the signs $\eps_{\bu_0}$ and $\eps_{\bu_1}$ in \cref{eq: mu from polygons}.
		For $\alpha = 1, \ldots, n$, let $\eps_\alpha$ (resp.\ $\eps_0$) be $1$ if the orientation of $C_\alpha$ (resp.\ $C_n$) coincides with the orientation induced by $\pd \Pi$ and $-1$ otherwise.
		The difference of the sign from
		$$
			\prod_{0 \le \alpha \le n} \eps_\alpha^{|x^\alpha|}
		$$
		is
		\begin{itemize}
			\item for \cref{subfig: ConcavePolygon_1a,subfig: ConcavePolygon_1c,fig: Concave Polygon 4}: $1$;
			\item for \cref{subfig: ConcavePolygon_1b,subfig: ConcavePolygon_1d,subfig: ConcavePolygon_3b}: $(-1)^{|x|}$;
			\item for \cref{subfig: ConcavePolygon_2a}: $(-1)^{|x^{k + j}| + |x|}$;
			\item for \cref{subfig: ConcavePolygon_2c,subfig: ConcavePolygon_2b,subfig: ConcavePolygon_2d,subfig: ConcavePolygon_3a}: $(-1)^{|x^k|}$;
			\item for \cref{fig: Concave Polygon 5}: $(-1)^{|x^0| + |x|}$,
		\end{itemize}
	\end{enumerate}

	Now multiply all up and use that by \cref{eq: |z^0| from |z^1| + ... + |z^n|} and \cref{eq: composition along 2-morphism} we have
	\begin{align*}
		|x| &= |x^{k + j}| + \ldots + |x^{k + 1}| + 2 - j;\\
		|\ph| &= |\ph^{k + j}| + \ldots + |\ph^{k + 1}|;\\
		|x^0| &= |x^n| + \ldots + |x^1| + 3 - n.
	\end{align*}
	Let $l$ be the number of vertices of $\Pi$ on the open interval of $\pd \Pi$ from $x^0$ to the concave corner, i.e.\ for \cref{fig: Concave Polygon 1,fig: Concave Polygon 2,fig: Concave Polygon 3} the concave corner is $x^{l + 1}$ and $l = 0$ for \cref{fig: Concave Polygon 4,fig: Concave Polygon 5}.
	We obtain that the overall sign is equal to
	$$
		(-1)^{\sum\limits_{0 < \alpha \le l} (|x^\alpha| - 1) + \sum\limits_{0 < \alpha < \beta \le n} |\ph^\alpha|(|x^\beta| - 1)}\prod_{0 \le \alpha \le n} \eps_\alpha^{|x^\alpha|};
	$$
	for \cref{fig: Concave Polygon 1,subfig: ConcavePolygon_3b,fig: Concave Polygon 4} and is opposite to that for \cref{fig: Concave Polygon 2,subfig: ConcavePolygon_3a,fig: Concave Polygon 5}.
	Now it is clear that the summands corresponding to $\bt$ and $\iota(\bt)$ cancel out for each $\bt \in T$.
\end{proof}

\subsection{Lagrangian branes}
\label{subs: Lagrangian branes}
Before defining Lagrangian branes in \cref{subsubs: Lagrangian branes} we recall what bounding cochains in curved $A_\infty$-algebras are in \cref{subsubs: bounding cochains}.

\subsubsection{Bounding cochains}
\label{subsubs: bounding cochains}
Let $\A$ be a nilpotent $A_\infty$-algebra.
Given any elements ${\delta_0, \ldots, \delta_n \in \A^1}$ we can define a map
$\mu^n_{\delta_n, \ldots, \delta_0} \colon \A^{\otimes n} \to \A$
by
\begin{equation} \label{eq: mu corrected by one bounding cochain}
\begin{aligned}
	&\mu^n_{\delta_n, \ldots, \delta_0}(a^n, \ldots, a^1) = \\
	&\sum_{k_0, \ldots, k_n \ge 0} \mu^{n + k_n + \ldots + k_0}(\underbrace{\delta_n, \ldots, \delta_n}_{k_n}, a^n, \ldots, a^2, \underbrace{\delta_1, \ldots, \delta_1}_{k_1}, a^1, \underbrace{\delta_0, \ldots, \delta_0}_{k_0}).
\end{aligned}
\end{equation}
The sum in \cref{eq: mu corrected by one bounding cochain} is finite since $\A$ is nilpotent.
These maps satisfy a version of the $A_\infty$-relation:
$$
	\sum_{n \ge j \ge i \ge 0} (-1)^{\dagger}\mu^{n - j + 1 + i}_{\delta_n, \ldots, \delta_j, \delta_i, \ldots, \delta_0}\left(a^n, \ldots, a^{j + 1}, \mu^{j - i}_{\delta_j, \ldots, \delta_i}(a^{j}, \ldots, a^{i + 1}), a^i, \ldots, a^1\right) = 0.
$$
Specialising to $\delta_0, \ldots, \delta_n = \delta$ we get a curved $A_\infty$-algebra structure
\begin{equation} \label{eq: mu_delta}
\mu_\delta = \{\mu^n_{\delta, \ldots, \delta}\}_{n = 0}^\infty \text{ on }\A.
\end{equation}
Note that if $(\A, \mu)$ is structured, then $(\A, \mu_\delta)$ is also structured, with the same $Z$, action, and degree.

A \emph{bounding cochain} in $\A$ is an element $\delta \in \A^1$ such that $(\A, \mu_\delta)$ is uncurved, i.e.\ such that the \emph{Maurer--Cartan equation}
\begin{equation} \label{eq: Maurer--Cartan equation}
\mu_\delta^0() = \sum_{n = 0}^\infty \mu^n(\underbrace{\delta, \ldots, \delta}_n) = 0
\end{equation}
holds.

\subsubsection{Lagrangian branes}
\label{subsubs: Lagrangian branes}
\begin{defn} \label{d: Lagrangian brane}
A \emph{Lagrangian brane} in the weighted Liouville surface $(\S; \C)$ is the data of
\begin{itemize}
	\item a Lagrangian pre-brane $\L$ in $(\S; \C)$;
	\item a bounding cochain $\delta$ in the obstruction algebra $\A_\L$.
\end{itemize}
\end{defn}

Lagrangian branes will be the objects of the Fukaya category.
For a Lagrangian brane $\LL = (\L, \delta)$ denote by $(\A_{\LL}, \mu_{\LL})$ the uncurved $A_\infty$-algebra $(\A_\L, \mu_\delta)$.

\begin{remark} \label{r: obstruction algebra and curved endomorphism algebra}
	Classically (see Akaho--Joyce~\cite{AJ10}), given a Lagrangian pre-brane $\L$, one looks for bounding cochains in a much larger algebra, the curved $A_\infty$-endomorphism algebra of $\L$, which we will not define in our setting, but only explain what the underlying $\k$-module is.
	It is the direct sum of the following:
	\begin{enumerate}
		\item for each circle component $C$ of the domain $I$ the cohomology space $H^*(C; \res{\M}{C})$;
		\item for each self-intersection point $x$ of $\L$ the spaces $\C_{x}(\M(x_+), \M(x_-))[-|x|]$ and $\C_{x}(\M(x_-), \M(x_+))[|x| - 1]$;
		\item for each Reeb chord $y$ between the ends of $\L$ (including the constant Reeb chords) the space $\C_{\gamma(y_-)}(\C(y)\M(y_+), \M(y_-))[-|y|]$.
	\end{enumerate}

	It is clear that the obstruction algebra is a submodule of the above, and it can be shown that it is a curved $A_\infty$-subalgebra, see Haiden \cite{H21}*{Proposition 3.1} for the proof in a related setting.
	The benefit of working with the obstruction algebra instead of the curved endomorphism algebra is not only that the former is smaller, but also that one does not need to perturb the Lagrangian to define the $A_\infty$-operations.
\end{remark}

\begin{remark} \label{r: bounding cochains and augmentations}
	Recall the relation between obstruction algebra and Chekanov--Eliashberg algebra from \cref{r: obstruction algebra and Chekanov-Eliashberg algebra}.
	Assume for concreteness that $\k = \FF_2$, $\C$ is the trivial $2$-functor with fibre $\k$.
	Consider a Lagrangian pre-brane $\L$ with closed domain and local system of rank 1.
	Then $\A^1_{\L}$ can be identified with the degree $1$ maps from $\FF_2$ to $\A_{\L}$ and taking the dual of such map gives a degree $0$ map from $B\A_\L^\lor$ to $\FF_2$.
	We therefore get a map
	$$
	\A_{\L}^1 \to \Hom_{\FF_2}(B\A_\L^\lor, \FF_2), \quad \delta \mapsto \epsilon_\delta
	$$
	that establishes a bijection between bounding cochains in $\A_{\L}$ and augmentations of $B\A_\L^\lor$, i.e.\ maps of dg-algebras $(B\A_\L^\lor, \partial) \to (\FF_2, 0)$.

	The linearised contact cohomology $LCH^*_{\eps_\delta}(\L)$ is then identified with the cohomology of $(\A_\L, \mu_{\delta})$, see Etnyre--Ng~\cite{EN22}*{\Sec 4} for details.
\end{remark}

\section{\texorpdfstring{$A_\infty$}{A∞}-machinery}
\label{s: A_infty-machinery}
In this section we study $A_\infty$-pre-categories with enough homotopy equivalences and develop the algebraic machinery that allows us to pass from the geometric Fukaya $A_\infty$-pre-category to a strictly unital $A_\infty$-category over an arbitrary commutative ring $\k$.
The section is purely algebraic and is independent of the geometric material in \cref{s: Geometric preparations,s: Obstruction algebra}.

In \cref{subs: Homotopy equivalence in A_infty-pre-categories} we define homotopy equivalences in $A_\infty$-pre-categories and $A_\infty$-pre-categories with enough homotopy equivalences, and establish basic properties of homotopy equivalences in \cref{l: 2-out-of-6 property for homotopy equivalences,l: homotopy equivalent transverse objects,l: exact morphisms and contractible objects in A_infty-pre-category,l: homotopy equivalence can be tested pointwise}.
In \cref{subs: Triangulated A_infty-pre-categories} we define the triangulated envelope of an $A_\infty$-pre-category and strongly triangulated $A_\infty$-pre-categories, and prove that homotopy equivalences behave well under passing to the triangulated envelope in \cref{p: homotopy equivalence in triangulated A_infty-pre-category}.
In \cref{subs: From A_infty-pre-categories to A_infty-categories} we construct a strictly unital $A_\infty$-category spanned by a collection of objects admitting a system of strict units in \cref{p: A_infty category spanned by a collection of objects is well-defined} and give a sufficient condition for this construction to be canonical up to homotopy equivalence in \cref{p: W_SU is independent of choices}.

\subsection{Homotopy equivalence in \texorpdfstring{$A_\infty$}{A∞}-pre-categories}
\label{subs: Homotopy equivalence in A_infty-pre-categories}

In \cref{subsubs: A_infty-pre-categories} we recall $A_\infty$-pre-categories and in \cref{subsubs: Homotopy equivalences in A_infty-pre-category} we define homotopy equivalences in them and $A_\infty$-pre-categories with enough homotopy equivalences.
In \cref{subsubs: composition and 2-out-of-6} we prove that homotopy equivalences are closed under composition and satisfy the 2-out-of-6 property, and deduce that transverse homotopy equivalent objects are connected by a single homotopy equivalence.
In \cref{subsubs: exact morphisms in A_infty-pre-category} we characterise exact closed morphisms of degree zero in terms of pre- and post-composition.
In \cref{subsubs: homotopy equivalence can be tested pointwise} we prove that homotopy equivalence can be tested pointwise.
In \cref{subsubs: homotopy equivalence for strictly unital A_infty-category} we characterise homotopy equivalences in strictly unital $A_\infty$-categories in terms of strict units.

\subsubsection{\texorpdfstring{$A_\infty$}{A∞}-pre-categories}
\label{subsubs: A_infty-pre-categories}
We recall the definition of $A_\infty$-pre-category by Kontsevich--Soibelman \cite{KS01}*{\Sec 4.3}.
Let $\k$ be a unital commutative ring.
A \emph{(small $\ZZ$-graded $\k$-linear) $A_\infty$-pre-category} $\W$ is the following data.
\begin{itemize}
	\item A set of objects $\Ob \W$.
	\item For each $n \ge 1$ a subset $\Ob^n_{tr} \W$ of $(\Ob \W)^n$, called the set of \emph{transverse sequences},
	such that $\Ob^1_{tr} \W = \Ob \W$ and subsequences of transverse sequences are transverse.
	\item For each pair $(X_0, X_1) \in \Ob^2_{tr} \W$ a $\ZZ$-graded $\k$-module of morphisms $\W(X_0, X_1)$.
	\item For each $n \ge 1$ and a transverse sequence $(X_0, \ldots, X_n)$ a $\k$-linear map
	$$
		\bmu^n\colon \W(X_{n - 1}, X_n) \otimes \ldots \otimes \W(X_0, X_1) \to \W(X_0, X_n)
	$$
	of degree $2 - n$ satisfying the $A_\infty$-relation \cref{eq: A-infinity relations} for all chains of morphisms
	$X_0 \xrightarrow{a^1} X_1 \xrightarrow{a^2} \ldots \xrightarrow{a^n} X_n$.
\end{itemize}

Note that, unlike $A_\infty$-algebras, all $A_\infty$-pre-categories are uncurved, i.e.\ $\bmu^0 = 0$.

An $A_\infty$-category is an $A_\infty$-pre-category such that any finite sequence of objects is transverse.
Given a transverse sequence $\bX = (X_0, \ldots, X_n)$ in an $A_\infty$-pre-category $\W$, we define the \emph{directed} $A_\infty$-category $\res{\W}{\bX}$ by setting $\Ob \res{\W}{\bX} = \{X_0, \ldots, X_n\}$,
$$
\res{\W}{\bX}(X_i, X_j) = \begin{cases}
	\W(X_i, X_j), &\text{if } i < j;\\
	0, &\text{otherwise.}
\end{cases}
$$
and the $A_\infty$-operations given by that in $\W$.

Let $\W$ be an $A_\infty$-pre-category, $I$ a set and $\bX_I = \{X_i\}_{i \in I}$ a collection of objects in $\W$.
We say that $\bX_I$ is a \emph{transverse set} (this is a slight abuse of language since objects of $\W$ can be repeated in $\bX_I$) if for any finite sequence of pairwise distinct elements $i_0, \ldots, i_n \in I$ the sequence $(X_{i_0}, \ldots, X_{i_n})$ is transverse.
Unless stated otherwise, our standing assumption is that transverse sequences in $A_\infty$-pre-categories \emph{form} transverse sets, that is, any permutation of a transverse sequence is a transverse sequence.

\subsubsection{Homotopy equivalences in \texorpdfstring{$A_\infty$}{A∞}-pre-category}
\label{subsubs: Homotopy equivalences in A_infty-pre-category}
Let $\W$ be an $A_\infty$-pre-category.
The structure maps $\bmu^1$ make morphism spaces into cochain complexes by the \ref{eq: A_1-relation}-relation.
The \ref{eq: A_2-relation},\ref{eq: A_3-relation}-relations imply that
$a \circ b = (-1)^{|b|} \bmu^2(a, b)$
satisfies the graded Leibniz rule with respect to the differential
$da = (-1)^{|a|}\bmu^1(a)$
and descends to an associative (partially defined) composition on cohomology.
We will use the associated homological algebra language (closed, cohomologous, etc.) for morphisms in $\W$.
Denote the degree $k$ component of $\W(X_0, X_1)$ by $\W^k(X_0, X_1)$.

Recall our standing assumption that transverse sequences in $\W$ form transverse sets.
Recall from \cite{KS01}*{Definition 9} that a closed morphism $q^1 \in \W^0(X_0, X_1)$ is called a \emph{quasi-isomorphism} if for any transverse set $\{T, X_0, X_1\}$ the maps
$$
	(\W(T, X_0), d) \xrightarrow{q^1 \circ -} (\W(T, X_1), d) \text{ and } (\W(X_1, T), d) \xrightarrow{- \circ q^1} (\W(X_0, T), d)
$$
are quasi-isomorphisms.
Since we are working with arbitrary modules over an arbitrary commutative ring, quasi-isomorphisms do not necessarily have quasi-inverses.
We therefore introduce a more convenient notion of a homotopy equivalence in $A_\infty$-pre-category.

\begin{defn} \label{d: homotopy inverse}
Let $\{X_0, X_1\}$ be a transverse set in $\W$ and $q^{1} \in \W^0(X_0, X_1)$, $q^{0} \in \W^0(X_1, X_0)$ two closed morphisms.
We say that $q^{1}$ and $q^{0}$ are \emph{homotopy inverse} if for any transverse set $\{T, X_0, X_1\}$ and $i = 0, 1$ the map
$
\bmu^2\left(q^{i}, \bmu^2(q^{1 - i}, -)\right) = q^i \circ (q^{1 - i} \circ -)
$
is homotopic to $\id_{\W(T, X_i)}$ and
the map
$
\bmu^2\left(\bmu^2(-, q^{i}), q^{1 - i}\right) = (- \circ q^{i}) \circ q^{1 - i}
$
is homotopic to $\id_{\W(X_i, T)}$.

A degree zero closed morphism is called a \emph{homotopy equivalence} if it admits a homotopy inverse.
Two objects $X, X' \in \Ob \W$ are called \emph{homotopy equivalent} if for some $n \ge 0$ there exist a chain of homotopy equivalences
$$
X = X_0 \to X_1 \to \ldots \to X_n = X'.
$$
\end{defn}

Note that, in the notation of \cref{d: homotopy inverse}, the morphisms $\bmu^2(q^0, q^1)$ and $\bmu^2(q^1, q^0)$ are usually not defined since $(X_0, X_1, X_0)$ and $(X_1, X_0, X_1)$ are usually not transverse sequences.
Standard homological algebra implies that homotopy equivalences are quasi-isomorphisms.

The following definition is a natural adaptation of \cite{KS01}*{Definition 10} from quasi-isomorphisms to homotopy equivalences.

\begin{defn} \label{d: W has enough homotopy equivalences}
	We say that the $A_\infty$-pre-category $\W$ \emph{has enough homotopy equivalences} if transverse sequences in $\W$ form transverse sets and for any object $X \in \Ob \W$ and a finite set of finite transverse sets $\{S_i\}_{i \in I}$ there exists an object $X'$ and a homotopy equivalence $q\colon X \to X'$ such that the set $S_i \sqcup \{X'\}$ is transverse for each $i \in I$.
\end{defn}
Having enough homotopy equivalences is a substitute for being unital for $A_\infty$-pre-categories and makes homotopy equivalences satisfy the expected properties, see \cref{subs: Homotopy equivalence in A_infty-pre-categories}.

Note that having enough homotopy equivalences is a property of an $A_\infty$-pre-category rather than an additional structure.
Also observe that, in contrast with the homotopy unit framework of Fukaya--Oh--Ohta--Ono~\cite{FOOO}*{\Sec 3.3} we do not require any higher coherences from homotopy equivalences.
This works well in our setting since our $A_\infty$-pre-categories are triangulated by construction.

\subsubsection{}
\label{subsubs: composition and 2-out-of-6}

We prove that composition of homotopy equivalences is a homotopy equivalence whenever it is defined and that homotopy equivalences satisfy the 2-out-of-6 property.
Note that one does not follow from the other since we do not have identity morphisms, but
the 2-out-of-3 property follows from these.

\begin{lemma} \label{l: 2-out-of-6 property for homotopy equivalences}
	Let $\W$ be an $A_\infty$-pre-category with enough homotopy equivalences.
	Then
	\begin{enumerate}
		\item Composition of homotopy equivalences is a homotopy equivalence: if $\{X_0, X_1, X_2\}$ is a transverse set and $q^1 \in \W(X_0, X_1), \, q^2 \in \W(X_1, X_2)$ are homotopy equivalences, then $q^2 \circ q^1$ is a homotopy equivalence.
		\item Homotopy equivalences satisfy the 2-out-of-6 property: for a transverse set $\{X_0, X_1, X_2, X_3\}$ and closed morphisms $q^1 \in \W^0(X_0, X_1)$, $q^2 \in \W^0(X_1, X_2)$, $q^3 \in \W^0(X_2, X_3)$ if $q^2 \circ q^1$ and $q^3 \circ q^2$ are homotopy equivalences, then so are $q^1, q^2, q^3$.
	\end{enumerate}
\end{lemma}

As can be seen from the proof below, having enough homotopy equivalences is crucial for any homotopy argument in $A_\infty$-pre-categories.

\begin{proof}
	We prove (1).
	Denote $q^0 = q^2 \circ q^1$.
	First note that $|q^0| = |q^1| + |q^2| + (2 - 2) = 0$ and $\bmu^1(q^0) = 0$ by the $A_2$-relation.
	Let $p^1 \in \W(X_1, X_0)$ be a homotopy inverse to $q^1$ and $p^2 \in \W(X_2, X_1)$ a homotopy inverse to $q^2$.
	We show that
	$p^0 = p^1 \circ p^2$ is homotopy inverse to $q^0$.
	Let $T \in \Ob \W$ be such that $\{T, X_0, X_2\}$ is a transverse set.
	Since $\W$ has enough homotopy equivalences, there is an object $T'$ and homotopy inverse morphisms $e \in \W(T', T), \, e' \in \W(T, T')$ such that the sets $\{T', T, X_0, X_2\}$ and $\{T', X_0, X_1, X_2\}$ are transverse.

	The $A_3$-relation implies that the following diagram is homotopy commutative.
	$$
	\xymatrix@C=1pc{
		\W(T, X_0) \ar[rr]^-{q^0 \circ -} \ar[d]_{- \circ e}&& \W(T, X_2) \ar@<-1ex>[d]_{- \circ e} \ar[rr]^-{p^0 \circ -} && \W(T, X_0)\\
		\W(T', X_0) \ar[rr]^-{q^0 \circ -} \ar[rd]_-{q^1 \circ -} && \W(T', X_2) \ar@<-1ex>[u]_{- \circ e'} \ar[rr]^-{p^0 \circ -} \ar[rd]_-{p^2 \circ -} && \W(T', X_0) \ar[u]_{- \circ e'}\\
		& \W(T', X_1) \ar[ru]_-{q^2 \circ -} && \W(T', X_1) \ar[ru]_-{p^1 \circ -} &
	}
	$$
	Using that $(- \circ e) \circ e'$, $p^2 \circ (q^2 \circ -)$, $p^1 \circ (q^1 \circ -)$ are all homotopic to identity we conclude that so is $p^0 \circ (q^0 \circ -)$.
	The other three maps from \cref{d: homotopy inverse} are homotopic to identity for similar reasons.

	We prove (2).
	Let $p^{20} \in \W(X_2, X_0)$ be a homotopy inverse to $q^{02} = q^2 \circ q^1$ and $p^{31} \in \W(X_3, X_1)$ a homotopy inverse to $q^{13} = q^3 \circ q^2$.
	We show that $p^2_0 = q^1 \circ p^{20}$ and $p^2_3 = p^{31} \circ q^3$ are homotopy inverse to $q^2$.
	Let $T \in \Ob \W$ be such that $\{T, X_1, X_2\}$ is a transverse set.
	Since $\W$ has enough homotopy equivalences, there is an object $T'$ and homotopy inverse morphisms $e \in \W(T', T)$, $e' \in \W(T, T')$ such that the sets $\{T', T, X_1, X_2\}$ and $\{T', X_0, X_1, X_2, X_3\}$ are transverse.

	The $A_3$-relation implies that the following diagram is homotopy commutative.
	$$
	\xymatrix{
		\W(X_2, T) \ar[rr]^-{- \circ q^2} \ar[d]_{e' \circ -}&& \W(X_1, T) \ar@<-1ex>[d]_{e' \circ -} \ar[rr]^-{- \circ p^2_0} && \W(X_2, T)\\
		\W(X_2, T') \ar[rr]^-{- \circ q^2} \ar[rrd]_-{- \circ q^{02}} && \W(X_1, T') \ar@<-1ex>[u]_{e \circ -} \ar[rr]^-{- \circ p^2_0} \ar[d]_-{- \circ q^1} && \W(X_2, T') \ar[u]_{e \circ -}\\
		& & \W(X_0, T') \ar[rru]_-{- \circ p^{20}} &&
	}
	$$
	Therefore, $(- \circ q^2) \circ p^2_0$ is homotopic to $\id_{\W(X_2, T)}$.
	A similar diagram implies that $(- \circ p^2_3) \circ q^2$ is homotopic to $\id_{\W(X_1, T)}$.
	Hence $- \circ p^2_0$ and $- \circ p^2_3$ are homotopic as morphisms $\W(X_1, T) \to \W(X_2, T)$, $(- \circ q^2) \circ p^2_3$ is homotopic to $\id_{\W(X_2, T)}$, and $(- \circ p^2_0) \circ q^2$ is homotopic to $\id_{\W(X_1, T)}$.
	The other four maps from \cref{d: homotopy inverse} applied to $(q^2, p^2_0)$ and $(q^2, p^2_3)$ are homotopic to identity for similar reasons hence $p^2_0$ and $p^2_3$ are homotopy inverse to $q^2$.

	A similar argument shows that $p^{20} \circ q^2$ is homotopy inverse to $q^1$ and $q^2 \circ p^{31}$ is homotopy inverse to $q^3$.
\end{proof}

\begin{lemma} \label{l: homotopy equivalent transverse objects}
	Let $\W$ be an $A_\infty$-pre-category with enough homotopy equivalences.
	Then for any pair of homotopy equivalent objects $X_0, X_1$ such that $\{X_0, X_1\}$ is a transverse set there exists a homotopy equivalence ${q\colon X_0 \to X_1}$.
\end{lemma}

Note that \cref{l: homotopy equivalent transverse objects} implies that for any homotopy equivalent objects $X_0, X_1$ there exists an object $X_2$ and homotopy equivalences $X_0 \to X_2$ and $X_2 \to X_1$.

\begin{proof}
	Since $X_0$ and $X_1$ are homotopy equivalent, there exist a sequence $X_0 = Y_0, Y_1, \ldots, Y_n = X_1$ of objects in $\W$ and homotopy equivalences $q^i \in \W(Y_{i - 1}, Y_i)$ for $i = 1, \ldots, n$.
	We need to construct a homotopy equivalence $q \in \W(Y_0, Y_n)$ provided that $Y_0$ and $Y_n$ are transverse.
	As $\W$ has enough homotopy equivalences, we can take an object $Y'$ and homotopy inverse morphisms $q_0' \in \W(Y', Y_0)$ and $p_0' \in \W(Y_0, Y')$ such that the sets $\{Y', Y_0, Y_n\}$ and $\{Y', Y_{i - 1}, Y_i\}, \, i = 1, \ldots, n$ are transverse.
	Using \cref{l: 2-out-of-6 property for homotopy equivalences}\,(1) we obtain homotopy equivalences $q_i' \in \W(Y', Y_i)$ by $q_i' = q_{i} \circ q_{i - 1}'$ for $i = 1, \ldots, n$.
	Setting $q = q_n' \circ p_0'$ gives the desired homotopy equivalence from $Y_0 = X_0$ to $Y_n = X_1$.
\end{proof}
\subsubsection{}
\label{subsubs: exact morphisms in A_infty-pre-category}
We give criteria for a degree zero closed morphism to be exact.

\begin{lemma} \label{l: exact morphisms and contractible objects in A_infty-pre-category}
	Let $\W$ be an $A_\infty$-pre-category with enough homotopy equivalences and $q \in \W^0(X_0, X_1)$ a closed morphism.
	Then the following are equivalent.
	\begin{enumerate}
		\item $q$ is exact;
		\item for any transverse set $\{T, X_0, X_1\}$ the map $\W(T, X_0) \xrightarrow{q \circ -} \W(T, X_1)$ is null-homotopic;
		\item there exists a homotopy equivalence $e_0 \in \W(X_0', X_0)$ such that the set $\{X_0', X_0, X_1\}$ is transverse and $q \circ e_0 \in \W(X_0', X_1)$ is exact;
		\item for any transverse set $\{X_0, X_1, T\}$ the map $\W(X_1, T) \xrightarrow{- \circ q} \W(X_0, T)$ is null-homotopic;
		\item there exists a homotopy equivalence $e_1' \in \W(X_1, X_1')$ such that the set $\{X_0, X_1, X_1'\}$ is transverse and $e_1' \circ q \in \W(X_0, X_1')$ is exact.
	\end{enumerate}
\end{lemma}
\begin{proof}
	We show $(1) \Rightarrow (2) \Rightarrow (3) \Rightarrow (1)$.
	The argument for $(1) \Rightarrow (4) \Rightarrow (5) \Rightarrow (1)$ is similar.

	Assume $q = \bmu^1(a)$ for some $a \in \W^{-1}(X_0, X_1)$.
	Let $\{T, X_0, X_1\}$ be a transverse set.
	The $A_2$-relation then implies that
	$$
		q \circ b = \bmu^1(\bmu^2(a, b)) + \bmu^2(a, \bmu^1(b)) \text{ for any homogeneous } b \in \W(T, X_0).
	$$
	Hence $\W(T, X_0) \xrightarrow{q \circ -} \W(T, X_1)$ is indeed null-homotopic.

	That (2) implies (3) is clear since we can take $X_0'$ for $T$ and $\bmu^1(e_0) = 0$.

	Now assume (3).
	As $e_0$ is a homotopy equivalence, $\W(X_0, X_1) \xrightarrow{- \circ e_0} \W(X_0', X_1)$ is a homotopy equivalence.
	Since $q \circ e_0 \in \W(X_0', X_1)$ is exact, so is $q \in \W(X_0, X_1)$.
\end{proof}

\subsubsection{}
\label{subsubs: homotopy equivalence can be tested pointwise}
We prove that homotopy equivalence can be tested pointwise.
\begin{lemma} \label{l: homotopy equivalence can be tested pointwise}
	Let $\W$ be an $A_\infty$-pre-category with enough homotopy equivalences, and $X_0 \xrightarrow{q} X_1$ a closed morphism of degree zero in $\W$.
	Suppose that for any transverse set $\{T, X_0, X_1\}$ the maps $\W(T, X_0) \xrightarrow{q \circ -} \W(T, X_1)$ and $\W(X_1, T) \xrightarrow{- \circ q} \W(X_0, T)$ are homotopy equivalences.

	Then $q$ is a homotopy equivalence in $\W$.
\end{lemma}
\begin{proof}
	It suffices to construct two degree zero closed morphisms $p^0, p^1 \in \W(X_1, X_0)$ such that for any transverse set $\{T, X_0, X_1\}$ the maps $p^0 \circ (q \circ -)$, $(- \circ p^0) \circ q$, $q \circ (p^1 \circ -)$, and $(- \circ q) \circ p^1$ are homotopic to identity.
	Then the maps $p^0 \circ -$ and $p^1 \circ -$ are homotopic hence by \cref{l: exact morphisms and contractible objects in A_infty-pre-category} $p^0$ is cohomologous to $p^1$ and they are both homotopy inverse to $q$.

	To construct $p^0$ pick homotopy inverse morphisms $e \in \W(X_0', X_0)$ and $e' \in \W(X_0, X_0')$ such that $\{X_0', X_0, X_1\}$ is a transverse set.
	Let $\W(X_0, X_0') \xrightarrow{\alpha} \W(X_1, X_0')$ be a homotopy inverse to $- \circ q$.
	Set $p^0 = e \circ \alpha(e') \in \W(X_1, X_0)$.
	Clearly, $p^0$ is closed of degree zero.
	Let $\{T, X_0, X_1\}$ be a transverse set.
	As in the proof of \cref{l: 2-out-of-6 property for homotopy equivalences} we can assume that $\{T, X_0, X_1, X_0'\}$ is transverse.
	Note that $\alpha(e') \circ q$ is cohomologous to $e'$ since $\alpha$ is homotopy inverse to $- \circ q$.
	The $A_3$-relation and \cref{l: exact morphisms and contractible objects in A_infty-pre-category} now imply that the diagram
	$$\xymatrix{
		\W(T, X_0) \ar[r]^-{q \circ -} \ar[rrd]_-{e' \circ -}& \W(T, X_1) \ar[rr]^-{p^0 \circ -} \ar[rd]^-{\alpha(e') \circ -} && \W(T, X_0)\\
		&& \W(T, X_0') \ar[ru]_-{e \circ -}
	}$$
	is homotopy commutative.
	As $e$ is homotopy inverse to $e'$ we conclude that $p^0 \circ (q \circ -)$ is homotopic to identity.
	A similar argument shows that $(- \circ p^0) \circ q$ is homotopic to identity.

	To construct $p^1$ pick homotopy inverse morphisms $e \in \W(X_1', X_1)$ and $e' \in \W(X_1, X_1')$ such that $\{X_1', X_0, X_1\}$ is a transverse set.
	Let $\W(X_1', X_1) \xrightarrow{\alpha} \W(X_1', X_0)$ be a homotopy inverse to $q \circ -$.
	Set $p^1 = \alpha(e) \circ e' \in \W(X_1, X_0)$.
	Clearly, $p^1$ is closed of degree zero.
	Let $\{T, X_0, X_1\}$ be a transverse set.
	As in the proof of \cref{l: 2-out-of-6 property for homotopy equivalences} we can assume that $\{T, X_0, X_1, X_1'\}$ is transverse.
	Note that $q \circ \alpha(e)$ is cohomologous to $e$ since $\alpha$ is homotopy inverse to $q \circ -$.
	The $A_3$-relation and \cref{l: exact morphisms and contractible objects in A_infty-pre-category} now imply that the diagram
	$$\xymatrix{
		\W(T, X_1) \ar[rr]^-{p^1 \circ -} \ar[rd]_-{e' \circ -}&&\W(T, X_0) \ar[r]^-{q \circ -}  & \W(T, X_1)\\
		& \W(T, X_1') \ar[ru]^-{\alpha(e) \circ -} \ar[rru]_-{e \circ -}
	}$$
	is homotopy commutative.
	As $e$ is homotopy inverse to $e'$ we conclude that $q \circ (p^1 \circ -)$ is homotopic to identity.
	A similar argument shows that $(- \circ q) \circ p^1$ is homotopic to identity.

	We have therefore constructed the desired $p^0, p^1 \in \W(X_1, X_0)$.
\end{proof}

\subsubsection{}
\label{subsubs: homotopy equivalence for strictly unital A_infty-category}
Recall that an $A_\infty$-category $\C$ is called \emph{strictly unital} if for each $X \in \Ob \C$ there exists a closed morphism $1_X \in \C^0(X, X)$ such that $1_X \circ a = a$ for any $a \in \C(Y, X)$, $a \circ 1_X = a$ for any $a \in \C(X, Y)$, and $\bmu^n(\ldots, 1_X, \ldots) = 0$ for $n \ge 3$.
Note that strict units $1_X, \, X \in \Ob \C$ are unique if they exist.

Observe that strict units in $\C$ are homotopy equivalences and that they certify that $\C$ has enough homotopy equivalences.
The following characterisation of homotopy equivalences in strictly unital $A_\infty$-categories indicates that homotopy equivalences are the right notion of invertible morphism in an $A_\infty$-pre-category.

\begin{lemma} \label{l: characterisation of homotopy equivalences in a strictly unital A_infty-category}
	Let $\C$ be a strictly unital $A_\infty$-category.
	Let $q \in \C^0(X_0, X_1)$ be a closed morphism.
	\begin{enumerate}
		\item If $p \in \C^0(X_1, X_0)$ is homotopy inverse to $q$, then $p \circ q$ is cohomologous to $1_{X_0}$ and $q \circ p$ is cohomologous to $1_{X_1}$.
		\item If there exists a closed morphism $p^1 \in \C^0(X_1, X_0)$ such that $q \circ p^1$ is cohomologous to $1_{X_1}$ and a closed morphism $p^0 \in \C^0(X_1, X_0)$ such that $p^0 \circ q$ is cohomologous to $1_{X_0}$, then $q$ is a homotopy equivalence.
	\end{enumerate}
\end{lemma}
\begin{proof}
	(1) follows immediately from \cref{l: exact morphisms and contractible objects in A_infty-pre-category}.

	To prove (2) note that by \cref{l: exact morphisms and contractible objects in A_infty-pre-category} $p^0, q, p^1$ satisfy the assumption of \cref{l: 2-out-of-6 property for homotopy equivalences}(2) hence $q$ is a homotopy equivalence.
\end{proof}

\subsection{Triangulated \texorpdfstring{$A_\infty$}{A∞}-pre-categories}
\label{subs: Triangulated A_infty-pre-categories}
In \cref{subsubs: twisted complexes} we define the triangulated envelope $\Tw \W$ of an $A_\infty$-pre-category $\W$.
In \cref{subsubs: A_infty-functors} we recall $A_\infty$-functors between $A_\infty$-pre-categories.
In \cref{subsubs: triangulated A_infty-pre-category} we define strongly triangulated $A_\infty$-pre-categories, i.e.\ those admitting a strict totalisation functor $\Tw \W \to \W$, and prove in \cref{p: homotopy equivalence in triangulated A_infty-pre-category} that for such $\W$ with enough homotopy equivalences homotopy equivalences behave well under passing to $\Tw \W$; in particular, any twisted complex is homotopy equivalent to an object of $\W$.
In \cref{subsubs: Shifts and cones} we define shifts and cones in $\Tw \W$ and in strongly triangulated $A_\infty$-pre-categories and prove that upper triangular morphisms of twisted complexes with homotopy equivalences on the diagonal are homotopy equivalences in \cref{l: upper triangular morphisms of twisted complexes}.

\subsubsection{Triangulated envelope}
\label{subsubs: twisted complexes}
The triangulated envelope of an $A_\infty$-pre-category is defined similarly to the case of $A_\infty$-categories, and the construction (in the case of dg-categories) goes back to Bondal--Kapranov~\cite{BK91}.
We follow Efimov~\cite{E11}*{\Sec 4} with slight modifications taking into account that transverse sequences form transverse sets in our $A_\infty$-pre-categories.

Let $\W$ be a $\ZZ$-graded $\k$-linear $A_\infty$-pre-category where transverse sequences form transverse sets.
The \emph{triangulated envelope} $\Tw \W$ of $\W$ is the $\ZZ$-graded $\k$-linear $A_\infty$-pre-category defined as follows.

\subsubsection*{Objects}
An object of $\Tw \W$ is a twisted complex, i.e.\ a pair $(X, \delta)$ where
\begin{itemize}
	\item $X = \bigoplus_{i \in I} X_i[m_i]$ is a formal direct sum of formal shifts of objects in $\W$ such that $I$ is finite, $\{X_i\}_{i \in I}$ is a transverse set in $\W$, and $m_i \in \ZZ$ for $i \in I$;
	\item $\delta$ is an element of degree $1$ in the graded $\k$-module
	$$
	\bigoplus_{i_0, i_1 \in I\colon i_0 \ne i_1} \W(X_{i_0}, X_{i_1})[m_{i_1} - m_{i_0}]
	$$
	and there exists a partial order $<$ on $I$ such that
	$$
		\delta \in \A_{X}^< = \bigoplus_{i_0, i_1 \in I\colon i_0 < i_1} \W(X_{i_0}, X_{i_1})[m_{i_1} - m_{i_0}]
	$$
	is a bounding cochain in the nilpotent (uncurved) $A_\infty$-algebra $\A_{X}^<$ where the $A_\infty$-operations in $\A_{X}^<$ are given by that in $\W$ with a Koszul sign rule,
	\begin{equation} \label{eq: mu with Koszul sign}
	\mu^n(a_n[i_n], \ldots, a_1[i_1]) = (-1)^{\sum_{1 \le \alpha < \beta \le n} i_\alpha(|a_\beta| - 1)}\bmu^n(a_n, \ldots, a_1)[i_1 + \ldots + i_n]
	\end{equation}
	for homogeneous morphisms $a_1, \ldots, a_n$ in $\W$ and $i_1, \ldots, i_n \in \ZZ$.
\end{itemize}
We have
$$
	\delta = \sum_{i_0, i_1 \in I} \delta_{i_0i_1}[m_{i_1} - m_{i_0}] \text{ where } \delta_{i_0i_1} \in \W^{m_{i_1} - m_{i_0} + 1}(X_{i_0}, X_{i_1}),
$$
the morphisms $\delta_{i_0i_1}$ are called \emph{components} of $\delta$.
Note that if $\delta_{i_0i_1} \ne 0$, then $i_0 < i_1$.
The partial order on $I$ is not part of the data, but we can always consider the weakest partial order such that $i_0 < i_1$ whenever $\delta_{i_0i_1} \ne 0$.

\subsubsection*{Transverse sequences}
A sequence of twisted complexes $((X^0, \delta^0), \ldots, (X^n, \delta^n))$ with $X^j = \bigoplus_{i \in I^j} X^j_{i}[m^j_i]$ for $j = 0, \ldots, n$ is transverse if and only if
$$
\bigsqcup_{j = 0}^n \bigsqcup_{i\in I^j} \{X^j_i\}
$$
is a transverse set in $\W$.

\subsubsection*{Morphisms}
For a transverse pair of twisted complexes $((X^0, \delta^0), (X^1, \delta^1))$ with
${X^j = \bigoplus_{i \in I^j} X^j_{i}[m^j_i]}$ for $j = 0, 1$ the morphism space is the graded $\k$-module
\begin{equation} \label{eq: morphism space between twisted complexes}
\Tw \W((X^0, \delta^0), (X^1, \delta^1)) = \bigoplus_{
i^0 \in I^0,
i^1 \in I^1}
\W(X^0_{i^0}, X^1_{i^1})[m^1_{i^1} - m^0_{i^0}].
\end{equation}

For a morphism $a \in \Tw \W((X^0, \delta^0), (X^1, \delta^1))$ of degree $d$ we have
$$
a = \sum_{i^0 \in I^0,
i^1 \in I^1}
a_{i^0i^1}[m^1_{i^1} - m^0_{i^0}]
\text{ where }
a_{i^0i^1} \in \W^{m^1_{i^1} - m^0_{i^0} + d}(X^0_{i^0}, X^1_{i^1}),
$$
the morphisms $a_{i^0i^1}$ are called \emph{components} of $a$.

\subsubsection*{Structure maps}
For a transverse sequence $((X^0, \delta^0), \ldots, (X^n, \delta^n))$ with $X^j = \bigoplus_{i \in I^j} X^j_{i}[m^j_i]$ for $j = 0, \ldots, n$, the $A_\infty$-operation $\bmu^n$ is given by $\mu_\delta$ from \cref{eq: mu_delta} in the $A_\infty$-algebra $\A_X^<$ where $\delta = \delta^0 + \ldots + \delta^n$ is the bounding cochain in $\A_X^<$ with $X = X^0 \oplus \ldots \oplus X^n$ and a partial order $<$ on $I^0 \sqcup \ldots \sqcup I^n$ that includes partial orders on each $I^j$ and satisfies $i^{j} < i^{j'}$ for any $j < j'$ and $i^j \in I^j$, $i^{j'} \in I^{j'}$.
\medskip

It is instructive to think of $X[m]$ for $X \in \W$ and $m \in \ZZ$ as a formal tensor product $\k[m] \otimes X$.
The Koszul sign in \cref{eq: mu with Koszul sign} comes from this point of view and the fact that $\mu^n$ acts from the right.

Note that $\W$ is naturally a full subcategory of $\Tw \W$ by identifying $X$ and the twisted complex $(X[0], 0)$ for each $X \in \Ob \W$.
Observe that transverse sequences in $\Tw \W$ form transverse sets.

\subsubsection{\texorpdfstring{$A_\infty$}{A∞}-functors}
\label{subsubs: A_infty-functors}
Let $\W$ and $\V$ be two $A_\infty$-pre-categories.
An \emph{$A_\infty$-functor} ${\F\colon \W \to \V}$ is the following data.
\begin{itemize}
	\item A map $\F^0\colon \Ob \W \to \Ob \V$ such that for any transverse sequence $(X_0, \ldots, X_n)$ in $\W$ the sequence $(\F^0(X_0), \ldots, \F^0(X_n))$ is transverse in $\V$.
	\item For each $n \ge 1$ and a transverse sequence $(X_0, \ldots, X_n)$ in $\W$ a $\k$-linear map
	$$
		\F^n\colon \W(X_{n - 1}, X_n) \otimes \ldots \otimes \W(X_0, X_1) \to \V(\F^0(X_0), \F^0(X_n))
	$$
	of degree $1 - n$ satisfying the $A_n$-functor-relation
	\begin{equation} \label{eq: A-infinity functor relations}
	\begin{aligned}
		&\sum_{n \ge j > i \ge 0} (-1)^{\dagger}\F^{n - j + 1 + i}\left(a^n, \ldots, a^{j + 1}, \bmu_\W^{j - i}(a^{j}, \ldots, a^{i + 1}), a^i, \ldots, a^1\right) = \\
		&\sum_{r > 0} \; \sum_{n = i_r > \ldots > i_0 = 0} \bmu^r_\V\left(\F^{n - i_{r - 1}}(a^n, \ldots, a^{i_{r - 1} + 1}), \ldots, \F^{i_1}(a^{i_1}, \ldots, a^1)\right).
	\end{aligned}
	\end{equation}
where $\dagger = |a^i| + \ldots + |a^1| - i$ for any $n \ge 1$ and any chain of homogeneous morphisms
$X_0 \xrightarrow{a^1} X_1 \xrightarrow{a^2} \ldots \xrightarrow{a^n} X_n$
in $\W$.
\end{itemize}

The $A_\infty$-functor $\F$ is called \emph{strict} if $\F^n = 0$ for $n \ge 2$.
For strict $A_\infty$-functors the $A_\infty$-relations \cref{eq: A-infinity functor relations} boils down to
$$
\F^1\left(\bmu^n_{\W}(a^n, \ldots, a^1)\right) = \bmu^n_{\V}\left(\F^1(a^n), \ldots, \F^1(a^1)\right)
$$
for any $n \ge 1$ and any chain of morphisms
$X_0 \xrightarrow{a^1} X_1 \xrightarrow{a^2} \ldots \xrightarrow{a^n} X_n$
in $\W$.

\subsubsection{Triangulated \texorpdfstring{$A_\infty$}{A∞}-pre-category}
\label{subsubs: triangulated A_infty-pre-category}
We call an $A_\infty$-pre-category $\W$ \emph{strongly triangulated} (cf.\ Drinfeld~\cite{D04}*{\Sec2.4}) if there exists a strict $A_\infty$-functor
${\T\colon \Tw \W \to \W}$, called the \emph{totalisation functor},
such that
\begin{enumerate}
		\item A sequence $(\bX^0, \ldots, \bX^n)$ of objects in $\Tw \W$ is transverse if $(\T^0\bX^0, \ldots, \T^0\bX^n)$ is a transverse sequence in $\W$;
		\item
		$\T^1\colon \Tw \W(\bX^0, \bX^1) \to \W(\T^0(\bX^0), \T^0(\bX^1))$
		is an isomorphism of graded modules for each transverse pair $(\bX^0, \bX^1)$ in $\Tw \W$;
		\item $\T^0(X) = X$ for each $X \in \Ob \W$;
		\item $\T^1$ is identity on $\Tw \W(X^0, X^1) = \W(X^0, X^1)$ for each transverse pair $(X^0, X^1)$ in $\W$.
\end{enumerate}

It is not hard to show that $\Tw \W$ is strongly triangulated for any $A_\infty$-pre-category $\W$.

If $\W$ has enough homotopy equivalences, then we can define $\W$ to be \emph{triangulated} if any twisted complex over $\W$ is homotopy equivalent to an object of $\W$.
The next statement shows that this alternative definition is implied by ours.

\begin{prop} \label{p: homotopy equivalence in triangulated A_infty-pre-category}
	Let $\W$ be a strongly triangulated $A_\infty$-pre-category with enough homotopy equivalences.
	Then
	two twisted complexes $\bX^0, \bX^1 \in \Ob \Tw \W$ are homotopy equivalent in $\Tw \W$ if and only if $X^0 = \T^0(\bX^0)$ is homotopy equivalent to $X^1 = \T^0(\bX^1)$ in $\W$.
	In particular,
	\begin{enumerate}
		\item homotopy equivalences between objects of $\W$ remain homotopy equivalences in $\Tw \W$;
		\item $\Tw \W$ has enough homotopy equivalences;
		\item any $\bX^0 \in \Ob \Tw \W$ is homotopy equivalent to $\T^0(\bX^0) \in \Ob \W$; and
		\item homotopy equivalence of twisted complexes can be tested only on objects of $\W$.
	\end{enumerate}
\end{prop}

\begin{remark}
	If $\W$ has enough homotopy equivalences, but not strongly triangulated, then we do not know how to prove that homotopy inverse morphisms in $\W$ remain homotopy inverse in $\Tw \W$ or that $\Tw \W$ has enough homotopy equivalences.
	The issue is that we did not impose any conditions on $\bmu^n(\ldots, q, \ldots), \, n \ge 3$ for homotopy equivalences $q$.
	This can probably be fixed in a way similar to the rather involved homotopy unit framework of \cite{FOOO}*{\Sec 3.3}, but we do not need it since the Fukaya categories that we construct geometrically are strongly triangulated.
\end{remark}

\begin{proof}[Proof of \cref{p: homotopy equivalence in triangulated A_infty-pre-category}]
	Since homotopy equivalence is defined as having a chain of homotopy equivalences between objects, it is enough to show that $\T^1$ and $(\T^1)^{-1}$ map homotopy inverse morphisms to homotopy inverse morphisms.
	That a strict surjective $A_\infty$-functor which is isomorphism on morphism spaces preserves homotopy inverses is clear.

	It remains to show that if $e^1 \in \W(X^0, X^1)$ and $e^0 \in \W(X^1, X^0)$ are homotopy inverse morphisms in $\W$, then
	$\be^1 = \T_{\bX^0, \bX^1}^{-1}(e^1) \in \Tw \W(\bX^0, \bX^1)$
	and
	$\be^0 = \T_{\bX^1, \bX^0}^{-1}(e^0) \in \Tw \W(\bX^1, \bX^0)$
	are homotopy inverse morphisms in $\Tw \W$.
	To see that this is indeed the case let $\{\bT, \bX^0, \bX^1\}$ be a transverse set in $\Tw \W$.
	Denote $T = \T^0(\bT)$.
	If $H$ is a homotopy of $\W(T, X^0)$ certifying that $e^0 \circ (e^1 \circ -)$ is homotopic to $\id_{\W(T, X^0)}$, then $\T^{-1}_{\bT, \bX^0}H\T_{\bT, \bX^0}$ certifies that $\be^0 \circ (\be^1 \circ -)$ is homotopic to $\id_{\Tw \W(\bT, \bX^0)}$.
	The argument for the other three homotopies is analogous.

\end{proof}

Now suppose that an $A_\infty$-pre-category $\W$ with enough homotopy equivalences has two strict totalisation functors $\T_0, \T_1\colon \Tw \W \to \W$.
We make the following observation.

\begin{lemma} \label{l: comparison of two totalisation functors}
	Let $\bX^0$ and $\bX^1$ be two twisted complexes over $\W$.
	Then $\bX^0$ is homotopy equivalent to $\bX^1$ in $\Tw \W$ if and only if $\T^0_0(\bX^0)$ is homotopy equivalent to $\T^0_1(\bX^1)$ in $\W$.
\end{lemma}
\begin{proof}
	Applying \cref{p: homotopy equivalence in triangulated A_infty-pre-category} repeatedly we obtain two chains of homotopy equivalences
	$$
		\T^0_1(\bX^0) \cong \bX^0 \cong \T^0_0(\bX^0)
	$$
	and
	$$
		\T^0_1(\bX^1) \cong \bX^1 \cong \T^0_0(\bX^1),
	$$
	from which the claim clearly follows.
\end{proof}

\subsubsection{Shifts and cones}
\label{subsubs: Shifts and cones}
Let $\W$ be an $A_\infty$-pre-category where transverse sequences form transverse sets.
For a twisted complex $\bX = (\bigoplus_{i \in I} X_i[m_i], \delta)$ over $\W$ and $k \in \ZZ$ we define the \emph{shift}
$$
\bX[k] = \Big(\bigoplus_{i \in I} X_i[m_i + k], \delta\Big)
$$
where $\delta$ has the same components as before.
This is indeed a twisted complex since the Koszul sign in \cref{eq: mu with Koszul sign} only depends on the differences $m_{i_1} - m_{i_0}$, so $\A_{X[k]}^< = \A_X^<$.
The same observation applied to \cref{eq: morphism space between twisted complexes} and to the definition of $\bmu$ in $\Tw \W$ gives that $[k]$ with identity on morphisms is a strict autoequivalence of $\Tw \W$.

Now let $\bX^j = (\bigoplus_{i \in I^j} X^j_i[m^j_i], \delta^j)$, $j = 0, 1$ be a transverse pair of twisted complexes over $\W$ and $a \in \Tw \W^0(\bX^0, \bX^1)$ a closed morphism.
We define the \emph{cone} of $a$ to be the twisted complex
$$
\Cone(a) = \Big(\bigoplus_{i \in I^0} X^0_i[m^0_i + 1] \oplus \bigoplus_{i \in I^1} X^1_i[m^1_i], \, \delta^a \Big)
\text{ where } \delta^a_{i^0 i^1} = \delta^0_{i^0 i^1} + \delta^1_{i^0 i^1} + (-1)^{m^1_{i^1}} a_{i^0 i^1}.
$$
Extending the orders on $I^0$ and $I^1$ to $I^0 \sqcup I^1$ so that $I^0 < I^1$ we see that $\delta^a \in \A^{<}_{\Cone(a)}$.
A straightforward computation then shows that $\delta^a$ is a bounding cochain.
Note that if $\bX^1$ is an object of $\W$, then the sign in $\delta^a$ is trivial.

If $\W$ is strongly triangulated with a totalisation functor $\T$, then we define the \emph{shift} $X[k] = \T^0(X[k])$ of $X \in \Ob \W$ and the \emph{cone} $\Cone(a) = \T^0(\Cone(a))$ of a closed morphism $a \in \W^0(X^0, X^1)$, where on the right-hand sides $X$ and $a$ are viewed as an object and a morphism in $\Tw \W$.
In particular, $\W$ is closed under shifts and cones, i.e.\ triangulated in the sense of the \nameref{s: Introduction}.
If moreover $\W$ has enough homotopy equivalences, then by \cref{p: homotopy equivalence in triangulated A_infty-pre-category} the shift and the cone in $\W$ are homotopy equivalent in $\Tw \W$ to the corresponding twisted complexes.

We finish with a criterion for a morphism of twisted complexes to be a homotopy equivalence.

\begin{lemma} \label{l: upper triangular morphisms of twisted complexes}
	Let $\W$ be a strongly triangulated $A_\infty$-pre-category with enough homotopy equivalences and
	$\bX = (\bigoplus_{i \in I} X_i[m_i], \delta^{\bX})$ and $\bY = (\bigoplus_{i \in I} Y_i[m_i], \delta^{\bY})$
	a transverse pair of twisted complexes over $\W$ with respect to the same partial order $<$ on $I$.
	Let $q \in \Tw \W^0(\bX, \bY)$ be a closed morphism which is \emph{upper triangular}, i.e.\ $q_{i_0i_1} \ne 0$ implies $i_0 \le i_1$, and such that $q_{ii} \in \W^0(X_i, Y_i)$ is a homotopy equivalence in $\W$ for each $i \in I$.
	Then $q$ is a homotopy equivalence in $\Tw \W$.
\end{lemma}
\begin{proof}
	By \cref{p: homotopy equivalence in triangulated A_infty-pre-category} $\Tw \W$ has enough homotopy equivalences, so by \cref{l: homotopy equivalence can be tested pointwise} it suffices to show that for any transverse set $\{\bT, \bX, \bY\}$ in $\Tw \W$ the maps
	$\Tw \W(\bT, \bX) \xrightarrow{q \circ -} \Tw \W(\bT, \bY)$
	and
	$\Tw \W(\bY, \bT) \xrightarrow{- \circ q} \Tw \W(\bX, \bT)$
	are homotopy equivalences of complexes.
	We treat $q \circ -$, the argument for $- \circ q$ is analogous with the roles of sources and targets exchanged.

	Identify $I$ with $\{1, \ldots, n\}$ with the order on the latter extending $<$.
	For $r = 1, \ldots, n + 1$ let $F^r\Tw \W(\bT, \bX)$ be the sum of the summands of \cref{eq: morphism space between twisted complexes} with target $X_s$ for $s \ge r$, and similarly for $\bY$.
	Since the components of $\delta^{\bX}$, $\delta^{\bY}$, and $q$ do not decrease the index, the definition of $\bmu$ in $\Tw \W$ shows that $F^r$ are subcomplexes, which are direct summands as graded modules, and that $q \circ -$ preserves them.
	Moreover, $F^r/F^{r + 1}$ is the complex $\Tw \W(\bT, X_r[m_r])$ and $q \circ -$ induces
	$$
		\Tw \W(\bT, (X_r[m_r], 0)) \xrightarrow{q_{rr} \circ -} \Tw \W(\bT, (Y_r[m_r], 0))
	$$
	on it.
	Since the shift is a strict autoequivalence of $\Tw \W$, this map coincides with
	$\Tw \W(\bT[-m_r], X_r) \xrightarrow{q_{rr} \circ -} \Tw \W(\bT[-m_r], Y_r)$,
	which is a homotopy equivalence since $q_{rr}$ is a homotopy equivalence in $\Tw \W$ by \cref{p: homotopy equivalence in triangulated A_infty-pre-category}.

	It remains to recall that a map of complexes of $\k$-modules is a homotopy equivalence if and only if its cone is contractible, and that a complex with a finite filtration by subcomplexes which are direct summands as graded modules is contractible if all the associated graded pieces are.
	Applying this to the cone of $q \circ -$ with the induced filtration we conclude that $q \circ -$ is a homotopy equivalence.
\end{proof}

\subsection{From \texorpdfstring{$A_\infty$}{A∞}-pre-categories to \texorpdfstring{$A_\infty$}{A∞}-categories}
\label{subs: From A_infty-pre-categories to A_infty-categories}

Let $\W$ be an $A_\infty$-pre-category with enough homotopy equivalences.
We construct strictly unital $A_\infty$-categories spanned by collections of objects in $\W$.
The construction requires a collection of objects $\bX^0 = \{X^0_i\}_{i \in I}$ to admit a system of strict units $\bX^\bullet$, which we define in \cref{subsubs: system of strict units}.
In \cref{subsubs: A_infty-category spanned by a collection of objects} we construct the $A_\infty$-category $\W_{\bX^\bullet}$ with $\Ob \W_{\bX^\bullet} = \bX^0$, and in \cref{subsubs: A_infty category spanned by a collection of objects is well-defined} we prove that it is strictly unital and that two objects of $\bX^0$ are homotopy equivalent in $\W_{\bX^\bullet}$ if and only if they are homotopy equivalent in $\W$, see \cref{p: A_infty category spanned by a collection of objects is well-defined}.
In \cref{subsubs: Invariance} we introduce saturated collections $\SU$ of systems of strict units and prove that they give a strictly unital $A_\infty$-category $\W_{\SU}$, well-defined up to homotopy equivalence, see \cref{p: W_SU is independent of choices}.
In \cref{subsubs: Induced triangulated structure} we prove that $\W_{\SU}$ is triangulated if $\W$ is strongly triangulated, see \cref{l: W_SU is triangulated}.

\subsubsection{System of strict units}
\label{subsubs: system of strict units}
A \emph{system of strict units} on a collection $\{X_i^0\}_{i \in I}$ of objects in $\W$ is the data of
\begin{itemize}
	\item a transverse set $\{X^j_i\}^{j \in \ZZ_{\ge 0}}$ in $\W$ for each $i \in I$;
	\item homotopy equivalences $e_i^{j_0, j_1} \in \W(X_i^{j_0}, X_i^{j_1})$ for each $i \in I$ and $j_0 \ne j_1$
\end{itemize}
such that
\begin{enumerate}
	\item $e_i^{j_1, j_2} \circ e_i^{j_0, j_1} = e_i^{j_0, j_2}$
	for any $i \in I$ and pairwise distinct $j_0, j_1, j_2 \ge 0$;
	\item for any $i_0, i_1 \in I$ there exists $N > 0$ such that for any
	${j_0, j_1, j_1' \in \{0\} \cup \ZZ_{\ge N}}$
	satisfying
	$j_0 < \min\{j_1, j_1'\} < \max\{j_1, j_1'\}$
	the set
	$\{X_{i_0}^{j_0}, X_{i_1}^{j_1}, X_{i_1}^{j_1'}\}$
	is transverse and the map
	$$
		e_{i_1}^{j_1', j_1} \circ (e_{i_1}^{j_1, j_1'} \circ -)
		\text{ is the identity on }
		\W(X_{i_0}^{j_0}, X_{i_1}^{j_1});
	$$
	\item for any $i_0, i_1 \in I$ there exists $N > 0$ such that for any
	${j_0, j_0', j_1 \in \{0\} \cup \ZZ_{\ge N}}$
	satisfying
	$\min\{j_0, j_0'\} < \max\{j_0, j_0'\} < j_1$
	the set
	$\{X_{i_0}^{j_0}, X_{i_0}^{j_0'}, X_{i_1}^{j_1}\}$
	is transverse and the map
	$$
		(- \circ e_{i_0}^{j_0', j_0}) \circ e_{i_0}^{j_0, j_0'}
		\text{ is the identity on }
		\W(X_{i_0}^{j_0}, X_{i_1}^{j_1});
	$$
	\item for any $n \ge 3$, $0 \le l \le n$, and
	$i_0, \ldots, i_{l - 2}, i_{l - 1} = i = i_{l}, i_{l + 1}, \ldots, i_n \in I$
	there exists $N > 0$ such that for any
	$j_0, \ldots, j_n \in \{0\} \cup \ZZ_{\ge N}$
	satisfying
	$j_0 < \ldots < j_{l - 2} < \min\{j_{l - 1}, j_l\} < \max\{j_{l - 1}, j_l\} < j_{l + 1} < \ldots < j_n$
	the set
	$\{X_{i_0}^{j_0}, \ldots, X_{i_n}^{j_n}\}$
	is transverse and the map
	$$
	\xymatrix{
		\makebox[.5\textwidth][l]{$\W(X_{i_{n - 1}}^{j_{n - 1}}, X_{i_n}^{j_n}) \ldots \W(X_{i}^{j_{l}}, X_{i_{l + 1}}^{j_{l + 1}}) \otimes \W(X_{i_{l - 2}}^{j_{l - 2}}, X_i^{j_{l - 1}}) \ldots \W(X_{i_0}^{j_0}, X_{i_1}^{j_1})$} \ar[d]^-{\bmu^n(\ldots, e_i^{j_{l - 1}, j_l}, \ldots)}\\
		\W(X_{i_0}^{j_0}, X_{i_n}^{j_n})
	}
	$$
	vanishes.
\end{enumerate}

\subsubsection{\texorpdfstring{$A_\infty$}{A∞}-category spanned by a collection of objects}
\label{subsubs: A_infty-category spanned by a collection of objects}
Let
$$
\bX^\bullet = \left\{\left(\{X_i^j\}^{j = 0, \ldots, \infty}, \{e_i^{j_0, j_1}\}^{j_0 \ne j_1}\right)\right\}_{i \in I}
$$
be a system of strict units on a collection $\bX^0 = \{X^0_i\}_{i \in I}$ of objects in $\W$.
We define an $A_\infty$-category $\W_{\bX^\bullet}$ \emph{spanned by $\bX^0$} as follows.

\subsubsection*{Objects}
$\Ob \W_{\bX^\bullet} = \bX^0$.

\subsubsection*{Morphisms}
Fix $X_{i_0}^0, X_{i_1}^0 \in \bX^0$.
The morphism space is given by
$$
\W_{\bX^\bullet}(X_{i_0}^0, X_{i_1}^0) = \colim \W(X^0_{i_0}, X^{j_1}_{i_1})
$$
where the colimit is taken with respect to the directed system
$$
	e_{i_1}^{j_1, j_1'} \circ - \colon \W(X^0_{i_0}, X^{j_1}_{i_1}) \to \W(X^0_{i_0}, X^{j_1'}_{i_1})
$$
for $0 \ll j_1 < j_1'$.
Properties (1) and (4) in \cref{subsubs: system of strict units} and the $A_3$-relation imply that the above is indeed a directed system.
Property (2) implies that the natural map
$
\W(X_{i_0}^{0}, X_{i_1}^{j_1}) \to \W_{\bX^\bullet}(X_{i_0}^0, X_{i_1}^0)
$
is an isomorphism of graded $\k$-modules for $j_1 \gg 0$.
Combining this with Property (3), for some $N > 0$ we obtain isomorphisms
$$
\ph^{j_0, j_1}_{i_0, i_1}\colon \W_{\bX^\bullet}(X_{i_0}^0, X_{i_1}^0) \to \W(X_{i_0}^{j_0}, X_{i_1}^{j_1})
$$
for any $j_0, j_1 \in \{0\} \cup \ZZ_{\ge N}$ satisfying $j_0 < j_1$.

\subsubsection*{\texorpdfstring{$A_\infty$}{A∞}-operations}
For $X_{i_0}^0, \ldots, X_{i_n}^0 \in \bX^0$ the $A_\infty$-operation
$$
\bmu^n_{\W_{\bX^\bullet}}\colon \W_{\bX^\bullet}(X_{i_{n - 1}}^0, X_{i_n}^0) \otimes \ldots \otimes \W_{\bX^\bullet}(X_{i_0}^0, X_{i_1}^0) \to \W_{\bX^\bullet}(X_{i_0}^0, X_{i_n}^0)
$$
is determined by
$$
\ph^{j_0, j_n}_{i_0, i_n} \circ \bmu^n_{\W_{\bX^\bullet}} = \bmu^n_{\W} \circ \left(\ph^{j_{n - 1}, j_n}_{i_{n - 1}, i_n} \otimes \ldots \otimes \ph^{j_0, j_1}_{i_0, i_1}\right)
$$
for sufficiently large $N > 0$ and any $j_0, \ldots, j_n \in \{0\} \cup \ZZ_{\ge N}$ satisfying ${j_0 < \ldots < j_n}$.
We check that this definition makes sense in \cref{p: A_infty category spanned by a collection of objects is well-defined}.

\subsubsection{}
\label{subsubs: A_infty category spanned by a collection of objects is well-defined}
We check that the definition of an $A_\infty$-category spanned by a collection of objects in \cref{subsubs: A_infty-category spanned by a collection of objects} makes sense and prove that homotopy equivalence in $\W$ and $\W_{\bX^\bullet}$ is the same.

\begin{prop} \label{p: A_infty category spanned by a collection of objects is well-defined}
	Let
	$\bX^\bullet = \left\{\left(\{X_i^j\}^{j = 0, \ldots, \infty}, \{e_i^{j_0, j_1}\}^{j_0 \ne j_1}\right)\right\}_{i \in I}$
	be a system of strict units on a collection of objects ${\bX^0 = \{X_i^0\}_{i \in I}}$ in an $A_\infty$-pre-category $\W$ with enough homotopy equivalences.

	Then $\W_{\bX^\bullet}$ is a strictly unital $A_\infty$-category with the strict unit $e_i \in \W_{\bX^\bullet}(X_i^0, X_i^0)$ defined by $\ph^{j_0, j_1}_{i, i}(e_i) = e^{j_0, j_1}_i$ for sufficiently large $N > 0$ and $j_0, j_1 \in \{0\} \cup \ZZ_{\ge N}$ satisfying $j_0 < j_1$.

	Moreover, $X^0_{i_0}$ and $X^0_{i_1}$ are homotopy equivalent in $\W$ if and only if they are homotopy equivalent in $\W_{\bX^\bullet}$.
\end{prop}
\begin{proof}
	That $\W_{\bX^\bullet}$ is a strictly unital $A_\infty$-category follows from the statement below.
	\begin{claim}
		Let $n \ge 1$ and $i_0, \ldots, i_n \in I$.
		Then there exists $N > 0$ such that for any $j_0, \ldots, j_n \in \{0\} \cup \ZZ_{\ge N}$ satisfying ${j_0 < \ldots < j_n}$ the following holds.
		\begin{enumerate}
			\item if
			$j_0' \in \{0\} \cup \ZZ_{\ge N}$
			is such that
			$j_0 < j_0' < j_1$,
			then
			$$
				\bmu^n_{\W}\left(a_n, \ldots, a_2, a_1 \circ e_{i_0}^{j_0', j_0}\right) =
				\bmu^n_{\W}\left(a_n, \ldots, a_1\right) \circ e_{i_0}^{j_0', j_0}
			$$
			for any sequence of homogeneous morphisms
			$X_{i_0}^{j_0} \xrightarrow{a_1} \ldots \xrightarrow{a_n} X_{i_n}^{j_n}$
			in $\W$;
			\item if
			$j_l' \in \{0\} \cup \ZZ_{\ge N}$
			is such that
			$j_l < j_l' < j_{l + 1}$,
			then
			$$
				\bmu^n_{\W}\left(a_n, \ldots, a_{l + 2}, a_{l + 1} \circ e_{i_l}^{j_l', j_l}, e_{i_l}^{j_l, j_l'} \circ a_{l}, a_{l - 1}, \ldots, a_1\right) =
				\bmu^n_{\W}\left(a_n, \ldots, a_1\right)
			$$
			for any sequence of homogeneous morphisms
			$X_{i_0}^{j_0} \xrightarrow{a_1} \ldots \xrightarrow{a_n} X_{i_n}^{j_n}$
			in $\W$;
			\item if
			$j_n' \in \{0\} \cup \ZZ_{\ge N}$
			is such that
			$j_n < j_n'$,
			then
			$$
				\bmu^n_{\W}\left(e_{i_n}^{j_n, j_n'} \circ a_n, a_{n - 1}\ldots, a_1\right) =
				e_{i_n}^{j_n, j_n'} \circ\bmu^n_{\W}\left(a_n, \ldots, a_1\right)
			$$
			for any sequence of homogeneous morphisms
			$X_{i_0}^{j_0} \xrightarrow{a_1} \ldots \xrightarrow{a_n} X_{i_n}^{j_n}$
			in $\W$;
		\end{enumerate}
	\end{claim}

	The Claim implies that for $n \ge 1$ and $i_0, \ldots, i_n \in I$ there exists $N > 0$ such that
	$$
		\left(\ph^{j_0, j_n}_{i_0, i_n}\right)^{-1} \circ \bmu^n_{\W} \circ \left(\ph^{j_{n - 1}, j_n}_{i_{n - 1}, i_n} \otimes \ldots \otimes \ph^{j_0, j_1}_{i_0, i_1}\right)
	$$
	does not depend on the choice of $j_0, \ldots, j_n \in \{0\} \cup \ZZ_{\ge N}$ satisfying ${j_0 < \ldots < j_n}$.
	Indeed, if we have two such choices $j_0' < \ldots < j_n'$ and $j_0'' < \ldots < j_n''$, then they give the same result as $j_0 < \ldots < j_n$ for $j_0 > \max\{j_n', j_n'', N\}$.
	Therefore, $\bmu_{\W_{\bX^\bullet}}$ is well-defined and satisfies the $A_\infty$-relations.

	Property (1) in \cref{subsubs: system of strict units} and the definition of $\ph$ imply that
	$\left(\ph^{j_0, j_1}_{i, i}\right)^{-1}\left(e^{j_0, j_1}_i\right)$
	does not depend on the choice of $j_0, j_1 \in \{0\} \cup \ZZ_{\ge N}$ satisfying $j_0 < j_1$ for sufficiently large $N > 0$.
	Since $e^{j_0, j_1}_i$ is a homotopy equivalence, $e_i$ is a closed morphism of degree zero.
	Properties (2) and (3) imply that $e_i$ satisfies the part of the definition of strict unit involving $\circ$, and Property (4) implies that $e_i$ satisfies the remaining part of the definition, involving $\bmu^n$ for $n \ge 3$.

	Suppose $X^0_{i_0}$ and $X^0_{i_1}$ are homotopy equivalent in $\W$.
	Then by \cref{l: homotopy equivalent transverse objects} there exists a homotopy equivalence $X^{j_0}_{i_0} \xrightarrow{q^{j_0, j_1}} X^{j_1}_{i_1}$ in $\W$ for $0 \ll j_0 < j_1$ hence
	$q = \left(\ph_{i_0, i_1}^{j_0, j_1}\right)^{-1}(q^{j_0, j_1})$
	is a homotopy equivalence in $\W_{\bX^\bullet}$ by \cref{l: homotopy equivalence can be tested pointwise} since $q \circ -$ (resp.\ $- \circ q$) on $\W_{\bX^\bullet}(X_{i_2}^0, X_{i_0}^0)$ (resp.\ on $\W_{\bX^\bullet}(X_{i_1}^0, X_{i_2}^0)$) can be computed as $\ph_{i_0, i_1}^{j_0, j_1}(q) \circ -$ (resp.\ $- \circ \ph_{i_0, i_1}^{j_0, j_1}(q)$) on $\W(X_{i_2}^{j_2}, X_{i_0}^{j_0})$ (resp.\ $\W(X_{i_1}^{j_1}, X_{i_2}^{j_2})$) for $0 \ll j_2 < j_0 < j_1$ (resp.\ $0 \ll j_0 < j_1 < j_2$).

	Now suppose that $X^0_{i_0}$ and $X^0_{i_1}$ are homotopy equivalent in $\W_{\bX^\bullet}$.
	Then by \cref{l: characterisation of homotopy equivalences in a strictly unital A_infty-category} there exist $j_0 < j_1 < j_0' < j_1'$ and closed morphisms of degree zero
	$X^{j_0}_{i_0} \xrightarrow{q^{j_0, j_1}} X^{j_1}_{i_1} \xrightarrow{q^{j_1, j_0'}} X^{j_0'}_{i_0} \xrightarrow{q^{j_0', j_1'}} X^{j_1'}_{i_1}$
	such that $q^{j_1, j_0'} \circ q^{j_0, j_1}$ is cohomologous to $e_{i_0}^{j_0, j_0'}$ and $q^{j_0', j_1'} \circ q^{j_1, j_0'}$ is cohomologous to $e_{i_1}^{j_1, j_1'}$.
	By \cref{l: exact morphisms and contractible objects in A_infty-pre-category}\,(1) and \cref{l: 2-out-of-6 property for homotopy equivalences} $q^{j_0, j_1}$ is a homotopy equivalence and hence $X^0_{i_0}$ and $X^0_{i_1}$ are homotopy equivalent in $\W$.
\end{proof}

\begin{proof}[Proof of Claim]
	Take $N > 0$ such that Properties (2) and (3) in \cref{subsubs: system of strict units} are satisfied for any $i, i' \in \{i_0, \ldots, i_n\}$ and Property (4) is satisfied for any tuple $(i_0', \ldots, i_{n+2}')$ where $i_0', \ldots, i_{n + 2}' \in \{i_0, \ldots, i_n\}$.

	(1) of the Claim is equivalent to
	$$
		\bmu^n\left(a_n, \ldots, a_2, \bmu^2(a_1, e_{i_0}^{j_0', j_0})\right) + (-1)^{0 - 1}
		\bmu^2\left(\bmu^n\left(a_n, \ldots, a_1\right), e_{i_0}^{j_0', j_0}\right) = 0,
	$$
	which is the $A_{n + 1}$-relation since $\bmu^1(e_{i_0}^{j_0', j_0}) = 0$ and all other terms vanish by Property (4).

	(3) of the Claim is equivalent to
	$$
	\begin{aligned}
		&(-1)^{|a_n|}\bmu^n\left(\bmu^2(e_{i_n}^{j_n, j_n'}, a_n), a_{n - 1} \ldots, a_1\right) +\\
		&(-1)^{|a_n| + \ldots + |a_1| - n - 1}\bmu^2(e_{i_n}^{j_n, j_n'}, \bmu^n\left(a_n, \ldots, a_1\right)) = 0,
	\end{aligned}
	$$
	which is a multiple of the $A_{n + 1}$-relation since $\bmu^1(e_{i_n}^{j_n', j_n}) = 0$ and all other terms vanish by Property (4).

	To prove (2) of the Claim consider the $A_{n + 1}$-relation for the sequence of morphisms
	$X_{i_0}^{j_0} \xrightarrow{a_1} \ldots \xrightarrow{a_{l}} X_{i_{l}}^{j_{l}} \xrightarrow{e_{i_l}^{j_l, j_l'}} X_{i_{l}}^{j_{l}'} \xrightarrow{\bmu^2(a_{l + 1}, e_{i_l}^{j_l', j_l})} X_{i_{l + 1}}^{j_{l + 1}} \xrightarrow{a_{l + 2}} \ldots \xrightarrow{a_n} X_{i_n}^{j_n}$
	and use that $e_{i_l}^{j_l, j_l'}$ is closed and Property (4) to get
	$$
	\begin{aligned}
		&(-1)^{|a_{l - 1}| + \ldots + |a_1| + l}\bmu^n\left(a_n, \ldots, a_{l + 2}, \bmu^2(a_{l + 1}, e_{i_l}^{j_l', j_l}), \bmu^2(e_{i_l}^{j_l, j_l'}, a_{l}), a_{l - 1}, \ldots, a_1\right) +\\
		&(-1)^{|a_l| + \ldots + |a_1| + l - 1}\bmu^{n}\left(a_n, \ldots, a_{l + 2}, \bmu^2\left(\bmu^2(a_{l + 1}, e_{i_l}^{j_l', j_l}), e_{i_l}^{j_l, j_l'}\right), a_l, \ldots, a_1\right) = 0.
	\end{aligned}
	$$
	Using Property (3) we obtain
	$$
	\begin{aligned}
		&\bmu^n_{\W}\left(a_n, \ldots, a_{l + 2}, a_{l + 1} \circ e_{i_l}^{j_l', j_l}, e_{i_l}^{j_l, j_l'} \circ a_{l}, a_{l - 1}, \ldots, a_1\right) =\\
		&(-1)^{|a_l|}\bmu^n\left(a_n, \ldots, a_{l + 2}, \bmu^2(a_{l + 1}, e_{i_l}^{j_l', j_l}), \bmu^2(e_{i_l}^{j_l, j_l'}, a_{l}), a_{l - 1}, \ldots, a_1\right)=\\
		&\bmu^{n}\left(a_n, \ldots, a_{l + 2}, \bmu^2\left(\bmu^2(a_{l + 1}, e_{i_l}^{j_l', j_l}), e_{i_l}^{j_l, j_l'}\right), a_l, \ldots, a_1\right) =
		\bmu^n_{\W}\left(a_n, \ldots, a_1\right)
	\end{aligned}
	$$
	as required.
\end{proof}

\subsubsection{Invariance}
\label{subsubs: Invariance}
A strict $A_\infty$-functor $\F$ is called a \emph{strict equivalence} if $\F^0$ is a bijection and $\F^1$ is an isomorphism of graded $\k$-modules on each morphism space in $\W$.
Note that this is equivalent to $\F$ having an inverse \emph{strict} $A_\infty$-functor.
A strict equivalence from $\W$ to itself is called a \emph{strict autoequivalence} of $\W$.

Let
$$
\bX^\bullet = \left\{\left(\{X_i^j\}^{j = 0, \ldots, \infty}, \{e_i^{j_0, j_1}\}^{j_0 \ne j_1}\right)\right\}_{i \in I}
$$
be a system of strict units in $\W$ and $\I$ a strict autoequivalence of $\W$.
Then
$$
	\I_*\bX^\bullet = \left\{\left(\{\I^0 X_i^j\}^{j = 0, \ldots, \infty}, \{\I^1 e_i^{j_0, j_1}\}^{j_0 \ne j_1}\right)\right\}_{i \in I}
$$
is also a system of strict units in $\W$ and $\I$ induces a strict equivalence of strictly unital $A_\infty$-categories ${\W_{\bX^\bullet} \to \W_{\I_*\bX^\bullet}}$.

We call a system of strict units $\bX^\bullet$ \emph{full} if any $X \in \Ob \W$ is homotopy equivalent in $\W$ to some $X^0 \in \bX^0$.
\begin{defn} \label{d: saturated collection of systems of strict units}
A collection $\SU$ of systems of strict units in $\W$ is called \emph{saturated} if
\begin{enumerate}
	\item $\SU$ contains a full system of strict units;
	\item for any two systems of strict units $\bX^\bullet_0, \bX^\bullet_1 \in \SU$ there exists a strict autoequivalence $\I$ of $\W$ such that $\I_*\bX^\bullet_0 \sqcup \bX^\bullet_1 \in \SU$.
\end{enumerate}
\end{defn}

\begin{prop} \label{p: W_SU is independent of choices}
	Let $\SU$ be a saturated collection of systems of strict units in $\W$.
	Then the strictly unital $A_\infty$-categories $\W_{\bX^\bullet}$ for full $\bX^\bullet \in \SU$ are all homotopy equivalent.
\end{prop}
We denote (the homotopy equivalence class of) any of the $\W_{\bX^\bullet}$ with $\bX^\bullet \in \SU$ full by $\W_{\SU}$; such $\bX^\bullet$ exists by \cref{d: saturated collection of systems of strict units}\,(1).
Since $\bX^\bullet$ is full, by \cref{p: A_infty category spanned by a collection of objects is well-defined} every object of $\W$ is homotopy equivalent in $\W$ to an object of $\W_{\SU}$, and two objects of $\W_{\SU}$ are homotopy equivalent in $\W$ if and only if they are homotopy equivalent in $\W_{\SU}$.
Moreover, for any $\bX^\bullet \in \SU$, applying \cref{d: saturated collection of systems of strict units}\,(2) to $\bX^\bullet$ and a full $\bY^\bullet \in \SU$ we get a strict autoequivalence $\I$ of $\W$ such that $\I_*\bX^\bullet \sqcup \bY^\bullet \in \SU$.
The latter system of strict units is full, so $\I$ induces a strict full embedding
$\W_{\bX^\bullet} \hookrightarrow \W_{\I_*\bX^\bullet \sqcup \bY^\bullet} = \W_{\SU}$.

The homotopy equivalence of strictly unital $A_\infty$-categories in \cref{p: W_SU is independent of choices} is defined as follows.
If $\C$ and $\D$ are two strictly unital $A_\infty$-categories, then an $A_\infty$-functor $\F\colon \C \to \D$ is called \emph{strictly unital} if $\F^1(1_X) = 1_{\F^0(X)}$ and $\F^n(\ldots, 1_X, \ldots) = 0$ for $n \ge 2$ for any $X \in \Ob \C$.
A strictly unital $A_\infty$-functor $\F\colon \C \to \D$ is called a \emph{homotopy equivalence} if
\begin{itemize}
	\item For each $X_0, X_1 \in \Ob \C$ the map ${\F^1\colon \C(X_0, X_1) \to \D(\F^0(X_0), \F^0(X_1))}$ is a homotopy equivalence; and
	\item each $Y \in \Ob \D$ is homotopy equivalent to $\F^0(X)$ for some $X \in \Ob \C$.
\end{itemize}
Two strictly unital $A_\infty$-categories $\C$ and $\D$ are called \emph{homotopy equivalent} if there exists a finite sequence of strictly unital $A_\infty$-categories $\C_0, \ldots, \C_n$, $\C_0 = \C$, $\C_n = \D$ such that for any $i = 1, \ldots, n$ there exists a homotopy equivalence from $\C_{i - 1}$ to $\C_i$ or vice versa.

\begin{proof}[Proof of \cref{p: W_SU is independent of choices}]
	Consider two full systems of strict units $\bX^\bullet_0, \bX^\bullet_1 \in \SU$.
	Pick a strict autoequivalence $\I$ of $\W$ such that $\I_*\bX^\bullet_0 \sqcup \bX^\bullet_1 \in \SU$ and set $\V = \W_{\I_*\bX^\bullet_0 \sqcup \bX^\bullet_1}$.
	Then the inclusion $\W_{\bX^\bullet_1} \hookrightarrow \V$ and the functor $\W_{\bX^\bullet_0} \to \V$ induced by $\I$ are strict full embeddings.
	Let $X$ be an object of $\V$.
	Since $\bX^\bullet_1$ and $\bX^\bullet_0$ are full (and $\I$ is a strict autoequivalence of $\W$), we see that $X$ is homotopy equivalent in $\W$ to an object of $\bX^0_1$ and to an object of $\I^0\bX^0_0$, hence also in $\V$ by \cref{p: A_infty category spanned by a collection of objects is well-defined}.
	Thus the inclusions above are homotopy equivalences of strictly unital $A_\infty$-categories hence $\W_{\bX^\bullet_0}$ and $\W_{\bX^\bullet_1}$ are homotopy equivalent.
\end{proof}

\subsubsection{Induced triangulated structure}
\label{subsubs: Induced triangulated structure}
We show that triangulated structure on $\W$ induces a triangulated structure on $\W_{\SU}$.

\begin{lemma} \label{l: W_SU is triangulated}
	Let $\W$ be a strongly triangulated $A_\infty$-pre-category with enough homotopy equivalences and $\SU$ a saturated collection of systems of strict units on $\W$.
	Then $\W_{\SU}$ is triangulated.
\end{lemma}
Recall that a strictly unital $A_\infty$-category $\C$ is triangulated if any twisted complex over $\C$ is homotopy equivalent to an object of $\C$.
It is clear that being triangulated is invariant under homotopy equivalence of strictly unital $A_\infty$-categories.

\begin{proof}
	Fix a full system of strict units $\bX^\bullet \in \SU$.
	Let
	$\bY = (\bigoplus_{\alpha \in A} X^0_{i_\alpha}[m_\alpha], \delta)$
	be a twisted complex over $\W_{\bX^\bullet}$.
	We need to show $\bY$ is homotopy equivalent to an object of $\W_{\bX^\bullet}$.
	Choose $j_\alpha \gg 0$ for each $\alpha \in A$ such that $j_{\alpha_0} < j_{\alpha_1}$ and
	$
		\ph^{j_{\alpha_0}, j_{\alpha_1}}_{i_{\alpha_0}, i_{\alpha_1}}
	$
	is an isomorphism whenever $\alpha_0 < \alpha_1$, and all $\bmu_{\W_{\bX^\bullet}}$ on components of $\delta$ can be computed in
	$$
		\bigoplus_{\alpha_0 < \alpha_1} \W(X^{j_{\alpha_0}}_{i_{\alpha_0}}, X^{j_{\alpha_1}}_{i_{\alpha_1}})[m_{\alpha_1} - m_{\alpha_0}]
	$$
	via
	$\ph^{j_{\alpha_0}, j_{\alpha_1}}_{i_{\alpha_0}, i_{\alpha_1}}$
	and Koszul sign.
	We then obtain a twisted complex
	$\bY' = (\bigoplus_{\alpha \in A} X^{j_\alpha}_{i_\alpha}[m_\alpha], \delta')$
	over $\W$ where
	$\delta'_{{\alpha_0\alpha_1}} = \ph^{j_{\alpha_0}, j_{\alpha_1}}_{i_{\alpha_0}, i_{\alpha_1}}(\delta_{\alpha_0\alpha_1})$.
	Let $X^0_k \in \Ob \W_{\bX^\bullet}$ be an object homotopy equivalent to $\T\bY'$ in $\W$, which exists since $\bX^\bullet$ is full.
	We prove that $\bY$ is homotopy equivalent to $X^0_k$ in $\Tw \W_{\bX^\bullet}$.

	By (4) in \cref{subsubs: system of strict units} we can assume that
	$\{X^{j_\alpha}_{i_\alpha}\}_{\alpha \in A} \sqcup \{X^0_k,X^j_k\}$
	is a transverse set for some $j \gg \max\{j_\alpha\}_{\alpha \in A}$.
	Define a twisted complex
	$\bY'' = (\bigoplus_{\alpha \in A} X^{j'_\alpha}_{i_\alpha}[m_\alpha], \delta'')$
	analogously to $\bY'$, but now with all $j'_\alpha \gg j$.
	Again by (4) in \cref{subsubs: system of strict units} we can assume that
	$\{X^{j_\alpha}_{i_\alpha}\}_{\alpha \in A} \sqcup \{X^0_k,X^j_k\} \sqcup \{X^{j'_\alpha}_{i_\alpha}\}_{\alpha \in A}$
	is a transverse set.
	By \cref{l: homotopy equivalent transverse objects,p: homotopy equivalence in triangulated A_infty-pre-category} we obtain a homotopy equivalence $q \in \Tw \W^0(\bY', X^j_k)$.
	Since $q$ is a homotopy equivalence, there is a morphism $p \in \Tw \W^0(X^0_k, \bY')$ such that $q \circ p$ is cohomologous to $e_k^{0j}$ and a morphism $p' \in \Tw \W^0(X^j_k, \bY'')$ such that $p' \circ q$ is cohomologous to
	$e_{\bY'\bY''} = \sum_{\alpha \in A} e^{j_\alpha j'_\alpha}_{i_\alpha}$.
	All of these morphisms define the corresponding morphisms in $\Tw \W_{\bX^\bullet}$ via $\ph^{-1}$, which is assembled from the inverses to $\ph^{j_0j_1}_{i_0i_1}$.
	This yields $\ph^{-1}(q) \circ \ph^{-1}(p)$ cohomologous to $1_{X^0_k}$ and $\ph^{-1}(p') \circ \ph^{-1}(q)$ cohomologous to $1_{\bY}$.
	Therefore, $\ph^{-1}(q) \in \Tw \W_{\bX^\bullet}^0(\bY, X^0_k)$ is a homotopy equivalence by \cref{l: 2-out-of-6 property for homotopy equivalences}.
\end{proof}

\section{Geometric Fukaya category}
\label{s: Geometric Fukaya category}
In this section we define the geometric Fukaya $A_\infty$-pre-category, show that it satisfies the assumptions of the machinery of \cref{s: A_infty-machinery}, and prove \cref{t: A_infty-categories of arc systems,t: geometric Fukaya A_infty-pre-category is good,t: geometric Fukaya strictly unital A_infty-category,t: topological Fukaya category inside geometric}.

In \cref{subs: Fukaya A_infty-pre-category} we define the geometric Fukaya $A_\infty$-pre-category $\Wgp(\S; \C)$, whose $A_\infty$-operations are computed in obstruction algebras of composite Lagrangian branes, and prove its invariance properties.
In \cref{subs: Transverse perturbations} we perturb morphisms of local systems on a Lagrangian brane to morphisms between the brane and its transverse perturbation and compute $A_\infty$-operations involving them.
Using these perturbed morphisms, we prove that $\Wgp(\S; \C)$ has enough homotopy equivalences in \cref{subs: Homotopy equivalences in geometric Fukaya category}, construct a saturated collection of systems of strict units on it in \cref{subs: Strict units in geometric Fukaya category}, and prove that it is strongly triangulated in \cref{subs: Triangulated structure on geometric Fukaya A_infty-pre-category}, see \cref{p: W(S) has enough homotopy equivalences,p: saturated collection of systems of strict units on Fukaya category,p: W(S) is triangulated}.
In \cref{subs: Embedding of topological Fukaya category} we identify the $A_\infty$-category of an arc brane system with the $A_\infty$-category of the corresponding tame collection of Lagrangian branes and prove that full arc brane systems generate the same subcategory.
In \cref{subs: Proofs of the main theorems} we deduce \cref{t: A_infty-categories of arc systems,t: geometric Fukaya A_infty-pre-category is good,t: geometric Fukaya strictly unital A_infty-category,t: topological Fukaya category inside geometric}.

\subsection{Fukaya \texorpdfstring{$A_\infty$}{A∞}-pre-category}
\label{subs: Fukaya A_infty-pre-category}
Fix a weighted Liouville surface $(\S = (S, Q, \eta, \lambda); \C)$.

In \cref{subsubs: Transverse Lagrangians} we define transverse tuples of Lagrangians and their composite Lagrangians, obtained by shifting the potentials vertically, and show that a suitable shift makes the actions of any given finitely many intersection points and Reeb chords positive in \cref{l: make actions positive by vertical shift}.
In \cref{subsubs: Fukaya category} we define the geometric Fukaya $A_\infty$-pre-category $\Wgp(\S; \C)$, whose $A_\infty$-operations are computed in the obstruction algebras of composite Lagrangian branes.
In \cref{subsubs: Well-definedness of geometric Fukaya category} we prove that these operations do not depend on the choice of the shift and that $\Wgp(\S; \C)$ is an $A_\infty$-pre-category, see \cref{p: W is well-defined}.
In \cref{subsubs: Invariance of Fukaya category} we prove that graded exact symplectomorphisms induce strict equivalences of geometric Fukaya $A_\infty$-pre-categories, and deduce that $\Wgp(\S; \C)$ is invariant under isotopies of $\eta$ and $\C$ and does not depend on the choice of $\lambda$ up to strict equivalence.

\subsubsection{Transverse Lagrangians}
\label{subsubs: Transverse Lagrangians}
Let $L_i = (I_i, \gamma_i, P_i), \, i = 0, \ldots, n$ be Lagrangians in $\S$ such that
$$
\gamma = \bigsqcup_{i = 0}^n \gamma_i \colon I = \bigsqcup_{i = 0}^n I_i \to S
$$
has at worst transverse double points and $\res{\gamma}{\pd I}$ is an embedding.
In this case we say that $L_0, \ldots, L_n$ are \emph{transverse}.
Note that $\gamma$ is an immersion satisfying conditions \cref{item: boundary of gamma,item: gamma avoid D and Lambda,item: gamma has at worst transcerse double points,item: ends of gamma are distinct} from \cref{d: Lagrangian}.

Given a vector $h = (h_0, \ldots, h_n) \in \RR^{n + 1}$ define
$$
P_h \colon I \to \RR, \, P_h(t) = P_i(t) + h_i \text{ if } t \in I_i.
$$
The data $L_h = L_h(L_0, \ldots, L_n) = (I, \gamma, P_h)$ satisfies all conditions on being a Lagrangian in \cref{d: Lagrangian} except that $P_h$ can fail to satisfy conditions \cref{item: P at double point,item: P for a Reeb chord}.
Call the vector $h \in \RR^{n + 1}$ \emph{admissible} for $L_0, \ldots, L_n$ if $L_h$ is a Lagrangian and let $\H(L_0, \ldots, L_n)$ be the set of admissible shift vectors.
For $h \in \H(L_0, \ldots, L_n)$ we call the Lagrangian $L_h$ the \emph{composite Lagrangian} of $L_0, \ldots, L_n$ with \emph{shift vector} $h$.

For $i, j \in \{0, \ldots n\}, i \ne j$ let $Z(L_i, L_j)$ be the union of the set $X(L_i, L_j)$ of intersection points of $\gamma_i$ and $\gamma_j$ and the set $Y(L_i, L_j)$ of Reeb chords from an end of $\gamma_i$ to an end of $\gamma_j$.
Note that $Z(L_i, L_j)$ is infinite if $L_i$ and $L_j$ both have ends on a component $C$ of $\pd S$ with $C \cap \Lambda = \varnothing$.
The set $Z(L_i, L_j)$ will give rise to the morphism space in the Fukaya category between Lagrangian branes supported on $L_i$ and $L_j$.

For $z \in X(L_i, L_j)$ let $z_+ \in I_i$ and $z_- \in I_j$ be such that $\gamma_i(z_+) = z = \gamma_j(z_-)$, and for $z \in Y(L_i, L_j)$ let $z_+ \in I_i$ and $z_- \in I_j$ be such that $z$ goes from $\gamma_i(z_+)$ to $\gamma_j(z_-)$.

Given a (not necessarily admissible) shift vector $h$ and $z \in Z(L_i, L_j)$ for $i \ne j$, we define the action
\begin{equation} \label{eq: def of A_h(z)}
A_h(z) = \begin{cases}
	P_h(z_+) - P_h(z_-), &\text{if } z \text{ is an intersection point};\\
	P_h(z_+) - \int_{z} \lambda - P_h(z_-), &\text{if } z \text{ is a Reeb chord}.
\end{cases}
\end{equation}
Note that $h$ being admissible is equivalent to the non-vanishing of $A_h$ for all $z$, but, unlike the action in \cref{eq: def of A(x) for double point}, $A_h(z)$ can be negative for intersection points.
Observe that
$$
A_h(z) - A_0(z) = h_i - h_j \text{ for any } z \in Z(L_i, L_j).
$$

We therefore have a characterisation of admissible shift vectors,
$$
\H(L_0, \ldots, L_n) = \left\{h \in \RR^{n + 1} \colon h_j - h_i \ne A_0(z) \text{ for all} \quad \begin{aligned}
&i, j \in \{0, \ldots n\}, i \ne j\\
&\text{and } z \in Z(L_i, L_j)
\end{aligned}\right\}.
$$
That is, $\H(L_0, \ldots, L_n)$ is the complement of a hyperplane arrangement in $\RR^{n + 1}$.
The number of hyperplanes in the arrangement can be infinite, but the arrangement is discrete since ${\{z \in Z(L_i, L_j) \colon A_0(z) > B\}}$ is finite for all $i, j \in \{0, \ldots n\}, i \ne j$ and $B \in \RR$ by \cref{l: Y^>C is finite}.

Let $L_h$ be a composite Lagrangian of $L_0, \ldots, L_n$.
Fix $i, j \in \{0, \ldots n\}, i \ne j$.
Then $A_h$ and $A_{L_h}$ are two action functions defined on $Z(L_i, L_j)$.
Note that
$$
	A_{L_h}(z) = \begin{cases}
		|A_h(z)|, &\text{if } z \in X(L_i, L_j);\\
		A_h(z), &\text{if } z \in Y(L_i, L_j).
	\end{cases}
$$
We finish this discussion with the following easy but useful lemma.

\begin{lemma} \label{l: make actions positive by vertical shift}
	Given finite sets $Z_{ij} \subset Z(L_i, L_j), \, 0 \le i < j \le n$ there exists an admissible shift vector $h \in \RR^{n + 1}$ such that for any $0 \le i < j \le n$ we have
	\begin{equation} \label{eq: inequalities on A_h(z)}
		A_h(z) > 0 \text{ for any } z \in Z_{ij} \text{ and } A_h(z) < 0 \text{ for any } z \in Z(L_j, L_i)
	\end{equation}
\end{lemma}
\begin{proof}
	Recall that $A_h(z) = A_0(z) + h_i - h_j$ for $z \in Z(L_i, L_j)$ and that ${\{z \in Z(L_i, L_j) \colon A_0(z) > 0\}}$ is finite for all $i, j \in \{0, \ldots n\}, i \ne j$.
	Therefore, for each $0 \le i < j \le n$ there exists
	$$
		M_{ij} = \min\big(\{0\} \cup \{A_0(z) \mid z \in Z_{ij}\} \cup \{-A_0(z) \mid z \in Z(L_j, L_i)\}\big) \in \RR.
	$$
	Observe that if $h_j - h_i < M_{ij}$ for each $0 \le i < j \le n$, then \cref{eq: inequalities on A_h(z)} holds.
	Since the left-hand sides of these inequalities are given by linear combinations with non-negative coefficients of $n$ linearly independent functionals $h \mapsto h_i - h_{i - 1}, \, i = 1, \ldots, n$, these inequalities define a non-empty open subset of $\RR^{n + 1}$.
	The lemma follows since $\H(L_0, \ldots, L_n)$ is dense in $\RR^{n + 1}$.
\end{proof}

If we are given a grading $\gr_i$ for each $L_i$, then we get a grading $\gr$ on $\gamma$ and hence can define the degree $|z|$ as in \cref{subsubs: Degrees of double points and Reeb chords}, but with the new definition of $z_{\pm}$.

Let
$\L_i = (I_i, \gamma_i, P_i, \gr_i, \M_i), \, i = 0, \ldots, n$
be Lagrangian pre-branes such that their underlying Lagrangians are transverse and $h \in \RR^{n+1}$ an admissible shift vector, i.e.\ $L = L_h(L_0, \ldots, L_n)$ is a Lagrangian.
Then
${\L = (L, \gr, \sqcup \M_i) = \L_h(\L_0, \ldots, \L_n)}$
is a Lagrangian pre-brane, so the obstruction algebra $(\A_\L, \mu_\L)$ is defined.
For any $i =0, \ldots, n$ and $z \in Z_{L_i}$ we have $\A_\L(z) = \A_{\L_i}(z)$,
and if $z^0, \ldots, z^k \in Z_L$ are such that $z^j \in Z_{L_{i_j}}$ for each $j = 0, \ldots, k$, then
$$
\mu^{z^0; z^k, \ldots, z^1}_{\L} = \begin{cases}
	\mu^{z^0; z^k, \ldots, z^1}_{i}, & \text{if } i_0 = \ldots = i_k = i;\\
	0, & \text{otherwise}.
\end{cases}
$$
In other words, $\bigoplus_{i = 0}^n (\A_{\L_i}, \mu_i)$ is an $A_\infty$-subalgebra of $(\A_\L, \mu_\L)$.

If $\delta_i \in \A_{\L_i}$, $i = 0, \ldots, n$ are bounding cochains, i.e.\ $\LL_i = (\L_i, \delta_i)$ are Lagrangian branes, then we obtain a bounding cochain
$
\delta = \delta_0 + \ldots + \delta_n \in \A_{\L}.
$
Indeed, $\delta$ lies in the $A_\infty$-subalgebra $\bigoplus_{i = 0}^n (\A_{\L_i}, \mu_i)$ and hence
$\mu^0_\delta() = \mu^0_{\delta_0}() + \ldots + \mu^0_{\delta_n}() = 0$.
We call the resulting Lagrangian brane $\LL_h(\LL_0, \ldots, \LL_n) = (\L, \delta)$ the \emph{composite} Lagrangian brane of $\LL_0, \ldots, \LL_n$ with the shift vector $h$.
The $A_\infty$-operations in the Fukaya category will be computed in $A_\infty$-algebras $(\A_{\LL}, \mu_{\LL})$ for composite Lagrangian branes $\LL$.

\subsubsection{Definition of Fukaya \texorpdfstring{$A_\infty$}{A∞}-pre-category}
\label{subsubs: Fukaya category}
We define the geometric Fukaya $A_\infty$-pre-category $\W = \Wgp(\S; \C)$ of a graded Liouville surface ${\S = (S, Q, \eta, \lambda)}$ with coefficients in the graded system of dg-categories $\C\colon \P^+_2(S, Q) \to \dgcat_\k$ as a $\ZZ$-graded $\k$-linear $A_\infty$-pre-category.

\subsubsection*{Objects}
An object of $\W$ is a Lagrangian brane $\LL = (I, \gamma, P, \gr, \M, \delta)$ in $(\S; \C)$.

\subsubsection*{Transverse sequences}
A sequence $(\LL_0, \ldots, \LL_n) \in (\Ob \W)^{n + 1}$ is transverse if the underlying Lagrangians $L_0, \ldots, L_n$ are transverse.

\subsubsection*{Morphisms}
For a transverse pair of objects $(\LL_0, \LL_1)$ let $Z = Z(\LL_0, \LL_1)$ be the set of intersection points between $\gamma_0$ and $\gamma_1$ and Reeb chords from an end of $\gamma_0$ to an end of $\gamma_1$.
The morphism space is the graded $\k$-module
\begin{equation} \label{eq: morphism space of W}
	\W(\LL_0, \LL_1) = \bigoplus_{z \in Z} \W(z)[-|z|]
\end{equation}
where
\begin{equation} \label{eq: W(z)}
	\W(z) = \C_{\gamma_1(z_-)}(\C(z)\M_0(z_+), \M_1(z_-))
\end{equation}

\subsubsection*{\texorpdfstring{$A_\infty$}{A∞}-operations}
For $n \ge 1$, a transverse sequence $(\LL_0, \ldots, \LL_n)$, and intersection points or Reeb chords
$$
z^0 \in Z(\LL_0, \LL_n), \, z^i \in Z(\LL_{i - 1}, \LL_i), \, i = 1, \ldots, n
$$
the entry
\begin{equation} \label{eq: tilde mu^n_z^0 z^n ... z^1}
	\bmu^{z^0; z^n, \ldots, z^1}\colon \W(z^n) \otimes \ldots \otimes \W(z^1) \to \W(z^0)
\end{equation}
of the structure map $\bmu^n$ is defined as follows.
Pick an admissible shift vector $h \in \H(L_0, \ldots, L_n)$ such that $A_h(z^0), \ldots, A_h(z^n) > 0$ and define
\begin{equation} \label{eq: bmu from mu}
	\bmu^{z^0; z^n, \ldots, z^1} = \mu^{z^0; z^n, \ldots, z^1}_{\LL}
\end{equation}
where $\LL = \LL_h(\LL_0, \ldots, \LL_n)$ is the composite Lagrangian brane of $\LL_0, \ldots, \LL_n$ with the shift vector $h$.

See \cref{p: W is well-defined} for the proof that the choices required to define the structure maps can be made and do not affect the result and that we get a $\ZZ$-graded $\k$-linear $A_\infty$-pre-category where transverse sequences form transverse sets.

\begin{remark} \label{r: relation to Haiden}
	A related approach to defining Fukaya category of surfaces is outlined by Haiden~\cite{H21}.
	The setup there is different in that $\k$ is a field, $D = \varnothing$, the system of dg-coefficients is constant, and the Lagrangians are only allowed to end in a fixed finite set of points on $\pd S$, so Reeb chords on the boundary do not enter the construction.
	The resulting category is called an infinitesimal Fukaya category of a surface.

	The idea of using a version of Chekanov--Eliashberg algebra as an obstruction algebra also appears in Chantraine--Dimitroglou Rizell--Ghiggini--Golovko~\cite{CRGG22}.
\end{remark}

\begin{remark}
	If $S$ is a disc, $Q = \varnothing$, $\C$ is a constant system equal to a field $\k$, and $L$ is a graded Lagrangian whose domain is closed, then the positive augmentation $A_\infty$-category $Aug_+(L)$ of Ng--Rutherford--Shende--Sivek--Zaslow~\cite{NRSSZ20} is defined.
	Modulo sign conventions, the $A_\infty$-operation $\bmu^n$ on a transverse sequence of Lagrangian branes $\LL_0, \ldots, \LL_n$ such that each $\LL_i$ is supported on a perturbation of $L$ can be computed in the augmentation category $Aug_+(L)$.
\end{remark}

\subsubsection{Well-definedness}
\label{subsubs: Well-definedness of geometric Fukaya category}
In this section we demonstrate the following.

\begin{prop} \label{p: W is well-defined}
	The geometric Fukaya $A_\infty$-pre-category $\W = \Wgp(\S; \C)$ defined in \cref{subsubs: Fukaya category} is a small $\ZZ$-graded $\k$-linear $A_\infty$-pre-category where transverse sequences form transverse sets.
\end{prop}

We defined the $A_\infty$-operations $\bmu$ in $\W$ by computing them in certain composite Lagrangian branes.
We first check that this computation does not depend on the choice of a composite Lagrangian brane.

\begin{lemma} \label{l: independence of higher products from shifts}
Let
$\LL_i = (I_i, \gamma_i, P_i, \gr_i, \M_i, \delta_i), \, i = 0, \ldots, n$
be Lagrangian branes, $h, h' \in \RR^{n + 1}$ two admissible shift vectors, and
${\LL = \LL_h = (I, \gamma, P, \gr, \M, \delta), \, \LL' = \LL_{h'} = (I, \gamma, P', \gr, \M, \delta)}$
the corresponding composite Lagrangian branes.
Then
\begin{enumerate}
	\item For any $z \in Z(\LL_i, \LL_j)$ if the actions $A_h(z)$ and $A_{h'}(z)$ are positive, then $A_{L}(z) = A_h(z)$, $A_{L'}(z) = A_{h'}(z)$, and $\A_{\LL}(z) = \A_{\LL'}(z)$.
	\item Let $z^0 \in Z(\LL_0, \LL_n)$ and $z^i \in Z(\LL_{i - 1}, \LL_i), \, i = 1, \ldots, n$ be such that $A_h(z^i) > 0$ and ${A_{h'}(z^i) > 0}$ for all $i = 1, \ldots, n$.
	Suppose that
	$$
		\mu_{\LL}^{z^0; z^n, \ldots, z^1} \colon \A_{\LL}(z^n) \otimes \ldots \otimes \A_{\LL}(z^1) \to \A_{\LL}(z^0)
	$$
	is non-zero.
	Then $A_h(z^0) > 0$, $A_{h'}(z^0) > 0$ and
	$\mu_{\LL}^{z^0; z^n, \ldots, z^1} = \mu_{\LL'}^{z^0; z^n, \ldots, z^1}$.
\end{enumerate}
\end{lemma}
\begin{proof}
	To show (1) first note that $A_h(z) \ne 0$ and $A_{h'}(z) \ne 0$ since $h$ and $h'$ are admissible.
	Positivity of these actions implies that the definitions of $z_+$ and $z_-$ are the same for $\LL$, $\LL'$ (these are given in \cref{subsubs: Double points and Reeb chords}) and for $(L_0, \ldots, L_n)$ (given in \cref{subsubs: Transverse Lagrangians}).
	Namely, we have $z_+ \in I_i$ and $z_- \in I_j$.
	This implies $A_{L}(z) = A_h(z)$, $A_{L'}(z) = A_{h'}(z)$ by comparing definitions in \cref{eq: def of A(x) for double point,eq: def of A(y),eq: def of A_h(z)}, and $\A_{\LL}(z) = \A_{\LL'}(z)$ by \cref{eq: def of A_L(z)}.

	We prove (2).
	Observe that
	$$
	A_h(z^i) = A_L(z^i) > 0 \text{ and } A_{h'}(z^i) = A_{L'}(z^i) > 0 \text{ for }i = 1, \ldots, n
	$$
	by (1) of this lemma and
	$$
	A_{L_i}(z) = A_h(z) = A_{h'}(z) = A_L(z) \text{ for any } z \in Z_{L_i},\ i = 0, \ldots, n
	$$
	by definition.
	Non-vanishing of $\mu_{\LL}^{z^0; z^n, \ldots, z^1}$ implies that there exists
	$$
		\bu \in \Delta^0(L; z^0; z^n_{j_n}, \ldots, z^n_1, z^n, z^{n - 1}_{j_{n - 1}}, \ldots, z^1_1, z^1, z^0_{j_0}, \ldots, z^0_{1})
	$$
	for some $j_0, \ldots, j_n \ge 0$ and $z^i_{j} \in Z_{L_i}$.

	Note that $z^0_+ \in I_0$ and $z^0_- \in I_n$ by any definition.
	For $L$ it follows from the existence of $\bu$, and for $L'$ this follows by considering $\bu$ as an element in $\Delta^0(L'; z^0; \ldots)$ and applying \cref{l: moduli spaces are finite}(2).
	Therefore, we have $A_h(z^0) = A_L(z^0) > 0$ and $A_{h'}(z^0) = A_{L'}(z^0) > 0$.
	Since any such $\bu$ can be viewed as an element of $\Delta^0(L'; z^0; \ldots)$ and vice versa, and as the definition of $\mu_{\LL}$ and $\mu_{\LL'}$ does not depend on the values of the actions, but only on their signs, we conclude
	$\mu_{\LL}^{z^0; z^n, \ldots, z^1} = \mu_{\LL'}^{z^0; z^n, \ldots, z^1}$.
\end{proof}

\begin{proof}[Proof of \cref{p: W is well-defined}]
	That subsets of transverse sets of Lagrangians are transverse is clear from the definition in \cref{subsubs: Transverse Lagrangians} hence transverse sequences form transverse sets in $\W$.
	$\W$ is $\k$-linear by construction.
	It is small since we chose a set of domains of curves in \cref{eq: I m n}.

	That the choice of an admissible shift vector required in the definition of $A_\infty$-operations in $\W$ from \cref{subsubs: Fukaya category} can be made follows from \cref{l: make actions positive by vertical shift}.
	By \cref{l: independence of higher products from shifts} the operation $\bmu^{z^0; z^n, \ldots, z^1}$ does not depend on this choice, and has the correct degree since the $A_\infty$-algebra $\A_{\LL}$ is $\ZZ$-graded.

	To prove that the $A_\infty$-relations hold in $\W$ consider a transverse sequence $(\LL_0, \ldots, \LL_n)$ of Lagrangian branes and a finite set of morphisms
	${\{a^k \in \W(\LL_{i_k}, \LL_{j_k})\}_{k \in K}}$
	such that $i_k < j_k$ for all $k \in K$.
	By \cref{l: make actions positive by vertical shift} there exists a composite Lagrangian brane $\LL = \LL_h(\LL_0, \ldots, \LL_n)$ such that $A_h(z) > 0$ whenever the $\W(z)$-component of one of the morphisms $a^k$ is non-zero.
	We call such composite Lagrangian brane $\LL$ \emph{adapted} to $\{a^k\}_{k \in K}$.
	Note that if $\LL$ is adapted to $\{a^k\}_{k \in K}$, then for any $l \ge 1$, $k_1, \ldots, k_l \in K$ we have
	\begin{equation} \label{eq: tmu_W = tmu_L}
		\mu^l_{\LL}(a^{k_l}, \ldots, a^{k_1}) = \begin{cases}
		\bmu^l_{\W}(a^{k_l}, \ldots, a^{k_1}), &\text{if } \bmu^l_{\W}(a^{k_l}, \ldots, a^{k_1}) \text{ is defined};\\
		0, &\text{otherwise}.
		\end{cases}
	\end{equation}

	That is, the $A_\infty$-relation on a given sequence $\LL_0 \xrightarrow{a^1} \ldots \xrightarrow{a^n} \LL_n$ of morphisms in $\W$ boils down to an $A_\infty$-relation in a composite Lagrangian brane $\LL$ and holds since $\A_{\LL}$ is an uncurved $A_\infty$-algebra.
\end{proof}

\subsubsection{Invariance}
\label{subsubs: Invariance of Fukaya category}
We prove two invariance statements for the geometric Fukaya $A_\infty$-pre-categories.
One is that a graded exact symplectomorphism induces a strict equivalence of the Fukaya categories, see \cref{l: I_bph is a strict equivalence}.
Another is that deforming the line field or the weighted system of dg-categories (\cref{l: invariance of W under isotopies}) or changing the Liouville form (\cref{cor: invariance of W of the choice of Liouville form}) gives a strictly equivalent Fukaya category.

Let $\Phi\colon (S_0, Q_0) \to (S_1, Q_1)$ be a diffeomorphism of marked surfaces (i.e.\ $\Phi$ is orientation preserving and $\Phi(Q_0) = Q_1$).
We have an induced 2-functor $\Phi_*\colon \P^+_2(S_0, Q_0) \to \P^+_2(S_1, Q_1)$.
Therefore, given a weighted system of dg-categories $\C\colon \P^+_2(S_1, Q_1) \to \dgcat_\k$ we obtain a weighted system of dg-categories $\Phi^*\C\colon \P^+_2(S_0, Q_0) \to \dgcat_\k$ by pre-composing $\C$ with $\Phi_*$.

Now let $\S_i = (S_i, Q_i, \eta_i, \lambda_i), \, i = 0, 1$ be two graded Liouville marked surfaces, $\bPhi = (\Phi, \tilde \Phi) \colon \S_0 \to \S_1$ a graded exact symplectomorphism, and $\C_i\colon \P^+_2(\S_i) \to \dgcat_\k$ two weighted systems of dg-categories such that $\Phi^*\C_1 = \C_0$.
We construct a strict equivalence of $A_\infty$-pre-categories
$
\I_{\bPhi}\colon \Wgp(\S_0; \C_0) \to \Wgp(\S_1; \C_1).
$
First note that we have a bijection between Lagrangian pre-branes in $(\S_0; \C_0)$ and $(\S_1; \C_1)$ given by
$$
	\bPhi_*\L = (I, \Phi \circ \gamma, P - g \circ \gamma, \Phi_*(\gr) \cdot (\Phi\gamma)^*(\tilde \Phi)^{-1}, \M) \text{ for } \L = (I, \gamma, P, \gr, \M)
$$
where $g \in C^{\infty}(S_0, \RR)$ is such that $\Phi^*\lambda_1 - \lambda_0 = dg$ and $g(q_0) = 0$ for some fixed basepoint $q_0 \in S_0$.
Secondly, observe that $\Phi$ induces a na\"ive isomorphism $\F_{\bPhi}\colon \A_{\L} \to \A_{\bPhi_*\L}$ for any Lagrangian pre-brane $\L$ in $(\S_0; \C_0)$.
Indeed, for $z \in Z_{\L}$ and $\Phi(z) \in Z_{\bPhi_*\L}$ we have $|z| = |\Phi(z)|$ and
\begin{equation} \label{eq: A_L(z) = A_phL(ph(z))}
	\A_{\L}(z) = (\C_0)_{\gamma(z_-)}(\C_0(z)\M(z_+), \M(z_-)) = (\C_1)_{\Phi\gamma(z_-)}(\C_1(\Phi(z))\M(z_+), \M(z_-)) = \A_{\bPhi_* \L}(\Phi(z)),
\end{equation}
and $\Phi$ gives a bijection between the polygons contributing to $\A_{\L}$ and to $\A_{\bPhi_* \L}$ via \cref{eq: mu from polygons}.

Therefore, we obtain a bijection
$$
\I_{\bPhi}^0 \colon \Ob \Wgp(\S_0; \C_0) \to \Ob \Wgp(\S_1; \C_1), \; (\L, \delta) \mapsto (\bPhi_*\L, \F_{\bPhi}(\delta)).
$$

Note that $\bPhi_*$ preserves transverse sets.
Given two transverse Lagrangian branes $\LL_0, \LL_1$ in $(\S_0; \C_0)$ and $z \in Z(\LL_0, \LL_1)$ we have
$$
\I^1_{\bPhi}(z) = \id \colon \Wgp(\S_0; \C_0)(z) \to \Wgp(\S_1; \C_1)(\Phi(z))
$$
via the same computation as in \cref{eq: A_L(z) = A_phL(ph(z))}.
Since $A_\infty$-operations in the Fukaya categories are computed via $A_\infty$-operations in obstruction algebras, we obtain
\begin{lemma} \label{l: I_bph is a strict equivalence}
$\I_{\bPhi}\colon \Wgp(\S_0; \C_0) \to \Wgp(\S_1; \C_1)$
is a strict equivalence of $A_\infty$-pre-categories.
\end{lemma}

We now pass to the other invariance.
Let $(\S = (S, Q, \eta, \lambda); \C)$ be a weighted Liouville surface.
\begin{lemma} \label{l: invariance of W under isotopies}
	Let $\Phi_s \in \Diff(S, Q), \; s \in [0, 1]$ be a family of diffeomorphisms with $\Phi_0 = \id$.
	Then
	\begin{enumerate}
		\item $\Wgp((S, Q, (\Phi_1)_*\eta, \lambda); \C)$ is strictly equivalent to $\Wgp((S, Q, \eta, \lambda); \C)$.
		\item $\Wgp(\S; \Phi_1^*\C)$ is strictly equivalent to $\Wgp(\S; \C)$.
	\end{enumerate}
\end{lemma}
\begin{proof}
	To show (1) observe that the identity functor works since the line field only enters the constructions to determine the degrees, the degrees are defined as differences in gradings, and a grading with respect to $(\Phi_1)_*\eta$ can be obtained as concatenation of the grading with respect to $\eta$ and $s \mapsto (\Phi_{1 - s})_*\eta$.

	The proof of (2) is similar to the proof of \cref{p: isotopy induces isomorphism of obstruction algebras}.
	Namely, given a Lagrangian pre-brane $\L = (I, \gamma, P, \gr, \M)$ in $(\S; \C)$ we define the Lagrangian pre-brane $\L_1 = (I, \gamma, P, \gr, \M_1)$ in $(\S; \Phi_1^*\C)$ by setting $\M_1(p) = \C(s \mapsto \Phi_s \gamma(p))\M(p)$ for $p \in I$ and for a path $a\colon [0, 1] \to I$ the transition map is
	$$
		\M_1(a) = \C(\Phi_1\gamma a; s \mapsto \Phi_s \gamma a(1), \gamma a, s \mapsto \Phi_{1 - s}\gamma a(0); \sigma_a)(\id_{\M_1(a(1))}, \M(a), \id_{\M(a(0))})
	$$
	where $\sigma_a$ is the 2-morphism in $\P_2^+(\S)$ given by $\Phi_s \gamma a(t)$.

	Then we define a na\"ive isomorphism $(\id, \F)\colon \A_{\L} \to \A_{\L_1}$ by
	\begin{equation} \label{eq: F^z for isotopy of C}
	\begin{aligned}
		\F^z(\ph) = \C\big(\Phi_1 z&; s \mapsto \Phi_s \gamma(z_-)&,& z&,& s \mapsto \Phi_{1 - s} \gamma(z_+)&;& [(s, t) \mapsto \Phi_s z(t)]\big)\\
		&\big(\id_{\M_1(z_-)}&,& \ph&,& \id_{\M(z_+)}&\big)&
	\end{aligned}
	\end{equation}
	for $z \in Z_{\L}$ and
	$\ph \in \A_{\L}(z) = \C_{\gamma(z_-)}(\C(z)\M(z_+), \M(z_-))$.
	The family $\Phi_s$ then gives the isotopies of polygons contributing via \cref{eq: mu from polygons} that certify that $(\id, \F)$ is a na\"ive isomorphism.

	Therefore, we obtain a bijection
	$$
		\I^0 \colon \Ob \Wgp(\S; \C) \to \Ob \Wgp(\S; \Phi_1^*\C), \; (\L, \delta) \mapsto (\L_1, \F(\delta)).
	$$
	Given two Lagrangian branes $\LL_0, \LL_1$ in $(\S; \C)$ and $z \in Z(\LL_0, \LL_1)$ we have
	$
	{\I^1(z) \colon \Wgp(\S; \C)(z) \to \Wgp(\S; \Phi_1^*\C)(z)}
	$
	via \cref{eq: F^z for isotopy of C}.
	Since $A_\infty$-operations in the Fukaya $A_\infty$-pre-categories are computed via $A_\infty$-operations in obstruction algebras, we conclude that $\I\colon \Wgp(\S; \C) \to \Wgp(\S; \Phi_1^*\C)$ is a strict equivalence of $A_\infty$-pre-categories.
\end{proof}

\begin{cor} \label{cor: invariance of W of the choice of Liouville form}
	Let $\lambda'$ be another Liouville form on $S$ compatible with the marking $Q$.
	Then $\Wgp((S, Q, \eta, \lambda); \C)$ is strictly equivalent to $\Wgp((S, Q, \eta, \lambda'); \C)$.
\end{cor}
\begin{proof}
	By \cref{cor: Liouville structures on the same marked surface are symplectomorphic} there exists a family $\Phi_s \in \Diff(S, Q), \; s \in [0, 1]$ of diffeomorphisms such that $\Phi_0 = \id$ and $\Phi_1$ is an exact symplectomorphism from $(S, Q, \lambda)$ to $(S, Q, \lambda')$.
	By \cref{l: I_bph is a strict equivalence} $\Wgp((S, Q, \eta, \lambda); \C)$ is strictly equivalent to $\Wgp((S, Q, (\Phi_1)_*\eta, \lambda'); (\Phi_1^{-1})^*\C)$.
	By \cref{l: invariance of W under isotopies} the latter category is strictly equivalent to $\Wgp((S, Q, \eta, \lambda'); \C)$.
\end{proof}

\subsection{Transverse perturbations}
\label{subs: Transverse perturbations}
Fix a weighted Liouville surface $(\S = (S, Q, \eta, \lambda); \C\colon \P^+_2(\S) \to \dgcat_\k)$.
Denote by $\W = \Wgp(\S; \C)$ the geometric Fukaya $A_\infty$-pre-category.

In \cref{subsubs: Transverse perturbations} we define transverse perturbations of Lagrangian branes and in \cref{subsubs: Morphisms of local systems} morphisms of local systems compatible with the brane structures.
In \cref{subsubs: Transverse double of a Lagrangian,subsubs: Lagrangians containing a transverse double} we describe the intersection points and Reeb chords between a Lagrangian and its transverse perturbation and the polygons bounded by a composite Lagrangian containing them, see \cref{l: polygons bounded by Lagrangian containing a transverse double}.
In \cref{subsubs: Perturbed morphisms} we perturb a morphism $f$ of local systems on a Lagrangian brane to a morphism in $\W$ and compute its differential.
In \cref{subsubs: Composition with perturbed morphisms} we prove that composition with the perturbation of $f$ is given by $f$ and that all higher products involving it vanish.

\subsubsection{Transverse perturbations}
\label{subsubs: Transverse perturbations}
Let $L_0 = (I, \gamma_0, P_0)$ be a Lagrangian in $\S$.
Recall from \cref{l: neighbourhood of L in Lag_I(S)} that a neighbourhood of $L_0$ in $\Lag_I(\S)$ is parametrised by a neighbourhood of $\{0\} \times \{\id\}$ in $C^\infty(I, \RR) \times \Diff(I)$ via $(F, \ph) \mapsto L(F, \ph)$.

\begin{defn} \label{d: transverse perturbation}
A Lagrangian $L_1 = (I, \gamma_1, P_1)$ is called a \emph{transverse perturbation} of $L_0$ if $L_1 = L(F, \id)$ for some $F \in C^\infty(I, \RR)$ and $L_s = L(sF, \id)$ is transverse to $L_0$ for any $s \in (0, 1]$.
\end{defn}

Note that $L(F, \id)$ is transverse to $L_0$ if and only if $F'$ has only simple zeros and is nowhere zero on $\pd I$.
In particular, it is enough for $L_1$ to be transverse to $L_0$ in \cref{d: transverse perturbation}.
Also note that if $L_1$ is a transverse perturbation of $L_0$, then $L_0$ is a transverse perturbation of $L_1$ after choosing an appropriate Weinstein chart for $L_1$.

Let $L_1$ be a transverse perturbation of $L_0$.
Recall from \cref{subsubs: isotopy induces isomorphism of obstruction algebras} that if we have a grading $\gr_0$ or a local system $\M_0$ on $L_0$, then we can transport them from $L_0$ to $L_1$ along the isotopy $L_s$ from \cref{d: transverse perturbation}.
We therefore get the notion of \emph{transverse perturbation} for graded Lagrangians and Lagrangian pre-branes.

Let $\LL_0 = (\L_0, \delta_0)$ be a Lagrangian brane in $(\S; \C)$
and $\L_1$ a transverse perturbation of the Lagrangian pre-brane $\L_0$.
Then, by \cref{p: isotopy induces isomorphism of obstruction algebras} the obstruction algebras $\A_{\L_0}$ and $\A_{\L_1}$ are na\"ively isomorphic.
In particular, we get a bounding cochain $\delta_1 \in \A_{\L_1}$ corresponding to the bounding cochain $\delta_0 \in \A_{\L_0}$.
We call the resulting Lagrangian brane $\LL_1 = (\L_1, \delta_1)$ a \emph{transverse perturbation} of $\LL_0$.

\subsubsection{Morphisms of local systems}
\label{subsubs: Morphisms of local systems}
Let $L_0 = (I, \gamma_0, P_0)$ be a Lagrangian in $\S$ and $\M^0_0$, $\M^1_0$ two local systems of objects over $\gamma_0^* \C$ on $I$.
A collection of morphisms
$$
f(t) \in \C_{\gamma_0(t)}(\M^0_0(t), \M^1_0(t)) \text{ for each } t \in I
$$
is called a \emph{morphism of local systems} ${f\colon \M^0_0 \to \M^1_0}$ if it commutes with transition maps, i.e.\
\begin{equation} \label{eq: f commutes with transition maps}
	\M^1_0(a) \cdot \C(\gamma_0 a)f(a(0)) = f(a(1)) \cdot \M^0_0(a) \text{ for each path } a\colon [0, 1] \to I.
\end{equation}
Note that \cref{eq: f commutes with transition maps} implies that a morphism $f\colon \M^0_0 \to \M^1_0$ is determined by $\{f(t_C)\}_{C \in \pi_0 I}$ for any choice of $t_C \in C$.
If $C \in \pi_0 I$ is a segment, then \cref{eq: f commutes with transition maps} poses no restriction on $f(t_C)$.
If $C \in \pi_0 I$ is a circle, then \cref{eq: f commutes with transition maps} for the component $C$ is equivalent to
$$
\M^1_0(a) \cdot \C(\gamma_0 a)f(t_C) = f(t_C) \cdot \M^0_0(a)
$$
for one path $a\colon [0, 1] \to C$ that parametrises $C$ with $a(0) = a(1) = t_C$.

If $\LL_0^i = (L_0, \gr_0, \M^i_0, \delta^i_0)$ are Lagrangian branes, then a morphism of local systems $f\colon \M^0_0 \to \M^1_0$ is \emph{compatible with the brane structures} if
\begin{equation} \label{eq: compatibility of f with the brane structure}
	f(z_-) \cdot (\delta^0_0)_z = (-1)^{|f|(1 - |z|)}(\delta^1_0)_z \cdot \C(z)f(z_+) \text{ for any } z \in Z_{L_0}
\end{equation}
where $|f| = m$ if $|f(t)| = m$ for all $t \in I$ and \cref{eq: compatibility of f with the brane structure} is understood to hold for each homogeneous component of $f$.

\subsubsection{Transverse double of a Lagrangian}
\label{subsubs: Transverse double of a Lagrangian}

We introduce some notation for intersection points and Reeb chords between a Lagrangian and its transverse perturbation.

\begin{figure}[!htbp]
	\centering
	\begin{subfigure}[t]{0.3\linewidth}
		\includegraphics[width=\linewidth]{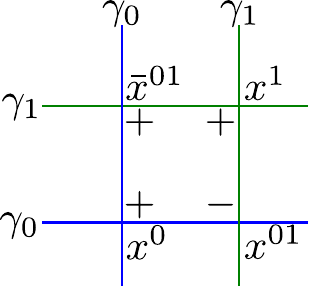}
		\caption{}
		\label{subfig: Neighbourhood of a perturbed self-intersection+}
	\end{subfigure}%
	\hspace{0.1\linewidth}%
	\begin{subfigure}[t]{0.3\linewidth}
		\includegraphics[width=\linewidth]{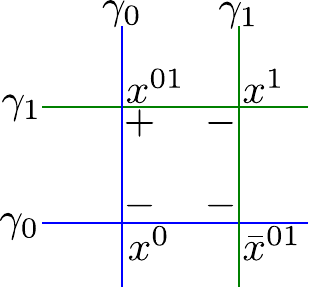}
		\caption{}
		\label{subfig: Neighbourhood of a perturbed self-intersection-}
	\end{subfigure}
	\caption{Neighbourhood of a perturbed self-intersection}
	\label{fig: Neighbourhood of a perturbed self-intersection}
\end{figure}

\begin{figure}[!htbp]
	\centering
	\begin{subfigure}[t]{0.25\linewidth}
		\includegraphics[width=\linewidth]{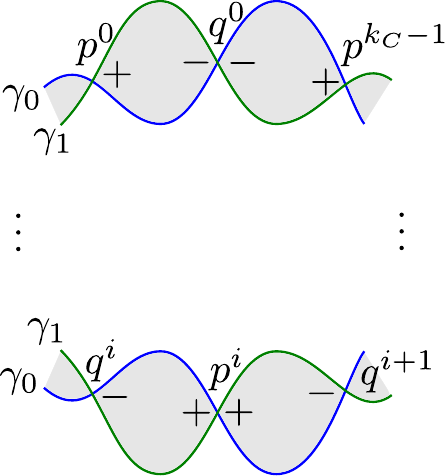}
		\caption{}
		\label{subfig: Pullback of a perturbed circle component of a Lagrangian along bigons}
	\end{subfigure}%
	\hfill%
	\begin{subfigure}[t]{0.68\linewidth}
		\includegraphics[width=\linewidth]{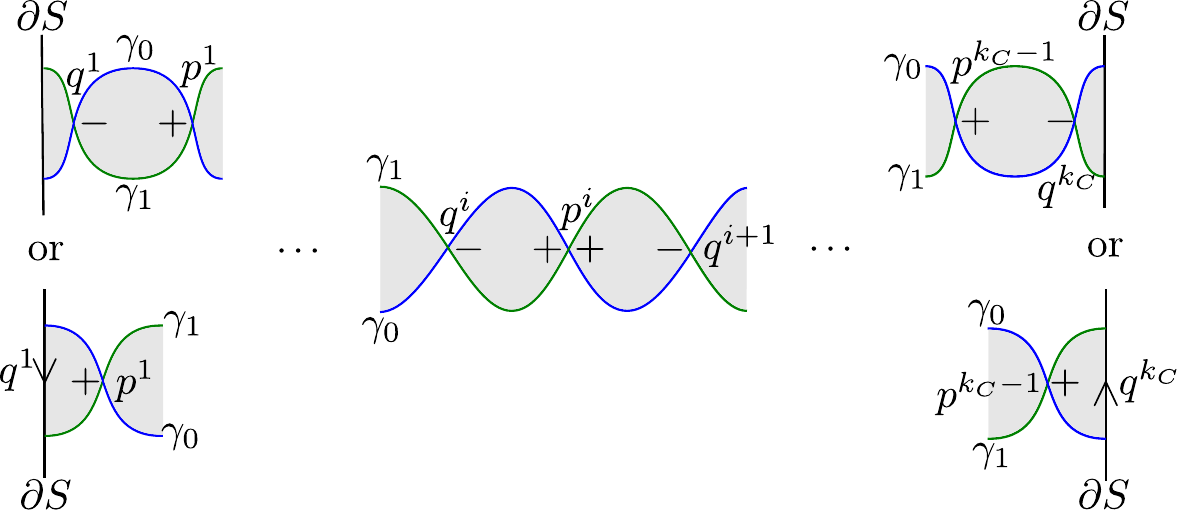}
		\caption{}
		\label{subfig: Pullback of a perturbed segment component of a Lagrangian along bigons}
	\end{subfigure}
	\caption{Pullbacks of perturbed components of the Lagrangian along bigons}
	\label{fig: Pullback of a perturbed component of a Lagrangian along bigons}
\end{figure}

\begin{defn} \label{d: transverse double}
A \emph{transverse double} of a Lagrangian $L_0 = (I, \gamma_0, P_0)$ is a composite Lagrangian $L_h = L_h(L_0, L_1)$ such that
\begin{enumerate}
\item $L_1 = L(F_1, \id)$ is a transverse perturbation of $L_0$ and $h \in \RR^2$ is an admissible shift vector for $(L_0, L(sF_1, \id))$ for any $s \in (0, 1]$.

\item For each double point $x^0$ of $L_0$, there is a double point $x^1 = x^1(x^0)$ of $L_1$, intersection points
$x^{01} = x^{01}(x^0), \bar x^{01} = \bar x^{01}(x^0) \in X(L_0, L_1)$
and a neighbourhood $U_{x^0}$ such that $\{U_{x^0}\}_{x^0 \in X_{L_0}}$ are pairwise disjoint and $L_h$ in $U_{x^0}$ is diffeomorphic to one of the pictures in \cref{fig: Neighbourhood of a perturbed self-intersection}.

\item For each connected component $C$ of $I$ that is a circle, there is an integer $k_C \ge 1$ such that the intersection points in
$Z(\res{L_0}{C}, \res{L_1}{C}) = X(\res{L_0}{C}, \res{L_1}{C})$
that are not of the form $x^{01}$ or $\bar x^{01}$ from (2) split into two sets, $Q_C = \{q^i\}_{i \in \ZZ/k_C\ZZ}$ and $P_C = \{p^i\}_{i \in \ZZ/k_C\ZZ}$.
These points satisfy $A_h(q^i) > 0$ and $A_h(p^i) > 0$ for each $i \in \ZZ/k_C\ZZ$ and
$$
\left|\Delta^{0}(L_h; p^i; q^j)\right| = \begin{cases}
	1,& \text{if } i = j \text{ or } i + 1 = j \text{ in } \ZZ/k_C\ZZ;\\
	0,& \text{otherwise}.
\end{cases}
$$
That is, there is an immersion
$i_C \colon T_C \subset \RR^2 \to S \setminus (\Lambda \cup D)$,
where $T_C$ is the shaded area in \cref{subfig: Pullback of a perturbed circle component of a Lagrangian along bigons}, that extends the immersions given by bigons from $\Delta^{0}(L_h; p^i; q^j)$ (that the numbering goes counter-clockwise is a condition on the indexing of $Q_C$ and $P_C$).

\item For each $t \in \pd I$ there exists a unique Reeb chord $y_t \in Z(L_0, L_1) \cup Z(L_1, L_0)$ that is isotoped to the constant path at $\gamma_1(t)$ by the isotopy from $L_0$ to $L_1$.
$A_h(y_t) > 0$ if $y_t \in Z(L_0, L_1)$ and $A_h(y) < 0$ for any $y \in Z(L_1, L_0)$.

\item For each connected component $C$ of $I$ that is a segment, there is an integer $k_C \ge 1$ such that the intersection points in $Z(\res{L_0}{C}, \res{L_1}{C})$ that are not of the form $x^{01}$ or $\bar x^{01}$ from (2) and the Reeb chords $y_t \in Z(\res{L_0}{C}, \res{L_1}{C})$ for $t \in \pd I$ split into two sets,
$Q_C = \{q^i\}_{i = 1}^{k_C}$ and $P_C = \{p^i\}_{i = 1}^{k_C - 1}$.
These points and chords satisfy $A_h(q^i) > 0$ and $A_h(p^i) > 0$ for each $i$ and
$$
\left|\Delta^{0}(L_h; p^i; q^j)\right| = \begin{cases}
	1,& \text{if } i = j \text{ or } i + 1 = j;\\
	0,& \text{otherwise}.
\end{cases}
$$
That is, there is an immersion
$i_C \colon T_C \subset \RR^2 \to S \setminus (\Lambda \cup D)$,
where $T_C$ is the closed shaded area in \cref{subfig: Pullback of a perturbed segment component of a Lagrangian along bigons}, that extends the immersions given by bigons from $\Delta^{0}(L_h; p^i; q^j)$.

\item The set $X_{L_h}$ of self-intersection points of $L_h$ satisfies
$$
X_{L_h} \sqcup \left(\{y_t \mid t \in \pd I\} \cap Z(L_0, L_1)\right) = \bigsqcup_{x^0 \in X_{L_0}} \{x^0, x^1, x^{01}, \bar x^{01}\} \sqcup \bigsqcup_{C \in \pi_0 I} Q_C \sqcup P_C.
$$
\end{enumerate}
\end{defn}

If $L_1$ is any transverse perturbation of $L_0$, then by \cref{l: make actions positive by vertical shift} we can pick an admissible shift vector $h \in \RR^2$ such that $L_h = L_h(L_0, L_1)$ is a transverse double of $L_0$.
Note that even if $L_1$ is $C^{\infty}$-close to $L_0$, we might not be able to choose the admissible vector $h$ close to $0 \in \RR^2$.
Namely, for $x^0 \in X_{L_0}$ the local model from (2) around this self-intersection implies $A_0(\bar x^{01}) + h_0 - h_1 = A_h(\bar x^{01}) > 0$, but
${A_0(\bar x^{01}) = P_0(\bar x^{01}_+) - P_1(\bar x^{01}_-)}$
is close to $P_0(x^0_-) - P_0(x^0_+) = -A_{L_0}(x^0)$ which is a negative number that depends only on $L_0$, not on the perturbation.

Suppose we have a grading $\gr_0$ on $L_0$.
Recall from \cref{subsubs: isotopy induces isomorphism of obstruction algebras} that then $\gr_0$ transports to a grading $\gr_1$ on $L_1$.
A straightforward computation using \cref{eq: |z^0| from |z^1| + ... + |z^n|} gives
\begin{enumerate}
	\item for each double point $x^0$ of $L_0$ we have $|x^1| = |x^{01}| = |x^0|$ and $|\bar x^{01}| = 1 - |x^0|$;
	\item for each connected component $C$ of $I$ we have $|q| = 0$ for any $q \in Q_C$ and $|p| = 1$ for any $p \in P_C$.
\end{enumerate}

\subsubsection{Lagrangians containing a transverse double}
\label{subsubs: Lagrangians containing a transverse double}
Let $L_0 = (I, \gamma_0, P_0)$ and $L = (J, \gamma, P)$ be two Lagrangians.
\begin{defn} \label{d: contains transverse double nicely}
	We say that $L$ \emph{contains a transverse double of $L_0$ nicely} if $L$ is the composite Lagrangian $L_{(h_0, h_1, h_2)}(L_0, L_1, L_2)$ such that
	\begin{enumerate}
		\item $L_1 = (I, \gamma_1, P_1) = L(F_1, \id)$ is a transverse perturbation of $L_0$ and $L_{(h_0, h_1)}(L_0, L_1)$ is a transverse double of $L_0$.
		\item For each $s \in (0, 1]$ the Lagrangian $L_{(h_0, h_1)}(L_0, L(sF_1, \id))$ is transverse to $L_2$.
	\end{enumerate}
\end{defn}

Note that a transverse double of $L_0$ contains itself nicely.
Observe that (2) in \cref{d: contains transverse double nicely} implies that each connected component $C$ of $I$ the closure (in $\RR^2$) of the pullback of $\gamma$ by the immersion $i_C$ to the interior of the shaded region $T_C$ in \cref{fig: Pullback of a perturbed component of a Lagrangian along bigons} is a collection of disjoint segments with one end on $\gamma_0$ and another end on $\gamma_1$.

We give a criterion for a polygon with convex corners bounded by $L$ to be contained in the shaded region $T_C$ from \cref{fig: Pullback of a perturbed component of a Lagrangian along bigons}.

\begin{figure}[!htbp]
	\centering
	\begin{subfigure}[t]{0.16\linewidth}
		\includegraphics[width=\linewidth]{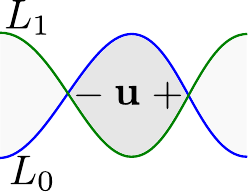}
		\caption{}
		\label{subfig: Polygons contained in T 1}
	\end{subfigure}%
	\hfill%
	\begin{subfigure}[t]{0.15\linewidth}
		\includegraphics[width=\linewidth]{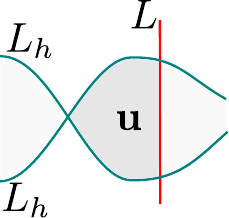}
		\caption{}
		\label{subfig: Polygons contained in T 2}
	\end{subfigure}%
	\hfill%
	\begin{subfigure}[t]{0.17\linewidth}
		\includegraphics[width=\linewidth]{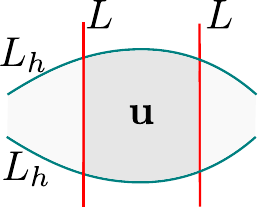}
		\caption{}
		\label{subfig: Polygons contained in T 3}
	\end{subfigure}%
	\hfill%
	\begin{subfigure}[t]{0.09\linewidth}
		\includegraphics[width=\linewidth]{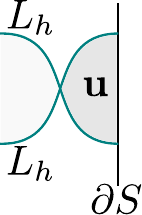}
		\caption{}
		\label{subfig: Polygons contained in T 4}
	\end{subfigure}%
	\hfill%
	\begin{subfigure}[t]{0.11\linewidth}
		\includegraphics[width=\linewidth]{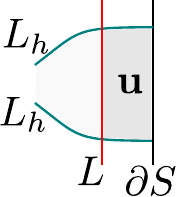}
		\caption{}
		\label{subfig: Polygons contained in T 5}
	\end{subfigure}%
	\hfill%
	\begin{subfigure}[t]{0.14\linewidth}
		\includegraphics[width=\linewidth]{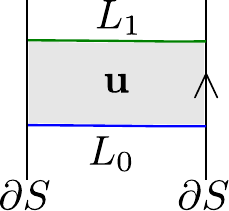}
		\caption{}
		\label{subfig: Polygons contained in T 6}
	\end{subfigure}%
	\caption{Polygons contained in $T_C$}
	\label{fig: Polygons contained in T}
\end{figure}

\begin{lemma} \label{l: polygons bounded by Lagrangian containing a transverse double}
	Let $L = (J, \gamma, P)$ be a Lagrangian containing a transverse double $L_h = L_h(L_0, L_1)$ of ${L_0 = (I, \gamma_0, P_0)}$ nicely, $C$ a connected component of $I$, $\bu = (\Pi, u)$ a polygon with convex corners bounded by $L$, and $v, v' \in \Pi$ two vertices of $\Pi$ such that $v = v'$ or $[v, v']$ is a side of $\Pi$.
	Suppose that $u(v) \in \Im \gamma_0$, $u(v') \in \Im \gamma_1$, and there is an embedding $\ph\colon U \to T_C$ where $U$ is a neighbourhood of $[v, v']$ in $\Pi$ such that $i_C \circ \ph = \res{u}{U}$.
	Then $\ph$ can be extended to an embedding $\ph \colon \Pi \to T_C$ such that $i_C \circ \ph = u$.
	That is, $\bu$ is one of the polygons in \cref{fig: Polygons contained in T}.
\end{lemma}

Note that convex polygons bounded by $L$ that have a negative corner at $Q_C$ or a positive corner at $P_C$ satisfy the assumption of \cref{l: polygons bounded by Lagrangian containing a transverse double}.
This lemma is similar to \cite{A08}*{Lemma 4.4}, our proof is somewhat different.

\begin{proof}
	Note that $i_C$ induces an orientation on $C$: it is the counter-clockwise orientation in \cref{subfig: Pullback of a perturbed circle component of a Lagrangian along bigons} if $C$ is a circle and the orientation from left to right in \cref{subfig: Pullback of a perturbed segment component of a Lagrangian along bigons} if $C$ is a segment.
	By abuse of notation we will not always distinguish between the objects in $S$ and their pullbacks to $\RR^2$ via $i_C$.

	Let the regions $R_1, \ldots, R_r \subset S$ be the closures of connected components of $S \setminus \Im \gamma$.
	Since $L$ contains a transverse double of $L_0$ nicely, the pullback of $\gamma$ to $T = T_C$ consists of the boundary of $T$ and a collection of disjoint segments with one end on $\gamma_0$ and another end on $\gamma_1$.
	These segments divide $T$ into regions each of which is a copy of some $R_i, \, i = 1, \ldots, r$ and is diffeomorphic to one of the darker shaded regions in \cref{fig: Polygons contained in T} (with an additional condition that there are no segments of $\gamma$ inside the region).
	Let $T_1, \ldots, T_l$ be these regions, listed in the direction of the orientation of $C$.
	When $C$ is a circle we think of the indices $1, \ldots, l$ as elements of $\ZZ/l\ZZ$, the indexing as above is not unique in this case, we pick any.

	Pulling $\gamma$ to $\Pi$ along $u$ we see that $\Pi$ is glued from copies of the regions $R_i, \, i = 1, \ldots, r$.
	To prove the lemma it suffices to show that these regions are $T_j, T_{j + 1}, \ldots, T_{j + m}$ or $T_j, T_{j - 1}, \ldots, T_{j - m}$ for some $m \ge 0$ and $\Pi$ is glued from them the same way as $T$.

	Let $R_{i_0}$ be the region that contains a neighbourhood of $[v, v']$ in $\Pi$.
	That $[v, v']$ is contained in one region is clear when $v = v'$ and follows from the existence of $\ph$ from the assumption when $[v, v']$ is a side of $\Pi$.
	The existence of $\ph$ also implies that $R_{i_0}$ has a copy $T_j \subset T$.
	Set $\sigma = 1$ if $T_j$ follows $[u(v), u(v')]$ in the direction of $C$, set $\sigma = -1$ otherwise.

	Let $v^{1}$ (resp.\ $v^{-1}$) be the vertex of $\Pi$ following (resp.\ preceding) the (possibly degenerate) segment $[v, v']$ in the counter-clockwise order on $\pd \Pi$.
	Let $w^1$ (resp.\ $w^{-1}$) be the vertex of $R_{i_0}$ that follows (resp.\ precedes) the segment $[v, v']$ in the counter-clockwise order on $\pd R_{i_0}$.
	If $w^1 \ne v^1$, then define $R_{i_1}$ to be the region that contains the corner at $w^1$ not contained in $R_{i_0}$.
	Observe that this region is the same as the region $T_{j + \sigma}$.
	Therefore, $w^{-1}$ is a vertex of $R_{i_1}$ and the regions $R_{i_0}$ and $R_{i_1}$ share the side $[w^{-1}, w^1]$.
	Let $w^2$ (resp.\ $w^{-2}$) be the vertex of $R_{i_1}$ that follows (resp.\ precedes) $w^1$ (resp.\ $w^{-1}$) in the counter-clockwise order on $\pd R_{i_1}$.
	If $w^2 \ne v^1$, then define $R_{i_2}$ to be the region that contains the corner at $w^2$ not contained in $R_{i_1}$.
	Observe that this region is the same as the region $T_{j + 2\sigma}$.
	Therefore, $w^{-2}$ is a vertex of $R_{i_2}$ and the regions $R_{i_1}$ and $R_{i_2}$ share the side $[w^{-2}, w^2]$.
	Continuing this process while we can we obtain a sequence of regions $R_{i_0}, \ldots, R_{i_m} \subset \Pi$ and the corresponding sequence of regions $T_j, \ldots, T_{j + m\sigma} \subset T$ for some $m \ge 0$.
	Since we cannot continue, $w^{m + 1} = v^1$.

	If $w^{-m - 1} = w^{m + 1}$, then $R_{i_0}, \ldots, R_{i_m}$ cover all $\Pi$ and we are done.
	In this case $\bu$ is a polygon from \cref{subfig: Polygons contained in T 1,subfig: Polygons contained in T 2,subfig: Polygons contained in T 4} with $u(w^{m + 1}) \in P_C \cup Q_C$.
	Suppose $w^{-m - 1} \ne w^{m + 1}$.
	Then $[w^{-m - 1}, w^{m + 1}]$ is a side of $R_{i_m}$ since $[u(w^{-m - 1}), u(w^{m + 1})]$ is a side of $T_{j + m\sigma}$.
	Moreover, as $v^1 = w^{m + 1}$, we have $[w^{-m - 1}, w^{m + 1}] \subset \pd \Pi$.
	Therefore, going from $v$ to $w^{m + 1}$ to $w^{-m - 1}$ to $v$ is a closed curve on $\pd \Pi$ hence $R_{i_0}, \ldots, R_{i_m}$ cover all $\Pi$ and we are done.
	If $[u(w^{-m - 1}), u(w^{m + 1})]$ is a sub-segment of $\gamma$, then $\bu$ is a polygon from \cref{subfig: Polygons contained in T 2,subfig: Polygons contained in T 3,subfig: Polygons contained in T 5} with $[u(w^{-m - 1}), u(w^{m + 1})]$ being a red side.
	If $[u(w^{-m - 1}), u(w^{m + 1})]$ is a sub-segment of $\pd S$, then $\bu$ is a polygon from \cref{subfig: Polygons contained in T 4,subfig: Polygons contained in T 5,subfig: Polygons contained in T 6} with $[u(w^{-m - 1}), u(w^{m + 1})]$ being a side lying on the boundary.

\end{proof}

\subsubsection{Perturbed morphisms}
\label{subsubs: Perturbed morphisms}

Let
$\LL_0^i = (L_0 = (I, \gamma_0, P_0), \gr_0, \M^i_0, \delta^i_0)$
be two Lagrangian branes and $f\colon \M^0_0 \to \M^1_0$ be a morphism of local systems compatible with the brane structures, i.e.\ \cref{eq: compatibility of f with the brane structure} holds.

Let
$\LL^1_1 = (L_1 = (I, \gamma_1, P_1), \gr_1, \M^1_1, \delta^1_1)$
be a transverse perturbation of $\LL_0^1$.
Recall from \cref{subsubs: Transverse double of a Lagrangian} that for each component $C$ of $I$ we have finite sets ${Q_C, P_C \subset Z(\res{L_0}{C}, \res{L_1}{C})}$ with $|z| = 0$ if $z \in Q_C$ and $|z| = 1$ if $z \in P_C$.
For any $q \in Q_C$ we have $\C(q)\M_0^1(q_+) = \M_1^1(q_-)$ hence we have a morphism
$\C(q)f(q_+) \in \W(q) = \C_{\gamma_1(q_-)}(\C(q)\M_0^0(q_+), \M_1^1(q_-))$.
We define the \emph{perturbation of $f$ to $\W(\LL_0^0, \LL_1^1)$} as
\begin{equation} \label{eq: f_LL_0LL_1}
	f_{\LL^0_0, \LL^1_1} = \sum_{C \in \pi_0 I} \sum_{q \in Q_C} \C(q)f(q_+) \in \W(\LL_0^0, \LL_1^1).
\end{equation}
\begin{lemma} \label{l: f_LL_0LL_1 is closed}
	Let Lagrangian branes $\LL_0^0$, $\LL_0^1$, $\LL_1^1$ be as above.
	Suppose that $f\colon \M^0_0 \to \M^1_0$ is a morphism of local systems compatible with the brane structures, homogeneous of degree $m$.
	Then
	$|f_{\LL^0_0, \LL^1_1}| = m$ and
	\begin{equation} \label{eq: bmu(f)}
		\bmu^1(f_{\LL^0_0, \LL^1_1}) = (-1)^m(df)_{\LL^0_0, \LL^1_1}.
	\end{equation}
\end{lemma}

\begin{figure}[!htbp]
	\centering
	\begin{subfigure}[t]{0.22\linewidth}
		\includegraphics[width=\linewidth]{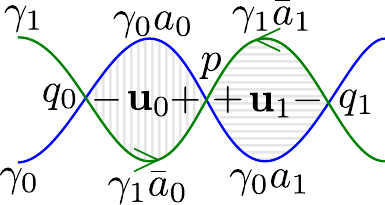}
		\caption{}
		\label{subfig: Polygons contributing to mu^1(f_01) 1}
	\end{subfigure}%
	\hfill%
	\begin{subfigure}[t]{0.3\linewidth}
		\includegraphics[width=\linewidth]{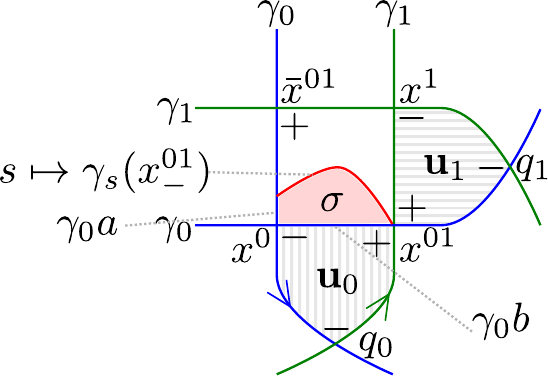}
		\caption{}
		\label{subfig: Polygons contributing to mu^1(f_01) 2}
	\end{subfigure}%
	\hfill%
	\begin{subfigure}[t]{0.13\linewidth}
		\includegraphics[width=\linewidth]{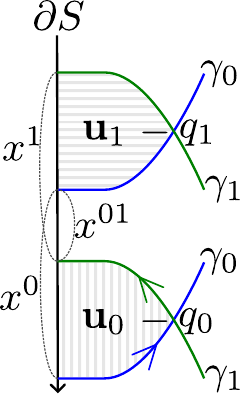}
		\caption{}
		\label{subfig: Polygons contributing to mu^1(f_01) 3}
	\end{subfigure}%
	\hfill%
	\begin{subfigure}[t]{0.3\linewidth}
		\includegraphics[width=\linewidth]{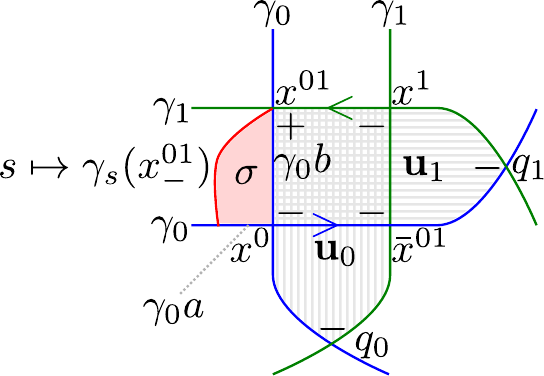}
		\caption{}
		\label{subfig: Polygons contributing to mu^1(f_01) 4}
	\end{subfigure}%
	\caption{Polygons contributing to $\bmu^1(f_{01})$}
	\label{fig: Polygons contributing to mu^1(f_01)}
\end{figure}

\begin{proof}
	Denote $f_{01} = f_{\LL^0_0, \LL^1_1}$.
	We have $|f_{01}| = m$ by construction since $|q| = 0$ for each $q \in Q_C,\, C \in \pi_0 I$ and $|f(t)| = m$ for $t \in I$.
	Let $L_h = L_h(L_0, L_1)$ be a transverse double of $L_0$, as in \cref{subsubs: Transverse double of a Lagrangian}.
	$\LL_h = \LL_h(\LL_0^0, \LL_1^1)$ is adapted to $f_{01}$, so $\bmu^1(f_{01}) = \mu^1_{\LL_h}(f_{01})$ by definition.
	By \cref{l: sufficiently self-wrapped model for Lagrangian,p: isotopy induces isomorphism of obstruction algebras} we can assume that all $q \in Q_C, \, C \in \pi_0 I$ are double points.

	The right-hand side of \cref{eq: bmu(f)} is the part of $\mu^1_{\LL_h}(f_{01})$ coming from \cref{eq: mu^z^0 z^0}.
	\cref{l: polygons bounded by Lagrangian containing a transverse double} applied to $L_h$ implies that the polygons contributing to $\mu^1_{\LL_h}(f_{01})$ via \cref{eq: mu from polygons} come in pairs $(\bu_0, \bu_1)$ as depicted in \cref{fig: Polygons contributing to mu^1(f_01)}.
	It therefore remains to prove the contributions from $\bu_0$ and $\bu_1$ cancel in each case.

	\emph{\cref{subfig: Polygons contributing to mu^1(f_01) 1}.}
	As $|p| = 1$ and $|q_0| = |q_1| = 0$, the contribution of $\bu_0$ and $\bu_1$ to $\mu^1_{\LL_h}(f_{01})$ have opposite sign in \cref{eq: mu from polygons}, so we need to prove $\mu^{\bu_0}(f(q_{0+})) = \mu^{\bu_1}(f(q_{1+}))$.

	For $j = 0, 1$ let $a_j$ be the path in $I$ from $p_+$ to $q_{j-}$ and $\bar a_j$ its reverse, as in \cref{subfig: Polygons contributing to mu^1(f_01) 1}.
	Then
	$$
		\mu^{\bu_j}(f(q_{j+})) =
		\C(p; \gamma_1\bar a_j, q_j, \gamma_0a_j; \bu_j)
		(\M_1^1(\bar a_j), f(q_{j+}), \M_0^0(a_j))
	$$
	and, by \cref{eq: def of M_1},
	$$
		\M_1^1(\bar a_j) = \C(\gamma_1\bar a_j; p, \gamma_0\bar a_j, q_j; \bu_j)(\id, \M_0^1(\bar a_j), \id).
	$$
	Since the orientations of $\bu_j$ in the above two expressions are opposite, \cref{l: properties of composition along a 2-path} and the fact that $f$ commutes with the transition maps of the local systems imply
	$$
		\mu^{\bu_j}(f(q_{j+})) = \M_0^1(\bar a_j) \cdot \C(\gamma_0 \bar a_j)f(q_{j+}) \cdot \C(\gamma_0 \bar a_j)\M_0^0(a_j) = f(p_+) \cdot \M_0^0(\bar a_j) \cdot \C(\gamma_0 \bar a_j)\M_0^0(a_j) = f(p_+)
	$$
	for $j = 0, 1$.

	\emph{\cref{subfig: Polygons contributing to mu^1(f_01) 2,subfig: Polygons contributing to mu^1(f_01) 3}.} We have $\bu_0 \in \Delta^0(L_h; x^{01}; q_0, x^0)$ and $\bu_1 \in \Delta^0(L_h; x^{01}; x^1, q_1)$ for some $x^0 \in X_{L_0}$, $q_0 \in Q_{C_0}$, $q_1 \in Q_{C_1}$, and $C_0, C_1 \in \pi_0 I$.
	By \cref{l: sufficiently self-wrapped model for Lagrangian,p: isotopy induces isomorphism of obstruction algebras} we can reduce \cref{subfig: Polygons contributing to mu^1(f_01) 3} to \cref{subfig: Polygons contributing to mu^1(f_01) 2}, so we only consider the latter case.

	Note that by the degree computations in \cref{subsubs: Transverse double of a Lagrangian} the components $(\delta^0_0)_{x^0}$ and $(\delta^1_1)_{x^1}$ both have degree $1 - |x^0| = 1 - |x^1|$.
	Comparing the signs in \cref{eq: mu from polygons} for $\bu_0$ and $\bu_1$, we see that we need to prove
	\begin{equation} \label{eq: mu^u_1 = mu^u_0 fig b}
		\mu^{\bu_0}(f(q_{0+}), (\delta_0^0)_{x^0}) = (-1)^{m(1 - |x^0|)}\mu^{\bu_1}((\delta_1^1)_{x^1}, f(q_{1+})).
	\end{equation}
	Expanding the definitions for $\mu^{\bu_i}$ from \cref{eq: def of mu^u}, $\M_1^1$ from \cref{eq: def of M_1} and $(\delta_1^1)_{x^1} = \F^{x^0}((\delta^1_0)_{x^0})$ from \cref{eq: F^z^0(ph)}, using that $f$ commutes with the transition maps of the local systems, the compatibility \cref{eq: compatibility of f with the brane structure} of $f$ with the brane structure at $z = x^0$, and \cref{l: properties of composition along a 2-path}, we obtain that both sides of \cref{eq: mu^u_1 = mu^u_0 fig b} are equal to
	\begin{equation} \label{eq: mu^u from sigma for f}
		\C(x^{01}; s \mapsto \gamma_{s}(x^{01}_-), \gamma_0a, x^0, \gamma_0b; \sigma)(\C(s \mapsto \gamma_{s}(x^{01}_-))f(x^{01}_-), \M_0^0(a), (\delta_0^0)_{x^0}, \M_0^0(b))
	\end{equation}
	in the notation of \cref{subfig: Polygons contributing to mu^1(f_01) 2} where $b$ is the path from $x^{01}_+$ to $x^0_+$ in $I$ and $a$ is the path from $x^0_-$ to $x^{01}_-$ in $I$.

	\emph{\cref{subfig: Polygons contributing to mu^1(f_01) 4}.} We have $\bu_0 \in \Delta^0(L_h; x^{01}; x^1, q_0)$ and $\bu_1 \in \Delta^0(L_h; x^{01}; q_1, x^0)$ for some $x^0 \in X_{L_0}$, $q_0 \in Q_{C_0}$, $q_1 \in Q_{C_1}$, and $C_0, C_1 \in \pi_0 I$.
	The signs in \cref{eq: mu from polygons} for $\bu_0$ and $\bu_1$ imply that we need to prove
	$$
	\mu^{\bu_0}((\delta_1^1)_{x^1}, f(q_{0+})) = (-1)^{m(1 - |x^0|)}\mu^{\bu_1}(f(q_{1+}), (\delta_0^0)_{x^0}).
	$$
	Arguing similarly to the case of \cref{subfig: Polygons contributing to mu^1(f_01) 2}, we conclude that both sides are equal to $(-1)^{m(1 - |x^0|)}$ times the expression in \cref{eq: mu^u from sigma for f} in the notation of \cref{subfig: Polygons contributing to mu^1(f_01) 4}.
\end{proof}

\subsubsection{Composition with perturbed morphisms}
\label{subsubs: Composition with perturbed morphisms}
Let $L_0$ and $T$ be Lagrangians in $\S$, $L_1$ a transverse perturbation of $L_0$ and $L = L_{(h_0, h_1, h_2)}(L_0, L_1, T)$ a composite Lagrangian containing the transverse double $L_{(h_0, h_1)}(L_0, L_1)$ of $L_0$ nicely (see \cref{d: contains transverse double nicely}).
Let $\LL^i_0 = (L_0, \gr_0, \M_0^i, \delta_0^i)$, $i = 0, 1$ and $\TT = (T, \gr_T, \M_T, \delta_T)$ be Lagrangian branes in $(\S; \C)$, $\LL_1^1$ the transverse perturbation of $\LL_0^1$ supported on $L_1$, and $\LL = \LL_{h_0, h_1, h_2}(\LL_0^0, \LL_1^1, \TT)$ the corresponding composite Lagrangian brane.
Finally, let $f\colon \M_0^0 \to \M_0^1$ be a morphism of local systems compatible with brane structures, homogeneous of degree $m$.

Note that $\LL$ is adapted to the perturbation $f_{01} = f_{\LL_0^0, \LL_1^1} \in \W(\LL_0^0, \LL_1^1)$.
We describe higher compositions $\bmu$ that can be computed as $\mu_{\LL}$ and have $f_{01}$ as one of the arguments.
The following is an immediate corollary of \cref{l: polygons bounded by Lagrangian containing a transverse double}.
\begin{lemma} \label{l: no higher products involving perturbed morphisms}
	Let $\LL$ and $f_{01}$ be as above.
	Then for any $n \ge 3$
	$$
		\mu^n_{\LL}\left(\ldots, f_{01}, \ldots\right) = 0.
	$$
\end{lemma}

Recall that $\mu^1_{\LL}(f_{01})$ is computed in \cref{l: f_LL_0LL_1 is closed}.
It remains to describe $\mu^2_{\LL}$ when one of the arguments is $f_{01}$.

Given $z^0 \in Z(\TT, \LL_0^0)$ (resp.\ $z^1 \in Z(\LL_1^1, \TT)$) we can pick $h_2$ above in such a way that $A_{(h_0, h_2)}(z^0) > 0$ (resp.\ $A_{(h_1, h_2)}(z^1) > 0$).
Tracking $z^0$ (resp.\ $z^1$) along the isotopy between $L_0$ and $L_1$ we obtain $z^1 \in Z(\TT, \LL_1^1)$ (resp.\ $z^0 \in Z(\LL_0^0, \TT)$).
The isotopy from $L_0$ to $L_1$ induces an isomorphism
$$
\F^{z^0}\colon \C_{\gamma_{0}(z^0_-)}(\C(z^0)\M_T(z^0_+), \M_0^1(z^0_-))\to \C_{\gamma_{1}(z^1_-)}(\C(z^1)\M_T(z^1_+), \M_1^1(z^1_-)) = \W(z^1)
$$
$$
\left(\text{resp.\ } \F^{z^1}\colon \W(z^1) = \C_{\gamma_T(z^1_-)}(\C(z^1)\M_1^1(z^1_+), \M_T(z^1_-)) \to \C_{\gamma_T(z^0_-)}(\C(z^0)\M_0^1(z^0_+), \M_T(z^0_-))\right)
$$
as in \cref{eq: F^z^0(ph)}.
We set
$$
f_*(\ph) = \F^{z^0}\left(f(z^0_-) \cdot \ph\right) \text{ for } \ph \in \W(z^0)$$
$$
\left(\text{resp.\ } f^*(\ph) = (-1)^{m|z^1|}\F^{z^1}(\ph) \cdot \C(z^0)f(z^0_+) \text{ for } \ph \in \W(z^1)\right).
$$
Combining these maps for all $z^0 \in Z(\TT, \LL_0^0)$ (resp.\ $z^1 \in Z(\LL_1^1, \TT)$) we obtain morphisms of graded modules of degree $m$
$$
f_* = f^{\LL_0^0, \LL_1^1}_{*}\colon \W(\TT, \LL_0^0) \to \W(\TT, \LL_1^1)\ (\text{resp.\ } f^* = f^{*}_{\LL_0^0, \LL_1^1}\colon \W(\LL_1^1, \TT) \to \W(\LL_0^0, \TT)).
$$
If $f$ is the identity, these are the isomorphisms induced by the isotopy from $L_0$ to $L_1$.

\begin{lemma} \label{l: composition with perturbed morphisms}
	Let $\LL_0^0$, $\LL_1^1$, and $\TT$ be Lagrangian branes as above.
	Then
	\begin{enumerate}
		\item
		$f_{01} \circ - = f_{*} \colon \W(\TT, \LL_0^0) \to \W(\TT, \LL_1^1)$;
		\item
		$- \circ f_{01} = f^{*} \colon \W(\LL_1^1, \TT) \to \W(\LL_0^0, \TT)$.
	\end{enumerate}
\end{lemma}
\begin{proof}
	We prove (1).
	Fix $z^0 \in Z(\TT, \LL_0^0)$ and $\ph \in \W(z^0)$, and pick $h_2$ such that $A_{(h_0, h_2)}(z^0) > 0$ and $A_{(h_1, h_2)}(z^s) > 0$ where $z^s$ is the trace of $z^0$ along the isotopy from $L_0$ to $L_1$.
	Let $C$ be the connected component of $I$ containing $z^0_-$.
	By \cref{l: polygons bounded by Lagrangian containing a transverse double} the only polygons contributing to $\mu^2_{\LL}(f_{01}, \ph)$ are the polygons $\bu \in \Delta^0(L; z^1; q, z^0)$ with $q \in Q_C$.
	By \cref{l: sufficiently self-wrapped model for Lagrangian,p: isotopy induces isomorphism of obstruction algebras} we can assume that $q$ is a double point, so $\bu$ is a polygon as in \cref{subfig: Polygons contained in T 2,subfig: Polygons contained in T 4}.

	We compute the sign of the contribution of $\bu$ in \cref{eq: mu from polygons}.
	The sign in front of the sum is $(-1)^{|\ph|}$.
	We have $|z^1| = |z^0|$ and the sides of $\bu$ lying on $L_0$ and on its parallel copy $L_1$ are traversed in opposite directions, so the product of the signs $\eps^{\bu}_i$ is $(-1)^{|z^0|}$.
	Multiplying all together, we obtain
	$$
		f_{01} \circ \ph = (-1)^{|\ph| + |z^0|}\mu^2_{\LL}(f_{01}, \ph) = \mu^{\bu}(\C(q)f(q_+), \ph).
	$$
	Expanding the definition of $\mu^{\bu}$ from \cref{eq: def of mu^u}, using that $f$ commutes with the transition maps of the local systems to move $f$ along the side of $\bu$ lying on $L_0$ from $q_+$ to $z^0_-$, and applying \cref{l: properties of composition along a 2-path} to the strip of $\bu$ between $L_0$ and $L_1$, we conclude
	$$
		\mu^{\bu}(\C(q)f(q_+), \ph) = \F^{z^0}\left(f(z^0_-) \cdot \ph\right) = f_*(\ph).
	$$

	The proof of (2) is analogous.
	For $z^1 \in Z(\LL_1^1, \TT)$ and $\ph \in \W(z^1)$ the only contributing polygons are $\bu \in \Delta^0(L; z^0; z^1, q)$ with $q \in Q_C$.
	The sign in front of the sum in \cref{eq: mu from polygons} is $(-1)^{(|z^1| - 1)m}$, and the product of the signs $\eps^{\bu}_i$ is $1$ since $\eps^{\bu}_0 = \eps^{\bu}_2$ and $|z^0| = |z^1|$.
	Hence
	$$
		\ph \circ f_{01} = (-1)^m\mu^2_{\LL}(\ph, f_{01}) = (-1)^{m|z^1|}\mu^{\bu}(\ph, \C(q)f(q_+)) = (-1)^{m|z^1|}\F^{z^1}(\ph) \cdot \C(z^0)f(z^0_+) = f^*(\ph).
	$$
\end{proof}

\begin{figure}[!htbp]
\centering
\begin{subfigure}[t]{0.269\linewidth}
	\includegraphics[width=\linewidth]{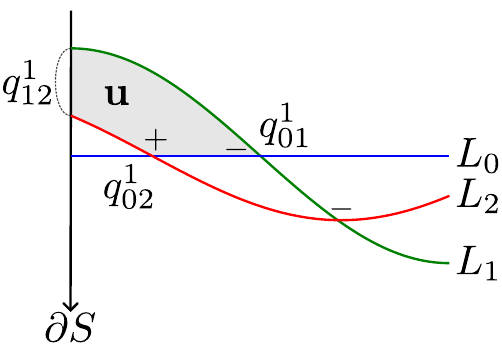}
	\caption{}
	\label{subfig: Pullback of L_2 by i_C a}
\end{subfigure}%
\hfill%
\begin{subfigure}[t]{0.385\linewidth}
	\includegraphics[width=\linewidth]{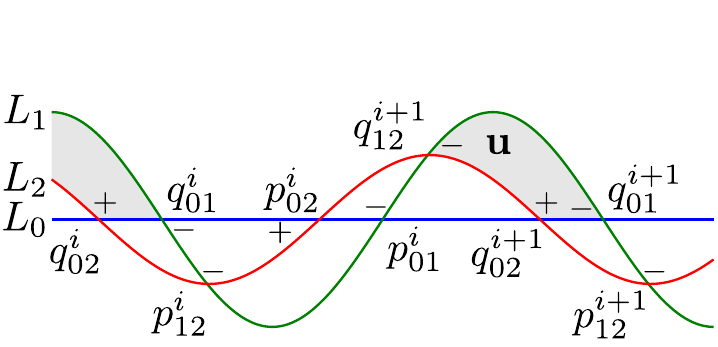}
	\caption{}
	\label{subfig: Pullback of L_2 by i_C b}
\end{subfigure}%
\hfill%
\begin{subfigure}[t]{0.296\linewidth}
	\includegraphics[width=\linewidth]{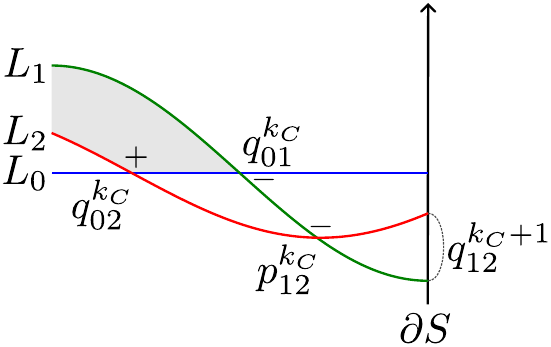}
	\caption{}
	\label{subfig: Pullback of L_2 by i_C c}
\end{subfigure}%

\vspace{1em}
\begin{subfigure}[t]{0.253\linewidth}
	\includegraphics[width=\linewidth]{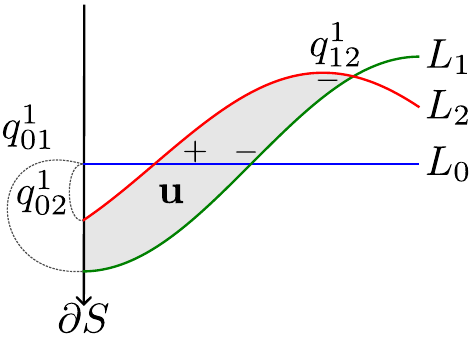}
	\caption{}
	\label{subfig: Pullback of L_2 by i_C d}
\end{subfigure}%
\hfill%
\begin{subfigure}[t]{0.250\linewidth}
	\includegraphics[width=\linewidth]{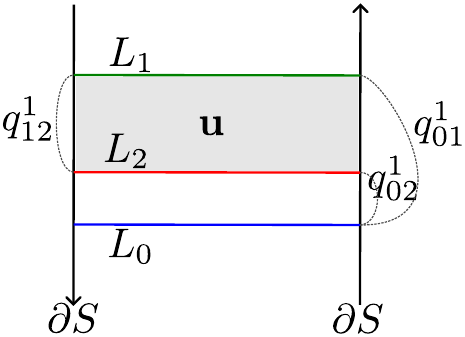}
	\caption{}
	\label{subfig: Pullback of L_2 by i_C e}
\end{subfigure}%
\hfill%
\begin{subfigure}[t]{0.256\linewidth}
	\includegraphics[width=\linewidth]{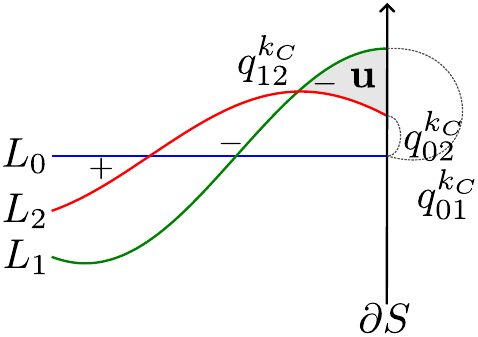}
	\caption{}
	\label{subfig: Pullback of L_2 by i_C f}
\end{subfigure}%
\caption{$L_2$ such that $L$ contains $L_{(h_0, h_2)}(L_0, L_2)$ and $L_{(h_1, h_2)}(L_1, L_2)$ nicely.}
\label{fig: Pullback of L_2 by i_C}
\end{figure}

We finish this section with a useful corollary of \cref{l: composition with perturbed morphisms}.
Let $L_0$ be a Lagrangian and $L_1 = L(F_1, \id)$ be its transverse perturbation.
Then we can find a Lagrangian $L_2 = L(F_2, \id)$ such that the composite Lagrangian $L = L_{(h_0, h_1, h_2)}(L_0, L_1, L_2)$ contains the transverse doubles $L_{(h_0, h_2)}(L_0, L_2)$ of $L_0$ and $L_{(h_1, h_2)}(L_1, L_2)$ of $L_1$ nicely as in \cref{fig: Pullback of L_2 by i_C}.
Note that $L_2$ can be chosen arbitrarily close in $\Lag_I(\S)$ to either $L_0$ or $L_1$.

Now let $\LL_0^i = (L_0, \gr_0, \M_0^i, \delta_0^i)$, $i = 0, 1, 2$ be three Lagrangian branes and $f\colon \M_0^0 \to \M_0^1$, $g\colon \M_0^1 \to \M_0^2$ two morphisms of local systems compatible with the brane structures.
Let $\LL_1^1$ be the transverse perturbation of $\LL_0^1$ supported on $L_1$ and $\LL_2^2$ the transverse perturbation of $\LL_0^2$ supported on $L_2$.

\cref{l: composition with perturbed morphisms} together with the fact that the local systems $\M_1$ and $\M_2$ are both obtained from $\M_0$ by transporting along the cotangent fibres gives the following.

\begin{cor} \label{cor: perturbed f circ g}
	Let $\LL_0^0$, $\LL_1^1$, $\LL_2^2$, $f$, and $g$ be as above.
	Then
	$$
		g_{\LL_1^1, \LL_2^2} \circ f_{\LL_0^0, \LL_1^1} = (g \cdot f)_{\LL_0^0, \LL_2^2}.
	$$
\end{cor}

\subsection{Homotopy equivalences}
\label{subs: Homotopy equivalences in geometric Fukaya category}
Keep the notation $\S$, $\C$, and $\W = \Wgp(\S; \C)$ from before.
In this subsection we prove the following proposition and then show in \cref{l: homotopy equivalence of local system induces homotopy equivalence of branes} that homotopy equivalences of local systems induce homotopy equivalences of Lagrangian branes.

\begin{prop} \label{p: W(S) has enough homotopy equivalences}
$\W$ has enough homotopy equivalences (see \cref{d: W has enough homotopy equivalences}).

Moreover, if $\LL_0 = (L_0, \gr_0, \M_0, \delta_0)$ is a Lagrangian brane in $(\S; \C)$ and $L_1$ a Lagrangian in $\S$ isotopic to $L_0$, then there exists a Lagrangian brane $\LL_1$ supported on $L_1$ such that $\LL_0$ and $\LL_1$ are homotopy equivalent as objects of $\W$ (see \cref{d: homotopy inverse}).
\end{prop}

Let $\LL_0 = (L_0, \gr_0, \M_0, \delta_0)$ be a Lagrangian brane in $(\S; \C)$ and $\LL_1$ be a transverse perturbation of $\LL_0$ (see \cref{subsubs: Transverse perturbations}).
The morphism of local systems $\id\colon \M_0 \to \M_0$ is clearly compatible with the brane structure hence we get a morphism
$\id_{\LL_0, \LL_1} \in \W(\LL_0, \LL_1)$ by \cref{eq: f_LL_0LL_1}.
Recall that we can view the reverse of the isotopy from \cref{d: transverse perturbation} as defining $\LL_0$ as a transverse perturbation of $\LL_1$, hence we get $\id_{\LL_1, \LL_0} \in \W(\LL_1, \LL_0)$.

\begin{lemma} \label{l: e is a homotopy equivalence}
In the notation above $\id_{\LL_0, \LL_1}$ and $\id_{\LL_1, \LL_0}$ are homotopy inverse morphisms in $\W$ (see \cref{d: homotopy inverse}).
\end{lemma}
\begin{proof}
	Denote $e_{10} = \id_{\LL_1, \LL_0}$ and $e_{01} = \id_{\LL_0, \LL_1}$.
	By \cref{l: f_LL_0LL_1 is closed} $e_{10}$ and $e_{01}$ are closed morphisms of degree $0$.
	It is enough to prove that for any transverse set $\{\TT, \LL_0, \LL_1\}$ in $\W$ the map $e_{10} \circ (e_{01} \circ -)$ is homotopic to $\id_{\W(\TT, \LL_0)}$ and the map $(- \circ e_{10}) \circ e_{01}$ is homotopic to $\id_{\W(\LL_0, \TT)}$ since we can swap $\LL_0$ and $\LL_1$.
	We prove that $e_{10} \circ (e_{01} \circ -)$ is homotopic to $\id_{\W(\TT, \LL_0)}$, the other homotopy is similar.

	Fix the transverse set $\{\TT, \LL_0, \LL_1\}$.
	Applying the discussion before \cref{cor: perturbed f circ g} to $L_1$ as a transverse perturbation of $L_0$ and to $L_0$ as a transverse perturbation of $L_1$, we can choose a Lagrangian $L_2$ arbitrarily close to $L_0$ which is a transverse perturbation of both $L_0$ and $L_1$ such that \cref{cor: perturbed f circ g} applies to both triples $(\LL_0, \LL_1, \LL_2)$ and $(\LL_1, \LL_0, \LL_2)$, where $\LL_2$ is the transverse perturbation of $\LL_0$ (equivalently, of $\LL_1$) supported on $L_2$.
	Taking $L_2$ close enough to $L_0$, by \cref{l: composition with perturbed morphisms,cor: perturbed f circ g} we can moreover assume that $\{\TT, \LL_0, \LL_1, \LL_2\}$ is a transverse set and
	$$
	\begin{aligned}
		&e_{02} \circ - = \id_* \colon \W(\TT, \LL_0) \to \W(\TT, \LL_2) \text{ is an isomorphism};\\
		&e_{12} \circ e_{01} = e_{02} \text{ and } e_{02} \circ e_{10} = e_{12}
	\end{aligned}
	$$
	where $e_{12} = \id_{\LL_1, \LL_2}$, and $e_{02} = \id_{\LL_0, \LL_2}$ are defined by \cref{eq: f_LL_0LL_1}.
	The $A_\infty$-relations imply that
	\begin{enumerate}
		\item $e_{02} \circ - =  (e_{12} \circ e_{01}) \circ -$ is homotopic to $e_{12} \circ (e_{01} \circ -)$ as maps $\W(\TT, \LL_0) \to \W(\TT, \LL_2)$;
		\item $e_{12} \circ - = (e_{02} \circ e_{10}) \circ -$ is homotopic to $e_{02} \circ (e_{10} \circ -)$ as maps $\W(\TT, \LL_1) \to \W(\TT, \LL_2)$.
	\end{enumerate}

	Combining the two homotopies above we obtain that $e_{02} \circ -$ is homotopic to $e_{02} \circ (e_{10} \circ (e_{01} \circ - ))$ as maps $\W(\TT, \LL_0) \to \W(\TT, \LL_2)$.
	As $e_{02} \circ -\colon \W(\TT, \LL_0) \to \W(\TT, \LL_2)$ is an isomorphism, we conclude
	that
	$e_{10} \circ (e_{01} \circ - )\colon \W(\TT, \LL_0) \to \W(\TT, \LL_0)$
	is homotopic to $\id$ as desired.
\end{proof}

We are now ready to prove \cref{p: W(S) has enough homotopy equivalences}.

\begin{proof}[Proof of \cref{p: W(S) has enough homotopy equivalences}]
	We first check that the Fukaya category $\W$ has enough homotopy equivalences, see \cref{d: W has enough homotopy equivalences}.
	Let $\LL_0 = (L_0 = (I, \gamma_0, P_0), \gr_0, \M_0, \delta_0)$ be a Lagrangian brane and $\{\{\LL_j^k\}_{k = 0}^{l_j}\}_{j = 1}^m$ be a set of transverse sets of Lagrangian branes.
	For each $j = 1, \ldots, m$ let $\LL_j$ be a composite Lagrangian brane of $\LL_j^0, \ldots, \LL_j^{l_j}$ supported on the Lagrangian $L_j$.

	By \cref{l: e is a homotopy equivalence} it suffices to find a function $F \in C^\infty(I, \RR)$ such that $L(F, \id)$ is transverse to $L_0$ and each $L_j$.
	This condition defines an open dense subset in a neighbourhood of $0$ in $C^\infty(I, \RR)$ hence the desired $F$ exists.

	Now let
	$\LL_0 = (L_0 = (I, \gamma_0, P_0), \gr_0, \M_0, \delta_0)$
	be a Lagrangian brane and $L_1 = (I, \gamma_1, P_1)$ be a Lagrangian isotopic to $L_0$.
	By \cref{l: neighbourhood of L in Lag_I(S)} we can find a sequence of Lagrangians $L_0 = L_0', L_1', \ldots, L_n' = L_1$ such that each $L_{i + 1}'$ is a transverse perturbation of $L_{i}'$ up to a reparametrisation of the domain.
	It remains to prove that reparametrisation of the domain preserves homotopy equivalences.

	Namely, let
	$\LL_i = (L_i = (I_i, \gamma_i, P_i), \gr_i, \M_i, \delta_i)$,
	$i = 0, 1$ be two transverse Lagrangian branes, $q\in \W(\LL_{0}, \LL_1)$ and $p \in \W(\LL_1, \LL_0)$ homotopy inverse morphisms, and $\ph \in \Diff(I_1)$ an orientation-preserving diffeomorphism.
	Then we get a Lagrangian pre-brane
	$\L_1' = (I_1, \gamma_1 \circ \ph, P_1 \circ \ph, \ph^*\gr_1, \ph^*\M_1)$
	and the obstruction algebra does not change, $\A_{\L_1'} = \A_{\L_1}$, so we have a Lagrangian brane $\LL_1' = (\L_1', \delta_1)$.
	For any Lagrangian brane $\TT$ the set $\{\TT, \LL_1\}$ is transverse if and only if $\{\TT, \LL_1'\}$ is, and in this case $Z(\TT, \LL_1) = Z(\TT, \LL_1')$ with the same modules $\W(z)$, so $\W(\TT, \LL_1) = \W(\TT, \LL_1')$, and similarly $\W(\LL_1, \TT) = \W(\LL_1', \TT)$.
	Since reparametrising the domains of composite Lagrangian branes by $\ph$ does not change their obstruction algebras, these identifications commute with all $\bmu$.
	Therefore the morphisms $q' \in \W(\LL_0, \LL_1')$ and $p' \in \W(\LL_1', \LL_0)$ corresponding to $q$ and $p$ are homotopy inverse.
\end{proof}

We prove that homotopy equivalences of local systems induce homotopy equivalences of Lagrangian branes.

\begin{lemma} \label{l: homotopy equivalence of local system induces homotopy equivalence of branes}
	Let $\LL_0^i = (L_0, \gr_0, \M_0^i, \delta_0^i)$ be two Lagrangian branes and $f\colon \M_0^0 \to \M_0^1$, $g\colon \M_0^1 \to \M_0^0$ closed degree zero morphisms of local systems compatible with brane structures.
	Suppose that there exist morphisms of local systems compatible with brane structures $H_i\colon \M_0^i \to \M_0^i$ such that $f \cdot g - \id = dH_1$ and $g \cdot f - \id = dH_0$.
	Then $\LL_0^0$ and $\LL_0^1$ are homotopy equivalent in $\W$.
\end{lemma}
\begin{proof}
	We have already proved \cref{p: W(S) has enough homotopy equivalences}, so $\W$ has enough homotopy equivalences and we can use the machinery of \cref{subs: Homotopy equivalence in A_infty-pre-categories}.
	By \cref{cor: perturbed f circ g} there exist transverse perturbations $\LL_1^1$ of $\LL_0^1$, $\LL_2^0$ of $\LL_0^0$, and $\LL_3^1$ of $\LL_1^1$ such that $\{\LL_0^0, \LL_1^1, \LL_2^0, \LL_3^1\}$ is a transverse set and
	$$
		g_{\LL_1^1, \LL_2^0} \circ f_{\LL_0^0, \LL_1^1} = (g \cdot f)_{\LL_0^0, \LL_2^0} \text{ and } f_{\LL_2^0, \LL_3^1} \circ g_{\LL_1^1, \LL_2^0} = (f \cdot g)_{\LL_1^1, \LL_3^1}.
	$$
	By assumption and \cref{l: f_LL_0LL_1 is closed}
	$$
		(g \cdot f)_{\LL_0^0, \LL_2^0} - \id_{\LL_0^0, \LL_2^0} \text{ and } (f \cdot g)_{\LL_1^1, \LL_3^1} - \id_{\LL_1^1, \LL_3^1} \text{ are exact}.
	$$
	By \cref{l: exact morphisms and contractible objects in A_infty-pre-category} composition with $(g \cdot f)_{\LL_0^0, \LL_2^0}$ on either side is homotopic to composition with $\id_{\LL_0^0, \LL_2^0}$, and similarly for $(f \cdot g)_{\LL_1^1, \LL_3^1}$ and $\id_{\LL_1^1, \LL_3^1}$.
	Hence by \cref{l: e is a homotopy equivalence,l: homotopy equivalence can be tested pointwise}
	$$
		(g \cdot f)_{\LL_0^0, \LL_2^0} \text{ and } (f \cdot g)_{\LL_1^1, \LL_3^1} \text{ are homotopy equivalences}.
	$$
	Finally, $f_{\LL_0^0, \LL_1^1}$ is a homotopy equivalence by \cref{l: 2-out-of-6 property for homotopy equivalences}, $\LL_0^1$ is homotopy equivalent to $\LL_1^1$ by \cref{l: e is a homotopy equivalence}, hence $\LL_0^0$ and $\LL_0^1$ are homotopy equivalent.
\end{proof}

\subsection{Strict units}
\label{subs: Strict units in geometric Fukaya category}
Keep the notation $\S$, $\C$, and $\W = \Wgp(\S; \C)$ from before.
Recall from \nameref{s: Introduction} that a collection of Lagrangian branes is called tame if the underlying Lagrangians form a transverse no more than countable set.
In this subsection we prove the following.

\begin{prop} \label{p: saturated collection of systems of strict units on Fukaya category}
	There exists a saturated collection of systems of strict units $\SU$ (see \cref{d: saturated collection of systems of strict units}) on $\W$ such that for any tame collection $\bL^0$ of Lagrangian branes there exists a system of strict units $\bL^\bullet \in \SU$ on $\bL^0$.
\end{prop}

\cref{p: saturated collection of systems of strict units on Fukaya category,p: W_SU is independent of choices} then give a strictly unital $A_\infty$-category $\W_{\SU} = \Wgsu(\S; \C)$, well-defined up to homotopy equivalence.

In \cref{subsubs: systems of strict units in Fukaya category} we describe how to get systems of strict units in $\W$.
In \cref{subsubs: Saturated collection of systems of strict units for Fukaya category} we describe $\SU$ and prove \cref{p: saturated collection of systems of strict units on Fukaya category}.

\subsubsection{Systems of strict units}
\label{subsubs: systems of strict units in Fukaya category}
Let $L^0 = (I, \gamma^0, P^0)$ be a Lagrangian.
Then there exists a transverse set $\{L^j\}_{j = 0}^\infty$ such that
\begin{enumerate}
	\item if $j_0 \ne j_1$, then $L^{j_1}$ is a transverse perturbation of $L^{j_0}$;
	\item $L^j$ tends to $L^0$ in $\Lag_I(\S)$ as $j \to \infty$;
	\item for pairwise distinct $j_0, j_1, j_2 \ge 0$ there exists a composite Lagrangian $L = L_{h_0, h_1, h_2}(L^{j_0}, L^{j_1}, L^{j_2})$ that contains the transverse double $L_{h_0, h_1}(L^{j_0}, L^{j_1})$ of $L^{j_0}$ nicely;
	\item for any
	$j_0 < \ldots < j_{l - 2} < \min\{j_{l - 1}, j_l\} < \max\{j_{l - 1}, j_l\} < j_{l + 1} < \ldots < j_n$
	there exists a composite Lagrangian $L = L_{h_0, \ldots, h_n}(L^{j_0}, \ldots, L^{j_n})$ that contains the transverse double $L_{h_{l - 1}, h_{l}}(L^{j_{l - 1}}, L^{j_l})$ nicely.
\end{enumerate}

\begin{figure}[!htbp]
	\begin{subfigure}[t]{0.45\linewidth}
		\includegraphics[width=\linewidth, trim={0 0 307.37bp 0}, clip]{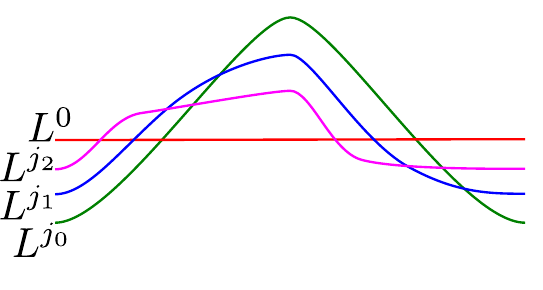}
		\caption{A circle component of $I$.}
		\label{subfig: system of strict units for a circle}
	\end{subfigure}%
	\hfill%
	\begin{subfigure}[t]{0.53\linewidth}
		\includegraphics[width=\linewidth, trim={260.79bp 0 0 0}, clip]{SystemOfStrictUnits.pdf}
		\caption{A segment component of $I$.}
		\label{subfig: system of strict units for a segment}
	\end{subfigure}
	\caption{$L^0, L^{j_0}, L^{j_1}, L^{j_2}$ for $0 < j_0 < j_1 < j_2$}
	\label{fig: system of strict units}
\end{figure}

Such a transverse set can be obtained, for example, by perturbing each connected component $C$ of $I$ that is a circle as in \cref{subfig: system of strict units for a circle} (the parts of $L^j$ that are not depicted run parallel to $L^0$) in a small enough neighbourhood of a point $\gamma(t)$ for $t \in C$ and by perturbing each connected component $C$ of $I$ that is a segment as in \cref{subfig: system of strict units for a segment} (the picture is of the pullback to a Weinstein neighbourhood of $\res{\gamma}{C}$).

Now let $\{L^0_k\}_{k \in K}$ be a transverse set of Lagrangians in $\S$ and $\{\LL_i^0\}_{i \in I}$ a collection of Lagrangian branes in $(\S; \C)$ such that each $\LL_i^0$ is supported on some $L^0_{\kappa(i)}$.
For each $k \in K$ pick a transverse set $\{L^j_k\}_{j = 0}^\infty$ as above and for each $i \in I$ and $j \ge 0$ let $\LL^j_{i}$ be the transverse perturbation of $\LL_i^0$ supported on $L^j_{\kappa(i)}$.
Note that $\LL_i^{j_1}$ is a transverse perturbation of $\LL_i^{j_0}$ for any $i \in I$ and $j_0 \ne j_1$ hence we have a homotopy equivalence
$e_i^{j_0, j_1} = \id_{\LL_i^{j_0}, \LL_i^{j_1}} \in \W(\LL_i^{j_0}, \LL_i^{j_1})$.

\begin{lemma} \label{l: system of strict units on a collection of Lagrangian branes}
	$\bL^\bullet = \left\{\left(\left\{\LL^j_i\right\}_{j = 0}^\infty, \left\{e_i^{j_0, j_1}\right\}^{j_0 \ne j_1}\right)\right\}_{i \in I}$
	is a system of strict units (see \cref{subsubs: system of strict units}) on the collection $\bL^0 = \{\LL_i^0\}_{i \in I}$ of objects in $\W$.
\end{lemma}
\begin{proof}
	We check properties (1)--(4) from \cref{subsubs: system of strict units}.
	Property (1) follows from (3) above and \cref{l: composition with perturbed morphisms}.

	To prove Property (2) fix $i_0, i_1 \in I$ and let $k_0 = \kappa(i_0)$ and $k_1 = \kappa(i_1)$.
	By \cref{l: composition with perturbed morphisms} it suffices to show that there exists $N > 0$ such that for any $j_0, j_1, j_1' \in \{0\} \cup \ZZ_{\ge N}$ satisfying $j_0 < j_1 < j_1'$ there exists a composite Lagrangian
	$L = L_{h_0, h_1, h_1'}(L_{k_0}^{j_0}, L_{k_1}^{j_1}, L_{k_1}^{j_1'})$
	that contains the transverse double $L_{h_1, h_1'}(L_{k_1}^{j_1}, L_{k_1}^{j_1'})$ nicely.
	If $k_0 = k_1$, then this follows from (3) above.
	If $k_0 \ne k_1$, then this is true since $L^0_{k_0}$ and $L^0_{k_1}$ are transverse, $L^j_{k_0} \xrightarrow{j \to \infty} L^0_{k_0}$, and $L^j_{k_1} \xrightarrow{j \to \infty} L^0_{k_1}$.
	Property (3) is proved analogously.

	Property (4) follows from \cref{l: no higher products involving perturbed morphisms} and conditions (2) and (4) above.
\end{proof}

\subsubsection{Saturated collection of systems of strict units}
\label{subsubs: Saturated collection of systems of strict units for Fukaya category}

\cref{p: characterisation of isotopy} and the fact that there are countably many classes of immersed curves up to a reparametrisation of the domain and the action of $\Diff(\S)$ imply that the number of connected components of the space $\Lag_I(\S)$ of all Lagrangians with given domain $I$ is countable.
The number of different domains is also countable.
Therefore, there exists a transverse set of Lagrangians $\{L_k^0\}_{k = 0}^\infty$ in $\S$ such that for any Lagrangian $L$ in $\S$ there exists a Lagrangian $L_k^0$ isotopic to $L$.

Let $\SU$ be the collection of all systems of strict units obtained as in \cref{subsubs: systems of strict units in Fukaya category} from a no more than countable transverse set of Lagrangians $\{L_k^0\}_{k = 0}^\infty$ and some set of Lagrangian branes each of which is supported on one of $L_k^0$.

\begin{proof}[Proof of \cref{p: saturated collection of systems of strict units on Fukaya category}]
The extra condition on $\SU$ in \cref{p: saturated collection of systems of strict units on Fukaya category} holds by construction, it remains to prove $\SU$ is saturated.

We prove $\SU$ satisfies (1) from \cref{d: saturated collection of systems of strict units}.
Let $\bL^\bullet \in \SU$ be a system of strict units on a collection $\LL^0 = \{\LL_i^0\}_{i \in I}$ of all Lagrangian brane structures on a transverse set of Lagrangians $\{L_k^0\}_{k = 0}^\infty$ such that any Lagrangian in $\S$ is isotopic to one of $L_k^0$.
Then by construction and \cref{p: W(S) has enough homotopy equivalences} for any Lagrangian brane $\LL$ there exists a Lagrangian brane $\LL_i^0, i \in I$ homotopy equivalent to $\LL$.

We prove $\SU$ satisfies (2) from \cref{d: saturated collection of systems of strict units}.
Recall the space $C^\infty_{[0, 1]}(\S, \RR)$ of Hamiltonian functions from \cref{subsubs: Hamiltonian diffeomorphisms} and note that the standard Frech{\'e}t structure makes it into a Baire space.
For $H \in C^\infty_{[0, 1]}(\S, \RR)$ and a Lagrangian $L = (I, \gamma, P)$ define
$$
H_*L = (I, \Phi_1^H \circ \gamma, P - g_H \circ \gamma) \text{ where } g_H = \int_0^1 \big(H_s + \lambda(X^H_s)\big)\circ\Phi^H_sds
$$
and $\Phi_1^H$ is the time-1 map of the Hamiltonian flow of $H$.
Note that $H_*L$ is a Lagrangian since $(\Phi_1^H)^*\lambda - \lambda = dg_H$.

Given any two transverse sets $\{L^i_k\}_{k = 0}^\infty$, $i = 0, 1$ there exists a Hamiltonian $H \in C^\infty_{[0, 1]}(\S, \RR)$ such that
$\{H_*L^0_k\}_{k = 0}^\infty \sqcup \{L^1_k\}_{k = 0}^\infty$
is transverse.
Indeed, for each
$0 \le k_0^0 < \ldots < k_a^0$ and $0 \le k_0^1 < \ldots < k_b^1$ the condition that
$\left\{H_*L^0_{k_0^0}, \ldots, H_*L^0_{k_a^0}, L^1_{k_0^1}, \ldots, L_{k_b^1}\right\}$
is transverse defines an open dense subset of $C^\infty_{[0, 1]}(\S, \RR)$.

If now $\bL^\bullet_i \in \SU$, $i = 0, 1$ are two systems of strict units supported on $\{L^i_k\}_{k = 0}^\infty$, $i = 0, 1$, then by \cref{l: I_bph is a strict equivalence,l: invariance of W under isotopies} we have a strict autoequivalence $\I_H$ of $\W$ such that $\I_{H*}\bL^\bullet_0 \sqcup \bL^\bullet_1 \in \SU$.
\end{proof}

\subsection{Triangulated structure}
\label{subs: Triangulated structure on geometric Fukaya A_infty-pre-category}
Keep the notation $\S$, $\C$, and $\W = \Wgp(\S; \C)$ from before.
In this subsection we prove the following statement.
\begin{prop} \label{p: W(S) is triangulated}
	The Fukaya category $\W$ is strongly triangulated (see \cref{subsubs: triangulated A_infty-pre-category}).
\end{prop}

This is the last ingredient we need to prove \cref{t: geometric Fukaya A_infty-pre-category is good,t: geometric Fukaya strictly unital A_infty-category} in \cref{subs: Proofs of the main theorems}.

We prove \cref{p: W(S) is triangulated} by definition, i.e.\ we construct a strict totalisation functor $\T\colon \Tw \W \to \W$ from the triangulated envelope defined in \cref{subsubs: twisted complexes}.
In fact, we construct two of them, a geometric and an algebraic totalisation functors $\T_g$ and $\T_a$.
$\T_g$ uses rotation of $\gr$ to define the shift and is defined for any $\C$, we construct it in \cref{subsubs: Shifts of Lagrangian pre-branes,subsubs: Geometric totalisation functor}.
To construct $\T_a$ we need $\C$ to have a shift.
Then we can use the shift in $\C$ to define the shift in $\Tw \Wgp(\S; \C)$ and obtain the strict totalisation functor $\T_a\colon\Tw \Wgp(\S; \C) \to \Wgp(\S; \C)$, we do this in \cref{subsubs: Algebraic totalisation functor}.
In \cref{subsubs: From algebraic to geometric cones} we prove that cones of degree zero closed morphisms in $\C$ give rise to cones over $\Wgp(\S; \C)$ in two equivalent ways, algebraic and geometric.

\subsubsection{Shifts of Lagrangian pre-branes}
\label{subsubs: Shifts of Lagrangian pre-branes}
Let $\L = (I, \gamma, P, \gr, \M)$ be a Lagrangian pre-brane in $(\S; \C)$ and $\bm\colon I \to \ZZ$ a locally constant function.
We define the \emph{geometric shift of $\L$ by $\bm$} to be
$$
\L[\bm] = (I, \gamma, P, \gr[\bm], \M)
$$
where $\gr[\bm](t)$ is $\gr(t)$ followed by a clockwise rotation of $\dot \gamma(t)$ by $\bm(t)\pi$ for $t \in I$.
We get the corresponding obstruction algebra $\A_{\L[\bm]}$.
Note that $|z|_{\L[\bm]} = |z|_{\L} + \bm(z_+) - \bm(z_-)$ for any $z \in Z = Z_\L$.
Unfortunately, the signs in $\A_{\L[\bm]}$ are different from the signs appearing for twisted complexes in \cref{subsubs: twisted complexes}.
We therefore define another algebra, the \emph{Koszul shift $\A_\L[\bm]$ of $\A_\L$ by $\bm$}.
Namely,
$$
\A_\L[\bm] = \bigoplus_{z \in Z} \A_\L(z)[-|z|_{\L} + \bm(z_-) - \bm(z_+)]
$$
as a graded module and $\mu_{\A_\L[\bm]}$ is given by the Koszul sign rule as in \cref{eq: mu with Koszul sign}, its components
$$
\mu^{z^0; z^n, \ldots, z^1}_{\A_\L[\bm]} \colon \A_\L(z^n) \otimes \ldots \otimes \A_\L(z^1) \to \A_\L(z^0)
$$
are
$$
\mu^{z^0; z^n, \ldots, z^1}_{\A_\L[\bm]}(\ph^n, \ldots, \ph^1) = (-1)^{\sum_{1 \le \alpha < \beta \le n} (\bm(z^\alpha_-) - \bm(z^\alpha_+))(|\ph^\beta| + |z^\beta|_{\L} - 1)}\mu^{z^0; z^n, \ldots, z^1}_{\A_\L}(\ph^n, \ldots, \ph^1).
$$
Note that $\A_{\L}[\bm]$ is a structured module and the degree of $\mu^{n}_{\A_\L[\bm]}$ is $2 - n$ since
$$
\sum_{\alpha = 1}^n (\bm(z^\alpha_-) - \bm(z^\alpha_+)) = \bm(z^0_-) - \bm(z^0_+)
$$
whenever $\mu^{z^0; z^n, \ldots, z^1}_{\A_\L} \ne 0$.
Observe that the summands of $\A_{\L[\bm]}$ and $\A_{\L}[\bm]$ corresponding to $z \in Z$ are the same graded module.

We now construct a na\"ive isomorphism $(\id_Z, \F_\L)\colon \A_{\L}[\bm] \to \A_{\L[\bm]}$, see \cref{subsubs: Naive isomorphisms} for the definition.
For $z \in Z$ we set
\begin{equation} \label{eq: F^z for geometric and Koszul shift}
	\F_\L^z\colon \A_\L(z) \to \A_\L(z), \; \ph \mapsto (-1)^{(\bm(z_-) - \bm(z_+))|\ph|}\nu(z)^{\bm(z_+)} \ph
\end{equation}
where
$$
\nu(z) = \begin{cases}
(-1)^{|z|_{\L} - 1}& \text{if } \left(\res{\frac{d\gamma}{dt}}{t = z_+}, \res{\frac{d\gamma}{dt}}{t = z_-}\right) \text{ is a positive basis of } T_zS;\\
(-1)^{|z|_{\L}}& \text{otherwise}
\end{cases}
$$
if $z$ is a double point and
$$
\nu(z) = \begin{cases}
(-1)^{|z|_{\L}}& \text{if } \res{\frac{d\gamma}{dt}}{t = z_+} \text{ and } \res{\frac{d\gamma}{dt}}{t = z_-} \text{ both point inside or outside } S;\\
(-1)^{|z|_{\L} - 1}& \text{otherwise}
\end{cases}
$$
if $z$ is a Reeb chord.

\begin{lemma} \label{l: naive isomorphism for geometric and Koszul shifts of obstruction algebra}
	$
		(\id_Z, \F_\L)\colon \A_{\L}[\bm] \to \A_{\L[\bm]}
	$
	defined above is a na\"ive isomorphism.
	In particular, $\A_\L[\bm]$ is a structured $A_\infty$-algebra.
\end{lemma}
\begin{proof}
	Since
	$|z|_{\A_{\L}[\bm]} = |z|_{\A_{\L[\bm]}} = |z|_{\L} + \bm(z_+) - \bm(z_-)$,
	and all $\F^z_\L$ are isomorphisms of graded modules, by
	\cref{l: characterisation of naive isomorphism} it remains to check that higher products are preserved.

	Fix $z^0, \ldots, z^n \in Z$ and $\ph^i \in \A_\L(z^i), \, i = 1, \ldots, n$.
	We need to check
	$$
		\F_\L^{z^0} \circ \mu^{z^0; z^n, \ldots, z^1}_{\A_{\L}[\bm]}(\ph^n, \ldots, \ph^1) = \mu^{z^0; z^n, \ldots, z^1}_{\A_{\L[\bm]}}(\F_\L^{z^n}(\ph^n), \ldots, \F_\L^{z^1}(\ph^1)).
	$$
	We compute the sign by which both sides differ from
	$\mu^{z^0; z^n, \ldots, z^1}_{\A_{\L}}(\ph^n, \ldots, \ph^1)$.
	If $n = 1$ and $z^0 = z^1$, then from the definition in \cref{eq: mu^z^0 z^0} we obtain that the sign for both sides is
	$$
		(-1)^{(\bm(z^0_-) - \bm(z^0_+))(|\ph^1| + 1)}\nu(z^0)^{\bm(z^0_+)}.
	$$
	To compute the sign in all other cases, since $\nu$ does not change when we replace a Reeb chord by a double point, by \cref{l: sufficiently self-wrapped model for Lagrangian,p: isotopy induces isomorphism of obstruction algebras} we can assume that $\L$ is sufficiently self-wrapped.
	Then $z^0, \ldots, z^n$ are all double points.
	Fix a polygon $\bu$ contributing to $\mu^{z^0; z^n, \ldots, z^1}_{\A_{\L}}$ by \cref{eq: mu from polygons}.

	We can assume that all relevant components of $I$ are oriented as sides of $\bu$ since changing the orientation of one of the sides only changes the contribution of $\nu$ and of $\eps$ from \cref{d: orientation induced by u} to the signs, and they change by the same factor on both sides.
	A straightforward computation then shows that in this case the sign for both sides is
	$$
		(-1)^{\bm(z^0_-) |z^0|_{\L} + \sum_{0 \le \alpha < \beta \le n} (\bm(z^\alpha_-) - \bm(z^\alpha_+))(|\ph^\beta| + |z^\beta|_\L - 1)}.
	$$
\end{proof}

\subsubsection{Geometric totalisation functor}
\label{subsubs: Geometric totalisation functor}
We construct the strict \emph{geometric totalisation functor} $\T_g\colon \Tw \W \to \W$.
Let
$\bX = (\bigoplus_{\alpha \in A} \LL_\alpha[m_\alpha], \delta_{\bX})$
be a twisted complex over $\W$ with
$\LL_\alpha = (\L_\alpha = (I_\alpha, \gamma_\alpha, P_\alpha, \gr_\alpha, \M_\alpha), \delta_\alpha)$.
We define the Lagrangian brane $\T_g^0(\bX)$ as follows.
First, identify $A$ with $\{0, \ldots, n\}$ in such a way that ${\delta_{\bX}}_{\alpha_0\alpha_1} \ne 0$ implies $\alpha_0 < \alpha_1$.
Secondly, by \cref{l: make actions positive by vertical shift} there exists a composite Lagrangian brane
$(\L, \delta) = \LL_{(0, h_1,\ldots, h_n)}(\LL_0, \ldots, \LL_n)$
that is adapted to $\delta_{\bX}$.
Then, by definition of a twisted complex, $\delta_{\bX} + \delta$ is a bounding cochain in the Koszul shift $\A_{\L}[\bm]$ where
\begin{equation} \label{eq: bm from twisted complex}
\bm\colon \bigsqcup_{\alpha = 0}^n I_\alpha \to \ZZ, \; \res{\bm}{I_\alpha}= m_\alpha.
\end{equation}
We set $\T_g^0(\bX) = (\L[\bm], \F_\L(\delta_{\bX} + \delta))$ where $\F_\L\colon \A_{\L}[\bm] \to \A_{\L[\bm]}$ is the na\"ive isomorphism from \cref{l: naive isomorphism for geometric and Koszul shifts of obstruction algebra}.

Now let
$\bX^i = (\bigoplus_{\alpha \in A^i} \LL^i_\alpha[m^i_\alpha], \delta_{\bX^i}), \, i = 0, 1$
be two twisted complexes over $\W$ with
$\LL_\alpha^i = (\L_\alpha^i = (I_\alpha^i, \gamma_\alpha^i, P_\alpha^i, \gr_\alpha^i, \M_\alpha^i), \delta_\alpha^i)$.
We define $\T_g^1\colon \Tw \W(\bX^0, \bX^1) \to \W(\LL^0, \LL^1)$ where
${\LL^i = (\L^i = (I^i, \gamma^i, P^i, \gr^i, \M^i), \delta^i) = \T_g^0(\bX^i), \, i = 0, 1}$.

Unwrapping the definitions we see that
\begin{equation} \label{eq: Tw W(bX^0 bX^1) = W(LL^0 LL^1)}
	\Tw \W(\bX^0, \bX^1) = \bigoplus_{\substack{\alpha^0 \in A^0\\\alpha^1 \in A^1}}\bigoplus_{z \in Z(\LL^0_{\alpha^0}, \LL^1_{\alpha^1})} \W(z)[-|z| + m^1_{\alpha^1} - m^0_{\alpha^0}] = \W(\LL^0, \LL^1)
\end{equation}
where
$$
\W(z) = \C_{\gamma^1_{\alpha^1}(z_-)}(\C(z)\M^0_{\alpha^0}(z_+), \M^1_{\alpha^1}(z_-)).
$$
Fix $z \in Z(\LL^0_{\alpha^0}, \LL^1_{\alpha^1})$.
By \cref{l: make actions positive by vertical shift} we can find a composite Lagrangian brane $(\L, \delta) = \LL_{h}(\LL^0_{\alpha^0}, \LL^1_{\alpha^1})$ such that $A_{h}(z) > 0$.
We set
$$
\res{\T_g^1}{\W(z)} = \F^z_\L
$$
where $\F_\L\colon \A_{\L}[\bm] \to \A_{\L[\bm]}$ is the na\"ive isomorphism from \cref{l: naive isomorphism for geometric and Koszul shifts of obstruction algebra} for
$$
\bm\colon I^0_{\alpha^0} \sqcup I^1_{\alpha^1} \to \ZZ, \; \res{\bm}{I^i_{\alpha^i}} = m^i_{\alpha^i}.
$$

The following lemma implies \cref{p: W(S) is triangulated}.
\begin{lemma}\label{l: properties of T_g for W}~
	\begin{enumerate}
		\item $\T_g^0(\LL) = \LL$ for any Lagrangian brane $\LL$.
		\item $\T_g^1$ is identity on $\Tw \W(\LL^0, \LL^1) = \W(\LL^0, \LL^1)$ for each transverse pair $(\LL^0, \LL^1)$ in $\W$.
		\item
		$\T_g^1\colon \Tw \W(\bX^0, \bX^1) \to \W(\T_g^0(\bX^0), \T_g^0(\bX^1))$
		is an isomorphism of graded modules for each transverse pair $(\bX^0, \bX^1)$ in $\Tw \W$.
		\item $\T_g\colon \Tw \W \to \W$ is a strict $A_\infty$-functor.
\end{enumerate}
\end{lemma}
\begin{proof}
	(1) follows from our choice of shift for $\LL_0$ in the definition of $\T_g^0$. (2) holds since the sign in \cref{eq: F^z for geometric and Koszul shift} is $1$ when $\bm = 0$.
	As $\F^z$ from \cref{l: naive isomorphism for geometric and Koszul shifts of obstruction algebra} is an isomorphism, we obtain (3).
	It remains to prove (4).

	Let
	$\{\bX^i = (\bigoplus_{\alpha \in A^i} \LL^i_\alpha[m^i_\alpha], \delta_{\bX^i}) \mid i = 0, \ldots, n\}$
	be a transverse set of twisted complexes over $\W$ with
	${\LL_\alpha^i = (\L_\alpha^i = (I_\alpha^i, \gamma_\alpha^i, P_\alpha^i, \gr_\alpha^i, \M_\alpha^i), \delta_\alpha^i)}$.
	Then
	$\{\LL^i = (\L^i, \delta^i) = \T_g^0\bX^i \mid i = 0, \ldots, n\}$
	is a transverse sequence since both transversality conditions are equivalent to
	$$
		\{L^i_{\alpha^i} \mid i = 0, \ldots, n;\, \alpha^i \in A^i\}
	$$
	being a transverse set.
	Let $a^i \in \Tw \W(\bX^{i - 1}, \bX^{i}), \, i = 1, \ldots, n$.
	We need to prove
	\begin{equation} \label{eq: Tot^1 bmu = bmu Tot^1}
		\T_g^1 \bmu^n_{\Tw \W}(a^n, \ldots, a^1) = \bmu^n_{\W}\left(\T_g^1 a^n, \ldots, \T_g^1a^1\right).
	\end{equation}
	By \cref{l: make actions positive by vertical shift} there exists a composite Lagrangian pre-brane
	$\L = \L_h(\L^i_{\alpha^i})_{\alpha^i \in A^i}^{i = 0, \ldots, n}$
	adapted to $a^n, \ldots, a^1$ and all $\delta^{\bX^i}$.
	As before,
	$$
	\delta = \sum_{i = 0}^n \left(\delta_{\bX^i} + \sum_{\alpha \in A^i} \delta^i_\alpha\right)
	$$
	is a bounding cochain in $\A_{\L}[\bm]$ where
	$$
	\bm\colon \bigsqcup_{i = 0}^n \bigsqcup_{\alpha^i \in A^i} I^i_{\alpha^i} \to \ZZ, \; \res{\bm}{I^i_{\alpha^i}} = m^i_{\alpha^i}.
	$$
	By \cref{l: independence of higher products from shifts} the left-hand side of \cref{eq: Tot^1 bmu = bmu Tot^1} is equal to $\F_\L \circ \mu^n_{\delta}(a^n, \ldots, a^1)$, while the right-hand side is equal to $\mu^n_{\F_\L(\delta)}(\F_\L(a^n), \ldots, \F_\L(a^1))$.
	Since $\F_\L$ is a na\"ive morphism, these are equal.
\end{proof}

\subsubsection{Algebraic totalisation functor}
\label{subsubs: Algebraic totalisation functor}
$\Tw$ defines a strict 2-functor $\dgcat_\k \to \dgcat_\k$.
Let $\Tw \C\colon \P^+_2(\S) \to \dgcat_\k$ be the weighted system of dg-categories given by ${\Tw} \circ \C$.
$\Tw$ has a subfunctor $\Sigma$, called \emph{the shift enlargement}, consisting of twisted complexes of the form $(X[m], 0)$ for $m \in \ZZ$.

Let $\L = (I, \gamma, P, \gr, \M)$ be a Lagrangian pre-brane in $(\S; \Sigma \C)$ and $\bm\colon I \to \ZZ$ a locally constant function.
We define the \emph{algebraic shift of $\L$ by $\bm$} to be
$$
\L\{\bm\} = (I, \gamma, P, \gr, \M[\bm])
$$
where $\M[\bm](t) = \M(t)[\bm(t)]$ for $t \in I$ and $\M[\bm](a) = \M(a)$ for any path $a\colon [0, 1] \to I$ (this makes sense since $\bm(a(0)) = \bm(a(1))$ as $\bm$ is locally constant).
A straightforward computation shows that the algebraic shift is the same as the Koszul shift on the obstruction algebra, ${\A_{\L\{\bm\}} = \A_{\L}[\bm]}$.

We construct the \emph{algebraic strict totalisation functor} $\T_a\colon \Tw \W \to \W$ where $\W = \Wgp(\S; \Sigma \C)$.
The construction is similar to the construction in \cref{subsubs: Geometric totalisation functor} but simpler since the Koszul and algebraic shifts coincide.
Let
$\bX = (\bigoplus_{\alpha \in A} \LL_\alpha[m_\alpha], \delta_{\bX})$
be a twisted complex over $\W$ with
${\LL_\alpha = (\L_\alpha = (I_\alpha, \gamma_\alpha, P_\alpha, \gr_\alpha, \M_\alpha), \delta_\alpha)}$.
We define the Lagrangian brane $\T_a^0(\bX)$ as follows.
Identify $A$ with $\{0, \ldots, n\}$ in such a way that ${\delta_{\bX}}_{\alpha_0\alpha_1} \ne 0$ implies $\alpha_0 < \alpha_1$ and take a composite Lagrangian brane
$(\L, \delta) = \LL_{(0, h_1,\ldots, h_n)}(\LL_0, \ldots, \LL_n)$
that is adapted to $\delta_{\bX}$.
Then, by definition of a twisted complex, $\delta_{\bX} + \delta$ is a bounding cochain in the Koszul shift $\A_{\L}[\bm] = \A_{\L\{\bm\}}$ where $\bm$ is as in \cref{eq: bm from twisted complex}.
We set $\T_a^0(\bX) = (\L\{\bm\}, \delta_{\bX} + \delta)$.

Now let
$\bX^i = (\bigoplus_{\alpha \in A^i} \LL^i_\alpha[m^i_\alpha], \delta_{\bX^i}), \, i = 0, 1$
be two twisted complexes over $\W$ with
$\LL_\alpha^i = (\L_\alpha^i = (I_\alpha^i, \gamma_\alpha^i, P_\alpha^i, \gr_\alpha^i, \M_\alpha^i), \delta_\alpha^i)$.
We define $\T_a^1\colon \Tw \W(\bX^0, \bX^1) \to \W(\LL^0, \LL^1)$ where
${\LL^i = (\L^i = (I^i, \gamma^i, P^i, \gr^i, \M^i), \delta^i) = \T_a^0(\bX^i), \, i = 0, 1}$.
The formula \cref{eq: Tw W(bX^0 bX^1) = W(LL^0 LL^1)} still holds and we set $\res{\T_a^1}{\W(z)} = \id$.
Arguing analogously to \cref{l: properties of T_g for W} we obtain the following.

\begin{lemma}\label{l: properties of T_a for W}~
	\begin{enumerate}
		\item $\T_a^0(\LL) = \LL$ for any Lagrangian brane $\LL$.
		\item $\T_a^1$ is identity on $\Tw \W(\LL^0, \LL^1) = \W(\LL^0, \LL^1)$ for each transverse pair $(\LL^0, \LL^1)$ in $\W$.
		\item
		$\T_a^1\colon \Tw \W(\bX^0, \bX^1) \to \W(\T_a^0(\bX^0), \T_a^0(\bX^1))$
		is an isomorphism of graded modules for each transverse pair $(\bX^0, \bX^1)$ in $\Tw \W$.
		\item $\T_a\colon \Tw \W \to \W$ is a strict $A_\infty$-functor.
\end{enumerate}
\end{lemma}

We therefore have two strict totalisation functors
$\T_g, \T_a\colon \Tw \Wgp(\S; \Sigma \C) \to \Wgp(\S; \Sigma \C)$
and \cref{l: comparison of two totalisation functors} applies.
In particular, if $\LL$ is a Lagrangian brane and $m \in \ZZ$, then the geometric shift $\T^0_g (\LL[m])$ is homotopy equivalent to the algebraic shift $\T^0_a (\LL[m])$.

\subsubsection{From algebraic to geometric cones}
\label{subsubs: From algebraic to geometric cones}
Let $\LL_i = (L, \M_i, \delta_i), \, i = 0, 1$ be two Lagrangian branes in $(\S; \C)$ supported on the same graded Lagrangian $L = (I, \gamma, P, \gr)$ and $f\colon \M_0 \to \M_1$ a morphism of local systems compatible with brane structures, closed of degree zero.

To $f$ we associate an \emph{algebraic cone} $\CC_a(f)$ in $\Wgp(\S; \Tw \C)$ and a \emph{geometric cone} $\bC_g(f)$ in $\Tw \Wgp(\S; \C)$ and prove that they are homotopy equivalent in $\Tw \Wgp(\S; \Tw \C)$.

The algebraic cone is the Lagrangian brane $\CC_a(f) = (L, \Cone(f), \delta_{\CC_a(f)})$ in $(\S; \Tw \C)$ where
$$
\Cone(f)(t) = (\M_0(t)[1] \oplus \M_1(t), f(t)[-1]) \in \Tw \C_{\gamma(t)} \text{ for any } t \in I
$$
with the transition maps $\M_0(a) \oplus \M_1(a)$ for any path $a\colon [0, 1] \to I$, and $\delta_{\CC_a(f)} = \delta_0 \oplus \delta_1$.
The compatibility of $f$ with the brane structure implies that $\delta_{\CC_a(f)}$ is indeed a bounding cochain.

To construct the geometric cone we choose a transverse perturbation $\LL_1^1$ of $\LL_1$ and set
${\bC_g(f) =(\LL_0[1] \oplus \LL_1^1, f_{\LL_0, \LL_1^1}[-1])}$.
\cref{l: f_LL_0LL_1 is closed} implies that ${f_{\LL_0, \LL_1^1} \in \Wgp(\S; \C)(\LL_0, \LL_1^1)}$ is a closed morphism of degree $0$ and hence $\bC_g(f)$ is indeed a twisted complex over $\Wgp(\S; \C)$.

\begin{lemma} \label{l: algebraic and geometric cones}
	Let $\LL_0$, $\LL_1$, and $f$ be as above.
	Then the algebraic cone $\CC_a(f)$ is homotopy equivalent to the geometric cone $\bC_g(f)$ in $\W = \Tw \Wgp(\S; \Tw \C)$.
\end{lemma}
\begin{proof}
	Let
	$\left\{\left(\left\{\LL^j_i\right\}_{j = 0}^\infty, \left\{e_i^{j_0, j_1}\right\}^{j_0 \ne j_1}\right)\right\}_{i = 0, 1, a}$
	be a system of strict units as in \cref{l: system of strict units on a collection of Lagrangian branes} with $\LL_0^0 = \LL_0$, $\LL_0^1 = \LL_1$, $\LL_a^0 = \CC_a(f)$, and $\LL_1^1$ as in the construction of $\bC_g(f)$.
	Let $\iota^0\colon \M_0 \to \Cone(f)$ and $\iota^1\colon \M_1 \to \Cone(f)$ be the inclusions, and $p^0\colon \Cone(f) \to \M_0$ and $p^1\colon \Cone(f) \to \M_1$ the projections.
	These are morphisms of local systems compatible with bounding cochains.
	$\iota^1$ and $p^1$ have degree $0$, $|\iota^0| = -1$, and $|p^0| = 1$.
	$\iota^1$ and $p^0$ are closed, $d\iota^0 = \iota^1 \cdot f$, and $d p^1 = f \cdot p^0$.

	Consider the following chain of morphisms in $\W$
	$$
		\left(\begin{aligned}
		\LL_0^0[1] \oplus \LL_1^1\\
		f_{\LL_0^0, \LL_1^1}[-1]
		\end{aligned}
		\right)
		\xrightarrow{\begin{pmatrix} \iota^0_{\LL_0^0, \LL_a^2}[-1] & \iota^1_{\LL_1^1, \LL_a^2}\end{pmatrix}}
		\LL_a^2
		\xrightarrow{\begin{pmatrix} {p^0_{\LL_a^2, \LL_0^3}[1]}\\\\ {p^1_{\LL_a^2, \LL_1^4}}\end{pmatrix}}
		\left(\begin{aligned}
		\LL_0^3[1] \oplus \LL_1^4\\
		f_{\LL_0^3, \LL_1^4}[-1]
		\end{aligned}
		\right)
		\xrightarrow{\begin{pmatrix} \iota^0_{\LL_0^3, \LL_a^5}[-1] & \iota^1_{\LL_1^4, \LL_a^5}\end{pmatrix}}
		\LL_a^5.
	$$
	By \cref{l: f_LL_0LL_1 is closed,l: composition with perturbed morphisms,cor: perturbed f circ g} these three morphisms are closed of degree 0.
	By \cref{l: composition with perturbed morphisms,l: no higher products involving perturbed morphisms} the composition of the first two morphisms is
	$$
		\begin{pmatrix}
			\id_{\LL_0^0, \LL_0^3} & 0\\
			0 & \id_{\LL_1^1, \LL_1^4}
		\end{pmatrix}
	$$
	and of the second two morphisms is $\id_{\LL_a^2, \LL_a^5}$.
	The former is a homotopy equivalence by \cref{l: e is a homotopy equivalence,l: upper triangular morphisms of twisted complexes} applied to $\Wgp(\S; \Tw \C)$, which is strongly triangulated with enough homotopy equivalences by \cref{p: W(S) is triangulated,p: W(S) has enough homotopy equivalences}, and the latter is a homotopy equivalence by \cref{l: e is a homotopy equivalence}.
	As $\W$ has enough homotopy equivalences by \cref{p: homotopy equivalence in triangulated A_infty-pre-category}, by \cref{l: 2-out-of-6 property for homotopy equivalences} these three morphisms are homotopy equivalences.
	In particular, the first morphism is a homotopy equivalence from $\bC_g(f)$ to a Lagrangian brane homotopy equivalent to $\CC_a(f)$ hence $\bC_g(f)$ and $\CC_a(f)$ are homotopy equivalent.
\end{proof}

\subsection{Embedding of topological Fukaya category}
\label{subs: Embedding of topological Fukaya category}
Fix a weighted Liouville surface
${(\S = (S, Q, \eta, \lambda); \C\colon \P^+_2(\S) \to \dgcat_\k)}$.

In \cref{subsubs: Arc brane systems and tame collections of Lagrangian branes} we prove that the $A_\infty$-category $\V_\bGamma$ of an arc brane system defined in \cref{subs: Fukaya category from arc systems} is strictly equivalent to the $A_\infty$-category $\W_{\bGamma^\bullet}$ of the corresponding tame collection of Lagrangian branes $\bGamma^0$.
In \cref{subsubs: From disc sequence to homotopy equivalence} we prove that disc sequences induce morphisms of twisted complexes in $\V_\bGamma$.
In \cref{subsubs: Image of the topological Fukaya category} we prove that the subcategories generated by full arc brane systems coincide, see \cref{l: any arc brane is generated by a full arc brane system}.

\subsubsection{Arc brane systems and tame collections of Lagrangian branes}
\label{subsubs: Arc brane systems and tame collections of Lagrangian branes}
Note that an arc brane $(\gamma, \gr, M)$ from \cref{subsubs: Arc systems and boundary paths} gives a Lagrangian brane $([0, 1], \gamma, P, \gr, \M, 0)$ from \cref{d: Lagrangian brane} via choosing arbitrary $P\colon[0, 1] \to \RR$ satisfying $\gamma^*\lambda + dP = 0$ and turning $M$ into a local system on $[0, 1]$ by transition functors of $\gamma^*\C$.

Let $\bGamma = \{\LL_i\}_{i \in I}$ be an arc brane system in $(\S; \C)$.
Then the corresponding collection of Lagrangian branes $\bGamma^0 = \{\LL_i^0\}_{i \in I}$ is tame and hence there exists a system of strict units
$$
\bL^\bullet = \left\{\left(\left\{\LL^j_i\right\}_{j = 0}^\infty, \left\{e_i^{j_0, j_1}\right\}^{j_0 \ne j_1}\right)\right\}_{i \in I}\in \SU.
$$

We construct a strict equivalence $\F$ from $\V_\bGamma$ (which we have not yet proved to satisfy $A_\infty$-relations) into the strictly unital $A_\infty$-category $\W_{\bL^\bullet}$.
On objects it is $\F^0(\LL_i) = \LL_{i}^0$.
To define $\F$ on morphisms note that $\LL^j_i$ is given by \cref{subfig: system of strict units for a segment} and for $i_0, i_1 \in I$ and $j \gg 0$ the underlying curves of $\LL^0_{i_0}$ and $\LL^j_{i_1}$ do not intersect and a Reeb chord $y$ from $\LL^0_{i_0}$ to $\LL^j_{i_1}$ gives a boundary path $\tilde y$ from $\LL_{i_0}$ to $\LL_{i_1}$ of the same degree, and the correspondence is a bijection  after we have identified the two constant boundary paths.
Moreover, we have an isomorphism of dg-modules
$$
\F^{y}\colon \V_\bGamma(\tilde y) = \C_{\tilde y(1)}\left(\C(\tilde y)(M_{i_0}^{y_+}), M_{i_1}^{y_-}\right) \xrightarrow{\C(s_y)} \C_{y(1)}\left(\C(y)(M_{i_0}^{y_+}), \C(s_y)(M_{i_1}^{y_-})\right) = \W(y)
$$
where $s_y$ is the short segment of $\pd S$ from $\tilde y(1)$ to $y(1)$, so that $\C(y) = \C(s_y) \circ \C(\tilde y)$.
The maps above assemble into an isomorphism of graded modules
$$
\F^1\colon \V_\bGamma(\LL_{i_0}, \LL_{i_1}) \to \W(\LL_{i_0}^0, \LL_{i_1}^j) = \W_{\bL^\bullet}(\LL_{i_0}^0, \LL_{i_1}^0)
$$

\begin{lemma} \label{l: A_Gamma embeds to W_bX}
	$\F$ defined above is a strict full equivalence $\V_\bGamma \to \W_{\bL^\bullet}$.
	In particular, $\V_\bGamma$ is a strictly unital $A_\infty$-category.
\end{lemma}
\begin{proof}
	We need to check that for $i_0, \ldots, i_n \in I$ and $0 = j_0 \ll j_1 \ll \ldots \ll j_n$
	$$
		\F^{y^0} \circ \bmu_{\V_\bGamma}^{\tilde y^0; \tilde y^n, \ldots, \tilde y^1}(\ph^n, \ldots, \ph^1) = \bmu_{\W}^{y^0; y^n, \ldots, y^1}(\F^{y^n}(\ph^n), \ldots, \F^{y^1}(\ph^1))
	$$
	where
	$\ph^\alpha \in \V_\bGamma(\tilde y^\alpha) \subset \V_\bGamma(\LL_{i_{\alpha - 1}}, \LL_{i_\alpha})[|y^i|]$
	and
	$\F^{y^\alpha}(\ph^\alpha) \in \W(y^\alpha) \subset \W(\LL_{i_{\alpha - 1}}^{j_{\alpha - 1}}, \LL_{i_\alpha}^{j_\alpha})[|y^i|]$
	for $\alpha = 1, \ldots, n$ and
	$\F^{y^0}\colon \V_\bGamma(\tilde y^0) \subset \V_\bGamma(\LL_{i_{0}}, \LL_{i_n})[|y^0|] \to \W(y^0) \subset \W(\LL_{i_{0}}^{j_{0}}, \LL_{i_n}^{j_n})[|y^0|]$.
	This is indeed the case as \cref{eq: mu for y^0 = y^1} corresponds to \cref{eq: mu^z^0 z^0}, \cref{eq: mu for y^0 = y^2y^1} corresponds to \cref{eq: mu from polygons} for a degenerate polygon, while \cref{eq: mu for a polygon with y^n = y^0z,eq: mu for a polygon with y^1 = zy^0} are particular cases of \cref{eq: mu from polygons} for a non-degenerate polygon.
\end{proof}

\subsubsection{From disc sequence to homotopy equivalence}
\label{subsubs: From disc sequence to homotopy equivalence}
Let $\bGamma$ be an arc brane system, $n \ge 2$, $\LL_i = (\gamma_i, \gr_i, M_i)$, $i = 1, \ldots, n$ arc branes from $\bGamma$, $(y^n, \ldots, y^1)$ a disc sequence bounding a region cut out by $\Gamma_\bGamma$ with boundary path $y^i$ from $\gamma_{i - 1}$ to $\gamma_i$ (the index is mod $n$), and $\ph^i \in \V_\bGamma(y^i)$, $i = 1, \ldots, n$ closed morphisms.
Suppose $\sum_{i = 1}^n (|\ph^i| + |y^i|) = n - 2$ and set
$\psi^i = \bmu^{c^i; y^{i + n}, \ldots, y^{i + 1}}(\ph^{i + n}, \ldots, \ph^{i + 1}) \in \V_\bGamma(\LL_i, \LL_i)$
where $c^i$ is the constant boundary path at $\gamma_i$.
A straightforward computation similar to \cite{HKK17}*{p. 275} gives the following.

\begin{lemma} \label{l: morphisms from a disc sequence}
	In the setup above,
	$$
		\bL = \left(\bigoplus_{i = 1}^{n - 1} \LL_i\left[\sum_{1 \le j \le i}\left(|\ph^j| + |y^j|\right) - i + 1\right], \sum_{2 \le i < n}\ph^i\left[|\ph^i| + |y^i| - 1\right]\right)
	$$
	is a twisted complex.
	The morphisms
	$$
	\bph^1 = \ph^1\left[|\ph^1| + |y^1|\right]\colon \LL_n \to \bL \text{ and } \bph^n = \ph^n\left[|\ph^n| + |y^n|\right]\colon \bL \to \LL_n
	$$
	are closed of degree zero.
	Moreover,
	$$
	\bmu^2(\bph^1, \bph^n) = (-1)^{\dagger}\diag(\psi^1, \ldots, \psi^{n - 1}) \text{ and }\bmu^2(\bph^n, \bph^1) = (-1)^{\dagger}\psi^n
	$$
	where
	$\dagger = |\ph^1| + |y^1| - 1 + \#\{i \in \{1, \ldots, n\}\colon |\ph^i| + |y^i| \text{ is even}\}/2$.
\end{lemma}

Note that \cref{l: morphisms from a disc sequence} implies that $\bL$ is homotopy equivalent to $\LL_n$ if all $H^*\ph^i$ are invertible and the weight of the disc bounded by $(y^n, \ldots, y^1)$ is an invertible morphism of dg-functors.

\subsubsection{Image of the topological Fukaya category}
\label{subsubs: Image of the topological Fukaya category}
Let $\bL^\bullet \in \SU$ be a system of strict units on a full tame collection of Lagrangian branes $\bL^0$ and $\bGamma$ a full arc brane system such that the corresponding tame collection of Lagrangian branes $\bGamma^0$ is contained in $\bL^0$.
Consider the strict full embedding $\V_\bGamma \hookrightarrow \W_{\bL^\bullet}$ given in \cref{subsubs: Arc brane systems and tame collections of Lagrangian branes}.
Denote by $\W_\bGamma$ the smallest full $A_\infty$-subcategory of $\W_{\bL^\bullet}$ containing the image of $\V_\bGamma$ and closed under taking shifts, cones, and homotopy equivalent objects.

\begin{figure}[!htbp]
	\centering
	\includegraphics[width=0.45\linewidth]{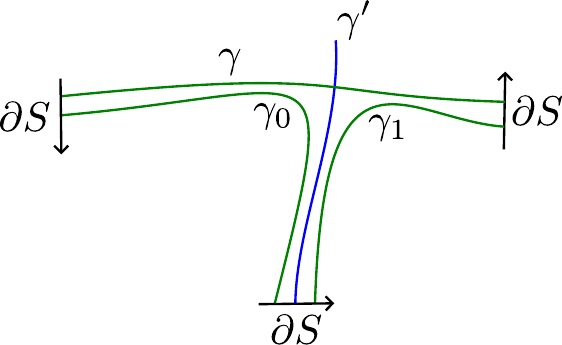}
	\caption{Replacing $\gamma$ with $\gamma_0$ and $\gamma_1$.}
	\label{fig: arc replacement}
\end{figure}

\begin{lemma} \label{l: any arc brane is generated by a full arc brane system}
	 Let $\LL$ be a Lagrangian brane corresponding to an arc brane.
	 Then $\LL$ is homotopy equivalent in $\Wgp(\S; \C)$ to an object of $\W_\bGamma$.
\end{lemma}
\begin{proof}
	By \cref{p: W_SU is independent of choices,p: saturated collection of systems of strict units on Fukaya category} we can assume that $\W_{\bL^\bullet}$ contains all Lagrangian branes supported on any no more than countable transverse set of Lagrangians.
	Replacing $\LL$ with an isotopic brane we can obtain a brane $\LL'$ which is transverse to all branes from $\bGamma$.
	By \cref{p: W(S) has enough homotopy equivalences} $\LL'$ is homotopy equivalent to $\LL$ in $\Wgp(\S; \C)$.
	So we need to prove $\LL' \in \bL^0$ is an object of $\W_{\bGamma}$.

	As $\W_{\bL^\bullet}$ is triangulated by \cref{l: W_SU is triangulated,p: W(S) is triangulated}, any twisted complex on objects of $\W_{\bGamma}$ is homotopy equivalent to an object of $\W_{\bGamma}$.
	By \cref{l: homotopy equivalence of local system induces homotopy equivalence of branes,l: properties of T_a for W,l: algebraic and geometric cones} any Lagrangian brane on an arc from $\Gamma$ lies in $\W_{\bGamma}$, so we can assume that $\bGamma$ contains all arc branes supported on an arc from $\Gamma$.
	We can also assume that $\bGamma$ is supported on a finite full arc system $\Gamma$ since any full arc system has a finite full arc subsystem.

	We prove that $\LL'$ lies in $\W_{\bGamma}$ by induction on the number of intersection points between its underlying curve $\gamma'$ and arcs from $\Gamma$.
	If $\LL'$ does not intersect any of the arcs in $\bGamma$, then we are done by \cref{l: morphisms from a disc sequence}.
	If $t \in [0, 1]$ is the smallest with the property that $\gamma'(t)$ lies on some $\gamma \in \Gamma$, then consider the arcs $\gamma_0$ and $\gamma_1$ that run close to $\Im \gamma \cup \Im \gamma'$ as in \cref{fig: arc replacement}.
	Again by \cref{l: morphisms from a disc sequence} all Lagrangian branes supported on $\gamma_0$ and $\gamma_1$ lie in $\W_{\bGamma}$ so we can replace $\Gamma$ with another full arc system $(\Gamma \setminus \{\gamma\}) \cup \{\gamma_0, \gamma_1\}$ and reduce the number of intersections of $\gamma'$ with arcs from $\Gamma$ by one.
\end{proof}

\subsection{Proofs of the main theorems}
\label{subs: Proofs of the main theorems}
We collect the results of this section to prove \cref{t: geometric Fukaya A_infty-pre-category is good,t: geometric Fukaya strictly unital A_infty-category,t: topological Fukaya category inside geometric,t: A_infty-categories of arc systems}.

\begin{proof}[Proof of \cref{t: geometric Fukaya A_infty-pre-category is good}]
	(1) follows from \cref{p: W is well-defined,p: W(S) is triangulated,p: W(S) has enough homotopy equivalences}. (2) follows from \cref{p: W(S) has enough homotopy equivalences}. (3) follows from \cref{l: invariance of W under isotopies,cor: invariance of W of the choice of Liouville form}.
\end{proof}

\begin{proof}[Proof of \cref{t: geometric Fukaya strictly unital A_infty-category}]
	(1) follows from \cref{p: isotopy induces isomorphism of obstruction algebras}, definition of $\W_{\bL^\bullet}$ in \cref{subsubs: A_infty-category spanned by a collection of objects} and the construction of $\SU$ since $L^j_i$ tends to $L^0_i$ as $j \to \infty$. The homotopy equivalence part follows from \cref{p: A_infty category spanned by a collection of objects is well-defined}.

	(2) is immediate from construction and (3) is part of the definition of a saturated collection of systems of strict units. (4) follows formally, see \cref{subsubs: Invariance}. (5) follows from \cref{l: W_SU is triangulated}.
\end{proof}

\begin{proof}[Proof of \cref{t: topological Fukaya category inside geometric}]
	(1) is \cref{l: A_Gamma embeds to W_bX}. (2) follows from \cref{l: any arc brane is generated by a full arc brane system} since given any two full arc brane systems $\bGamma_0, \bGamma_1$ satisfying the assumption we have $\bGamma_{1 - i}^0 \subset \Ob \W_{\bGamma_{i}}$ for $i = 0, 1$.
\end{proof}

\begin{proof}[Proof of \cref{t: A_infty-categories of arc systems}]
	(1) follows from \cref{l: A_Gamma embeds to W_bX}. (2) follows from \cref{t: topological Fukaya category inside geometric}\,(2) since $\W_{\bL^\bullet}$ is triangulated for a full tame collection $\bL^0$ by \cref{t: geometric Fukaya strictly unital A_infty-category}\,(5) and for any two arc brane systems $\bGamma_0, \bGamma_1$ we can embed $\V_{\bGamma_0}$ and $\V_{\bGamma_1}$ in the same $\W_{\bL^\bullet}$ for a full tame collection $\bL^0$ by \cref{t: geometric Fukaya strictly unital A_infty-category}\,(3).
\end{proof}
\end{document}